\documentclass[reqno]{amsproc}
\usepackage{amsmath}
\usepackage{amsbsy}
\usepackage{amsfonts}
\usepackage{amssymb}
\usepackage{graphics}
\usepackage{graphicx}
\usepackage{array}
\usepackage{epsfig}
\usepackage[dvipdfm]{hyperref}
\usepackage[all]{xy}%

\newtheorem{theorem}{Theorem}[section]
\newtheorem{corollary}[theorem]{Corollary}

\newtheorem{lemma}[theorem]{Lemma}
\newtheorem{conjecture}[theorem]{Conjecture}
\newtheorem{proposition}[theorem]{Proposition}

\theoremstyle{definition}
\newtheorem{definition}[theorem]{Definition}

\newtheorem{example}[theorem]{Example}
\newtheorem{examples}[theorem]{Examples}
\newtheorem{exercise}[theorem]{Exercise}
\newtheorem{application}[theorem]{Application}
\newtheorem{remark}[theorem]{Remark}
\newtheorem{note}[theorem]{Note}
\numberwithin{equation}{section} 
\begin{document}
\author[D.I. Dais]{Dimitrios I. Dais}
\address{University of Crete, Department of Mathematics, Division Algebra and Geometry,
Knossos Avenue, P.O. Box 2208, GR-71409, Heraklion, Crete, Greece}
\email{ddais@math.uoc.gr}
\author[M. Henk]{Martin Henk}
\address{Technical University Otto von Guericke, Institut f\"{u}r Algebra und
Geometrie, PSF 4120, D-39016 Magdeburg, Germany}
\email{henk@math.uni-magdeburg.de}
\author[G.M. Ziegler]{G\"{u}nter M. Ziegler}
\address{TU-Berlin, Institut f\"{u}r Mathematik, MA 6.2, Stra\ss e des 17 Juni 136, D-10623
Berlin, Germany} \email{ziegler@math.tu-berlin.de}
\date{}
\title[Crepant Resolutions of Gorenstein AQ-Singularities]{On the Existence of Crepant Resolutions\\of Gorenstein Abelian Quotient Singularities\\in Dimensions $\geq{4}$}

\begin{abstract}
For which finite subgroups $G$ of SL$(r,\mathbb{C}),$ $r\geq4,$
are there \textit{crepant} \textit{desingularizations }of the
quotient space $\mathbb{C}^{r}/G$? A complete answer to this
question (also known as \textquotedblleft Existence
Problem\textquotedblright\ for such desingularizations) would
classify all those groups for which the high-dimensional versions
of \textit{McKay correspondence }are valid. In the paper we
consider this question in the case of \textit{abelian} finite
subgroups of SL$(r,\mathbb{C})$ by using techniques from toric and
discrete geometry. We give two necessary existence conditions,
involving the Hilbert basis elements of the cone supporting the
junior simplex, and an Upper Bound Theorem, respectively.
Moreover, to the known series of Gorenstein abelian quotient
singularities admitting projective, crepant resolutions (which are
briefly recapitulated) we add a new series of non-c.i. cyclic
quotient singularities having this property.
\end{abstract}

\maketitle


\section{Introduction\label{INTRO}}

\noindent{}The \textit{McKay correspondence} can be understood as
a \textquotedblleft bridge\textquotedblright\ between the
irreducible representations (or, dually, the conjugacy classes) of
finite subgroups $G$ of the special linear group
SL$(r,\mathbb{C})$ and the (co)homology of
$\ \widehat{X}$ 's, for any crepant desingularization $\widehat{X}%
\longrightarrow X$ of $X=\mathbb{C}^{r}/G.$ (Here, \textit{crepant} simply
means that the canonical divisor $K_{\widehat{X}}$ of $\widehat{X}$ is
trivial.) Before we are going to focus on the constant companion of this
correspondence (i.e., the so-called \textquotedblleft Existence
Problem\textquotedblright, whenever $r\geq4$), let us briefly recall some
basic facts about quotient singularities and summarize both classical and
recent results concerning it.\medskip

\noindent{}$\bullet$ \textbf{Quotient singularities}. Let $G$ be a finite
subgroup of GL$\left(  r,\mathbb{C}\right)  $ which is \textit{small}, i.e.,
with no pseudoreflections, acting linearly on $\mathbb{C}^{r}$, and let
\[
\varpi:\mathbb{C}^{r}\longrightarrow\mathbb{C}^{r}/G=\text{Spec}%
(\mathbb{C}\left[  \mathfrak{x}_{1},\ldots,\mathfrak{x}_{r}\right]  ^{G})
\]
be the quotient map. Denote by $\left(  \mathbb{C}^{r}/G,\left[
\mathbf{0}\right]  \right)  $ the corresponding \textit{quotient singularity}
as germ at $\left[  \mathbf{0}\right]  :=\varpi\left(  \mathbf{0}\right)  $.
Quotient singularities are known to be normal and Cohen-Macaulay (see \cite[p.
129]{Lamotke}, \cite[Proposition 13]{Hochster-Eagon}, and \cite[Thm.
3.2]{Stanley1}).

\begin{proposition}
[Singular locus]\label{SLOC}If $G$ is a small finite subgroup of
\emph{GL}$\left(  r,\mathbb{C}\right)  $, then the singular locus of the
entire geometric quotient space $\mathbb{C}^{r}/G$ equals
\[
\text{\emph{Sing}}\left(  \mathbb{C}^{r}/G\right)  =\varpi\left(  \left\{
\mathbf{z}\in\mathbb{C}^{r}\ \left\vert \ G_{\mathbf{z}}\neq\left\{
\text{\emph{Id}}\right\}  \right.  \right\}  \right)
\]
where $G_{\mathbf{z}}:=\left\{  g\in G\ \left\vert \ g\cdot\mathbf{z=z}%
\right.  \right\}  $ is the isotropy group of $\mathbf{z=}\left(  z_{1}%
,\ldots,z_{r}\right)  \in\mathbb{C}^{r}$.
\end{proposition}

\begin{theorem}
[{Prill's isomorphism criterion, \cite[Thm. 2]{Prill}}]\label{prill-th}Let
$G_{1}$, $G_{2}$ be two small finite subgroups of \emph{GL}$\left(
r,\mathbb{C}\right)  ,$ $r\geq2.$ Then there exists an analytic isomorphism
$\left(  \mathbb{C}^{r}/G_{1},\left[  \mathbf{0}\right]  \right)  \cong\left(
\mathbb{C}^{r}/G_{2},\left[  \mathbf{0}\right]  \right)  $ \emph{(}i.e.,
$\mathcal{O}_{\mathbb{C}^{r}/G_{i},\left[  \mathbf{0}\right]  }=\mathbb{C}%
\{\mathfrak{x}_{1},\ldots,\mathfrak{x}_{r}\}^{G_{i}},$ $i=1,2,$ are isomorphic
as local $\mathbb{C}$-algebras\emph{)} if and only if $G_{1}$ and $G_{2}$ are
conjugate to each other within \emph{GL}$\left(  r,\mathbb{C}\right)  $.
\end{theorem}

\noindent{}Hence, the classification of quotient singularities is reduced to
the classification (up to conjugacy) of the small finite subgroups of
GL$\left(  r,\mathbb{C}\right)  $, or equivalently, to those of U$\left(
r,\mathbb{C}\right)  $, as it follows from the next Lemma.

\begin{lemma}
[{\cite[Lemma 4.2.15 (i), p. 82.]{Springer}}]\label{SU(r)}If $G$ is a finite
subgroup of \emph{GL}$\left(  r,\mathbb{C}\right)  $ \emph{(}resp., of
\emph{SL}$\left(  r,\mathbb{C}\right)  $\emph{)}, then $G$ is conjugate in
\emph{GL}$\left(  r,\mathbb{C}\right)  $ \emph{(}resp., in \emph{SL}$\left(
r,\mathbb{C}\right)  $\emph{)} to a finite subgroup of \emph{U}$\left(
r,\mathbb{C}\right)  $ \emph{(}resp., of \emph{SU}$\left(  r,\mathbb{C}%
\right)  $\emph{).}
\end{lemma}

\noindent Subclasses of quotient singularities $\left(  \mathbb{C}%
^{r}/G,\left[  \mathbf{0}\right]  \right)  $ of special theoretical value are
those dictated by the hierachy of the local Noetherian rings $\mathcal{O}%
_{\mathbb{C}^{r}/G,\left[  \mathbf{0}\right]  }$ (cf. \cite[\S VI.3]{Kunz}):%
\[
\fbox{$%
\begin{array}
[c]{c}%
\ \\
\text{(complete intersection (\textquotedblleft c.i\textquotedblright))
}\Longrightarrow\text{(Gorenstein) }\Longrightarrow\text{(Cohen-Macaulay) }\\
\
\end{array}
$}%
\]
\smallskip

\noindent$\bullet$ \textbf{Gorenstein quotient singularities}. These are
characterized as follows:

\begin{theorem}
[cf. \cite{Hin, Wat1}]\label{Gor}$\left(  \mathbb{C}^{r}/G,\left[
\mathbf{0}\right]  \right)  $ is a Gorenstein quotient singularity if and only
if\emph{\footnote{Note that every finite subgroup of SL$(r,\mathbb{C})$ is
small (cf. \cite[p. 503]{Stanley1}).}} $G\subset$ \emph{SL}$\left(
r,\mathbb{C}\right)  .$
\end{theorem}

\begin{note}
[Classification for small $r$]The finite subgroups of SL$\left(
r,\mathbb{C}\right)  $ have been completely classified (up to
conjugacy) only for $r$ small.\smallskip\ \newline(i) For $r=2$
the classification appears already in Klein's book
\textquotedblleft Vorlesungen \"{u}ber das
Ikosaeder\textquotedblright\ \cite{Klein} in 1884. (See also
\cite[Thm. on p. 35]{Lamotke} and \cite[\S 4.4]{Springer}). The
list contains the cyclic groups $\mathbf{Cyc}_{k}$ of order
$k\geq2,$ the binary dihedral groups $\mathbf{Dih}_{n-2}$ of order
$4(n-2),$ $n\geq4,$ and the binary tetrahedral, octahedral and
icosahedral groups $\mathbf{T},\mathbf{O},\mathbf{I},$ having
orders $24,48$ and $120,$ respectively.\smallskip\ \newline(ii)
For $r=3$ the main part of the classification goes back to works
of Jordan, Klein, Gordan, Maschke, Valentiner, Wiman, Gerbaldi,
and Blichfeldt, written around the end of the 19th century. For
the corresponding list of groups we refer to Blichfeldt's book
\cite{Blichfeldt}, as well as to Yau \& Yu \cite{Yau-Yu}
(containing some useful updates).\smallskip\ \newline(iii) For
$r=4$ see Hanany \& He \cite{Hanany-He}. In the group theory
literature there are lots of scattered results concerning this
topic for $r\in\{5,\ldots,10\}$ though they are far from giving
complete classifications. One of the main problems seems to be the
appearance of \textquotedblleft individual\textquotedblright\
groups. For instance, already for $r=6$ we meet the complex
irreducible representation of the Hall-Janko sporadic simple group
of order $604800$ within SL$(6,\mathbb{C}).$ Bigger $r$'s lead to
more \textquotedblleft exotic\textquotedblright\ groups. On the
other hand, there exist many families of finite subgroups of
SL$\left( r,\mathbb{C}\right) $ (e.g., the abelian ones,
dihedral-like groups, imprimitive groups etc) being present in
\textit{all} dimensions.
\end{note}

\begin{remark}
Let $g$ be an element of a finite subgroup $G$ of SL$\left(  r,\mathbb{C}%
\right)  $\emph{. }By Lemma \ref{SU(r)} there exists an $h\in$ SL$\left(
r,\mathbb{C}\right)  $ such that $hGh^{-1}\subset$ SU$\left(  r,\mathbb{C}%
\right)  $\emph{. }Since\emph{\ }$hgh^{-1}$ is a unitary matrix, it is known
(see, e.g., \cite[Thm. 10.2, p. 208]{Nering}) that\emph{\smallskip}%
\newline(i)\emph{\ }$hgh^{-1}$ is unitary similar to a diagonal
matrix,\smallskip\newline(ii)\emph{\ }the diagonal entries of this matrix are
the eigenvalues of\emph{\ }$hgh^{-1}$, and\smallskip\newline(iii)\emph{\ }the
eigenvalues of\emph{\ }$hgh^{-1}$\emph{\ }have absolute value\emph{\ }%
$1$.\smallskip\newline Thus, there is a suitable matrix $k\in$ U$\left(
r,\mathbb{C}\right)  $,\emph{\ }so that \emph{\ }%
\[
k\left(  hgh^{-1}\right)  k^{-1}=\left(  kh\right)  g\left(  kh\right)
^{-1}\in\text{SU}\left(  r,\mathbb{C}\right)
\]
is of the form%
\[
\left(  kh\right)  g\left(  kh\right)  ^{-1}=\text{\emph{\ }diag}(e^{2\pi
\sqrt{-1}\ \gamma_{1}},\ldots,e^{2\pi\sqrt{-1}\ \gamma_{r}}),
\]
for some $\gamma_{1},\ldots,\gamma_{r}\in\mathbb{Q\ \cap\ }\left[  0,1\right)
.$
\end{remark}

\begin{definition}
[\textquotedblleft Ages\textquotedblright\ and \textquotedblleft
heights\textquotedblright](i) The \textit{age}\footnote{In \cite[\S 6]%
{Bat-Arch} and \cite[\S 5]{BD} the age of an element $g\in G$ is called the
\textit{weight} of $g.$ Here, we shall adopt the terminology of
\cite{Ito-Reid}. To the word \textit{weight} we ascribe a different meaning.
(See below Definition \ref{deftype}.)}\textit{\ }of an element\emph{\ }$g\in
G$ is defined to be\emph{\ }the sum%
\begin{equation}
\fbox{$%
\begin{array}
[c]{ccc}
& \text{age}\left(  g\right)  :=\gamma_{1}+\gamma_{2}+\cdots+\gamma
_{r}.\smallskip &
\end{array}
$} \label{defage}%
\end{equation}
(ii) The \textit{height} ht$\left(  g\right)  $ of an element\emph{\ }$g\in G$
is defined to be\emph{\ }the rank%
\begin{equation}
\fbox{$%
\begin{array}
[c]{ccc}
& \text{ht}\left(  g\right)  :=\text{rank}(g-\text{Id}_{G}).\smallskip &
\end{array}
$} \label{defht}%
\end{equation}

\end{definition}

\begin{proposition}
[{\cite[Prop. 5.2.]{BD}}]\label{ht-ag}For every $g\in G$ we have
\begin{equation}
\emph{ht}\left(  g\right)  =\emph{ht}\left(  g^{-1}\right)  =\emph{age}\left(
g\right)  +\emph{age}\left(  g^{-1}\right)  . \label{HTAGEREL}%
\end{equation}

\end{proposition}

\begin{note}
\label{NOTEHTAGE}(i) Obviously$,$ $0\leq$ age$\left(  g\right)  \leq r-1,$
with age$(g)=0\Longleftrightarrow g=$ Id$_{G},$ and age$\left(  g_{1}\right)
=$ age$\left(  g_{2}\right)  $ for all pairs $\left(  g_{1},g_{2}\right)  \in
G\times G$ of group elements belonging to the same conjugacy class. Moreover,
$2\leq$ ht$\left(  g\right)  \leq r,$ for all $g\in G\smallsetminus\left\{
\text{Id}_{G}\right\}  .\smallskip$ \newline(ii) The group elements having age
$1$ (resp., age $i\geq2$) are usually called \textit{junior elements} (resp.,
\textit{senior elements}) of $G.$ (Correspondingly, by (i), we may speak of
\textit{junior} (resp., \textit{senior}) \textit{conjugacy classes}, or in
general of \textit{conjugacy classes of age} $i\in\{0,1,\ldots r-1\}$%
).$\smallskip$ \newline(iii) As it was pointed out by Ito \& Reid \cite[Thm.
1.3]{Ito-Reid}, the \textquotedblleft Tate twist\textquotedblright\
\[
G\left(  -1\right)  :=\text{ Hom}(\widehat{\mathbb{Z}}\left(  1\right)
,G),\ (\text{with }\ \widehat{\mathbb{Z}}\left(  1\right)  :=\underset
{\longleftarrow}{\text{ }\lim}\left(  \mathbb{Z\,}/\,d\,\mathbb{Z}\right)  ),
\]
(which is isomorphic to $G$ as long as one makes a concrete choice of roots of
unity) has a canonical grading inherited by the ages of its elements,
invariant under conjugacy in $G\left(  -1\right)  $ (or $G$), and it is
essentially used in Theorem \ref{IRTHM}.\smallskip\ \newline(iv) If $r$ is
even and $\mathbb{C}^{r}/G$ \textit{symplectic} (i.e., $G\subset$
Sp$(r,\mathbb{C})$), let
\[
\mathcal{F}_{\bullet}(\mathbb{C}[G])=\left\{  \left.  F_{k}(\mathbb{C}%
[G])\right\vert k\in\mathbb{Z}_{\geq0}\right\}
\]
denote the increasing filtration of the group algebra
$\mathbb{C}[G]$ defined
by setting%
\[
F_{k}(\mathbb{C}[G]):=\mathbb{C}\text{-span of}\mathbb{\ }\{\left.  g\in
G\ \right\vert \ \text{ht}(g)\leq k\},\ \ \forall k\in\mathbb{Z}_{\geq0}.
\]
$\mathcal{F}_{\bullet}(\mathbb{C}[G])$ is compatible with the algebra
structure on $\mathbb{C}[G].$ Using the induced filtration $\mathcal{F}%
_{\bullet}(\mathcal{Z}G)$ on the center $\mathcal{Z}G$ of $\mathbb{C}[G]$ (see
\cite{Alvarez, Ginz-Kal}) one determines the associated graded algebra%
\begin{equation}
\mathbf{gr}^{\mathcal{F}_{\bullet}}(\mathcal{Z}G):=\left\{  \left.
F_{k+1}(\mathcal{Z}G)/F_{k}(\mathcal{Z}G)\right\vert \ k\in\mathbb{Z}_{\geq
0}\right\}  \label{ASSGRA}%
\end{equation}
whose singificance is revealed in Theorem \ref{GKTHM}.
\end{note}

\noindent{}$\bullet$ \textbf{Quotient c.i.-singularities}. In the mid 1980's
Nakajima \& Watanabe \cite{NW}, and independently Gordeev \cite{Gordeev},
classified the quotient singularities which are complete intersections
(\textquotedblleft c.i's\textquotedblright) in \textit{all} dimensions. Even
to write down without further ado their group lists would demand several
pages. Instead, let us remind a previous result which constitutes the
foundation stone for their classification.\

\begin{theorem}
[Kac \& Watanabe \cite{Ka-Wa}]\label{KWTHM}If $\left(  \mathbb{C}%
^{r}/G,\left[  \mathbf{0}\right]  \right)  $ is a quotient c.i.- singularity,
then $G$ is generated by the set $\{\left.  g\in G\right\vert $ \emph{ht}%
$(g)\leq2\}.$
\end{theorem}

\noindent{}$\bullet$ \textbf{McKay correspondence in dimension }$r=2.$ In
dimension $2,$ the \textit{classical} \textit{McKay correspondence} exploits
the elegance of the invariant theory of finite subgroups of SL$\left(
2,\mathbb{C}\right)  $ and the uniqueness (and simple description) of the
\textit{minimal} desingularization of the quotient spaces $\mathbb{C}^{2}/G.$

\begin{theorem}
[{cf. \cite[II. 9-13]{Klein}, \cite[Ch. II, \S 8]{Lamotke}, \cite[\S 4.5]%
{Springer}}]\label{INVARC2}The quotient spaces $\mathbb{C}^{2}/G=\emph{Spec}%
(\mathbb{C}\left[  \mathfrak{x}_{1},\mathfrak{x}_{2}\right]  ^{G}),$ for $G$ a
finite subgroup of \emph{SL}$\left(  2,\mathbb{C}\right)  $, are minimally
embedded as hypersurfaces $\left\{  \left(  z_{1},z_{2},z_{3}\right)
\in\mathbb{C}^{3}\left\vert \varphi\left(  z_{1},z_{2},z_{3}\right)
=0\right.  \right\}  ,$ i.e.,
\[
\mathbb{C}\left[  \mathfrak{x}_{1},\mathfrak{x}_{2}\right]  ^{G}%
\cong\mathbb{C}\left[  z_{1},z_{2},z_{3}\right]  /\left(  \varphi\left(
z_{1},z_{2},z_{3}\right)  \right)  .
\]
\emph{(The normal form of the ideal generator is given in the 4th column of
Table \ref{Table1}.)}%
\begin{table}[h]
\setlength\extrarowheight{2pt}
\[%
\begin{tabular}
[c]{|l|l|l|l|}\hline
\emph{Nr.} & \emph{Groups} & \emph{Type} & $\varphi\left(  z_{1},z_{2}%
,z_{3}\right)  $\\\hline\hline
$1.$ & $\mathbf{Cyc}_{n+1},\ n\geq1$ & $\underset{}{A_{n}}$ & $z_{1}%
^{n+1}+z_{2}^{2}+z_{3}^{2}$\\\hline
$2.$ & $\mathbf{Dih}_{n-2},\ n\geq4$ & $\underset{}{D_{n}}$ & $z_{1}%
^{n-1}+z_{1}z_{2}^{2}+z_{3}^{2}$\\\hline
$3.$ & $\mathbf{T}$ & $\underset{}{E_{6}}$ & $z_{1}^{4}+z_{2}^{3}+z_{3}^{2}%
$\\\hline
$4.$ & $\mathbf{O}$ & $\underset{}{E_{7}}$ & $z_{1}^{3}z_{2}+z_{2}^{3}%
+z_{3}^{2}$\\\hline
$5.$ & $\mathbf{I}$ & $\underset{}{E_{8}}$ & $z_{1}^{5}+z_{2}^{3}+z_{3}^{2}%
$\\\hline
\end{tabular}
\ \ \ \ \ \ \ \ \
\]
\smallskip
\caption{}\label{Table1}
\end{table}
\setlength\extrarowheight{-2pt}
\end{theorem}

\begin{theorem}
[{cf. \cite{DuVal}, \cite[\S 5]{Durfee}, \cite[Thm. 2 on p. 152]{Lamotke}}]Let
$G$ be finite subgroup of \emph{SL}$\left(  2,\mathbb{C}\right)  $ and
$X=\mathbb{C}^{2}/G.$ The minimal \emph{(=} crepant\emph{)} resolution
\[
(\widehat{X},\text{\emph{Exc}}(f)\mathbf{)}\overset{f}{\mathbf{\longrightarrow
}}\text{ }\left(  X,\left[  \mathbf{0}\right]  \right)
\]
of the Gorenstein quotient singularity $\left(  X,\left[
\mathbf{0}\right] \right)  $ has exceptional set \emph{Exc}$(f)$
consisting of a configuration of rational smooth curves with
self-intersection number $-2$. The intersection form $\left(  \ ,\
\right)  :H_{2}(\widehat{X},\mathbb{Z})\times
H_{2}(\widehat{X},\mathbb{Z})\rightarrow\mathbb{Z}$ of
\emph{Exc}$(f)$ is negative definite, and therefore the dual
graphs \emph{DG}$(\emph{Exc}(f))$ of the irreducible components of
\emph{Exc}$(f)$ are exactly the Dynkin diagrams (of simply
connected complex Lie groups) of type $A$-$D$-$E$. \emph{(See
Tables \ref{Table1} and \ref{Table2}.)}
\end{theorem}

\begin{table}[h]
\epsfig{file=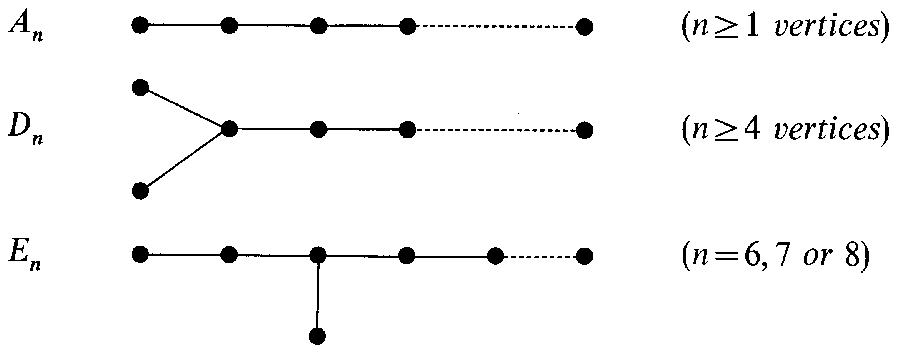, width=9.5cm, height=3.8cm}
\smallskip
\caption{}\label{Table2}
\end{table}

\noindent{}McKay \cite{McKay1, McKay2} established a remarkable connection
between the representation theory of the finite subgroups of SL$\left(
2,\mathbb{C}\right)  $ and the above Dynkin diagrams. This was the starting
point for Gonzalez-Sprinberg, Verdier \cite{GSV}, and Kn\"{o}rrer \cite{Kn},
to construct a purely geometric, direct correspondence
\begin{equation}
\fbox{$%
\begin{array}
[c]{ccc}
& \mathbf{McK}(G;f):\text{Irr}^{0}\left(  G\right)  \longrightarrow
\text{DG}(\text{Exc}\left(  f\right)  \smallskip) &
\end{array}
$} \label{McKay1}%
\end{equation}
\textquotedblleft of McKay-type\textquotedblright\ between the set
Irr$^{0}\left(  G\right)  $ of non-trivial irreducible representations of $G$
and DG$($Exc$\left(  f\right)  \smallskip)$ or, equivalently, between the
irreducible representations of $G$ and the members of the natural base of the
cohomology ring $H^{\bullet}(\widehat{X},\mathbb{Z)}$ (cf. \cite[\S 4]%
{Ito-Nak}). The bijection $\mathbf{McK}(G;f)$ induces an isomorphism between
Irr$^{0}\left(  G\right)  $ and the graph DG$($Exc$\left(  f\right)
\smallskip)$, i.e., the product of the images of two distinct elements of
Irr$^{0}\left(  G\right)  $ under $\mathbf{McK}(G;f)$ is mapped onto the
exceptional prime divisor corresponding to the \textquotedblleft
right\textquotedblright\ graph vertex. Brylinski \cite{Brylinski} constructed
subsequently a canonical \textquotedblleft dual\textquotedblright%
\ correspondence%
\begin{equation}
\fbox{$%
\begin{array}
[c]{ccc}
& \mathbf{McK}(G;f)^{\text{dual}}:\text{DG}(\text{Exc}\left(  f\right)
\smallskip)\longrightarrow\left\{
\begin{array}
[c]{c}%
\text{Non-trivial conjugacy }\\
\text{classes of }G
\end{array}
\right\}  &
\end{array}
$}\smallskip\label{McKay2}%
\end{equation}

\noindent{}relating the natural base of $H_{\bullet}(\widehat{X},\mathbb{Z)}$
to the set of conjugacy classes of $G.$ Finally, Ito \& Nakamura reinterpreted
(\ref{McKay1}) in terms of\thinspace\footnote{In general, for $r\geq2,$ by a
\textit{cluster} in $\mathbb{C}^{r}$ is meant a zero-dimensional subscheme
$Z\subset\mathbb{C}^{r},$ defined by an ideal $\mathcal{I}_{Z}\subset
\mathcal{O}_{\mathbb{C}^{r}},$ so that $\mathcal{O}_{Z}=\mathcal{O}%
_{\mathbb{C}^{r}}/\mathcal{I}_{Z}$ \ is a finite dimensional $\mathbb{C}%
$-vector space. If $G$ is a finite subgroup of SL$(r,\mathbb{C})$ of order
$l:=\left\vert G\right\vert ,$ then every $G$-invariant cluster ($G$%
-\textit{cluster}, for short) in $\mathbb{C}^{r}$ has global sections
$H^{0}(Z,\mathcal{O}_{Z})$ isomorphic (as $\mathbb{C}[G]$-module) to the
regular representation of $G.$ Consider the quasiprojective Hilbert scheme
Hilb$^{l}(\mathbb{C}^{r})$ parametrizing all clusters $Z$ of degree
dim$(\mathcal{O}_{Z})=l$ (cf. \cite[Lemma 5.1]{Reid-McK2}), the Hilbert-Chow
morphism Hilb$^{l}(\mathbb{C}^{r})\rightarrow$ Sym$^{l}(\mathbb{C}^{r}),$ as
well as the unique irreducible component $G$-Hilb$(\mathbb{C}^{r})$ of the
fixed locus $($Hilb$^{l}(\mathbb{C}^{r}))^{G}$ containing a general orbit of
$G$ on $\mathbb{C}^{r}.$ Then the Hilbert-Chow morphism induces a proper
birational morphism $G$-Hilb$(\mathbb{C}^{r})\rightarrow\mathbb{C}^{r}/G,$ and
$G$-Hilb$(\mathbb{C}^{r})$ parametrizes all $G$-clusters in $\mathbb{C}^{r}.$
(That's why it is shortly called \textit{the} \textit{Hilbert scheme of}
$G$-\textit{clusters.})} $G$-Hilb$(\mathbb{C}^{2})$:

\begin{theorem}
[{Hilbert scheme interpretation \cite[Thm. 10.4, p. 190]{Ito-Nak} }]%
$\widehat{X}$ is isomorphic to $G$-\emph{Hilb}$(\mathbb{C}^{2})$\emph{ }and
\emph{(\ref{McKay1}) }can be viewed as the bijection
\[
\text{\emph{Irr}}^{0}\left(  G\right)  \ni\rho\longmapsto D_{\rho}:=\left\{
\left.  I\in G\text{-}\emph{Hilb}(\mathbb{C}^{2})\right\vert V(I)\supset
V(\rho)\right\}  \in\text{\emph{DG}}(\text{\emph{Exc}}\left(  f\right)
\smallskip),
\]
where $V(\rho)$ is the $G$-module corresponding to $\rho,$ $\mathfrak{m}%
=\mathfrak{m}_{\mathbb{C}^{2}}$ \emph{(resp.,} $\mathfrak{m}_{X}$\emph{)} the
maximal ideal of $\mathbb{C}^{2}$ \emph{(resp.,} of $X$\emph{)} at the origin,
$\mathfrak{n}:=\mathfrak{m}_{X}\mathcal{O}_{\mathbb{C}^{2}},$ and
$V(I):=I/(\mathfrak{m}I+\mathfrak{n})$.
\end{theorem}

\begin{remark}
\label{2PROPERTIES}From the above mentioned results, the following attributes
of the quotient space $X=\mathbb{C}^{2}/G$ and of\emph{\ }the desingularizing
crepant morphism \emph{\ }$f:\widehat{X}\longrightarrow X$ are worth
recording:\smallskip\emph{\newline}(i) $X$ is always minimally embedded as a
\textit{hypersurface} in $\mathbb{C}^{3}$.\smallskip\newline(ii) The
desingularizing crepant morphism $f$ \textit{exists} for all groups $G$ of
Table \ref{Table1}, and is uniquely determined over $X$ up to \textit{isomorphism}%
.\smallskip\ \newline(iii) The singularity\emph{\ }$\left[
\mathbf{0}\right] \in X$ is \textit{isolated} (and consequently
all exceptional prime divisors w.r.t. $f$ are
compact).\newline(iv) $f$ is a \textit{projective} birational
morphism (i.e., $\widehat{X}$ is always a quasiprojective complex
variety), and can be decomposed into a finite sequence of blow-ups
of
\textit{points.\smallskip\ }\newline(v) $\widehat{X}\cong G$-Hilb$(\mathbb{C}%
^{2})$.\smallskip\medskip\newline\textit{None}\emph{\ }of (i)-(iv)
\textquotedblleft survives\textquotedblright\ \textit{in general} in higher
dimensions (and, as we explain below, $G$-Hilb$(\mathbb{C}^{r})$ is \textit{in
general} a good choice for an $\widehat{X}$ \textit{only for} $r=3$).
Nevertheless, passing to dimensions $r\geq3,$ it is useful to keep in mind
under which additional conditions the one or the other property (or a
reasonably \textit{weakened version} thereof) is preserved.
\end{remark}

\noindent{}$\bullet$ \textbf{McKay correspondence in dimension }$r=3.$ Based
on the classification table \cite{Yau-Yu} of the finite subgroups of
SL$(3,\mathbb{C}),$ Ito \cite{Ito1, Ito2}, Markushevich \cite{Mark-Comm,
Markushevich}, and Roan \cite{Roan} provided (by a case-by-case thorough
examination) a constructive proof of the following:

\begin{theorem}
[Existence Theorem in Dimension 3]\label{MARIR}All three-dimensional
Gorenstein quotient singularities possess crepant resolutions.
\end{theorem}

\noindent{}The resolution morphisms are unique only up to
\textquotedblleft isomorphism in codimension
$1$\textquotedblright\ (i.e., up to a finite number of canonical
\textit{flops}, cf. \cite[\S 6.4]{Kollar-Mori}), and to win
projectivity\footnote{According to a result of Kawamata \& Matsuki
\cite{Kawamata-Matsuki}, there is only a finite number of
projective crepant desingularizations of $\mathbb{C}^{3}/G.$} one
has to make particular choices (leading to \textit{smooth}
\textquotedblleft minimal models\textquotedblright). Furthermore,
for any such desingularization $f:\widehat{X}\longrightarrow
X=\mathbb{C}^{3}/G,$ $H^{\bullet}(\widehat{X},\mathbb{Z})$ is a
free $\mathbb{Z}$-module of rank equal to the number of conjugacy
classes of $G.$

\begin{theorem}
[{McKay Correspondence over $\mathbb{Q}$ in dimension $3;$ \cite[Prop.
5.6]{BD}, \cite[1.5-1.6]{Ito-Reid}}]For any crepant desingularization
$f:\widehat{X}\longrightarrow X=\mathbb{C}^{3}/G$ there are canonical
one-to-one correspondences%
\begin{equation}
\left\{  \text{\emph{conjugacy classes of }}G\text{ \emph{of age} }i\right\}
\longleftrightarrow\{\text{\emph{a basis of }}H^{2i}(\widehat{X}%
,\mathbb{Q})\}. \label{MCKAY3}%
\end{equation}

\end{theorem}

\noindent{}Besides, in analogy to the two-dimensional case, it turns out that,
among all possible projective crepant resolutions of $X,$ the Hilbert scheme
$G$-Hilb$(\mathbb{C}^{3})$ of $G$-clusters is a \textit{distinguished choice}.
(See Nakamura \cite{Nakamura}, and Craw \& Reid \cite{Craw-Reid} for the
abelian case, and \cite{BKR, GNS1, GNS2} for the non-abelian case.) There are
also several articles devoted to \textit{intrinsic interpretations} of
(\ref{MCKAY3}) for $\widehat{X}=G$-Hilb$(\mathbb{C}^{3})$ in terms of the
relevant ideals (see, e.g., \cite{CRAW-Ishii} for a GIT- and
\cite{Ito-Nakajima} for a $K$-theoretic description). More recently, Craw
\cite[Thm. 1.1]{Craw}, working with a natural base of the cohomology ring of
$\widehat{X}$ with \textit{integer} coefficients, succeeded in establishing an
explicit $3$-\textit{dimensional version of }(\ref{McKay1}) in the abelian
case.\medskip\

\noindent{}$\bullet$ \textbf{McKay correspondence in higher dimensions}$.$ In
lack of space we recall only a few highlights and refer to the survey articles
\cite{Fu3, Loeser, Reid-McK2} for further reading.

There are many obstructions in generalizing McKay correspondence in dimensions
$r\geq4,$ beginning with the Existence Problem\ (see comments below). But even
if one \textit{assumes} the existence of crepant desingularizations
$f:\widehat{X}\longrightarrow X$ of a given $X=\mathbb{C}^{r}/G,$ it is not
-as yet- clear if there might be a \textit{direct} analogue of (\ref{McKay1})
or (\ref{McKay2}) \textit{over} $\mathbb{Z}$ (cf. \cite[\S 1]{Reid-McK1}).
Most of the known results use homology and cohomology with coefficients taken
from $\mathbb{Q}$ or $\mathbb{C}.$

\begin{theorem}
[{Batyrev \cite[Thm 8.4]{Bat-Arch}, Denef \& Loeser \cite[Corollary
5.3]{Denef-Loeser}}]\label{BATYRTHM}If $G$ is a finite subgroup of
\emph{SL}$\left(  r,\mathbb{C}\right)  ,$ then for any crepant
desingularization $\widehat{X}\longrightarrow X$ of $X=\mathbb{C}^{r}/G$ we
have
\begin{equation}
\dim_{\mathbb{Q}}H^{2i}(\widehat{X},\mathbb{Q})=\sharp\left\{
\text{\emph{conjugacy classes of }}G\text{ \emph{having age} }i\right\}  ,
\label{BATFORMULA}%
\end{equation}
whereas the odd dimensional cohomology groups of $\widehat{X}$ are trivial. In
particular, the Euler number $\chi(\widehat{X})$ of $\widehat{X}$ equals the
number of the conjugacy classes of $G.$
\end{theorem}

\begin{theorem}
[{Ito-Reid Correspondence, \cite[Thm. 1.4]{Ito-Reid}}]\label{IRTHM}If $\left(
\mathbb{C}^{r}/G,\left[  \mathbf{0}\right]  \right)  $ is a Gorenstein
quotient singularity, then there is a canonical one-to-one correspondence
between the junior conjugacy classes in $G\left(  -1\right)  $ \emph{(or }$G,$
\emph{cf. \ref{NOTEHTAGE} (iii))} and the crepant discrete valuations of
$\mathbb{C}^{r}/G$.
\end{theorem}

\noindent{}Passing to Borel-Moore homology $H_{\bullet}^{\text{BM}}%
(\widehat{X},\mathbb{Q})$ (the dual of cohomology $H_{c}^{\bullet}(\widehat
{X},\mathbb{Q})$ with compact supports) for which the notion of fundamental
class of an algebraic cycle is well-defined, Theorem \ref{IRTHM} indicates
what the exact expectation for a \textit{high-dimensional McKay
correspondence} \textit{over} $\mathbb{Q}$ ought to be.

\begin{conjecture}
[{\cite{Reid-McK2}, \cite[2.8]{Kaledin1}}]\label{McKbig}If $G$ is a finite
subgroup of \emph{SL}$\left(  r,\mathbb{C}\right)  ,$ and $\widehat
{X}\longrightarrow X$ a crepant desingularization of $X=\mathbb{C}^{r}/G,$
then there is a canonical one-to-one correspondence%
\[
\fbox{$%
\begin{array}
[c]{ccc}
&  & \\
& \left\{
\begin{array}
[c]{c}%
\text{\emph{conjugacy classes of }}G\\
\text{ \emph{having age} }i
\end{array}
\right\}  \ni\lbrack g]\longleftrightarrow\text{\textsc{cl}}\left(
Z_{g}\right)  \in H_{2(r-i)}^{\text{\emph{BM}}}(\widehat{X},\mathbb{Q}), & \\
&  &
\end{array}
$}%
\]
mapping $[g]$ onto the fundamental class of the algebraic cycle $Z_{g},$ where
$Z_{g}$ denotes the Zariski closure of the center of the monomial valuation
\emph{(}of the function field of $X$\emph{)} corresponding to $g.$
\end{conjecture}

\begin{theorem}
[{Kaledin \cite[2.9]{Kaledin1}}]\label{KALEDINTHM}Conjecture
\emph{\ref{McKbig}} is true for symplectic $X$'s.
\end{theorem}

\noindent{}In fact, in the symplectic case, working with coefficients from
$\mathbb{C},$ it is also possible to confirm a \textit{multiplicative version}
of the high-dimensional McKay correspondence.

\begin{theorem}
[{Ginzburg \& Kaledin \cite[Thm. 1.2]{Ginz-Kal}, \cite[Thm. 2.4]{Kaledin3}}%
]\label{GKTHM}For symplectic $X$'s, there is a canonical graded algebra
isomorphism%
\begin{equation}
\fbox{$%
\begin{array}
[c]{ccc}
& H^{\bullet}(\widehat{X},\mathbb{C})\cong\mathbf{gr}^{\mathcal{F}_{\bullet}%
}(\mathcal{Z}G), &
\end{array}
$} \label{GKISO}%
\end{equation}
with $\mathbf{gr}^{\mathcal{F}_{\bullet}}(\mathcal{Z}G)$ as defined in
\emph{(\ref{ASSGRA})}.
\end{theorem}

\noindent{}For further information about symplectic quotient singularities and
their relation to McKay correspondence, the reader is referred to \cite{B-Kal,
Boissiere, Fu1, Fu2, Fu3, Fu4, Fu-Na, Kaledin2}.\medskip

\noindent{}{}$\bullet$ \textbf{The failure of} $G$-Hilb$(\mathbb{C}^{r})$
\textbf{in the role of an} $\widehat{X}$ \ \textbf{for }$r\geq4.$ Another
serious problem occurring in dimensions $r\geq4$ is that, in most of the
cases, and though it has plenty of nice properties, the Hilbert scheme of
$G$-clusters is no more suitable for our purposes. There are namely only a few
exceptional examples in which $G$-Hilb$(\mathbb{C}^{r})$ serves as a crepant
desingularization of $X=\mathbb{C}^{r}/G$, like the four-dimensional cyclic
quotient singularity of type $\frac{1}{15}(1,2,4,8)$ (cf. Note \ref{NOTEGPSS}
(ii)), the symplectic singularities $(\mathbb{C}^{r}/\mathfrak{S}%
_{r/2},\left[  \mathbf{0}\right]  ),$ $r\in2\mathbb{Z},$ (with the symmetric
group $\mathfrak{S}_{r/2}$ acting on $\mathbb{C}^{r}$ by permuting
coordinates) or $(\mathbb{C}^{r}/\Gamma\wr\mathfrak{S}_{r/2},\left[
\mathbf{0}\right]  ),$where $\Gamma$ denotes a finite subgroup of
SL$(2,\mathbb{C})$ and \textquotedblleft$\wr$\textquotedblright\ the wreath
product, and some others. (See Haiman \cite[\S 5]{Haiman}, Wang \cite{Wang},
and Kuznetsov \cite{Kuznetsov}.) In general, $G$-Hilb$(\mathbb{C}^{r})$
\textit{has} \textit{not necessarily} trivial canonical divisor, can be
\textit{singular}, or even \textit{non-normal} (see Note \ref{ROANREM} and
\cite[Corollary 1.6]{Craw-McG-Thomas}, respectively).\medskip

\noindent{}$\bullet$ \textbf{On the existence of crepant resolutions in
dimensions }$r\geq4.$ The presence of \textit{terminal}\footnote{For the
definition of terms like \textit{canonical} (resp., \textit{terminal})
\textit{singularities }(of \textit{index} $i\geq1$),\textit{ crepant divisor}
etc., see \cite{Matsuki, ReidC3F, ReidYPG}.} Gorenstein quotient singularities
in dimensions $r\geq4$ (for which there are no \textquotedblleft crepant
divisors\textquotedblright\ to pull out, cf. \cite{Mo-Ste}) means
automatically that, in contrast to what happens in dimension $2$ and $3,$
\textit{not all }Gorenstein quotient spaces $\mathbb{C}^{r}/G$ can be
desingularized by crepant birational morphisms.\smallskip

\noindent{}$\blacktriangleright$ \textbf{Existence Problem }(cf.
\cite[\S 4.5]{Ito-Reid}\ and \cite[\S 7]{Reid-McK2}): \textit{For
which} $G\subset$ SL$(r,\mathbb{C}),$ $r\geq4,$ \textit{do there
exist crepant (preferably projective) desingularizations of}
$\mathbb{C}^{r}/G?\smallskip$

\noindent{}Our first guess is that if we do not move too far away from the
hypersurface case (cf. Theorem \ref{INVARC2} and Remark \ref{2PROPERTIES}
(i)), then the existence of birational morphisms of this sort is indeed guaranteed.

\begin{conjecture}
[\cite{DHZ}]\label{CONJCIS}All quotient c.i.-singularities admit projective,
crepant resolutions in all dimensions.
\end{conjecture}

\noindent{}By (\ref{HTAGEREL}) and Theorem \ref{KWTHM} we see that for every
quotient c.i.-singularity $\left(  \mathbb{C}^{r}/G,\left[  \mathbf{0}\right]
\right)  $ the group $G$ \textit{is generated by its junior elements}. In view
of Theorem \ref{IRTHM}, we believe that this property is sufficient for the
existence of the desired desingularizations of $\mathbb{C}^{r}/G$. Conjecture
\ref{CONJCIS} is true whenever $G$ is abelian (see below Theorem
\ref{DHZTHM}). Moreover, the same assertion for all \textit{toric} (not
necessarily quotient) c.i.-singularities has been proven to be true in
\cite{DHaZ}.

Now going \textit{beyond} the \textquotedblleft c.i.'s\textquotedblright,\ and
putting the terminal ones aside, the remaining Gorenstein quotient
singularities have \textit{rarely }resolutions of this kind. Nonetheless, to
our surprise, the singularity series which do so, are \textit{not} negligible
as one would expect at first sight. (See below \S \ref{ONETWOPARSER} and
\S \ref{GPSERIES}.)\smallskip

Henceforth, we consider exclusively the Existence Problem for
Gorenstein \textit{\textbf{abelian}} quotient singularities. There
are several reasons to give priority to the abelian ones:

(i) Abelian finite subgroups $\ G$ of SL$(r,\mathbb{C})$ exist in
all dimensions $r,$ their conjugacy classes are singletons, and
their character groups are isomorphic to themselves. Hence,
letting them act linearly on $\mathbb{C}^{r},$ the age and the
height of any element $g\in G$ are determined by its
\textit{weights} appearing in the \textit{type }of $\left(
\mathbb{C}^{r}/G,\left[  \mathbf{0}\right]  \right)  $ (as long as
we fix eigencoordinates and generators; cf. \ref{deftype}).

(ii) For abelian $G$'s the Gorenstein quotient spaces
$\mathbb{C}^{r}/G$ can be treated by the toric machinery and,
particularly, by studying the properties of the so-called
\textit{junior simplices} $\mathfrak{s}_{G}$. Note that there is
no loss of generality when one works in the toric category because
the existence of an \textit{arbitrary} projective crepant
desingularization of $\mathbb{C}^{r}/G$ implies the existence of a
$\mathbb{T}_{N_{G}}$-\textit{equivariant} projective crepant
desingularization, where $\mathbb{T}_{N_{G}}:=\left(
\mathbb{C}^{\ast }\right)  ^{r}/G$ (cf. \cite[proof of Lemma
1]{Matsushita}). Moreover, the \textit{secondary polytope} of
$\mathfrak{s}_{G}$ describes conveniently the corresponding flops.

(iii) It is expected (cf. \cite[\S 4.6]{Ito-Reid}) that abelian
quotient singularities will be good candidates for proving both
Conjecture \ref{McKbig} and an analogue of Theorem
\ref{KALEDINTHM}, and for removing the restrictive hypothesis on
$\mathbb{C}^{r}/G$ (i.e., to be symplectic).

(iv) Given a \textit{non-abelian} group $G,$ it is also
conjectured (see \cite[\S 3]{Reid-McK0}) that the existence of
crepant desingularizations of $\mathbb{C}^{r}/G$ \ may be related
to the existence of such desingularizations for the quotients
$\mathbb{C}^{r}/H,$ for all \textit{maximal} cyclic (or abelian)
subgroups $H$ contained in $G.$

The present paper has been written trying to be self-contained and
partially expository. In particular, it includes more background
material than the average research paper has. The new results are
essentially in sections \ref{AQSTOR}, \ref{GPSERIES},
\ref{SECEXCRITERION} and \ref{ALGORITHM}, together with some parts
of \S\ref{CREPANTAQS} and of Appendices \ref{APPA} and
\ref{LATPJS}. (In \S\ref{CICASE}, \S\ref{FIRSTEXCR} , and
\S\ref{ONETWOPARSER} we summarize results from \cite{DHH,DHZ} and
from the unpublished manuscript \cite{DH}.)

More precisely, the paper is organized as follows: In \S
\ref{Prel} we recall fundamental notions from toric geometry and
introduce our notation. A detailed study of the abelian quotient
singularities as toric singularities (including various properties
of the junior and senior simplices of those which are Gorenstein)
is presented in \S \ref{AQSTOR}. In section \ref{CREPANTAQS} we
explain: (a) why the Existence Problem (in the abelian case) is
equivalent to the problem of finding junior simplices possessing
\textit{basic} (preferably \textit{coherent}) \textit{lattice
triangulations}, (b) how one can compute (\ref{BATFORMULA}) by
means of the Ehrhart polynomials of these simplices, and (c) why
two different maximal (partial or full) projective crepant
desingularizations of a Gorenstein abelian quotient space
$\mathbb{C}^{r}/G$ can be obtained from each other by a finite
succession of flops.

Wide classes of Gorenstein abelian quotient singularities
admitting projective, crepant resolutions are given in \S
\ref{CICASE}, \S \ref{ONETWOPARSER}, and \S \ref{GPSERIES}. (In
\S\ref{GPSERIES} we prove a long-standing conjecture concerning
the so-called \textit{GP-singularity series}, cf. \cite[\S
10]{DH}.)

On the other hand, to \textit{exclude} candidates for having
crepant resolutions (whenever the available lattice points in the
junior simplex are either \textquotedblleft strangely
located\textquotedblright\ or \textquotedblleft not
enough\textquotedblright\ to triangulate suitably), we apply two
\textit{necessary} existence conditions. The first of them (see
(\ref{HILBCON}) in \S \ref{FIRSTEXCR}) informs us that, provided
such a resolution is present, each of the Hilbert basis elements
of the cone supporting $\mathfrak{s}_{G}$ has to be either a
junior element or a vertex of $\mathfrak{s}_{G}.$\ The second one
(see (\ref{fircr2})-(\ref{criter2}) in \S \ref{SECEXCRITERION})
states that the existence of a crepant resolution implies the
boundedness of the acting group order from above by a number which
depends on the number of lattice points of $\mathfrak{s}_{G}$ and
of $\partial\mathfrak{s}_{G}$.

Next, combining our results, we outline an algorithm by means of
which it is possible to handle the Existence Problem (in the
abelian case); see \S \ref{ALGORITHM} (and especially Figure
\ref{Fig.7}). Only in the last two steps of this algorithm the
computational complexity grows rapidly. In particular, Step $5$
(involving the determination of maximal coherent triangulations of
$\mathfrak{s}_{G}$) is
added just for excluding some \textquotedblleft sporadic\textquotedblright%
\ counterexamples which happen to \textquotedblleft survive\textquotedblright%
\ after having used the above mentioned existence criteria.
(Computer programs like \texttt{Puntos} \cite{DELOERA2} or
\texttt{TOPCOM} \cite{RAMBAU} offer practical assistance in this
situation.)

Useful technical notions and results from the theory of subdivisions,
triangulations, lattice triangulations, and lattice point enumerators, are
presented separately in four appendices at the end of the paper. An
extensive\ part of Appendix \ref{APPA} is devoted to Upper Bound Theorems
(UBT's). Apart from the UBT \textit{for the facets of simplicial balls}
(Theorem \ref{SECUBB}), we conjecture the validity of a more effective UBT
\textit{for the facets of geometric simplex triangulations}, and give a proof
of it in dimension $3$ by means of PL-topological methods (see Theorem
\ref{BIGTH}). The \textit{coherence} of triangulations and certain
combinatorial properties of \textit{bistellar flips }belong to the topics
covered in Appendix \ref{APPCOH}. In Appendix \ref{APPB} we explain how one
passes from the coordinates of the $\mathbf{h}^{\ast}$-vector of a lattice
polytope $\ P$ (or, equivalently, from the coordinates of the $\mathbf{h}%
$-vector of any \textit{basic} triangulation $\mathcal{T}$ of $P$) to the
coefficients of its Ehrhart polynomial. Finally, in Appendix \ref{LATPJS} we
compute the coefficients of the Ehrhart polynomial of \textit{any junior
simplex} by making use of Mordell-Pommersheim and Diaz-Robins formulae.

\section{Toric Glossary\label{Prel}}

\noindent At first we recall some basic facts from the theory of toric
varieties. We mostly use the same notation as in \cite{DHH, DH, DHZ}. Our
standard references on toric geometry are the books of Oda \cite{Oda} and
Fulton \cite{Fulton}. \medskip

\noindent$\bullet$ \textbf{General notation}.\textsf{ }The \textit{linear
hull, }the\textit{\ affine hull}, the \textit{integral affine hull}, the
\textit{positive hull} and \textit{the convex hull} of a set $B$ of vectors of
$\mathbb{R}^{r}$, $r\geq1,$ will be denoted by lin$\left(  B\right)  $,
aff$\left(  B\right)  $, aff$_{\mathbb{Z}}\left(  B\right)  $, pos$\left(
B\right)  $ (or $\mathbb{R}_{\geq0}\,B$) and conv$\left(  B\right)  ,$
respectively. The \textit{dimension} dim$\left(  B\right)  $ of a
$B\subset\mathbb{R}^{r}$ is defined to be the dimension of its affine
hull.\medskip\ \newline$\bullet$ \textbf{Lattice determinants}.\textsf{ }Let
$N\cong\mathbb{Z}^{r}$ be a free $\mathbb{Z}$-module of rank $r\geq1$. $N$ can
be regarded as a \textit{lattice }in $N_{\mathbb{R}}:=N\otimes_{\mathbb{Z}%
}\mathbb{R}\cong\mathbb{R}^{r}$. (For fixed identification, we shall represent
the elements of $N_{\mathbb{R}}$ by column-vectors in $\mathbb{R}^{r}$). If
$\left\{  n_{1},\ldots,n_{r}\right\}  $ is a $\mathbb{Z}$-basis of $N$, then
\[
\text{det}\left(  N\right)  :=\left\vert \text{det}\left(  n_{1},\ldots
,n_{r}\right)  \right\vert
\]
is the \textit{lattice determinant}. An $n\in N$ is called \textit{primitive}
if conv$\left(  \left\{  \mathbf{0},n\right\}  \right)  \cap N$ contains no
other points except $\mathbf{0}$ and $n$.\smallskip

\noindent{}$\bullet$ \textbf{Cones}.\textsf{ }Let $N\cong\mathbb{Z}^{r}$ be as
above, $M:=$ Hom$_{\mathbb{Z}}\left(  N,\mathbb{Z}\right)  $ its dual lattice,
$N_{\mathbb{R}},M_{\mathbb{R}}$ their real scalar extensions, and
$\left\langle .,.\right\rangle :M_{\mathbb{R}}\times N_{\mathbb{R}}%
\rightarrow\mathbb{R}$ the natural $\mathbb{R}$-bilinear pairing. A subset
$\sigma$ of $N_{\mathbb{R}}$ is called \textit{strongly convex polyhedral cone
}(\textit{s.c.p. cone}, for short), if there exist $n_{1},\ldots,n_{k}\in
N_{\mathbb{R}}$, such that $\sigma=$ pos$\left(  \left\{  n_{1},\ldots
,n_{k}\right\}  \right)  $ and $\sigma\cap\left(  -\sigma\right)  =\left\{
\mathbf{0}\right\}  $. The \textit{dual cone} of such a $\sigma$ is
$\sigma^{\vee}:=\left\{  \mathbf{x}\in M_{\mathbb{R}}\ \left\vert
\ \left\langle \mathbf{x},\mathbf{y}\right\rangle \geq0,\ \forall
\mathbf{y},\ \mathbf{y}\in\sigma\right.  \right\}  .$ A subset $\tau$ of a
s.c.p. cone $\sigma$ is called a \textit{face} of $\sigma$ (notation:
$\tau\prec\sigma$), if $\tau=\left\{  \mathbf{y}\in\sigma\ \left\vert
\ \left\langle m_{0},\mathbf{y}\right\rangle =0\right.  \right\}  $, for some
$m_{0}\in\sigma^{\vee}$. A s.c.p. cone $\sigma=$ pos$\left(  \left\{
n_{1},\ldots,n_{k}\right\}  \right)  $ is called \textit{simplicial} (resp.,
\textit{rational}) if $n_{1},\ldots,n_{k}$ are $\mathbb{R}$-linearly
independent (resp., if $n_{1},\ldots,n_{k}\in N_{\mathbb{Q}}$, where
$N_{\mathbb{Q}}:=N\otimes_{\mathbb{Z}}\mathbb{Q}$).\medskip\newline$\bullet$
\textbf{Hilbert bases}.\textsf{ }If $\sigma\subset N_{\mathbb{R}}%
\cong\mathbb{R}^{r}$ is a rational s.c.p. cone, then $\sigma$ has $\mathbf{0}$
as its apex and the subsemigroup $\sigma\cap N$ of $N$ is a monoid.

\begin{proposition}
[Gordan's lemma]\label{Gorlem}$\sigma\cap N$ is finitely generated
as an additive semigroup, i.e. there exist
$n_{1},\ldots,n_{\nu}\in\sigma\cap N\ $\ such that\
\[
\ \sigma\cap N=\mathbb{Z}_{\geq0}\ n_{1}+\cdots+\mathbb{Z}_{\geq0}\ n_{\nu
}\ .
\]

\end{proposition}

\begin{proposition}
[{Minimal generating system, \cite[p. 233]{Schrijver}}]\label{MINGS}Among all
the systems of generators of $\sigma\cap N$, there is a unique system
$\mathbf{Hlb}_{N}\left(  \sigma\right)  $ of \emph{minimal cardinality},
namely\emph{:\smallskip}
\begin{equation}
\mathbf{Hlb}_{N}\left(  \sigma\right)  =\left\{  n\in\sigma\cap\left(
N\smallsetminus\left\{  \mathbf{0}\right\}  \right)  \ \left\vert \
\begin{array}
[c]{c}%
n\ \text{\emph{cannot be expressed as }}\\
\text{\emph{the sum of two other vectors }}\\
\text{\emph{belonging to\ } }\sigma\cap\left(  N\smallsetminus\left\{
\mathbf{0}\right\}  \right)
\end{array}
\right.  \right\}  . \label{Hilbbasis}%
\end{equation}

\end{proposition}

\begin{definition}
\emph{\ }$\mathbf{Hlb}_{N}\left(  \sigma\right)  $ is called\emph{\ }%
\textit{the Hilbert basis of }$\sigma$ w.r.t. $N.$
\end{definition}

\begin{theorem}
[Seb\"{o} \cite{SEBO}]\label{SEBOTHM}Given a rational s.c.p. cone
$\sigma\subset N_{\mathbb{R}}$ and an element $n\in\sigma\cap N,$ it is
co$\,\mathcal{NP}$-complete to decide whether $n$ is contained in
$\mathbf{Hlb}_{N}\left(  \sigma\right)  .$
\end{theorem}

\begin{remark}
\label{HBREM}Seb\"{o}'s Theorem shows the difficulty of deciding whether an
integral vector is additively reducible. In general, at least $r+\left\lfloor
\frac{r}{6}\right\rfloor $ elements of $\mathbf{Hlb}_{N}\left(  \sigma\right)
$ are needed to write an $n\in\sigma\cap N$ as non-negative integer linear
combination of elements of $\mathbf{Hlb}_{N}\left(  \sigma\right)  $ (see
\cite{BrunsETAL}). For an algorithm computing $\mathbf{Hlb}_{N}\left(
\sigma\right)  $ by the determination of the Graver basis of a suitable
integer matrix we refer to Sturmfels \cite[Algorithm 13.2, p. 128]{Sturmfels}.
Another efficient algorithm (which relies on a project-and-lift approach,
without making use of additional variables, and is implemented in the computer
program \texttt{MLP}) is due to Hemmecke \cite{Hemmecke}.
\end{remark}

\noindent$\bullet$ \textbf{Affine toric varieties}.\textsf{ }For a lattice
$N\cong\mathbb{Z}^{r}$ having $M$ as its dual, we define an $r$-dimensional
\textit{algebraic torus }$\mathbb{T}_{N}\cong\left(  \mathbb{C}^{\ast}\right)
^{r}$ by setting
\[
\mathbb{T}_{N}:=\text{Hom}_{\mathbb{Z}}\left(  M,\mathbb{C}^{\ast}\right)
=N\otimes_{\mathbb{Z}}\mathbb{C}^{\ast}.
\]
\ We identify $M$ with the character group of $\mathbb{T}_{N}$ and $N$ with
the group of $1$-parameter subgroups of $\mathbb{T}_{N}$. Let $\sigma$ be a
rational s.c.p. cone with
\[
M\cap\sigma^{\vee}=\mathbb{Z}_{\geq0}\ m_{1}+\mathbb{Z}_{\geq0}\ m_{2}%
+\cdots+\mathbb{Z}_{\geq0}\ m_{k}.
\]
To the finitely generated, normal, monoidal $\mathbb{C}$-subalgebra
$\mathbb{C}\left[  M\cap\sigma^{\vee}\right]  $ of $\mathbb{C}\left[
M\right]  $ we associate an affine complex variety
\[
U_{\sigma}:=\text{Spec}\left(  \mathbb{C}\left[  M\cap\sigma^{\vee}\right]
\right)
\]
endowed with a canonical $\mathbb{T}_{N}$-action. The analytic structure
induced on $U_{\sigma}$ is independent of the semigroup generators $\left\{
m_{1},\ldots,m_{k}\right\}  .$ Moreover, $\sharp\left(  \mathbf{Hlb}%
_{M}\left(  \sigma^{\vee}\right)  \right)  \ \left(  \leq k\right)  $ is
nothing but the \textit{embedding dimension} of $U_{\sigma}$, i.e. the minimal
number of generators of the maximal ideal of the local $\mathbb{C}$-algebra
$\mathcal{O}_{U_{\sigma},\mathbf{0}}$ (cf. \cite[1.2-1.3]{Oda}).\medskip

\noindent$\bullet$ \textbf{Fans}. A \textit{fan }w.r.t.\textit{\ }%
$N\cong\mathbb{Z}^{r}$ is a finite collection $\Delta$ of rational s.c.p.
cones in $N_{\mathbb{R}}$, such that the faces of any member belongs to
$\Delta$ and such that the intersection of any two members is a face of each
of them. We denote by $\left\vert \Delta\right\vert $ the support, and by
$\Delta\left(  i\right)  $ the set of $i$-dimensional cones of $\Delta.$ If
$\varrho$ is a \textit{ray} of $\Delta,$ i.e, if $\varrho\in\Delta\left(
1\right)  ,$ then there exists a unique primitive vector $n\left(
\varrho\right)  \in N\cap\varrho$ with $\varrho=\mathbb{R}_{\geq0}\ n\left(
\varrho\right)  $ and each cone $\sigma\in\Delta$ can be therefore written as
$\sigma=\sum_{\varrho\in\Delta\left(  1\right)  ,\,\varrho\prec\sigma
}\ \mathbb{R}_{\geq0}\ n\left(  \varrho\right)  .$ The set Gen$\left(
\sigma\right)  :=\left\{  n\left(  \varrho\right)  \ \left\vert \ \varrho
\in\Delta\left(  1\right)  ,\varrho\prec\sigma\right.  \right\}  $ is called
the\textit{\ set of minimal generators }of $\sigma$. For $\Delta$ itself one
defines analogously
\[
\text{Gen}\left(  \Delta\right)  :=\bigcup_{\sigma\in\Delta}\text{ Gen}\left(
\sigma\right)  .\medskip
\]
$\bullet$ \textbf{General toric varieties}. The \textit{toric variety
X}$\left(  N,\Delta\right)  $ associated to a fan\textit{\ }$\Delta$ w.r.t.
the lattice\textit{\ }$N$ is by definition the identification space $X\left(
N,\Delta\right)  :=\bigcup_{\sigma\in\Delta}\ U_{\sigma}\ /\ \sim,$ with
$U_{\sigma_{1}}\ni u_{1}\sim u_{2}\in U_{\sigma_{2}}$ if and only if there is
a $\tau\in\Delta$ such that $\tau\prec\sigma_{1}\cap\sigma_{2}$ and
$u_{1}=u_{2}$ within $U_{\tau}$. $X\left(  N,\Delta\right)  $ is called
\textit{simplicial }if all cones of $\Delta$ are simplicial. If $X\left(
N,\Delta\right)  $ is $r$-dimensional, then its topological Euler number
$\chi(X\left(  N,\Delta\right)  )$ equals
\begin{equation}
\chi(X\left(  N,\Delta\right)  )=\sharp\Delta\left(  r\right)  .\text{ (See
\cite[p. 59]{Fulton}.)}\label{EULERCHAR}%
\end{equation}
$X\left(  N,\Delta\right)  $ is also equipped with a canonical $\mathbb{T}%
_{N}$-action which is compatible with the above mentioned $\mathbb{T}_{N}%
$-actions on $U_{\sigma}$'s. The orbits w.r.t. this action are parametrized by
the set of all cones belonging to $\Delta$. For a $\tau\in\Delta$, we denote
by orb$\left(  \tau\right)  $ (resp., by $\mathbf{V}\left(  \tau\right)  $)
the orbit (resp., the closure of the orbit) which is associated to $\tau$. If
$\tau\in\Delta$, then $\mathbf{V}\left(  \tau\right)  =X\left(  N\left(
\tau\right)  ,\text{ Star}\left(  \tau;\Delta\right)  \right)  $ is itself a
toric variety w.r.t.
\[
N\left(  \tau\right)  :=N/N_{\tau}\ ,\ \ \text{Star}\left(  \tau
;\Delta\right)  :=\left\{  \overline{\sigma}\ \left\vert \ \sigma\in
\Delta,\ \tau\prec\sigma\right.  \right\}  ,
\]
where $N_{\tau}$ denotes the sublattice of $N$ generated (as subgroup) by the
intersection $N\cap$ lin$\left(  \mathbb{\tau}\right)  ,$ and $\overline
{\sigma}\ =\left(  \sigma+\left(  N_{\tau}\right)  _{\mathbb{R}}\right)
/\left(  N_{\tau}\right)  _{\mathbb{R}}$ is the image of $\sigma$ in $N\left(
\tau\right)  _{\mathbb{R}}=N_{\mathbb{R}}/\left(  N_{\tau}\right)
_{\mathbb{R}}$.\medskip

\noindent$\bullet$ \textbf{Equivariant maps}. A \textit{map of fans\ }%
$\varpi:\left(  N^{\prime},\Delta^{\prime}\right)  \rightarrow\left(
N,\Delta\right)  $ is a $\mathbb{Z}$-linear homomorphism $\varpi:N^{\prime
}\rightarrow N$ whose scalar extension $\varpi=\varpi_{\mathbb{R}%
}:N_{\mathbb{R}}^{\prime}\rightarrow N_{\mathbb{R}}$ satisfies the property:
$\forall\sigma^{\prime}\in\Delta^{\prime}$ $\exists\ \sigma\in\Delta\ \ $
with\ \ $\varpi\left(  \sigma^{\prime}\right)  \subset\sigma\,.$ Note that the
dual of the homomorphism $\varpi\otimes_{\mathbb{Z}}$id$_{\mathbb{C}^{\ast}%
}:\mathbb{T}_{N^{\prime}}=N^{\prime}\otimes_{\mathbb{Z}}\mathbb{C}^{\ast
}\rightarrow\mathbb{T}_{N}=N\otimes_{\mathbb{Z}}\mathbb{C}^{\ast}$ induces an
\textit{equivariant holomorphic map }$\varpi_{\ast}:X\left(  N^{\prime}%
,\Delta^{\prime}\right)  \rightarrow X\left(  N,\Delta\right)  $. This map
is\textit{\ proper} if and only if $\varpi^{-1}\left(  \left\vert
\Delta\right\vert \right)  =\left\vert \Delta^{\prime}\right\vert .$ In
particular, if $N=N^{\prime}$ and $\Delta^{\prime}$ is a refinement of
$\Delta$, then id$_{\ast}:X\left(  N,\Delta^{\prime}\right)  \rightarrow
X\left(  N,\Delta\right)  $ is \textit{proper}\emph{\ }and \textit{birational
}(cf. \cite[1.15 and 1.18]{Oda}).\medskip

\noindent$\bullet$ \textbf{Basic cones and desingularization}. Let
$N\cong\mathbb{Z}^{r}$ be a lattice of rank $r$ and $\sigma\subset
N_{\mathbb{R}}$ a simplicial, rational s.c.p. cone of dimension $k\leq r$.
$\sigma$ can be obviously written as $\sigma=\varrho_{1}+\cdots+\varrho_{k}$,
for distinct $1$-dimensional cones $\varrho_{1},\ldots,\varrho_{k}$. Let
\[
\mathbf{Par}\left(  \sigma\right)  :=\left\{  \mathbf{y}\in\left(  N_{\sigma
}\right)  _{\mathbb{R}}\ \left\vert \ \mathbf{y}=\sum_{j=1}^{k}\ \varepsilon
_{j}\ n\left(  \varrho_{j}\right)  ,\ \text{with\ \ }0\leq\varepsilon
_{j}<1,\ \forall j,\ 1\leq j\leq k\right.  \right\}
\]
be the \textit{fundamental }(\textit{half-open})\textit{\
parallelotope }associated to\textit{\ }$\sigma$. The
\textit{multiplicity} mult$\left( \sigma;N\right)  $ of $\sigma$
with respect to $N$ is defined to be $\text{mult}\left(
\sigma;N\right)  :=\sharp\left( \mathbf{Par}\left( \sigma\right)
\cap N_{\sigma}\right).$ As it turns out,
\begin{equation}
\text{mult}\left(  \sigma;N\right)  =\text{Vol}_{N_{\sigma}}\left(
\mathbf{Par}\left(  \sigma\right)  \right)  , \label{VOLFORMULA}%
\end{equation}
where Vol$\left(  \mathbf{Par}\left(  \sigma\right)  \right)  $ denotes the
usual volume (Lebesgue measure) of $\mathbf{Par}\left(  \sigma\right)  ,$ and
\[
\text{Vol}_{N_{\sigma}}\left(  \mathbf{Par}\left(  \sigma\right)  \right)
:=\tfrac{\text{Vol}\left(  \mathbf{Par}\left(  \sigma\right)  \right)
}{\text{det}\left(  N_{\sigma}\right)  }=\tfrac{\text{det}\left(
\mathbb{Z}\,n\left(  \varrho_{1}\right)  \oplus\cdots\oplus\mathbb{Z}%
\,n\left(  \varrho_{k}\right)  \right)  }{\text{det}\left(  N_{\sigma}\right)
}%
\]
its relative volume. If mult$\left(  \sigma;N\right)  =1$, then $\sigma$ is
called a \textit{basic cone} \textit{w.r.t.} $N$.

\begin{proposition}
[{\cite[Thm. 1.10 and Prop. 1.25]{Oda}}]\label{SMCR}The affine toric variety
$U_{\sigma}$ is $\mathbb{Q}$-factorial \emph{(}resp., smooth\emph{)} if and
only if $\sigma$ is simplicial \emph{(}resp., basic \textit{w.r.t.}
$N$\emph{).} Correspondingly, a toric variety $X\left(  N,\Delta\right)  $ is
$\mathbb{Q}$-factorial \emph{(}resp., smooth\emph{)} if and only if it is
simplicial (resp., simplicial and each s.c.p. cone $\sigma\in\Delta$ is
basic\emph{).}
\end{proposition}

\noindent By Carath\'{e}odory's theorem concerning convex polyhedral cones one
can choose a refinement $\Delta^{\prime}$ of any given fan $\Delta$, so that
$\Delta^{\prime}$ becomes simplicial. Since further subdivisions of
$\Delta^{\prime}$ reduce the multiplicities of its cones, we may arrive (after
finitely many subdivisions) at a fan $\widetilde{\Delta}$ having only basic cones.

\begin{theorem}
[Existence of Desingularizations]\label{THMTORDES}For every toric variety
$X\left(  N,\Delta\right)  $ there exists a refinement $\widetilde{\Delta}$ of
$\Delta$ consisting of exclusively basic cones w.r.t. $N$, i.e., such that
$f=\emph{id}_{\ast}:X(N,\widetilde{\Delta})\longrightarrow X\left(
N,\Delta\right)  $ is a desingularization.
\end{theorem}

\section{AQS as Toric Singularities\label{AQSTOR}}

\noindent{}Abelian quotient singularities (\textit{AQS}, for short) can be
directly investigated by means of the theory of toric varieties. If $G$ is a
finite abelian subgroup of GL$\left(  r,\mathbb{C}\right)  $, then $\left(
\mathbb{C}^{\ast}\right)  ^{r}/G$ is automatically an algebraic torus embedded
in $\mathbb{C}^{r}/G$.\medskip

\noindent$\bullet$ \textbf{General notation}. For $n\in\mathbb{N}$,
$m\in\mathbb{Z}$, we denote by $\left[  m\right]  _{n}$ the (uniquely
determined) integer for which $0\leq$ $\left[  m\right]  _{n}<n$, $m\equiv$
$\left[  m\right]  _{n}\left(  \text{mod }n\right)  $. If $x\in\mathbb{Q}$, we
define $\left\lceil x\right\rceil $ (resp., $\left\lfloor x\right\rfloor $) to
be the least integer number $\geq x$ (resp., the greatest integer number $\leq
x$). \textquotedblleft gcd\textquotedblright\ and \textquotedblleft
lcm\textquotedblright\ will be abbreviations for greatest common divisor and
least common multiple. Furthermore, for $n\in\mathbb{Z}_{\geq2}$, we denote by
$\zeta_{n}:=e^{\frac{2\pi\sqrt{-1}}{n}}$ the \textquotedblleft
first\textquotedblright\ $n$-th primitive root of unity.\medskip

\noindent$\bullet$ \textbf{The equivalence relation \textquotedblleft%
}$\backsim$\textbf{\textquotedblright}. \ For $\left(  q,r\right)  \in\left(
\mathbb{Z}_{\geq2}\right)  ^{2}$ we define\footnote{The symbol $\widehat
{\alpha_{i}}$ means here that $\alpha_{i}$ is omitted.}
\[
\Lambda\left(  q;r\right)  :=\left\{  \left(  \alpha_{1},..,\alpha_{r}\right)
\in\left\{  0,1,2,..,q-1\right\}  ^{r}\left\vert
\begin{array}
[c]{c}%
\ \text{gcd}(q,\alpha_{1},..,\widehat{\alpha_{i}},..,\alpha_{r})=1,\\
\text{for all }i,\ \ 1\leq i\leq r
\end{array}
\right.  \right\}
\]
and for $\left(  \left(  \alpha_{1},\ldots,\alpha_{r}\right)  ,\left(
\alpha_{1}^{\prime},\ldots,\alpha_{r}^{\prime}\right)  \right)  \in
\Lambda\left(  q;r\right)  \times\Lambda\left(  q;r\right)  $ the relation
\begin{equation}
\left(  \alpha_{1},\ldots,\alpha_{r}\right)  \backsim\left(  \alpha
_{1}^{\prime},\ldots,\alpha_{r}^{\prime}\right)  :\Longleftrightarrow\left\{
\begin{array}
[c]{l}%
\text{there exists a permutation }\\
\phi:\left\{  1,\ldots,r\right\}  \rightarrow\left\{  1,\ldots,r\right\}  \\
\text{and an integer }\lambda\text{, }1\leq\lambda\leq l-1\text{, }\\
\text{with gcd}\left(  \lambda,l\right)  =1\text{, such that}\\
\alpha_{\phi\left(  i\right)  }^{\prime}=\left[  \lambda\cdot\alpha
_{i}\right]  _{l}\text{, }\forall i,\ 1\leq i\leq r
\end{array}
\right\}  .\label{permute}%
\end{equation}
It is easy to see that $\backsim$ is an equivalence relation. We denote by
\[
\overline{\Lambda}\left(  q;r\right)  :=\Lambda\left(  q;r\right)  /\backsim
\]
the corresponding set of equivalent classes determined by \textquotedblleft%
$\backsim$\textquotedblright.\medskip\newline$\bullet$ \textbf{The
\textquotedblleft type\textquotedblright\ of an AQS}. \ Let $G$ be
a finite, small, abelian subgroup of GL$\left( r,\mathbb{C}\right)
$, $\ r\geq2$, with order $l=\left\vert G\right\vert \geq2$.
Consider as \textquotedblleft starting point\textquotedblright\ a
maximal decomposition of $G$ (viewed as an abstract group) into a
direct product of cyclic groups of orders, say,
$q_{1},\ldots,q_{\kappa\text{ }}$:
\begin{equation}
\fbox{$G\cong\left(  \mathbb{Z}/q_{1}\mathbb{Z}\right)  \times\cdots
\times\left(  \mathbb{Z}/q_{\kappa}\mathbb{Z}\right)  $}\label{finaldec}%
\end{equation}
and let
\[
\fbox{$\text{exp}\left(  G\right)  :=\text{ lcm}\left(  q_{1},\ldots
,q_{\kappa}\right)  $}%
\]
be the \textit{exponent} of $G$ and $\xi_{\mu}:=$ exp$\left(  G\right)  \cdot
q_{\mu}^{-1},$ for all $\mu\in\{1,\ldots,\kappa\}.$ Since $G$ is small, it is
easy to prove that $G$ possesses at most $r-1$ generators. Therefore we may
assume, from now on, that $\kappa\leq$ min$(r-1,\left\lfloor \frac{l}%
{2}\right\rfloor ).$ Choose after this fixing in advance of isomorphism
(\ref{finaldec}) suitable coordinates on $\mathbb{C}^{r}$ to diagonalize the
action of each factor $\left(  \mathbb{Z}/q_{\mu}\mathbb{Z}\right)  $ on
$\mathbb{C}^{r}$. According to Theorem \ref{prill-th}, and since every
representation of a finite cyclic group in a $\mathbb{C}$-vector space is the
direct sum of the one-dimensional representations, the action
\begin{equation}
\left(  \mathbb{Z}/q_{\mu}\mathbb{Z}\right)  \times\mathbb{C}^{r}\ni\left(
g_{\mu},\left(  z_{1},\ldots,z_{r}\right)  \right)  \mapsto\text{ diag}%
(\zeta_{q_{\mu}}^{\alpha_{\mu,1}}\ z_{1},\ldots,\zeta_{q_{\mu}}^{\alpha
_{\mu,r}}\ z_{r})\in\mathbb{C}^{r}\label{diagonal}%
\end{equation}
can be uniquely determined by the choice of a generator $g_{\mu}:=$
diag$\left(  \zeta_{q_{\mu}}^{\alpha_{\mu,1}},\ldots,\zeta_{q_{\mu}}%
^{\alpha_{\mu,r}}\right)  $ of each cyclic factor, i.e., by the choice of an
$r$-tuple $\left(  \alpha_{\mu,1},\ldots,\alpha_{\mu,r}\right)  \in
\Lambda\left(  q_{\mu};r\right)  $ as a representative of an equivalence class
within $\overline{\Lambda}\left(  q_{\mu};r\right)  $ w.r.t. $\backsim$ which
is defined by (\ref{permute}). (Any two $r$-tuples from $\Lambda\left(
q_{\mu};r\right)  $ belong to the same equivalence class w.r.t. $\backsim$ if
and only if the two associated representations of $\mathbb{Z}/q_{\mu
}\mathbb{Z}$ within GL$\left(  r,\mathbb{C}\right)  $, which correspond to
these two (probably different) generators of $\mathbb{Z}/q_{\mu}\mathbb{Z}$,
are conjugate to each other, i.e., if and only if the corresponding quotient
spaces are isomorphic; cf. \cite[Lemma 2, p. 296]{Fujiki}). Each element $g\in
G$, identified with the diagonalization of its image under (\ref{finaldec}),
is of the form
\begin{equation}
\fbox{$g=\text{diag}\left(  \zeta_{\text{exp}\left(  G\right)  }^{\delta
_{1}\left(  j_{1},\ldots,j_{\kappa}\right)  },\ldots,\zeta_{\text{exp}\left(
G\right)  }^{\delta_{r}\left(  j_{1},\ldots,j_{\kappa}\right)  }\right)  $%
}\label{exp-diag}%
\end{equation}
induced by a unique $\kappa$-tuple $\left(  j_{1},\ldots,j_{\kappa}\right)
\in\left\{  0,1,\ldots,q_{1}-1\right\}  \times\cdots\times\left\{
0,1,\ldots,q_{\kappa}-1\right\}  $, where
\begin{equation}
\delta_{i}\left(  j_{1},\ldots,j_{\kappa}\right)  :=\left[
{\textstyle\sum\limits_{\mu=1}^{\kappa}} \
j_{\mu}\cdot\xi_{\mu}\cdot\alpha_{\mu,i}\right]
_{\text{exp}\left(
G\right)  }\ ,\ \forall i\in\{1,\ldots,r\}.\label{deltas}%
\end{equation}

\begin{definition}
\label{deftype}For a $G$ having decomposition (\ref{finaldec}) and for a
given, predeterminated choice of the representation of the generators of its
factors within GL$\left(  r,\mathbb{C}\right)  $, as above in (\ref{diagonal}%
), (\ref{exp-diag}), (\ref{deltas}), we say that the AQS $\left(
\mathbb{C}^{r}/G,\left[  \mathbf{0}\right]  \right)  $ \textit{is of type}
\begin{equation}
\fbox{$%
\begin{array}
[c]{ccc}
&  & \\
& \frac{1}{q_{1}}\left(  \alpha_{1,1},\alpha_{1,2},\ldots,\alpha_{1,r}\right)
\times\ \cdots\ \cdots\ \times\frac{1}{q_{\kappa}}\left(  \alpha_{\kappa
,1},\alpha_{\kappa,2},\ldots,\alpha_{\kappa,r}\right)  & \\
&  &
\end{array}
$} \label{typeAQS}%
\end{equation}
(and view all the entries $\alpha_{i,j}$ as its \textit{weights}.) If
$G$\ happens to be cyclic, then we fix a choice of a generator of
$G\cong\mathbb{Z}/l\mathbb{Z}$ by making use of suitable exponents of
$\zeta_{l}$\ (i.e., by shortening of (\ref{finaldec}), in order to have just
one factor, and by suitable diagonalization), and we omit the first subscript
index of each of these regarded weights. In this case, we simply say that
$\left(  \mathbb{C}^{r}/G,\left[  \mathbf{0}\right]  \right)  $ is a cyclic
quotient singularity\textit{\ }(\textit{CQS}, for short)\textit{\ of type}%
\begin{equation}
\fbox{$\dfrac{1}{l}\left(  \alpha_{1},\alpha_{2},\ldots,\alpha_{r}\right)
\emph{.}$} \label{typeCQS}%
\end{equation}

\end{definition}

\noindent{}$\bullet$ \textbf{Abelian quotient spaces as toric varieties}. Let
$G$ be a finite, small, \textit{abelian} subgroup of GL$\left(  r,\mathbb{C}%
\right)  $, $r\geq2$, having order $l=\left\vert G\right\vert \geq2$, and let
\[
\left\{  e_{1}=\left(  1,0,\ldots,0,0\right)  ^{\intercal},\ldots
,e_{r}=\left(  0,0,\ldots,0,1\right)  ^{\intercal}\right\}
\]
denote the standard basis of $\mathbb{R}^{r}$, $N_{0}:=\mathbb{Z}^{r}%
=\sum_{i=1}^{r}\mathbb{Z}e_{i}$ the standard rectangular lattice, $M_{0}$ its
dual, and\smallskip\ $\mathbb{T}_{N_{0}}:=$ Spec$\left(  \mathbb{C}\left[
\mathfrak{x}_{1}^{\pm1},\ldots,\mathfrak{x}_{r}^{\pm1}\right]  \right)
=\left(  \mathbb{C}^{\ast}\right)  ^{r}.$ Clearly,
\[
\mathbb{T}_{N_{G}}:=\text{Spec}(\mathbb{C}\left[  \mathfrak{x}_{1}^{\pm
1},\ldots,\mathfrak{x}_{r}^{\pm1}\right]  ^{G})=\left(  \mathbb{C}^{\ast
}\right)  ^{r}/G
\]
is an $r$-dimensional algebraic torus with $N_{G}$ as its $1$-parameter group
and with $M_{G}$ as its group of characters. Using the map
\[
\left(  N_{0}\right)  _{\mathbb{R}}\ni\left(  y_{1},\ldots,y_{r}\right)
^{\intercal}=\mathbf{y\longmapsto}\text{ }\theta\left(  \mathbf{y}\right)
:=(e^{(2\pi\sqrt{-1})y_{1}},\ldots,e^{(2\pi\sqrt{-1})y_{r}})^{\intercal}%
\in\mathbb{T}_{N_{0}}%
\]
and the injection $\mathbf{\iota}:\mathbb{T}_{N_{0}}\hookrightarrow$
GL$\left(  r,\mathbb{C}\right)  $ defined by
\[
\mathbb{T}_{N_{0}}\ni\left(  t_{1},\ldots,t_{r}\right)  ^{\intercal
}\longmapsto\text{diag}\left(  t_{1},\ldots,t_{r}\right)  \in\text{ GL}\left(
r,\mathbb{C}\right)  ,
\]
we have obviously $N_{G}=\left(  \mathbf{\iota\circ}\text{ }\theta\right)
^{-1}\left(  G\right)  $ with det$\left(  N_{G}\right)  =\tfrac{1}{l}.$ In
fact, using (\ref{exp-diag}) and (\ref{deltas}), we get
\[
N_{G}=N_{0}+\sum_{\mu=1}^{\kappa}\,\mathbb{Z\,}(\frac{1}{q_{\mu}}\left(
\alpha_{\mu,1},\ldots,\alpha_{\mu,r}\right)  ^{\intercal})
\]%
\begin{equation}
=N_{0}+\sum_{\left(  j_{1},\ldots,j_{\kappa}\right)  \in\left\{
0,1,\ldots,q_{1}-1\right\}  \times\cdots\times\left\{  0,1,\ldots,q_{\kappa
}-1\right\}  }\,\mathbb{Z\,}\left(  \tfrac{\delta_{1}\left(  j_{1}%
,\ldots,j_{\kappa}\right)  }{\text{exp}\left(  G\right)  },\ldots
,\tfrac{\delta_{r}\left(  j_{1},\ldots,j_{\kappa}\right)  }{\text{exp}\left(
G\right)  }\right)  ^{\intercal} \label{ENG}%
\end{equation}
and
\[
M_{G}=\left\{  m\in M_{0}\,\left\vert
\begin{array}
[c]{c}%
\mathfrak{x}^{m}=\,\mathfrak{x}_{1}^{\mu_{1}}\,\cdots\,\mathfrak{x}_{r}%
^{\mu_{r}}\text{ \ \ is a }G\text{-invariant }\\
\text{Laurent monomial (}m=\left(  \mu_{1},\ldots,\mu_{r}\right)  \text{)}%
\end{array}
\right.  \right\}  \,\,\text{(with det}\left(  M_{G}\right)  =l\text{)}.
\]
\newline If we define
\[
\fbox{$%
\begin{array}
[c]{ccc}
& \sigma_{0}:=\ \text{pos}\left(  \left\{  e_{1},..,e_{r}\right\}  \right)  &
\end{array}
$}%
\]
to be the $r$-dimensional positive orthant, and $\Delta_{G}$ to be the
fan\smallskip%
\[
\Delta_{G}:=\left\{  \sigma_{0}\text{ together with its faces}\right\}  ,
\]
then by the exact sequence $0\rightarrow G\cong N_{G}/N_{0}\rightarrow
\mathbb{T}_{N_{0}}\rightarrow\mathbb{T}_{N_{G}}\rightarrow0$ induced by the
canonical duality pairing
\[
M_{0}/M_{G}\times N_{G}/N_{0}\rightarrow\mathbb{Q}/\mathbb{Z\hookrightarrow
C}^{\ast}\smallskip
\]
(cf. \cite[ p. 34]{Fulton} and \cite[pp. 22-23]{Oda}), we get $\mathbb{C}%
^{r}=X\left(  N_{0},\Delta_{G}\right)  \rightarrow X\left(  N_{G},\Delta
_{G}\right)  $ as projection map, where
\[
\fbox{$X\left(  N_{G},\Delta_{G}\right)  =U_{\sigma_{0}}{}=\mathbb{C}^{r}/G=$
Spec$\left(  \mathbb{C}\left[  \mathfrak{x}_{1},\ldots,\mathfrak{x}%
_{r}\right]  ^{G}\right)  \hookleftarrow\mathbb{T}_{N_{G}}$}%
\]
Formally, we identify $\left[  \mathbf{0}\right]  $ with orb$\left(
\sigma_{0}\right)  $. Moreover, the singular locus of $X\left(  N_{G}%
,\Delta_{G}\right)  $ can be written (by Propositions \ref{SLOC} and
\ref{SMCR}) as the union
\[
\text{Sing}\left(  X\left(  N_{G},\Delta_{G}\right)  \right)  =\text{
}\left\{  \text{orb}\left(  \sigma_{0}\right)  \right\}  \bigcup\left\{
\text{Spec}\left(  \mathbb{C}\left[  \overline{\sigma}_{0}^{\vee}\cap
M_{G}\left(  \tau\right)  \right]  \right)  \left\vert
\begin{array}
[c]{c}%
\,\tau\precneqq\sigma_{0}\text{,}\\
\text{dim}\left(  \tau\right)  \geq2,\text{and \thinspace}\\
\text{mult}\left(  \tau;N_{G}\right)  \geq2
\end{array}
\right.  \right\}  .
\]

\begin{proposition}
\label{isolat}For an AQS $\left(  \mathbb{C}^{r}/G,\left[  \mathbf{0}\right]
\right)  $ the following are equivalent\emph{:\smallskip}\newline\emph{(i)}
\emph{Sing}$\left(  X\left(  N_{G},\Delta_{G}\right)  \right)  =$ $\left\{
\text{\emph{orb}}\left(  \sigma_{0}\right)  \right\}  ,$ i.e., $\emph{orb}%
\left(  \sigma_{0}\right)  =\left[  \mathbf{0}\right]  \in X\left(
N_{G},\Delta_{G}\right)  $ is isolated.\smallskip\newline\emph{(ii)} For all
$\tau,\tau\precneqq\sigma_{0},$ with \emph{dim}$\left(  \tau\right)  \geq2$,
we have \emph{mult}$\left(  \tau\right)  =1.$
\end{proposition}

\begin{definition}
\label{MSCDEF}The \textit{splitting codimension} of\emph{\ }$U_{\sigma_{0}%
}=\mathbb{C}^{r}/G$\emph{\ }is defined to be the number
\[
\text{splcod}\left(  U_{\sigma_{0}}\right)  :=\text{min}\left\{
\mathfrak{k}\in\left\{  2,\ldots,r\right\}  \ \left\vert \
\begin{array}
[c]{c}%
U_{\sigma_{0}}\cong U_{\tau}\times\mathbb{C}^{r-\mathfrak{k}},\ \text{s.t.}\\
\tau\preceq\sigma_{0}\text{\emph{, }dim}\left(  \tau\right)  =\mathfrak{k}\\
\text{and Sing}\left(  U_{\tau}\right)  \neq\varnothing
\end{array}
\right.  \right\}
\]
If splcod$\left(  U_{\sigma_{0}}\right)  =r$,\ then orb$\left(  \sigma
_{0}\right)  $\emph{\ }is called an\emph{\ }\textit{msc-singularity}, i.e., a
singularity having the maximum splitting codimension.\smallskip
\end{definition}

\begin{proposition}
\ \label{non-msc}For an AQS $\left(  \mathbb{C}^{r}/G,\left[  \mathbf{0}%
\right]  \right)  $ of type \emph{(\ref{typeAQS})} the following conditions
are equivalent\emph{:\smallskip}\newline\emph{(i)} \emph{orb}$\left(
\sigma_{0}\right)  $ is a non-msc-singularity.\smallskip\newline\emph{(ii)}
\emph{splcod}$\left(  \emph{orb}\left(  \sigma_{0}\right)  ;U_{\sigma_{0}%
}\right)  =\mathfrak{k}_{0}$ with $2\leq\mathfrak{k}_{0}\leq r-1.\smallskip
$\newline\emph{(iii)} There exists a subfamily $\left\{  y_{1},y_{2}%
,\ldots,y_{r-\mathfrak{k}_{0}}\right\}  \subset\left\{  1,\ldots,r\right\}  $,
such that
\[
\delta_{y_{1}}\left(  j_{1},\ldots,j_{\kappa}\right)  =\cdots=\delta
_{y_{r-\mathfrak{k}_{0}}}\left(  j_{1},\ldots,j_{\kappa}\right)  =0\,,
\]
for all $\left(  j_{1},\ldots,j_{\kappa}\right)  \in\left\{  0,1,2,\ldots
,q_{1}-1\right\}  \times\cdots\times\left\{  0,1,2,\ldots,q_{\kappa
}-1\right\}  \,;\smallskip$\newline$\left\{  y_{1},y_{2},\ldots
,y_{r-\mathfrak{k}_{0}}\right\}  $ is, in addition, the largest subfamily of
$\left\{  1,\ldots,r\right\}  $ having this property.
\end{proposition}

\begin{theorem}
[Gorenstein condition]\label{Gor-prop} Let $\left(  \mathbb{C}^{r}/G,\left[
\mathbf{0}\right]  \right)  $ be an AQS of type \emph{(\ref{typeAQS}).} Then
the following conditions are equivalent:\smallskip\newline\emph{(i)} $X\left(
N_{G},\Delta_{G}\right)  =\mathbb{C}^{r}/G$ is Gorenstein.\smallskip
\newline$\smallskip$\emph{(ii)} $%
{\textstyle\sum\limits_{i=1}^{r}} \alpha_{\mu,i}\equiv\left(  0\
\text{\emph{mod} }q_{\mu}\right)  $, for all
$\mu,\,1\leq\mu\leq\kappa.$\newline$\smallskip$\emph{(iii)}
$\left\langle \left(  1,1,\ldots.1,1\right)  ,n\right\rangle
\geq1$, for all $n,\,n\in \sigma_{0}\cap\left(
N_{G}\smallsetminus\left\{  \mathbf{0}\right\}  \right)
.\smallskip$\newline\emph{(iv) }$\left(  X\left(
N_{G},\Delta_{G}\right) ,\text{\emph{orb}}\left(
\sigma_{0}\right)  \right)  $ is a canonical singularity of index
$1.$
\end{theorem}

\noindent{}\textsc{Proof}. It follows from Theorem \ref{Gor} and \cite[Thm.
3.1, p. 292]{ReidC3F}.\hfill{}$\square\medskip$

\noindent$\bullet$ \textbf{Junior and senior elements}. Let $\left(
\mathbb{C}^{r}/G,\left[  \mathbf{0}\right]  \right)  $ be a Gorenstein AQS of
type (\ref{typeAQS}). For all $i\in\{1,\ldots,r-1\}$\emph{, }we denote by
\[
\fbox{$\mathcal{H}_{i}:=\left\{  \left(  x_{1},\ldots,x_{r}\right)
^{\intercal}\in\mathbb{R}^{r}\,\left\vert \,x_{1}+x_{2}+\cdots+x_{r}=i\right.
\right\}  $}%
\]
the affine hyperplane \textit{of level} $i,$ and by
\[
\fbox{$\mathfrak{s}_{G}^{[i]}:=\sigma_{0}\cap\mathcal{H}_{i}=$ conv$\left(
\left\{  ie_{1},\ldots,ie_{r}\right\}  \right)  $}%
\]
the $\left(  r-1\right)  $-dimensional lattice simplex \textit{of age} $i$
which lies on $\mathcal{H}_{i}$. (We adopt here the terminology of
\cite[\S 1-2]{Ito-Reid}). In particular, \textit{the junior simplex }is
defined to be
\[
\fbox{$\mathfrak{s}_{G}:=\mathfrak{s}_{G}^{[1]}=\sigma_{0}\cap\mathcal{H}%
_{1}=$ conv$\left(  \left\{  e_{1},\ldots,e_{r}\right\}  \right)  .$}%
\]
An element $g\in G$ as in (\ref{exp-diag}) is a \textit{junior element}
(resp., a \textit{senior element of age }$i$, in the sense of \ref{NOTEHTAGE}
(ii)) whenever%
\[
n_{g}\,\,\,\text{belongs to \thinspace}\mathfrak{s}_{G}\cap N_{G}%
\,\,\,\ \text{(resp., to \thinspace\thinspace}\mathfrak{s}_{G}^{[i]}%
\,\cap\,N_{G},\text{ for }i\in\{2,\ldots,r-1\}\text{)\thinspace,}%
\]
where%
\begin{equation}
\fbox{$%
\begin{array}
[c]{lll}
& n_{g}:=\left(  \tfrac{\delta_{1}\left(  j_{1},\ldots,j_{\kappa}\right)
}{\text{exp}\left(  G\right)  },\ldots,\tfrac{\delta_{r}\left(  j_{1}%
,\ldots,j_{\kappa}\right)  }{\text{exp}\left(  G\right)  }\right)
^{\intercal}. &
\end{array}
$} \label{LPTREPRG}%
\end{equation}
$\mathbf{Par}\left(  \sigma_{0}\right)  $ is nothing but the unit half-open
cube in $\mathbb{R}^{r}$, with
\begin{equation}%
\begin{array}
[c]{l}%
\left\{  \text{lattice points }n_{g}\text{ representing the junior elements
}g\text{ of }G\right\} \\
=\left(  \mathfrak{s}_{G}\cap\mathbf{Par}\left(  \sigma_{0}\right)  \cap
N_{G}\right)  =\left(  \mathfrak{s}_{G}\cap\left(  N_{G}\smallsetminus\left\{
e_{1},\ldots,e_{r}\right\}  \right)  \right)
\end{array}
\label{JUN}%
\end{equation}
and
\begin{equation}%
\begin{array}
[c]{l}%
\left\{  \text{lattice points }n_{g}\text{ representing the seniors }g\in
G\text{ whose age is }i\right\} \\
=(\mathfrak{s}_{G}^{[i]}\cap\mathbf{Par}\left(  \sigma_{0}\right)  \cap
N_{G})\subset(\mathfrak{s}_{G}^{[i]}\cap\left(  N_{G}\smallsetminus\left\{
ie_{1},\ldots,ie_{r}\right\}  \right)  )
\end{array}
\label{SEN}%
\end{equation}
for all $i\in\{2,\ldots,r-1\}$ (cf. (\ref{ENG})). Obviously,
\begin{equation}%
{\textstyle\sum\limits_{i=1}^{r-1}}
(\sharp(\mathfrak{s}_{G}^{[i]}\cap\mathbf{Par}\left(
\sigma_{0}\right)  \cap
N_{G}))=l-1. \label{cardinalities}%
\end{equation}

\begin{figure}[h]
\epsfig{file=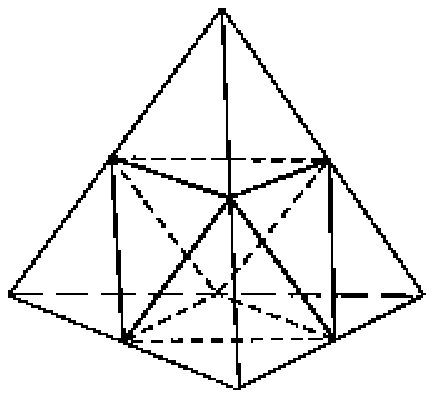}
\caption{}\label{Fig.1}
\end{figure}

\begin{note}
[Geometric interpetation via hypersimplices]If $i\in\{1,\ldots,r-1\},$ then
the $(r-1)$-dimensional polytope
\begin{equation}
\text{HypS}\left(  i,r\right)  :=\text{conv}\left(  \left\{  e_{\nu_{1}%
}+\cdots+e_{\nu_{i}}\,\left\vert \,1\leq\nu_{1}<\nu_{2}<\cdots<\nu_{i}\leq
r\right.  \right\}  \right)  \label{Hypers}%
\end{equation}%
\[
=\left\{  \left(  x_{1},\ldots,x_{r}\right)  ^{\intercal}\in\mathbb{R}%
^{r}\,\left\vert \,0\leq x_{j}\leq1\,\,\emph{\,}\text{for\thinspace
\thinspace\thinspace}1\leq j\leq r,\,\,\,\sum_{j=1}^{r}x_{j}=i\right.
\right\}
\]
has\emph{\ }$\binom{r}{i}$ vertices\emph{,} $2r$ facets for $i\in
\{2,\ldots,r-2\},$ and only $r$ facets for $i\in\{1,r-1\}$\emph{.
}It is the so-called\emph{\ }$\left(  i,r\right)
$-\textit{hypersimplex }and can be viewed as the convex hull of
the barycenters of the $\left(  i-1\right) $-dimensional faces of
the standard $\left(  r-1\right)  $-dimensional simplex (see
\cite[pp. 19-20]{Ziegler}). Figure \ref{Fig.1} shows Hyp$(2,4)$
(which is a regular octahedron)$.$ Note, in particular, that
HypS$\left( 1,r\right) $\emph{\ }and HypS$\left(  r-1,r\right)
$\emph{\ }are simplices\emph{.} The sets
$\mathfrak{s}_{G}^{[i]}\cap\mathbf{Par}\left(  \sigma_{0}\right) $
can be expressed in terms of\emph{\ }hypersimplices as follows:
\begin{equation}
\mathfrak{s}_{G}\cap\mathbf{Par}\left(  \sigma_{0}\right)  =\text{HypS}\left(
1,r\right)  \smallsetminus\left\{  e_{1},\ldots,e_{r}\right\}  =\mathfrak{s}%
_{G}\smallsetminus\left\{  e_{1},\ldots,e_{r}\right\}  , \label{Hyp1}%
\end{equation}
and for $i\in\{2,\ldots,r-1\},$ respectively,
\begin{equation}
\mathfrak{s}_{G}^{[i]}\cap\mathbf{Par}\left(  \sigma_{0}\right)
=\text{HypS}\left(  i,r\right)  \mathbb{r}\left\{  \mathbf{x}\in
\text{HypS}\left(  i,r\right)  \,\left\vert
\begin{array}
[c]{c}%
\emph{\,}\text{with\emph{\thinspace\thinspace}}x_{j}=1\\
\text{for all least one\emph{\ }}\\
j\in\{1,\ldots,r\}
\end{array}
\right.  \right\}  . \label{Hyp2}%
\end{equation}

\end{note}

\begin{remark}
The\emph{\ }height (\ref{defht}) of an element $g$ of $G$ (as in
(\ref{exp-diag})) equals\smallskip%
\begin{equation}
\fbox{ht$\left(  g\right)  =$ rank$\left(  g-\text{Id}_{G}\right)
=\sharp\left\{  i\,\left\vert \,1\leq i\leq r\text{ \thinspace with
\thinspace}\delta_{i}\left(  j_{1},\ldots,j_{\kappa}\right)  \neq0\right.
\right\}  $}\smallskip\label{height}%
\end{equation}
and specifies the dimension of the face of $\sigma_{0}$ on which $n_{g}\in
N_{G}$ lies.
\end{remark}

\begin{definition}
For $i\in\{1,\ldots,r-1\}$ and $k\in\{2,\ldots,r\}$ we define
\[
\fbox{$\mathfrak{B}_{G}\left(  i,k\right)  :=\left\{  g\in G\,\left\vert
\,\text{age}\left(  g\right)  =i,\,\text{ht}\left(  g\right)  =k\right.
\right\}  .$}%
\]

\end{definition}

\begin{lemma}
\label{smallK}$\mathfrak{B}_{G}\left(  i,k\right)  =\varnothing$ for $k\leq
i.$
\end{lemma}

\noindent\textsc{Proof}\textit{. }By Proposition \ref{ht-ag}, for each $g\in
G\smallsetminus\left\{  \text{Id}_{G}\right\}  $ we have ht$\left(  g\right)
>$ age$\left(  g\right)  $.\hfill{}$\square\medskip$\noindent\newline%
$G\smallsetminus\left\{  \text{Id}_{G}\right\}  $ can be decomposed as
follows:
\[
G\smallsetminus\left\{  \text{Id}_{G}\right\}  =\text{ JAI}\left(  G\right)
\overset{\bullet}{%
{\textstyle\bigcup} }\text{ SAI}\left(  G\right)  ,
\]
where%
\[
\text{JAI}\left(  G\right)  :=\,\overset{\bullet}{%
{\textstyle\bigcup}}\left\{  \mathfrak{B}_{G}\left(  i,k\right)
\,\left\vert \,1\leq i<k\leq r,\,k\neq2i\right.  \right\}  ,
\]
and%
\[
\text{SAI}\left(  G\right)  :=\ \overset{\bullet}{%
{\textstyle\bigcup}}\left\{  \mathfrak{B}_{G}\left(  i,2i\right)
\,\left\vert \,1\leq i \leq \left\lfloor \frac{r}{2}\right\rfloor
\right. \right\} \,.
\]
We call JAI$\left(  G\right)  $ (resp., SAI$\left(  G\right)  $) the set of
elements of $G\smallsetminus\left\{  \text{Id}_{G}\right\}  $ with inverses
\textit{of jumping age} (resp., \textit{of stationary age}), as we have:

\begin{lemma}
[\textquotedblleft Ping-Pong Lemma\textquotedblright]\label{pp-lem}\emph{(i)}
For $\,1\leq i<k\leq r,\,k\neq2i$, there is an one-to-one correspondence
\[
\text{\emph{JAI}}\left(  G\right)  \ni\mathfrak{B}_{G}\left(  i,k\right)  \ni
g\longmapsto g^{-1}\in\mathfrak{B}_{G}\left(  k-i,k\right)  \in
\text{\emph{JAI}}\left(  G\right)
\]
\emph{(ii)} If $k=2i$ and $g\in$ \emph{SAI}$\left(  G\right)  $, then
$g^{-1}\in$ \emph{SAI}$\left(  G\right)  .$
\end{lemma}

\noindent Note that JAI$\left(  G\right)  $ (resp. SAI$\left(  G\right)  $) is
expressed as the disjoint union of exactly $\left(  \binom r2-%
\genfrac{\lfloor}{\rfloor}{}{}{r}{2}%
\right)  $ ($\in2\,\mathbb{Z}$) sets (resp. of exactly $%
\genfrac{\lfloor}{\rfloor}{}{}{r}{2}%
$ sets). In addition, JAI$\left(  G\right)  $ consists of group elements the
cardinality of which is always an \textit{even }number.

\begin{definition}
To treat the group elements on each face of $\mathfrak{s}_{G}^{[i]}$'s
separately, for $1\leq i\leq r-1$, $2\leq k\leq r,$ and indices $1\leq\nu
_{1}<\nu_{2}<\cdots<\nu_{k}\leq r$\emph{, }we define
\[
\fbox{$%
\begin{array}
[c]{c}%
\mathfrak{s}_{G}^{[i]}\left(  \nu_{1},\nu_{2},\ldots,\nu_{k}\right)
:=\text{conv}\left(  \left\{  ie_{\nu_{1}},ie_{\nu_{2}},\ldots,ie_{\nu_{k}%
}\right\}  \right) \\
\\
\text{(}\mathfrak{s}_{G}^{[i]}=\mathfrak{s}_{G}^{[i]}\left(  1,2,\ldots
,r\right)  \text{\ )}%
\end{array}
$}%
\]
and
\[
\fbox{$%
\begin{array}
[c]{c}%
\mathfrak{B}_{G}\left(  i,k;\nu_{1},\nu_{2},\ldots,\nu_{k}\right)  :=\\
\\
\left\{  g\text{ as in (\ref{exp-diag})},\ g\in\mathfrak{B}_{G}\left(
i,k\right)  \,\left\vert \,n_{g}\in\mathfrak{s}_{G}^{[i]}\left(  \nu_{1}%
,\nu_{2},\ldots,\nu_{k}\right)  \cap\mathbf{Par}\left(  \sigma_{0}\right)
\cap N_{G}\right.  \right\}  .
\end{array}
$}%
\]

\end{definition}

\begin{lemma}
\label{interior}For any $g\in G$ the following conditions are
equivalent\smallskip\smallskip\emph{: }\newline\emph{(i) }$g\in\mathfrak{B}%
_{G}\left(  i,k;\nu_{1},\nu_{2},\ldots,\nu_{k}\right)  \emph{.\smallskip}%
$\newline\emph{(ii) }$n_{g}\in$ \emph{int}$(\mathfrak{s}_{G}^{[i]}\left(
\nu_{1},\nu_{2},\ldots,\nu_{k}\right)  )\cap\mathbf{Par}\left(  \sigma
_{0}\right)  \cap N_{G}.$
\end{lemma}

\noindent\textsc{Proof}\textit{. }The lattice points which belong to
\thinspace$\partial(\mathfrak{s}_{G}^{[i]}\left(  \nu_{1},\nu_{2},\ldots
,\nu_{k}\right)  )\cap\mathbf{Par}\left(  \sigma_{0}\right)  \cap N_{G}\,$
represent group elements with height $<k$.\hfill{}$\square$

\begin{lemma}
[\textquotedblleft Refined Ping-Pong Lemma\textquotedblright]\label{REF p-p}%
For $1\leq i<k\leq r,$ and indices $1\leq\nu_{1}<\nu_{2}<\cdots<\nu_{k}\leq
r$, there is an one-to-one correspondence
\[
\mathfrak{B}_{G}\left(  i,k;\nu_{1},\nu_{2},\ldots,\nu_{k}\right)  \ni
g\longmapsto g^{-1}\in\mathfrak{B}_{G}\left(  k-i,k;\nu_{1},\nu_{2},\ldots
,\nu_{k}\right)
\]

\end{lemma}

\noindent\textsc{Proof}\textit{. }As both $g$ and $g^{-1}$ (as diagonal
matrices) must have the entries which are equal to $1$ at the same positions,
the assertion can be varified by Lemma \ref{pp-lem} (i).\hfill{}$\square$

\begin{lemma}
\label{TWOI}For $1\leq i\leq r-1$ and $k=2i$ we have\emph{:}\newline\emph{(i)}
if $g\in\mathfrak{B}_{G}\left(  i,2i;\nu_{1},\nu_{2},\ldots,\nu_{2i}\right)
$, then $g^{-1}\in\mathfrak{B}_{G}\left(  i,2i;\nu_{1},\nu_{2},\ldots,\nu
_{2i}\right)  $, and\newline\emph{(ii)} if $l=\left\vert G\right\vert
\equiv1\left(  \text{\emph{mod }}2\right)  $, then $\mathfrak{B}_{G}\left(
i,2i;\nu_{1},\nu_{2},\ldots,\nu_{2i}\right)  $ consists of an even number of
group elements. Among them there are no elements $g$ with $g=g^{-1}.$
\end{lemma}

\noindent\textsc{Proof}\textit{.} For (i) we use the same argument as in the
proof of Lemma \ref{REF p-p}. Assertion (ii) can be shown easily as well. If
$\mathfrak{B}_{G}\left(  i,2i;\nu_{1},\nu_{2},\ldots,\nu_{2i}\right)  $ would
contain an element $g$ with $g=g^{-1}$, then this $g$ would have order $2$ and
therefore $l\equiv0\left(  \text{mod\emph{\ }}2\right)  $ by Lagrange Theorem.
This would contradict to our assumption.\hfill{}$\square$

\begin{proposition}
\label{isol2}For a Gorenstein AQS $\left(  \mathbb{C}^{r}/G,\left[
\mathbf{0}\right]  \right)  $ the following conditions are
equivalent\emph{:\smallskip}\newline\emph{(i)} $\emph{orb}\left(  \sigma
_{0}\right)  \in X\left(  N_{G},\Delta_{G}\right)  $ is isolated.\smallskip
\newline\emph{(ii)} $\bigcup_{i=1}^{r-1}(\partial\mathfrak{s}_{G}^{[i]}\cap
N_{G})=\varnothing$ \thinspace and $\,\bigcup_{i=1}^{r-1}($\emph{int}%
$(\mathfrak{s}_{G}^{[i]})\cap N_{G})\neq\varnothing.\smallskip$\newline%
\emph{(iii)} $\mathfrak{B}_{G}\left(  i,k;\nu_{1},\nu_{2},\ldots,\nu
_{k}\right)  =\varnothing$, for all $i\in\{1,\ldots,r-1\},k\in\{2,\ldots
,r-1\},$ and $1\leq\nu_{1}<\cdots<\nu_{k}\leq r,$ while $\mathfrak{B}%
_{G}\left(  i,r\right)  \neq\varnothing$ for at least one $i\in\{1,\ldots
,r-1\}.$
\end{proposition}

\noindent\textsc{Proof}\textit{. }It follows from Proposition \ref{isolat} and
Theorem \ref{Gor-prop}.\hfill{}$\square$

\begin{lemma}
\label{maxht}If there exists an element $g\in G$, such that \emph{ht}$\left(
g\right)  =r$, then \emph{orb}$\left(  \sigma_{0}\right)  $ is an
msc-singularity. For $G$ cyclic, the converse is also true.
\end{lemma}

\noindent{}\textsc{Proof}\textit{. }It follows from (\ref{height}) and the
definition of splitting codimension.\hfill{}$\square$

\begin{corollary}
\label{ISMSC}If $\emph{orb}\left(  \sigma_{0}\right)  \in X\left(
N_{G},\Delta_{G}\right)  $ is isolated, then \emph{orb}$\left(  \sigma
_{0}\right)  $ is an msc-singularity.
\end{corollary}

\begin{figure}[h]
\epsfig{file=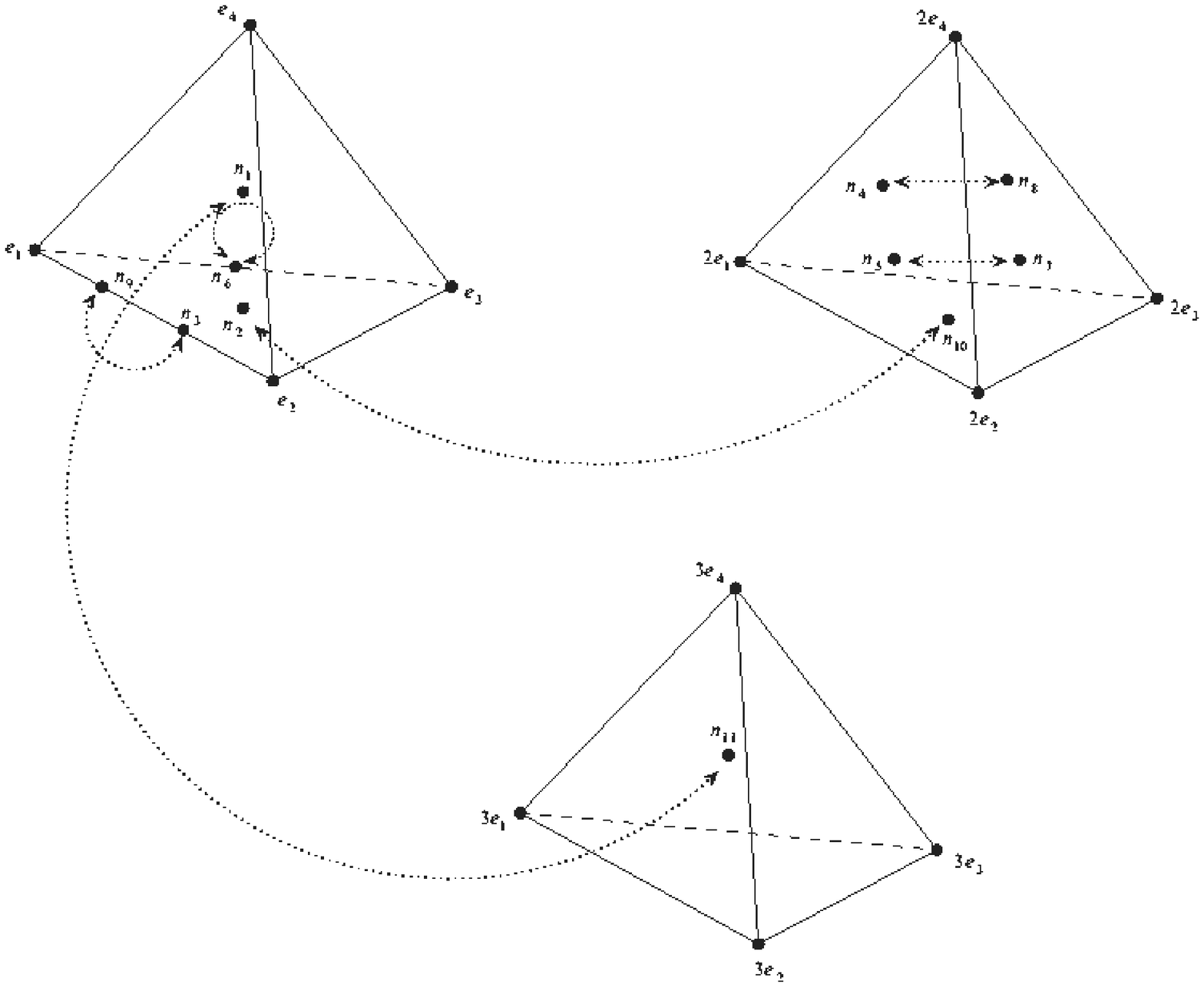, height=11cm, width=12.5cm}
\caption{}\label{Fig.2}
\end{figure}

\begin{example}
\label{Exam12}Let $\left(  \mathbb{C}^{4}/G,\left[  \mathbf{0}\right]
\right)  $ be the Gorenstein CQS\emph{\ }of type $\frac{1}{12}\left(
1,2,3,6\right)  $\emph{. }This is a non-isolated msc-singularity (by
Propositions \ref{isolat}, \ref{isol2} and Lemma \ref{maxht}). If\emph{\ }%
$G=\left\{  g_{0}=\text{Id}_{G},g_{1},g_{2},\ldots,g_{11}\right\}  $ and
$n_{i}:=n_{g_{i}}$ denotes the lattice point of\emph{\ }$N_{G}$
representing\emph{\ }$g_{i}$\emph{,} for $1\leq i\leq11,$\emph{\ }then, up to
reenumeration of indices, we find%
\smallskip
\setlength\extrarowheight{2pt}
\[%
\begin{tabular}
[c]{|c||c|c|}\hline
$i$ & $n_{i}$ & $n_{i+6}$\\\hline\hline
$0$ & $\underset{}{(0,0,0,0)^{\intercal}}$ & $\frac{1}{12}\left(
6,0,6,0\right)  ^{\intercal}$\\\hline
$1$ & $\underset{}{\frac{1}{12}\left(  1,2,3,6\right)  ^{\intercal}}$ &
$\frac{1}{12}\left(  7,2,9,6\right)  ^{\intercal}$\\\hline
$2$ & $\underset{}{\frac{1}{12}\left(  2,4,6,0\right)  ^{\intercal}}$ &
$\frac{1}{12}\left(  9,6,3,6\right)  ^{\intercal}$\\\hline
$3$ & $\underset{}{\frac{1}{12}\left(  4,8,0,0\right)  ^{\intercal}}$ &
$\frac{1}{12}\left(  8,4,0,0\right)  ^{\intercal}$\\\hline
$4$ & $\underset{}{\frac{1}{12}\left(  3,6,9,6\right)  ^{\intercal}}$ &
$\frac{1}{12}\left(  10,8,6,0\right)  ^{\intercal}$\\\hline
$5$ & $\underset{}{\frac{1}{12}\left(  5,10,3,6\right)  ^{\intercal}}$ &
$\frac{1}{12}\left(  11,10,9,6\right)  ^{\intercal}$\\\hline
\end{tabular}
\ \ \ \ \ \ \ \ \
\]
\setlength\extrarowheight{-2pt}
\smallskip
\emph{\ }Obviously,
$g_{1},g_{2},g_{3},g_{6},g_{9}$\emph{\ }are juniors,
$g_{4},g_{5},g_{7},g_{8},g_{10}$ are seniors of age $2$\emph{, }and $g_{11}%
$\emph{\ }is the only senior of age\emph{\ }$3$\emph{. }Furthermore,%
\[
\left\{
\begin{array}
[c]{l}%
\text{JAI}\left(  G\right)  =\mathfrak{B}_{G}\left(  1,3\right)
\overset{\bullet}{\,\bigcup\,}\mathfrak{B}_{G}\left(  1,4\right)
\overset{\bullet}{\,\bigcup\,}\mathfrak{B}_{G}\left(  2,3\right)
\overset{\bullet}{\,\bigcup\,}\mathfrak{B}_{G}\left(  3,4\right)  ,\\
\ \\
\text{ SAI}\left(  G\right)  =\mathfrak{B}_{G}\left(  1,2\right)
\overset{\bullet}{\,\bigcup\,}\mathfrak{B}_{G}\left(  2,4\right)  ,
\end{array}
\right.
\]
where
\begin{align*}
\mathfrak{B}_{G}\left(  1,3\right)   &  =\mathfrak{B}_{G}\left(
1,3;1,2,3\right)  =\left\{  g_{2}\right\}
,\,\,\mathfrak{B}_{G}\left(
1,4\right)  =\left\{  g_{1}\right\}  ,\,\,\\
\mathfrak{B}_{G}\left(  2,3\right)   &  =\mathfrak{B}_{G}\left(
2,3;1,2,3\right)  =\left\{  g_{10}\right\}
,\mathfrak{B}_{G}\left( 3,4\right)  =\left\{  g_{11}\right\}
,\,\,
\end{align*}%
\[
\mathfrak{B}_{G}\left(  1,2\right)  =\mathfrak{B}_{G}\left(
1,2;1,2\right)
\overset{\bullet}{\,\bigcup\,}\mathfrak{B}_{G}\left(
1,2;1,3\right) ,\,\,\mathfrak{B}_{G}\left(  1,2;1,2\right)
=\left\{  g_{3},g_{9}\right\}  ,
\]
and%
\[
\mathfrak{B}_{G}\left(  1,2;1,3\right)  =\left\{  g_{6}\right\}
,\,\,\,\mathfrak{B}_{G}\left(  2,4\right)  =\left\{  g_{4},g_{5},g_{7}%
,g_{8}\right\}  \,\,.
\]
Figure \ref{Fig.2} shows the location of the lattice points\emph{\ }$n_{1}%
,\ldots,n_{11}$ on the junior tetrahedron $\mathfrak{s}_{G}=\mathfrak{s}%
_{G}^{[1]}$\emph{\ }and on the two tetrahedra $\mathfrak{s}_{G}^{[2]}$ and
$\mathfrak{s}_{G}^{[3]}$\emph{, }containing the representatives of\emph{\ }the
junior and of the senior elements of ages\emph{\ }$2$\emph{\ }and $3$\emph{,
}respectively. (For aesthetic reasons, $\mathfrak{s}_{G}^{[2]}$ and
$\mathfrak{s}_{G}^{[3]}$ are scaled by $\frac{1}{2}$ and $\frac{1}{3},$
respectively.) The dotted lines (with arrows at their ends) indicate how the
refined Ping-Pong Lemma \ref{REF p-p} and Lemma \ref{TWOI} (i) are applied in
practice. Note that\emph{\ }$g_{12-i}=g_{i}^{-1}$\emph{, }for all $i$\emph{,}
$1\leq i\leq6.$
\end{example}

\noindent{}$\bullet$ \textbf{On the cardinality of} $\mathfrak{B}_{G}\left(
1,k\right)  .$ Let $\left(  \mathbb{C}^{r}/G,\left[  \mathbf{0}\right]
\right)  $ be a Gorenstein AQS of type (\ref{typeAQS}). The \textit{counting
}of the lattice points of $N_{G}$ representing all the elements of
$\mathfrak{B}_{G}\left(  i,k;\nu_{1},\nu_{2},\ldots,\nu_{k}\right)  $'s and
$\mathfrak{B}_{G}\left(  i,k\right)  $'s (with emphasis on $\mathfrak{B}%
_{G}\left(  1,k\right)  $) is of particular interest (cf.
\S \ref{SECEXCRITERION}). By Lemma \ref{interior} we get
\begin{equation}
\sharp\left(  \mathfrak{B}_{G}\left(  i,k;\nu_{1},\nu_{2},\ldots,\nu
_{k}\right)  \right)  =\sharp(\text{int}(\mathfrak{s}_{G}^{[i]}\left(  \nu
_{1},\nu_{2},\ldots,\nu_{k}\right)  \cap\mathbf{Par}\left(  \sigma_{0}\right)
\cap N_{G})), \label{N-BGIK}%
\end{equation}
and for $i=1,$%
\begin{equation}
\sharp\left(  \mathfrak{B}_{G}\left(  1,k;\nu_{1},\nu_{2},\ldots,\nu
_{k}\right)  \right)  \,\overset{\text{(\ref{JUN})}}{=}\,\sharp(\text{int}%
(\mathfrak{s}_{G}\left(  \nu_{1},\nu_{2},\ldots,\nu_{k}\right)  \cap N_{G})).
\label{N-BG1K}%
\end{equation}
Moreover,
\begin{align}
\sharp\left(  \mathfrak{B}_{G}\left(  i,k\right)  \right)   &  =\sum_{1\leq
\nu_{1}<\cdots<\nu_{k}\leq r\,}\,\sharp\left(  \mathfrak{B}_{G}\left(
i,k;\nu_{1},\nu_{2},\ldots,\nu_{k}\right)  \right) \label{BGIK}\\
& \nonumber\\
&  =\sum_{1\leq\nu_{1}<\cdots<\nu_{k}\leq r\,}\,\sharp(\text{int}%
(\mathfrak{s}_{G}^{[i]}\left(  \nu_{1},\nu_{2},\ldots,\nu_{k}\right)
\cap\mathbf{Par}\left(  \sigma_{0}\right)  \cap N_{G}))\nonumber
\end{align}
and for $i=1,$%
\begin{equation}
\fbox{$%
\begin{array}
[c]{c}%
\sharp\left(  \mathfrak{B}_{G}\left(  1,k\right)  \right)  =%
{\displaystyle\sum\limits_{1\leq\nu_{1}<\cdots<\nu_{k}\leq r\,}}
\,\sharp\left(  \mathfrak{B}_{G}\left(  1,k;\nu_{1},\nu_{2},\ldots,\nu
_{k}\right)  \right) \\
\\
=%
{\displaystyle\sum\limits_{1\leq\nu_{1}<\cdots<\nu_{k}\leq r\,}}
\,\,\sharp(\text{int}(\mathfrak{s}_{G}\left(
\nu_{1},\nu_{2},\ldots,\nu _{k}\right)  \cap N_{G}))
\end{array}
$} \label{BG1K}%
\end{equation}
The cardinality (\ref{N-BGIK}) can be written as follows:
\begin{equation}%
\begin{array}
[c]{c}%
\sharp\left(  \mathfrak{B}_{G}\left(  i,k;\nu_{1},\nu_{2},\ldots,\nu
_{k}\right)  \right)  =\\
\,\\
=\sharp\left\{  \lambda\in\lbrack0,l)\cap\mathbb{Z}\left\vert
\begin{array}
[c]{c}%
{\textstyle\sum\limits_{\rho=1}^{r}}
\ \delta_{\rho}\left(  j_{1},\ldots,j_{\kappa}\right)  =i\cdot\text{exp}%
\left(  G\right)  \text{ and}\\
\delta_{\rho}\left(  j_{1},\ldots,j_{\kappa}\right)  \ \left\{
\begin{array}
[c]{ll}%
\neq0, & \text{if }\rho\in\left\{  \nu_{1},\ldots,\nu_{k}\right\} \\
=0, & \text{if }\rho\notin\left\{  \nu_{1},\ldots,\nu_{k}\right\}
\end{array}
\right.  \ \\
\text{for all\ \ \ }\rho\in\left\{  1,2,\ldots,r\right\}
\end{array}
\right.  \right\}  .
\end{array}
\label{AB-FORM}%
\end{equation}

\begin{proposition}
Let $r\geq3$, $2\leq k\leq r-1$, $\left\{  \nu_{1},\nu_{2},\ldots,\nu
_{k}\right\}  \subset\left\{  1,2,\ldots,r\right\}  $ be a family of indices,
such that $1\leq\nu_{1}<\cdots<\nu_{k}\leq r$, and $\left\{  \nu_{1}^{\prime
},\nu_{2}^{\prime},\ldots,\nu_{r-k}^{\prime}\right\}  =\left\{  1,2,\ldots
,r\right\}  \smallsetminus\left\{  \nu_{1},\nu_{2},\ldots,\nu_{k}\right\}  $
be its complement. If $\left(  \mathbb{C}^{r}/G,\left[  \mathbf{0}\right]
\right)  $ is Gorenstein msc-CQS\emph{\ }of type \emph{(\ref{typeCQS}), }then
the number of group elements having fixed height $k$ and arbitrary age equals
\begin{align*}
\sum_{i=1}^{r-1}\sharp\left(  \mathfrak{B}_{G}\left(  i,k;\nu_{1},\nu
_{2},\ldots,\nu_{k}\right)  \right)   &  =\sum_{i=1}^{k-1}\sharp\left(
\mathfrak{B}_{G}\left(  i,k;\nu_{1},\nu_{2},\ldots,\nu_{k}\right)  \right) \\
&  =\left[  \text{\emph{gcd}}\left(  \alpha_{\nu_{1}^{\prime}},\ldots
,\alpha_{\nu_{r-k}^{\prime}},l\right)  -1\right]
\end{align*}%
\begin{equation}
-\sum_{\substack{1\leq\rho_{1}<\cdots<\rho_{y}\leq k\\\text{\emph{with \ }%
}1\leq y\leq k-2\ \text{\emph{and}}\\\left\{  \rho_{1},\ldots,\rho
_{y}\right\}  \subset\left\{  \nu_{1},\nu_{2},\ldots,\nu_{k}\right\}
}}\ \left[  \text{\emph{gcd}}\left(  \alpha_{\nu_{1}^{\prime}},..,\alpha
_{\nu_{r-k}^{\prime}},\alpha_{\rho_{1}},..,\alpha_{\rho_{y}},l\right)
-1\right]  \label{GCD-FORM}%
\end{equation}
\emph{(For} $k=2$ \emph{we simply omit this last sum.)}
\end{proposition}

\noindent\textsc{Proof}\textit{. }The first equality is clear by Lemma
\ref{smallK}. On the other hand,%
\[%
\begin{array}
[c]{l}%
{\textstyle\sum\limits_{i=1}^{r-1}}
\sharp(\mathfrak{s}_{G}^{[i]}\left(
\nu_{1},\nu_{2},\ldots,\nu_{k}\right)
\cap\mathbf{Par}\left(  \sigma_{0}\right)  \cap N_{G})\smallskip\\
=%
{\textstyle\sum\limits_{i=1}^{k-1}}
\sharp(\mathfrak{s}_{G}^{[i]}\left(
\nu_{1},\nu_{2},\ldots,\nu_{k}\right)
\cap\mathbf{Par}\left(  \sigma_{0}\right)  \cap N_{G})\smallskip\\
=\sharp\left\{  \lambda\in\mathbb{Z\ }\left\vert \ 1\leq\lambda\leq
l-1:\lambda\alpha_{\nu_{i}^{\prime}}\equiv0\left(  \text{mod }l\right)
,\ \forall i,\ 1\leq i\leq r-k\right.  \right\}  \smallskip\\
=\text{gcd}(\alpha_{\nu_{1}^{\prime}},\ldots,\alpha_{\nu_{r-k}^{\prime}},l)-1.
\end{array}
\]
The sum we subtract in (\ref{GCD-FORM}) is nothing but the evaluation of the
number
\[%
{\textstyle\sum\limits_{i=1}^{k-1}} \
\sharp\,(\partial\mathfrak{s}_{G}^{[i]}\left(
\nu_{1},\nu_{2},\ldots ,\nu_{k}\right)  \cap\mathbf{Par}\left(
\sigma_{0}\right)  \cap N_{G})
\]
of lattice points lying on the corresponding relative boundaries.\hfill
{}$\square$

\begin{remark}
\label{pprem}For arbitrary $r\geq3$\emph{, }the numbers $\sharp\,\left(
\mathfrak{B}_{G}\left(  1,2;\nu_{1},\nu_{2}\right)  \right)  $ and
\[
\sharp\,\left(  \mathfrak{B}_{G}\left(  1,3;\nu_{1},\nu_{2},\nu_{3}\right)
\right)  =\sharp\,\left(  \mathfrak{B}_{G}\left(  2,3;\nu_{1},\nu_{2},\nu
_{3}\right)  \right)
\]
can be determined by the formula (\ref{GCD-FORM}); the first of them directly
(because in this case the left-hand side of (\ref{GCD-FORM}), consists of only
one summand), and the latter just as the half of what we get from the
right-hand side of\emph{\ }(\ref{GCD-FORM}). For\emph{\ }$r\in\left\{
3,4\right\}  $\emph{, }as a byproduct of this formula and of the refined
Ping-Pong Lemma \ref{REF p-p}, one gets a simple method to count the number of
the lattice points lying in the relative interior of each proper face of the
junior simplex separately!
\end{remark}

\begin{note}
One can analogously compute $%
{\textstyle\sum\nolimits_{i=1}^{r-1}} \ \sharp\left(
\mathfrak{B}_{G}\left(  i,k;\nu_{1},\nu_{2},\ldots,\nu _{k}\right)
\right)  $ whenever the acting group $G$ is abelian but not
cyclic. In this case, the formula generalizing (\ref{GCD-FORM})
contains denumerants of restricted weighted vectorial partitions
instead of greatest common divisors.
\end{note}

\section{Crepant Resolutions of Gorenstein AQS\label{CREPANTAQS}}

\noindent{}Let $\left(  \mathbb{C}^{r}/G,\left[  \mathbf{0}\right]
\right) ,$ $r\geq3,$ be a Gorenstein AQS. In this section we
explain how the crepant (partial or full)
$\mathbb{T}_{N_{G}}$-equivariant desingularizations of its
underlying space are to be studied by means of \textit{lattice
triangulations} of the junior simplex $\mathfrak{s}_{G}.$

\begin{definition}
We denote the set of all lattice triangulations $\mathcal{T}$\emph{\ }of
$\mathfrak{s}_{G}$ w.r.t. $N_{G}$ (with vert$(\mathcal{T})\subseteq
\mathfrak{s}_{G}\cap N_{G},$ cf. \ref{LATTR}) by $\mathbf{LTR}_{N_{G}%
}\,\left(  \mathfrak{s}_{G}\right)  ,$ and define
\[%
\begin{array}
[c]{l}%
\mathbf{LTR}_{N_{G}}^{\text{{\scriptsize max}}}\,\left(  \mathfrak{s}%
_{G}\right)  :=\left\{  \mathcal{T}\in\mathbf{LTR}_{N_{G}}\left(
\mathfrak{s}_{G}\right)  \left\vert
\begin{array}
[c]{c}%
\mathcal{T}\text{\emph{\ \ }is a \textit{maximal triangulation}}\\
\text{ of\emph{\ }}\mathfrak{s}_{G},\text{ i.e., vert}(\mathcal{T}%
)=\mathfrak{s}_{G}\cap N_{G}%
\end{array}
\right.  \right\}  \smallskip,\smallskip\\
\mathbf{LTR}_{N_{G}}^{\text{{\scriptsize basic}}}\left(  \mathfrak{s}%
_{G}\right)  :=\left\{  \mathcal{T}\in\mathbf{LTR}_{N_{G}}%
^{\text{{\scriptsize max}}}\,\left(  \mathfrak{s}_{G}\right)  \ \left\vert
\begin{array}
[c]{c}%
\ \mathcal{T}\text{\emph{\ \ }is a \textit{basic} triangulation }\\
\text{of\emph{ \ }}\mathfrak{s}_{G}\text{ (see Definition \ref{BASICTR})}%
\end{array}
\right.  \right\}  .
\end{array}
\]
Adding the prefix $\mathbf{Coh}$- to any one of the above sets, we
mean the subset consisting of coherent triangulations (in the
sense of \ref{CohTriang}).\emph{ }The\emph{ }\textit{hierarchy of
lattice triangulations }of $\mathfrak{s}_{G}$\ is given by the
following inclusion
diagram:%
\[%
\begin{array}
[c]{ccccc}%
\mathbf{LTR}_{N_{G}}^{\text{{\scriptsize basic}}}\left(  \mathfrak{s}%
_{G}\right)  & \subset & \mathbf{LTR}_{N_{G}}^{\text{{\scriptsize max}}%
}\,\left(  \mathfrak{s}_{G}\right)  & \subset & \mathbf{LTR}_{N_{G}}\left(
\mathfrak{s}_{G}\right)  \smallskip\\
\bigcup &  & \bigcup &  & \bigcup\smallskip\\
\mathbf{Coh}\text{-}\mathbf{LTR}_{N_{G}}^{\text{{\scriptsize basic}}}\left(
\mathfrak{s}_{G}\right)  & \subset & \mathbf{Coh}\text{-}\mathbf{LTR}_{N_{G}%
}^{\text{{\scriptsize max}}}\,\left(  \mathfrak{s}_{G}\right)  & \subset &
\mathbf{Coh}\text{-}\mathbf{LTR}_{N_{G}}\left(  \mathfrak{s}_{G}\right)
\end{array}
\]

\end{definition}

\begin{note}
(i) There is a bijection between the triangulations belonging to
$\mathbf{LTR}_{N_{G}}\left(  \mathfrak{s}_{G}\right)  $ (resp., to
$\mathbf{Coh}$-$\mathbf{LTR}_{N_{G}}\left( \mathfrak{s}_{G}\right)
$) and the vertex set of the \textit{universal polytope}
Un$(\mathcal{V})$ (resp., of the \textit{secondary polytope}
Sec$(\mathcal{V})$) of $\mathfrak{s}_{G}$ w.r.t. the point
configuration $\mathcal{V}=\mathfrak{s}_{G}\cap N_{G};$ see
Appendix \ref{APPCOH}.\smallskip\ \newline(ii) There exist always
coherent maximal triangulations of $\mathfrak{s}_{G}$'s (see
\ref{NOTCOHEX}).\smallskip \
\newline(iii) For $r\geq3$ there are
$\mathfrak{s}_{G}$'s admitting maximal non-coherent
triangulations. \smallskip\ \newline(iv) For $r\geq4$ there are
lots of $\mathfrak{s}_{G}$'s admitting maximal non-basic
triangulations. \smallskip\ \newline (v) It is not known, as yet,
if there are $\mathfrak{s}_{G}$'s for which
$\mathbf{LTR}_{N_{G}}^{\text{{\scriptsize basic}}}\left(  \mathfrak{s}%
_{G}\right)  \neq\varnothing,$ whereas $\mathbf{Coh}$-$\mathbf{LTR}_{N_{G}%
}^{\text{{\scriptsize basic}}}\left(  \mathfrak{s}_{G}\right)
=\varnothing$ (see \ref{EXLATTTR} (iii)).\smallskip\ \newline(v)
An immediate consequence of Theorem \ref{TRCOR} is that the
\textquotedblleft Existence Problem\textquotedblright, as stated
in \S \ref{INTRO}, is in the abelian case \textit{equivalent} to
the following:
\end{note}

\noindent{}$\blacktriangleright$ \textbf{Existence Problem for
Gorenstein AQS}: For \textit{which abelian finite subgroups}
$G\subset$ SL$(r,\mathbb{C}),$
$r\geq4,$ \textit{do there exist triangulations }$\mathcal{T}\in\mathbf{LTR}%
_{N_{G}}^{\text{{\scriptsize basic}}}\left(
\mathfrak{s}_{G}\right)
$\textit{ }(\textit{and preferably }$\mathcal{T}\in\mathbf{Coh}$-$\mathbf{LTR}%
_{N_{G}}^{\text{{\scriptsize basic}}}\left(
\mathfrak{s}_{G}\right) $)?\medskip

\noindent{}This demonstrates one more example of the interplay between
algebraic and discrete geometry. The initial problem is fairly difficult, yet
once translated into discrete geometry, it can be treated by using familiar tools.

\begin{definition}
Identifying $\mathbb{C}^{r}/G$ with $X\left(  N_{G},\Delta_{G}\right)  $ as in
\S \ref{AQSTOR}, let
\[
\sigma_{\mathbf{s}}:=\text{pos}\left(  \mathbf{s}\right)  \subset\left(
N_{G}\right)  _{\mathbb{R}}\cong\mathbb{R}^{r}\
\]
denote the cone supporting a simplex\emph{ }$\mathbf{s}$ of a $\mathcal{T}%
\in\mathbf{LTR}_{N_{G}}\left(  \mathfrak{s}_{G}\right)  $. We define the
fan\emph{ }
\[
\widehat{\Delta_{G}}\left(  \mathcal{T}\right)  :=\left\{  \sigma_{\mathbf{s}%
}\ \left\vert \ \mathbf{s}\in\mathcal{T}\right.  \right\}
\]
and%
\[%
\begin{array}
[c]{ll}%
\mathbf{PCDES}\left(  X\left(  N_{G},\Delta_{G}\right)  \right)  := & \left\{
\begin{array}
[c]{c}%
\text{\textit{partial crepant} }\mathbb{T}_{N_{G}}\text{-equivariant }\\
\text{desingularizations of \ }X\left(  N_{G},\Delta_{G}\right) \\%
\begin{array}
[c]{c}%
\ \text{with overlying spaces having\emph{ }}\\
\text{\emph{\ }at worst\emph{ }}\mathbb{Q}\text{-factorial\emph{
}\textit{canonical} }\\
\text{singularities of index }1
\end{array}
\end{array}
\right\}  \ ,\smallskip\medskip\\
\mathbf{PCDES}^{\text{max}}\left(  X\left(  N_{G},\Delta_{G}\right)  \right)
:= & \left\{
\begin{array}
[c]{c}%
\text{\textit{partial crepant} }\mathbb{T}_{N_{G}}\text{-equivariant\emph{ }%
}\\
\text{desingularizations of \ }X\left(  N_{G},\Delta_{G}\right) \\%
\begin{array}
[c]{c}%
\ \text{with overlying spaces having\emph{ }}\\
\text{\emph{\ }at worst }\mathbb{Q}\text{-factorial \textit{terminal} \ }\\
\text{singularities\emph{ }of index\emph{ }}1
\end{array}
\end{array}
\right\}  \smallskip\ ,\medskip\\
\mathbf{CDES}\left(  X\left(  N_{G},\Delta_{G}\right)  \right)  := & \left\{
\begin{array}
[c]{c}%
\text{\textit{crepant} }\mathbb{T}_{N_{G}}\text{-equivariant (full)}\\
\text{\emph{\ \ }desingularizations of }X\left(  N_{G},\Delta_{G}\right)
\end{array}
\right\}  \ .\medskip\smallskip
\end{array}
\]
(Whenever we put the prefix $\mathbf{QP}$- in the front of any one
of them, we mean the corresponding subset of it consisting of
those desingularizations whose overlying spaces are\emph{
}\textit{quasiprojective}\emph{.})
\end{definition}

\begin{theorem}
[Desingularizing by triangulations]\label{TRCOR}Let $\left(  \mathbb{C}%
^{r}/G,\left[  \mathbf{0}\right]  \right)  $ be a Gorenstein AQS
\emph{(}$r\geq3$\emph{).} Then there exist one-to-one correspondences\emph{:}%
\smallskip%
\[
\fbox{$%
\begin{array}
[c]{ccc}
&  & \\
&
\begin{array}
[c]{ccc}%
\left(  \mathbf{Coh}\text{-}\right)  \mathbf{LTR}_{N_{G}}^{\emph{basic}%
}\left(  \mathfrak{s}_{G}\right)  \smallskip & \overset{\text{\emph{1:1}%
{\scriptsize \smallskip}}}{\longleftrightarrow} & \left(  \mathbf{QP}%
\text{-}\right)  \mathbf{CDES}\left(  X\left(  N_{G},\Delta_{G}\right)
\right) \\
\cap &  & \cap\\
\left(  \mathbf{Coh}\text{-}\right)  \mathbf{LTR}_{N_{G}}^{\emph{max}%
}\,\left(  \mathfrak{s}_{G}\right)  \smallskip & \overset{\text{\emph{1:1}%
{\scriptsize \smallskip}}}{\longleftrightarrow} & \left(  \mathbf{QP}%
\text{-}\right)  \mathbf{PCDES}^{\emph{max}}\left(  X\left(  N_{G},\Delta
_{G}\right)  \right) \\
\cap &  & \cap\\
\left(  \mathbf{Coh}\text{-}\right)  \mathbf{LTR}_{N_{G}}\left(
\mathfrak{s}_{G}\right)  & \overset{\text{\emph{1:1}{\scriptsize \smallskip}}%
}{\longleftrightarrow} & \left(  \mathbf{QP}\text{-}\right)  \mathbf{PCDES}%
\left(  X\left(  N_{G},\Delta_{G}\right)  \right)
\end{array}
& \\
&  &
\end{array}
$}\smallskip\smallskip
\]
which are realized by crepant $\mathbb{T}_{N_{G}}$-equivariant birational
morphisms of the form\smallskip\smallskip%
\begin{equation}
\fbox{$%
\begin{array}
[c]{ccc}
&  & \\
& f_{\mathcal{T}}=\text{\emph{id}}_{\ast}:X(N_{G},\text{ }\widehat{\Delta_{G}%
}\left(  \mathcal{T}\right)  )\longrightarrow X\left(  N_{G},\Delta_{G}\right)
& \\
&  &
\end{array}
$} \label{DESING}%
\end{equation}
induced by mapping
\[
\mathcal{T}\longmapsto\widehat{\Delta_{G}}\left(  \mathcal{T}\right)
,\ \ \ \ \ \widehat{\Delta_{G}}\left(  \mathcal{T}\right)  \longmapsto
X(N_{G},\text{ }\widehat{\Delta_{G}}\left(  \mathcal{T}\right)  ).
\]

\end{theorem}

\noindent\textsc{Sketch of proof}\textit{.} $X\left(  N_{G},\Delta_{G}\right)
$ is Gorenstein and has at most rational singularities, i.e., canonical
singularities of index $1$. Moreover, its dualizing sheaf is trivial. Let
\[
f=\text{id}_{\ast}:X(N_{G},\text{ }\widetilde{\Delta_{G}})\longrightarrow
X\left(  N_{G},\Delta_{G}\right)
\]
denote an arbitrary partial desingularization. Studying the behaviour of the
highest rational differentials on $X(N_{G},$ $\widetilde{\Delta_{G}})$ (see
\cite[\S 3]{ReidC3F} or \cite[Proposition 4.1]{DHZ}), one proves
\[
K_{X(N_{G},\text{ }\widetilde{\Delta_{G}})}=f^{\ast}\left(  K_{X\left(
N_{G},\Delta_{G}\right)  }\right)  -\sum_{\varrho\in\widetilde{\Delta_{G}%
}\left(  1\right)  }\ \left(  \left\langle \left(
1,\ldots,1\right) ,n\left(  \varrho\right)  \right\rangle
-1\right)  \ D_{n\left( \varrho\right)  }\ ,
\]
where
\[
D_{n\left(  \varrho\right)  }:=\mathbf{V}\left(  \mathbb{\varrho
}\right)  =\mathbf{V}\left(  \text{pos}\left(  \left\{  n\left(
\varrho\right)  \right\}  \right)  \right).
\]
Obviously, $f$ is crepant if and only if
Gen$(\widetilde{\Delta_{G})}\subset\mathcal{H}_{1},$ and since the
number of crepant exceptional prime divisors is independent of the
specific choice of $f$, the first and second $1$-$1$
correspondences (from below) are obvious by the
adjunction-theoretic definition of terminal (resp., canonical)
singularities. In particular, all $\mathbb{T}_{N_{G}}$-equivariant
partial crepant desingularizations of $X\left(
N_{G},\Delta_{G}\right)  $ of the form (\ref{DESING}) have
overlying spaces with at most $\mathbb{Q}$-factorial
singularities; and conversely, each partial
$\mathbb{T}_{N_{G}}$-equivariant
crepant desingularization with overlying space with at worst $\mathbb{Q}%
$-factorial singularities, has to be of this form. ($\mathbb{Q}$-factoriality
is here equivalent to the consideration only of triangulations instead of more
general polytopal subdivisions; cf. Proposition \ref{SMCR}. Furthermore, by
maximal triangulations we exhaust all crepant prime divisors). The top $1$-$1$
correspondence follows from the equivalence
\[
\mathcal{T}\ni\mathbf{s}\text{\ is a basic simplex w.r.t. }N_{G}%
\Longleftrightarrow\sigma_{\mathbf{s}}\in\widehat{\Delta_{G}}\left(
\mathcal{T}\right)  \text{ is a basic cone w.r.t. }N_{G}.
\]
To prove the $1$-$1$ correspondences after omitting the brackets, it suffices
to use the fact that the coherence of a triangulation $\mathcal{T}%
\in\mathbf{LTR}_{N_{G}}\left(  \mathfrak{s}_{G}\right)  $ implies
the existence of a strictly upper convex function defined on the
support of the entire fan $\widehat{\Delta_{G}}\left(
\mathcal{T}\right)  ,$ and then to apply ampleness criterion; cf.
\cite[Proposition 4.5, p. 211]{DHZ}.\hfill {}$\square$

\begin{remark}
Concerning the existence or non-existence of lattice triangulations
$\mathcal{T}\in\mathbf{LTR}_{N_{G}}^{\text{{\scriptsize basic}}}\left(
\mathfrak{s}_{G}\right)  $ for Gorenstein AQS $\left(  \mathbb{C}%
^{r}/G,\left[  \mathbf{0}\right]  \right)  $ of type (\ref{typeAQS}) we can
w.l.o.g. restrict ourselves to the class of \textit{msc}-AQS (as defined in
\ref{MSCDEF}), because the existence question for a non-msc-singularity is
obviously reduced to the same question for an msc-singularity of
\textit{strictly smaller} dimension. (The recognition of the those which are
msc-singularities follows from Proposition \ref{non-msc}.)
\end{remark}

\begin{note}
[Exceptional prime divisors]If $\mathcal{T}\in\mathbf{LTR}_{N_{G}%
}^{\text{{\scriptsize basic}}}\left(  \mathfrak{s}_{G}\right),$
then the exceptional prime divisors $D_{n\left(  \varrho\right) }$
 ($\varrho\in\widehat{\Delta_{G}}\left( \mathcal{T}\right)
(1)\mathbb{r}\Delta_{G}(1)$) w.r.t. $f_{\mathcal{T}}$ are
$(r-1)$-dimensional toric varieties whose topological Euler number
(\ref{EULERCHAR}) equals
\[
\chi(D_{n\left(  \varrho\right)  })=\sharp\{(r-1)\text{-dimensional simplices
of star}_{n\left(  \varrho\right)  }(\mathcal{T})\}.
\]
Moreover, $D_{n\left(  \varrho\right)  }$ is compact if and only if $n\left(
\varrho\right)  \in$ int$\left(  \mathfrak{s}_{G}\right)  .$ On the other
hand, if $\varrho\in\widehat{\Delta_{G}}\left(  \mathcal{T}\right)  \left(
1\right)  $ and $n\left(  \varrho\right)  \in$ int$\left(  \mathfrak{s}%
_{G}\left(  \nu_{1},\ldots,\nu_{k}\right)  \right)  \cap N_{G}$\emph{, }for
some $k,$ $2\leq k\leq r-1$\emph{, }and certain $1\leq\nu_{1}<\nu_{2}%
<\cdots<\nu_{k}\leq r$\emph{, }then the non-compact $D_{n\left(
\varrho\right)  }$\emph{\ }can be viewed as the total space of a fibration
$D_{n\left(  \varrho\right)  }\longrightarrow\mathbb{C}^{r-k}.$ The generic
fibers are isomorphic to the $\left(  k-1\right)  $-dimensional compact toric
variety associated to the star of $\varrho$ within
\[
\{\sigma\in\widehat{\Delta_{G}}\left(  \mathcal{T}\right)  \
\left\vert \ \sigma\prec\sigma_{\mathbf{s}},\
\mathbf{s}\in\mathfrak{s}_{G}\left(
\nu_{1},\ldots,\nu_{k}\right)\cap\mathcal{T}\}\right.
\]
In many cases, looking at the star$_{n\left(  \varrho\right)
}(\mathcal{T}),$ one can say more about the structure of
$D_{n\left(  \varrho\right)  }.$ (For concrete classes of
examples, see below Remark \ref{REMARKHYPERS} (i), Theorem
\ref{DESEXC}, and Remark \ref{WHITEPOINT} (ii).)
\end{note}

\noindent{}$\bullet$ \textbf{Cohomology dimensions}. \ Using (\ref{JUN}) and
(\ref{SEN}) we deduce from Theorem \ref{BATYRTHM} the following:

\begin{theorem}
\label{THMCOHDIMSIG}If $\left(  \mathbb{C}^{r}/G,\left[  \mathbf{0}\right]
\right)  $ is a Gorenstein AQS$,$ then for any crepant desingularization
$\widehat{X}\longrightarrow X$ of $X=\mathbb{C}^{r}/G$ we have $\dim
_{\mathbb{Q}}H^{0}(\widehat{X},\mathbb{Q})=1,$
\begin{equation}
\dim_{\mathbb{Q}}H^{2i}(\widehat{X},\mathbb{Q})=\sharp(\mathfrak{s}_{G}%
^{[i]}\cap\mathbf{Par}\left(  \sigma_{0}\right)  \cap N_{G}),\ \ \forall
i\in\{1,\ldots,r-1\}, \label{COHDIMSIG}%
\end{equation}
and the other cohomology groups of $\widehat{X}$ are trivial. In particular,
$\chi(\widehat{X})=\left\vert G\right\vert .$
\end{theorem}

\noindent{}To determine the cohomology dimensions (\ref{COHDIMSIG}) you may
exploit the description (\ref{Hyp1})-(\ref{Hyp2}) of $\mathfrak{s}_{G}%
^{[i]}\cap\mathbf{Par}\left(  \sigma_{0}\right)  $ in terms of hypersimplices
HypS$\left(  i,r\right)  $. But if you don't like to work directly with
hypersimplices, here is an alternative: Compute the coefficients of the
Ehrhart polynomial of the junior simplex $\mathfrak{s}_{G}$ by the formulae
given in Appendix \ref{LATPJS}, and then apply (\ref{COHDIMFORMULA}) instead
of (\ref{COHDIMSIG}).

\begin{theorem}
Maintaining the notation and the assumptions of \emph{\ref{THMCOHDIMSIG},} we
have%
\begin{equation}
\fbox{$%
\begin{array}
[c]{ccc}
& \dim_{\mathbb{Q}}H^{2i}(\widehat{X},\mathbb{Q})=\mathfrak{h}_{i}^{\ast
}(\mathfrak{s}_{G})=%
{\displaystyle\sum\limits_{j=0}^{r-1}} \left(
{\displaystyle\sum\limits_{\kappa=0}^{i}}
(-1)^{\kappa}\tbinom{r}{\kappa}\left(  i-\kappa\right)
^{j}\right) \mathbf{a}_{j}(\mathfrak{s}_{G}), &
\end{array}
$} \label{COHDIMFORMULA}%
\end{equation}
for all $i\in\{0,1,\ldots,r-1\},$ where by $\mathfrak{h}_{i}^{\ast
}(\mathfrak{s}_{G})$ is denoted the $i$-th component of the $\mathbf{h}^{\ast
}$-vector and by $\mathbf{a}_{j}(\mathfrak{s}_{G})$ the $j$-th coefficient of
the Ehrhart polynomial of $\mathfrak{s}_{G},$ respectively. \emph{(See
\ref{Ehr-Pol} and \ref{DEFHSTARVECTOR}.)}
\end{theorem}

\noindent{}\textsc{Proof}. If there exists a crepant desingularization
$\widehat{X}\longrightarrow X=\mathbb{C}^{r}/G,$ then there exists also a
$\mathbb{T}_{N_{G}}$-equivariant crepant desingularization (\ref{DESING})
induced by a triangulation $\mathcal{T}\in\mathbf{LTR}_{N_{G}}%
^{\text{{\scriptsize basic}}}\left(  \mathfrak{s}_{G}\right)  .$ Using Theorem
\ref{BATYRTHM}, \cite[Thm. 4.4]{BD} and Theorem \ref{BMcMTHM}, we get for all
$i\in\{0,1,\ldots,r-1\},$
\[
\dim_{\mathbb{Q}}H^{2i}(\widehat{X},\mathbb{Q})=\dim_{\mathbb{Q}}%
H^{2i}(X(N_{G},\widehat{\Delta_{G}}\left(  \mathcal{T}\right)  ),\mathbb{Q}%
)=\mathfrak{h}_{i}(\mathcal{T})=\mathfrak{h}_{i}^{\ast}(\mathfrak{s}_{G}),
\]
and it suffices to apply formula (\ref{H-A-FORMULA}) (for $P=\mathfrak{s}%
_{G},$ $d=r-1$) to obtain (\ref{COHDIMFORMULA}).\hfill{}$\square$

\begin{remark}
[A simple basicness criterion]\label{BASICCR}If $\left(  \mathbb{C}%
^{r}/G,\left[  \mathbf{0}\right]  \right)  $ is a Gorenstein AQS$,$ and
$\mathcal{T}\in\mathbf{LTR}_{N_{G}}\,\left(  \mathfrak{s}_{G}\right)  $, then
by (\ref{F-H-FORMULA}), (\ref{SUMOFHASTARS}) and (\ref{HAHASTAR}) we get%
\[
\mathfrak{f}_{r-1}(\mathcal{T})=\sharp\left\{
\begin{array}
[c]{c}%
(r-1)\text{-dim. }\\
\text{simplices of }\mathcal{T}%
\end{array}
\right\}  =\sum_{i=0}^{r-1}\mathfrak{h}_{i}(\mathcal{T})\leq\sum_{i=0}%
^{r-1}\mathfrak{h}_{i}^{\ast}(\mathfrak{s}_{G})=(r-1)!\,\text{Vol}%
(\mathfrak{s}_{G}),
\]
which implies%
\[
\mathfrak{f}_{r-1}(\mathcal{T})=\chi(X(N_{G},\widehat{\Delta}_{G}\left(
\mathcal{T}\right)  ))\leq(r-1)!\,\text{Vol}(\mathfrak{s}_{G})=\left\vert
G\right\vert =\tfrac{1}{\det(N_{G})},
\]
cf. (\ref{EULERCHAR}). This holds as equality:%
\begin{equation}
\mathfrak{f}_{r-1}(\mathcal{T})=\chi(X(N_{G},\widehat{\Delta}_{G}\left(
\mathcal{T}\right)  ))=(r-1)!\,\text{Vol}(\mathfrak{s}_{G})=\left\vert
G\right\vert \label{VOLEQUALITY}%
\end{equation}
\textit{if and only if} $\mathcal{T}\in\mathbf{LTR}_{N_{G}}%
^{\text{{\scriptsize basic}}}\left(  \mathfrak{s}_{G}\right)  $ (by Theorem
\ref{BMcMTHM}). \ In practice, having a concrete maximal triangulation
$\mathcal{T}$ of $\mathfrak{s}_{G}$ in hand, it suffices to compare
$\mathfrak{f}_{r-1}(\mathcal{T})$ with $\left\vert G\right\vert .$ If these
two numbers coincide, then $\mathcal{T}$ \ has to be basic.
\end{remark}

\noindent{}$\bullet$ \textbf{Flops}. If
$\mathcal{T},\mathcal{T}^{\prime}$ are two \textit{coherent}
lattice triangulations of $\mathfrak{s}_{G},$ are there
\textquotedblleft elementary operations\textquotedblright\ whose
repetitive use would geometrically describe how one can obtain $\mathcal{T}%
^{\prime}$ from $\mathcal{T}$ ? On the level of triangulations a satisfactory
answer is given by the \textit{bistellar flips}\footnote{One of the main
reasons for adding to the triangulations involved in the above formulation of
Existence Problem the phrase \textquotedblleft preferably
coherent\textquotedblright\ is their connection by a sequence of bistellar
flips. This does not hold in general for non-coherent triangulations. For
instance, Santos provided in \cite{Santos} a point set whose space of (all)
triangulations is bistellarly \textit{disconnected}.\smallskip} (as defined in
combinatorial topology). If, in addition, $\mathcal{T},\mathcal{T}^{\prime}$
are assumed to be maximal, this answer on the level of birational maps
connecting $X(N_{G},\widehat{\Delta}_{G}\left(  \mathcal{T}\right)  )$ with
$X(N_{G},\widehat{\Delta}_{G}\left(  \mathcal{T}^{\prime}\right)  )$ leads to
algebro-geometric \textit{flops}.

\begin{theorem}
[Bistellar flips, and flops]\label{FLOPTHM}\emph{(i)} If $\smallskip
\mathcal{T},\mathcal{T}^{\prime}\in\mathbf{Coh}$-$\mathbf{LTR}_{N_{G}}\left(
\mathfrak{s}_{G}\right)  ,$ then there exist finitely many circuits%
\[
\mathcal{C}_{1},\ldots,\mathcal{C}_{\kappa}\subset\mathfrak{s}_{G}\cap
N_{G},\ \ \text{and\ \ }\mathcal{T}_{1},\ldots,\mathcal{T}_{\kappa}%
\in\mathbf{Coh}\text{-}\mathbf{LTR}_{N_{G}}\left(  \mathfrak{s}_{G}\right)
\]
such that $\mathcal{T}_{i+1}=$ \emph{FL}$_{\mathcal{C}_{i}}\left(
\mathcal{T}_{i}\right)  $ \emph{(}i.e., such that $\mathcal{T}_{i+1}$ is the
\emph{bistellar flip} of $\mathcal{T}_{i}$ along $\mathcal{C}_{i},$ \emph{cf.
\ref{DEFFLIP})} for all $i\in\left\{  1,\ldots,\kappa-1\right\}  ,$ with
$\mathcal{T}_{1}=\mathcal{T}$ and $\mathcal{T}_{\kappa}=\mathcal{T}^{\prime}%
$.\medskip\newline\emph{(ii)} In particular, if $\mathcal{T},\mathcal{T}%
^{\prime}\in\mathbf{Coh}$-$\mathbf{LTR}_{N_{G}}^{max}\,\left(
\mathfrak{s}_{G}\right)  ,$ then the circuits
$\mathcal{C}_{1},\ldots
,\mathcal{C}_{\kappa}$ can be chosen in such a way that $\sharp(\mathcal{C}%
_{i})=r+1$ and \emph{dim}$($\emph{conv}$(\mathcal{C}_{i}))=r-1$ for all
$i\in\{1,\ldots,\kappa-1\}.$ Setting $X_{i}:=X(N_{G},\widehat{\Delta}%
_{G}\left(  \mathcal{T}_{i}\right)  ),$ $X:=X_{1},$ and $X^{\prime}%
:=X_{\kappa},$ we conclude that $X$ and $X^{\prime}$ can be obtained from each
other by a finite succession of flops\emph{\footnote{In the MMP-language (and
as long as one may work in the category of \textit{quasi}projective complex
varieties) we say that the \textquotedblleft minimal models\textquotedblright%
\ $X$ and $X^{\prime}$ are connected by a sequence of flops whose
\textquotedblleft termination\textquotedblright\ is due to our specific
setting; cf. \cite[Thm. 3.4.6, p. 158]{Matsuki}. In fact, these flops can be
conceived of as high-dimensional analogues of \ the original \textquotedblleft
Atiyah's flop\textquotedblright\ (see \cite[Example 3.4.3, p. 157]{Matsuki}).
}} which fit together into the following diagram\emph{:}

\[
\xymatrix @!0 @R=2.5pc @C=3.6pc {%
&  Y_{1}&  & \relax Y_{2}&  &  & &  \relax Y_{{\kappa}-1}&  &  &  & \\%
\relax X_{1} \ar@{<--}[rr]\ar[ur]^{\varphi _{1}} \ar@/_2pc/[rdrdrr]^{f_{\mathcal{T}_{1}}}&%
& \relax X_{2} \ar@{<--}[rr] \ar[ur]^{\varphi _{2}} \ar[ul]_{\vartheta _{1}}  \ar@/_1pc/[drdr]^{f_{\mathcal{T}_{2}}}&%
&X_{3} \ar[ul]_{\vartheta _{2}} \ar[dd]_{f_{\mathcal{T}_{3}}}%
\ar@{<--}[r]&...&X_{{\kappa}-1}\ar@{-->}[l] \ar@{<--}[rr] \ar[ur]^{\varphi _{\kappa -1}} \ar@/^1pc/[dldl]_{f_{\mathcal{T}_{\kappa -1}}}&%
&  X_{{\kappa}}\ar[ul]_{\vartheta _{\kappa -1}} \ar@/^2pc/[dldlll]_{f_{\mathcal{T}_{\kappa}}}& \\
&  &  &  &  &  &  &  &  & &  &  \\
&  &  &  &X(N_{G},\Delta _{G})&  &  &  &  & &  &}
\]
Here, by \emph{\textquotedblleft flops\textquotedblright} we mean the upper
triangles\ of the diagram, where both $\varphi_{i}$ and $\vartheta_{i}$ are
\emph{small} birational morphisms \emph{(}i.e., their exceptional loci have
codimension $\geq2$\emph{) }and $X_{i+1}\dashrightarrow X_{i}$ birational maps
which are isomorphisms in codimension $1.$
\end{theorem}

\noindent{}\textsc{Proof}. We shall use the notation and the terminology
introduced in Appendix \ref{APPCOH}.\smallskip\ \newline(i) Consider an edge
path $\overline{\mathbf{v}_{1}\mathbf{v}_{2}},$ $\overline{\mathbf{v}%
_{2}\mathbf{v}_{3}},\ldots,\overline{\mathbf{v}_{\kappa-1}\mathbf{v}_{\kappa}%
}$ on the polytope Sec$(\mathfrak{s}_{G}\cap N_{G})$ connecting $\mathbf{v}%
_{1}:=\mathbf{v}_{\mathcal{T}_{1}}$ with $\mathbf{v}_{\kappa}:=\mathbf{v}%
_{\mathcal{T}_{\kappa}}.$ By Theorem \ref{FLIPPING} one determines circuits
$\mathcal{C}_{1},\ldots,\mathcal{C}_{\kappa}\subset\mathfrak{s}_{G}\cap N_{G}$
such that $\mathbf{v}_{i}=\mathbf{v}_{\mathcal{T}_{i}}$ with $\mathcal{T}%
_{i+1}=$ FL$_{\mathcal{C}_{i}}\left(  \mathcal{T}_{i}\right)  ,$ $\forall
i\in\left\{  1,\ldots,\kappa-1\right\}  .\medskip$\newline(ii) If
$\mathcal{T}$ is a maximal triangulation of $\mathfrak{s}_{G},$ and
$\mathbf{s}_{1},\mathbf{s}_{2}$ two $(r-1)$-dimensional simplices of
$\mathcal{T}$ having $\mathbf{s}_{1}\cap\mathbf{s}_{2}$ as $(r-2)$-dimensional
common face, then $\left.  \mathcal{T}\right\vert _{\mathbf{s}_{1}%
\cup\mathbf{s}_{2}}$ is either the triangulation $\mathcal{Y}_{+}\left(
\mathcal{C}\right)  $ or the triangulation $\mathcal{Y}_{-}\left(
\mathcal{C}\right)  $ of conv$($vert$(\mathbf{s}_{1})\cup$vert$(\mathbf{s}%
_{2}))$ w.r.t. the circuit $\mathcal{C}=$ vert$(\mathbf{s}_{1})\cup
$vert$(\mathbf{s}_{2})$ with $\sharp(\mathcal{C})=r+1$ (cf. Lemma
\ref{CIRCUITTR}). To pass from $\mathcal{T}$ to another maximal triangulation
$\mathcal{T}^{\prime}$ it suffices to apply (i) for circuits $\mathcal{C}%
_{1},\ldots,\mathcal{C}_{\kappa}$ \textit{only of this kind}. (This follows
from results of Oda \& Park \cite[Corollary 3.9, Proposition 3.10, and Theorem
3.12, pp. 395-398]{OP}.) After having determined such $\mathcal{T}_{i}$'s
(with $\mathcal{T}_{i+1}=$ FL$_{\mathcal{C}_{i}}\left(  \mathcal{T}%
_{i}\right)  $ for all $i\in\{1,\ldots,\kappa-1\}$), it is enough to define
$Y_{i}$ to be the Gorenstein toric variety associated to the fan which
consists of the cones supporting the lattice polytopes of the polytopal
subdivision $\mathcal{T}_{i}\smallsetminus(\mathcal{Y}_{\odot}\left(
\mathcal{C}_{i}\right)  \ast\mathcal{T}^{\left[  \mathcal{C}_{i}\right]  }$ of
$\mathfrak{s}_{G}.$ By the birational morphism $\varphi_{i}$ we contract
$\mathbf{V}(\sigma_{\mathbf{s}_{i}}),$ where $\mathbf{s}_{i}$ denotes the
unique $(r-2)$-dimensional simplex of $\mathcal{Y}_{\odot}\left(
\mathcal{C}_{i}\right)  $ with int$(\mathbf{s}_{i})\subset$ int$(\mathcal{Y}%
_{\odot}\left(  \mathcal{C}_{i}\right)  ),$ and by $\vartheta_{i}$ we extract
$\mathbf{V}(\sigma_{\mathbf{t}_{i}}),$ where $\mathbf{t}_{i}$ denotes the
unique $(r-2)$-dimensional simplex of $\mathcal{Y}_{\boxtimes}\left(
\mathcal{C}_{i}\right)  $ with int$(\mathbf{t}_{i})\subset$ int$(\mathcal{Y}%
_{\boxtimes}\left(  \mathcal{C}_{i}\right)  ).$ The $\vartheta_{i}$'s are
non-divisorial extractions, because we do not introduce any new vertices in
the triangulation $\mathcal{T}_{i+1}.$\hfill{}$\square$

\begin{exercise}
Take again the example of CQS of type $\frac{1}{12}\left(  1,2,3,6\right)  $
as in \ref{Exam12}. Working with \texttt{Puntos} (cf. Note \ref{DELOERANOTE})
we find all $\mathcal{T}\in\mathbf{Coh}$-$\mathbf{LTR}_{N_{G}}%
^{\text{{\scriptsize max}}}\,\left(  \mathfrak{s}_{G}\right)  .$ These are
altogether $12$ triangulations: One of them has $9$ simplices, two have $10$
simplices, four have $11$ simplices, and the remaining five have $12$
simplices. The latter ones are necessarily the elements of the set
$\mathbf{Coh}$-$\mathbf{LTR}_{N_{G}}^{\text{{\scriptsize basic}}}\left(
\mathfrak{s}_{G}\right)  .$ The vertex sets of their simplices (in the
notation used in \ref{Exam12}) are recorded in the following list.{\small
\[%
\begin{tabular}
[c]{|c|c|c|c|c|}\hline
$\mathcal{T}_{1}$ & $\mathcal{T}_{2}$ & $\mathcal{T}_{3}$ & $\mathcal{T}_{4}$
& $\mathcal{T}_{5}$\\\hline\hline
$\left\{  e_{2},e_{4},n_{1},n_{3}\right\}  $ & $\left\{  e_{2},e_{4}%
,n_{1},n_{3}\right\}  $ & $\left\{  e_{2},e_{4},n_{1},n_{3}\right\}  $ &
$\left\{  e_{2},e_{4},n_{1},n_{3}\right\}  $ & $\left\{  e_{2},e_{4}%
,n_{1},n_{3}\right\}  $\\\hline
$\left\{  e_{3},n_{1},n_{2},n_{6}\right\}  $ & $\left\{  e_{3},n_{1}%
,n_{2},n_{6}\right\}  $ & $\left\{  e_{3},n_{1},n_{2},n_{9}\right\}  $ &
$\left\{  e_{3},n_{1},n_{2},n_{6}\right\}  $ & $\left\{  e_{3},n_{1}%
,n_{2},n_{6}\right\}  $\\\hline
$\left\{  e_{3},e_{4},n_{1},n_{6}\right\}  $ & $\left\{  e_{3},e_{4}%
,n_{1},n_{6}\right\}  $ & $\left\{  e_{3},e_{4},n_{1},n_{6}\right\}  $ &
$\left\{  e_{3},e_{4},n_{1},n_{6}\right\}  $ & $\left\{  e_{3},e_{4}%
,n_{1},n_{6}\right\}  $\\\hline
$\left\{  e_{2},e_{3},e_{4},n_{1}\right\}  $ & $\left\{  e_{2},e_{3}%
,e_{4},n_{1}\right\}  $ & $\left\{  e_{2},e_{3},e_{4},n_{1}\right\}  $ &
$\left\{  e_{2},e_{3},e_{4},n_{1}\right\}  $ & $\left\{  e_{2},e_{3}%
,e_{4},n_{1}\right\}  $\\\hline
$\left\{  e_{2},e_{3},n_{1},n_{2}\right\}  $ & $\left\{  e_{2},e_{3}%
,n_{1},n_{2}\right\}  $ & $\left\{  e_{2},e_{3},n_{1},n_{2}\right\}  $ &
$\left\{  e_{2},e_{3},n_{1},n_{2}\right\}  $ & $\left\{  e_{2},e_{3}%
,n_{1},n_{2}\right\}  $\\\hline
$\left\{  e_{2},n_{1},n_{2},n_{6}\right\}  $ & $\left\{  e_{2},n_{1}%
,n_{2},n_{6}\right\}  $ & $\left\{  e_{3},n_{1},n_{6},n_{9}\right\}  $ &
$\left\{  n_{1},n_{2},n_{6},n_{9}\right\}  $ & $\left\{  n_{1},n_{3}%
,n_{6},n_{9}\right\}  $\\\hline
$\left\{  e_{1},e_{4},n_{1},n_{9}\right\}  $ & $\left\{  e_{1},e_{4}%
,n_{1},n_{9}\right\}  $ & $\left\{  e_{1},e_{4},n_{1},n_{9}\right\}  $ &
$\left\{  e_{1},e_{4},n_{1},n_{9}\right\}  $ & $\left\{  e_{1},e_{4}%
,n_{1},n_{9}\right\}  $\\\hline
$\left\{  e_{4},n_{1},n_{3},n_{9}\right\}  $ & $\left\{  e_{4},n_{1}%
,n_{3},n_{9}\right\}  $ & $\left\{  e_{4},n_{1},n_{3},n_{9}\right\}  $ &
$\left\{  e_{4},n_{1},n_{3},n_{9}\right\}  $ & $\left\{  e_{4},n_{1}%
,n_{3},n_{9}\right\}  $\\\hline
$\left\{  e_{1},e_{4},n_{1},n_{6}\right\}  $ & $\left\{  e_{1},e_{4}%
,n_{1},n_{6}\right\}  $ & $\left\{  e_{1},e_{4},n_{1},n_{6}\right\}  $ &
$\left\{  e_{1},e_{4},n_{1},n_{6}\right\}  $ & $\left\{  e_{1},e_{4}%
,n_{1},n_{6}\right\}  $\\\hline
$\left\{  e_{2},n_{1},n_{3},n_{6}\right\}  $ & $\left\{  e_{2},n_{1}%
,n_{2},n_{3}\right\}  $ & $\left\{  e_{2},n_{1},n_{2},n_{3}\right\}  $ &
$\left\{  e_{2},n_{1},n_{2},n_{3}\right\}  $ & $\left\{  e_{2},n_{1}%
,n_{2},n_{3}\right\}  $\\\hline
$\left\{  n_{1},n_{3},n_{6},n_{9}\right\}  $ & $\left\{  n_{1},n_{2}%
,n_{3},n_{9}\right\}  $ & $\left\{  n_{1},n_{2},n_{3},n_{9}\right\}  $ &
$\left\{  n_{1},n_{2},n_{3},n_{9}\right\}  $ & $\left\{  n_{1},n_{2}%
,n_{3},n_{6}\right\}  $\\\hline
$\left\{  e_{1},n_{1},n_{6},n_{9}\right\}  $ & $\left\{  e_{1},n_{1}%
,n_{2},n_{9}\right\}  $ & $\left\{  e_{1},n_{1},n_{6},n_{9}\right\}  $ &
$\left\{  e_{1},n_{1},n_{6},n_{9}\right\}  $ & $\left\{  e_{1},n_{1}%
,n_{6},n_{9}\right\}  $\\\hline
\end{tabular}
\
\]
}For $1\leq i\leq j\leq5,$ $i\neq j,$ find sequences of flops connecting
$X(N_{G},\widehat{\Delta}_{G}\left(  \mathcal{T}_{i}\right)  )$ with
$X(N_{G},\widehat{\Delta}_{G}\left(  \mathcal{T}_{j}\right)  ),$ and
distinguish those possessing the \textit{smallest} number of flops.
\end{exercise}

\section{The C.I.-Case\label{CICASE}}

\noindent{}An evidence in support of Conjecture \ref{CONJCIS} is given by the following:

\begin{theorem}
[\cite{DHZ}]\label{DHZTHM}All abelian quotient c.i.-singularities admit
projective, crepant resolutions in all dimensions.
\end{theorem}

\noindent An extensive technical part of its proof is devoted to the rendering
of the original (purely algebraic) group classification of Watanabe
\cite{Wat2} into graph-theoretic terms and to a subsequent \textit{convenient}
description of the corresponding junior simplices. As it turns out, an AQS is
a c.i.-singularity if and only if the junior simplex $\mathfrak{s}_{G}$ is
(what we call) a \textit{Watanabe simplex w.r.t. }$N_{G}$. (In addition,
notice that every abelian quotient c.i.-msc-singularity of dimension $\geq3$
has to be \textit{non-cyclic}!)

\begin{definition}
Let $d$ be an integer $\geq0$\ and $N$ a free $\mathbb{Z}$-module of rank $d,$
regarded as a lattice within $N_{\mathbb{R}}\cong\mathbb{R}^{d}$\emph{. }The
\textit{Watanabe simplices} w.r.t. $N$ are the lattice\emph{\ }%
simplices\emph{\ }$\mathbf{s}$ (w.r.t. $N,$ of dimension $\leq d$) satisfying%
\[
\text{aff}_{\mathbb{Z}}\left(  \mathbf{s\cap}N\right)  =\text{aff}\left(
\mathbf{s}\right)  \cap N
\]
which are defined inductively (starting in dimension $0$) in the following
manner:\emph{\smallskip}\newline(i) Every $0$-dimensional lattice simplex
$\mathbf{s}=\left\{  n\right\}  $, $n\in N$, is a Watanabe simplex.\smallskip
\newline(ii) A lattice simplex $\mathbf{s}\subset N_{\mathbb{R}}$\emph{\ }of
dimension\emph{\ }$d^{\prime}$\emph{, }$1\leq d^{\prime}\leq d$, is a Watanabe
simplex if\medskip\newline$\bullet$ \textit{either }$\mathbf{s}=\mathbf{s}%
_{1}\ast\mathbf{s}_{2}$ (the \textit{join} of $\mathbf{s}_{1}$ and
$\mathbf{s}_{2}$), where $\mathbf{s}_{1}$, $\mathbf{s}_{2}$ are Watanabe
simplices of dimensions $d_{1}$, $d_{2}\geq0$, $d_{1}+d_{2}=d^{\prime}-1$,
with respect to\emph{\ }sublattices $N_{1}\subset$ aff$\left(  \mathbf{s}%
_{1}\right)  $\emph{,} $N_{2}\subset$ aff$\left(  \mathbf{s}_{2}\right)  $
of\emph{\ }$N$, such that aff$_{\mathbb{Z}}\left(  \mathbf{s\cap}N\right)  =$
aff$_{\mathbb{Z}}\left(  N_{1}\cup N_{2}\right)  $\emph{,\medskip}%
\newline$\bullet$ \textit{or }$\mathbf{s}$ is a lattice translate of some
dilation\emph{\ }$\lambda\,\mathbf{s}^{\prime}$, where $\lambda\in\mathbb{Z}%
$\emph{,} $\lambda\geq2$, and $\mathbf{s}^{\prime}$ is an $d^{\prime}%
$-dimensional Watanabe simplex with respect\emph{\ }to $N$.\medskip
\newline(These conditions are mutually exclusive; with this definition every
affine integral transformation that preserves $N$ also preserves the Watanabe
simplices w.r.t. $N$).
\end{definition}

\noindent Theorem \ref{DHZTHM} results from the following:

\begin{theorem}
\label{WATSIMPTR}All Watanabe simplices w.r.t. a lattice $N$ possess basic,
coherent triangulations.
\end{theorem}

\noindent{}To prove \ref{WATSIMPTR} it suffices to show that: (i) joins and
dilations of coherent triangulations of lattice polytopes remain coherent,
(ii) the join of two basic simplices is basic, and (iii) the dilation of a
basic simplex by a factor $k\in\mathbb{Z},$ $k\geq2,$ possesses a basic
triangulation (see \cite[Theorem 3.5, Lemma 3.6 and Proposition 6.1]{DHZ}).

\begin{figure}[h]
\epsfig{file=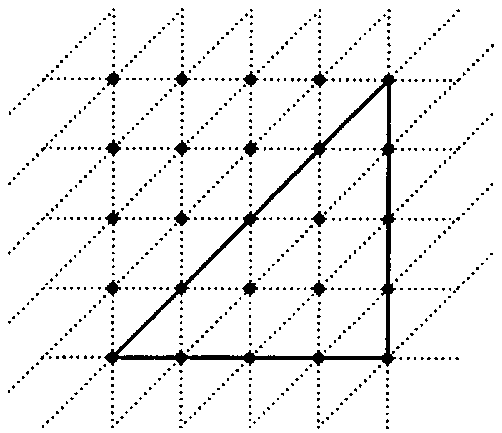}
\caption{}\label{Fig.3}
\end{figure}

\begin{example}
\label{HVS}Let\emph{\ }$\mathbb{H}_{d}$ denote the affine hyperplane
arrangement of type $\widetilde{\mathcal{A}}_{d}$\emph{\ }in $\mathbb{R}^{d}$
consisting of the union of hyperplanes
\[
\left\{  \!\left\{  \mathbf{x\in\,}\mathbb{R}^{d}\left\vert x_{i}%
=\kappa\right.  \right\}  \!,1\leq i\leq d,\kappa\in\mathbb{Z}\right\}
\cup\left\{  \!\left\{  \mathbf{x\in\,}\mathbb{R}^{d}\left\vert x_{i}%
-x_{j}=\lambda\right.  \right\}  \!,1\leq i<j\leq d,\lambda\in\mathbb{Z}%
\right\}  \!.
\]
$\mathbb{H}_{d}$ induces a basic triangulation $\mathbf{T}_{\mathbb{H}_{d}}$
(w.r.t. $\mathbb{Z}^{d}$) on the \textit{entire }space $\mathbb{R}^{d}.$ Let
Hvs$:\mathbb{R}\longrightarrow\mathbb{R}$ denote the Heaviside function%
\[
\text{Hvs}(x):=\left\{
\begin{array}
[c]{ll}%
x, & \text{if }x\geq0,\\
0, & \text{otherwise.}%
\end{array}
\right.
\]
The $\mathbf{T}_{\mathbb{H}_{d}}$-support function
\[
\psi_{\text{Hvs}}^{(d)}(\mathbf{x})\!:=\!-%
{\textstyle\sum\limits_{0\leq i<j\leq d}} \left\{
{\textstyle\sum\limits_{0\leq\kappa\leq x_{j}-x_{i}}}
\!\text{Hvs}(x_{i}-x_{j}-\kappa)+%
{\textstyle\sum\limits_{x_{j}-x_{i}\leq\kappa\leq0}}
\!\text{Hvs}(\kappa-x_{j}+x_{i})\!\right\}  \!,
\]
$\forall\mathbf{x}=(x_{1},\ldots,x_{d})\in\mathbb{R}^{d},$ with $x_{0}:=0,$ is
strictly upper convex. Thus, $\mathbf{T}_{\mathbb{H}_{d}}$ is also coherent.
Next, define
\begin{align*}
\mathbf{s}_{d}  &  :=\text{ conv}(\{\mathbf{0},e_{1},e_{1}+e_{2},\ldots
,e_{1}+\cdots+e_{d}\})\\
&  =\left\{  \mathbf{x\in\,}\mathbb{R}^{d}\,\left\vert \,0\leq x_{d}\leq
x_{d-1}\leq\cdots\leq x_{1}\leq1\right.  \right\}  ,
\end{align*}
and let $k$ be an integer $\geq2.$ Since the affine hulls of the facets of
$\mathbf{s}_{d}$ belong to $\mathbb{H}_{d},$ the restriction $\mathbf{T}%
\left(  d;k\right)  :=\left.  \mathbf{T}_{\mathbb{H}_{d}}\right\vert
_{k\mathbf{s}_{d}}$ of $\mathbf{T}_{\mathbb{H}_{d}}$ on $k\,\mathbf{s}_{d}$ is
a basic, coherent triangulation of $k\,\mathbf{s}_{d}$ w.r.t. $\mathbb{Z}%
^{d}.$ The triangulation $\mathbf{T}\left(  2;4\right)  $ of $4\,\mathbf{s}%
_{2}$ is depicted in Figure \ref{Fig.3}.
\end{example}

\begin{application}
[Hypersurface case]\label{HYPERSU}The underlying space of a
Gorenstein AQS $\left(  \mathbb{C}^{r}/G,\left[  \mathbf{0}\right]
\right)  ,$ $r\geq3,$ is embeddable as a hypersurface in
$\mathbb{C}^{r+1}$ if and only if $\ G$ is conjugate (within
SL$(r,\mathbb{C})$) to a group of the form
\[
G(r;k):=\left\langle \{\text{diag}(1,...,1,\underset{i\text{-th pos.}%
}{\underbrace{\zeta_{k}}},\underset{\left(  i+1\right)  \text{\emph{-}pos.}%
}{\underbrace{\zeta_{k}^{-1}}},1,...,1)\left\vert \ 1\leq i\leq
r-1\right. \}\right\rangle ,
\]
with $k\in\mathbb{Z},$ $k$ $\geq2,$ i.e., if and only if it is of type
\[
\frac{1}{k}\left(  1,1,0,\ldots,0\right)  \times\frac{1}{k}\left(
0,1,1,0,\ldots,0\right)  \times\ \cdots\ \cdots\ \times\frac{1}{k}\left(
0,0,\ldots,0,1,1\right)  .
\]
In this case, we may identify $\mathbb{C}^{r}/G(r;k)$ (or $\mathbb{C}^{r}/G,$
cf. Thm. \ref{prill-th}) with
\[
\left\{  \left(  z_{0},z_{1},\ldots,z_{r}\right)  \in\mathbb{C}^{r+1}%
\left\vert \ z_{0}^{k}=%
{\textstyle\prod\limits_{i=1}^{r}} z_{i}\right.  \right\}  ,
\]
and write the junior simplex $\mathfrak{s}_{G}$ as the dilation of a basic
simplex (w.r.t. $N_{G}$) by the factor $k$:
\[
\mathfrak{s}_{G}=k\ \text{conv}(\{\tfrac{1}{k}\ e_{1},\ldots,\tfrac{1}%
{k}\ e_{r}\}).
\]
There is an affine transformation $\mathbb{R}^{r}\longrightarrow
\mathbb{R}^{r-1}\times\{0\}\subset\mathbb{R}^{r}$ whose restriction on the
affine hull aff$(\mathfrak{s}_{G})$ of $\mathfrak{s}_{G}$ is a bijection, say
$\Xi,$ mapping the lattice aff$(\mathfrak{s}_{G})\cap N_{G}$ onto the standard
rectangular lattice $\mathbb{Z}^{r-1}=\mathbb{Z}^{r-1}\times\{0\}\subset
\mathbb{R}^{r-1}\times\{0\},$ and $\mathfrak{s}_{G}$ onto $\Xi\left(
\mathfrak{s}_{G}\right)  =k\,\mathbf{s}_{r-1}.$ Since $\mathbf{T}\left(
r-1;k\right)  $ (as defined in \ref{HVS}, with $d=r-1$) is a basic coherent
triangulation of $k\,\mathbf{s}_{r-1}$ w.r.t. $\mathbb{Z}^{r-1},$
\[
f_{\Xi^{-1}(\mathbf{T}\left(  r-1;k\right)  )}:X(N_{G},\widehat{\Delta}%
_{G}(\Xi^{-1}(\mathbf{T}\left(  r-1;k\right)  ))\longrightarrow X\left(
N_{G},\Delta_{G}\right)  =\mathbb{C}^{r}/G
\]
is a projective crepant desingularization of the quotient space $\mathbb{C}%
^{r}/G\cong\mathbb{C}^{r}/G(r;k).$
\end{application}

\begin{remark}
\label{REMARKHYPERS}If $G\subset$ SL$(r,\mathbb{C})$ is conjugate to $G(r;k),$
then\smallskip\ \newline(i) the star of any vertex of $\mathbf{T}\left(
r-1;k\right)  ,$ belonging to the interior of $k\,\mathbf{s}_{r-1},$ is
constructed as the image under an appropriate integral affine transformation
of the \textit{join} of the origin $\mathbf{0}\in\mathbb{R}^{r-1}$ with the
facets of the zonotope which is defined as the convex hull of the union of the
$[-1,0]$-cube and the $[0,1]$-cube; cf. \cite{DHaZ}. Hence, every compactly
supported exceptional prime divisor on $X(N_{G},\widehat{\Delta}_{G}(\Xi
^{-1}(\mathbf{T}\left(  r-1;k\right)  ))$ is a Fano manifold obtained by a
$\mathbb{T}_{N_{G}}$-equivariant projective crepant desingularization of a
toric Fano variety (with at worst Gorenstein singularities). This Fano variety
turns out to be a projective variety of degree $\tbinom{2r}{r}$ embedded in
$\mathbb{P}_{\mathbb{C}}^{r(r+1)}.$ \smallskip\newline(ii) Besides $\Xi
^{-1}(\mathbf{T}\left(  r-1;k\right)  $ there are lots of other basic,
coherent triangulations of $\mathfrak{s}_{G}$ corresponding to different
vertices of its secondary polytope. For instance, if $G$ is conjugate to
$G(4;2),$ \texttt{Puntos} \cite{DELOERA2} gives us $196$ maximal coherent
triangulations of $\mathfrak{s}_{G}.$ $192$ out of them are basic.
\smallskip\newline(iii) The non-trivial cohomology dimensions
(\ref{COHDIMFORMULA}) of any crepant desingularization $\widehat{X}$ of
$X=\mathbb{C}^{r}/G$ are equal to%
\[
\dim_{\mathbb{Q}}H^{2i}(\widehat{X},\mathbb{Q})=\sum_{j=0}^{i}\left(
-1\right)  ^{j}\tbinom{r}{j}\tbinom{k(i-j)+r-1}{r-1}.
\]

\end{remark}

\begin{note}
\label{ROANREM}(i) For $r=4,$ $k=2,$ Chiang \& Roan \cite[Thm. 4.1]%
{CHIANG-ROAN1}, \cite[\S 4]{CHIANG-ROAN2}, proved that the Hilbert scheme
$G(4;2)$-Hilb$(\mathbb{C}^{4})$ is a four-dimensional non-singular toric
variety with \textit{non-trivial canonical divisor}. The dualizing sheaf
$\omega_{G(4;2)\text{-Hilb}(\mathbb{C}^{4})}$ is $\cong\mathcal{O}%
_{G(4;2)\text{-Hilb}(\mathbb{C}^{4})}(\mathbb{P}_{\mathbb{C}}^{1}%
\times\mathbb{P}_{\mathbb{C}}^{1}\times\mathbb{P}_{\mathbb{C}}^{1}).$ There
are three different ways to blow down this divisor and pass to crepant
desingularizations of $\mathbb{C}^{4}/G(4;2),$ corresponding to the three
different projections $\mathbb{P}_{\mathbb{C}}^{1}\times\mathbb{P}%
_{\mathbb{C}}^{1}\times\mathbb{P}_{\mathbb{C}}^{1}\longrightarrow
\mathbb{P}_{\mathbb{C}}^{1}\times\mathbb{P}_{\mathbb{C}}^{1}.$ The first
blow-down leads to the crepant desingularization of $\mathbb{C}^{4}/G(4;2)$
described in \ref{HYPERSU}. The other two are obtained by flops, and belong to
the $192$ mentioned in \ref{REMARKHYPERS} (ii).\smallskip\ \newline(ii) For
$r=5,$ $k=2,$ the situation becomes worse. $G(5;2)$-Hilb$(\mathbb{C}^{5})$ is
a five-dimensional \textit{singular} toric variety with \textit{non-trivial
canonical divisor}. In this case, among all crepant $\mathbb{T}_{N_{G(5;2)}}%
$-equivariant desingularizations of $\mathbb{C}^{5}/G(5;2)$ there are only
$12$ dominated by $G(5;2)$-Hilb$(\mathbb{C}^{5})$ (see \cite[\S 5]%
{CHIANG-ROAN2}).
\end{note}

\section{First Existence Criterion via Hilbert basis\label{FIRSTEXCR}}

\noindent{}A \textit{necessary} condition for an \textit{arbitrary} Gorenstein
AQS $\left(  \mathbb{C}^{r}/G,\left[  \mathbf{0}\right]  \right)  $ to admit a
crepant resolution is described as follows:

\begin{theorem}
[First Necessary Existence Condition]\label{KILLER}\ \smallskip Let\emph{\ }%
$\left(  \mathbb{C}^{r}/G,\left[  \mathbf{0}\right]  \right)  $ be a
Gorenstein AQS. If $\mathfrak{s}_{G}$ has a basic triangulation, then
\begin{equation}
\fbox{$%
\begin{array}
[c]{l}
\end{array}
\mathbf{Hlb}_{N_{G}}\left(  \sigma_{0}\right)  =\mathfrak{s}_{G}\cap N_{G},%
\begin{array}
[c]{l}
\end{array}
$}\label{HILBCON}%
\end{equation}
i.e., each of the members of the Hilbert basis of $\sigma_{0}$ has to be
either a \emph{\textquotedblleft}junior\emph{\textquotedblright} element or a
vertex of $\mathfrak{s}_{G}$.
\end{theorem}

\noindent\textsc{Proof}. The inclusion \textquotedblleft$\supseteq
$\textquotedblright\ is always true (without any further assumption about the
existence or non-existence of such a triangulation) and is obvious by the
definition of Hilbert basis. Now if there were an element $n\in\mathbf{Hlb}%
_{N_{G}}\left(  \sigma_{0}\right)  \smallsetminus\left(  \mathfrak{s}_{G}\cap
N_{G}\right)  $, then by Lemma \ref{Gorlem} this would be written as a
non-negative integer linear combination
\[
\ n=\lambda_{1}n_{1}+\cdots+\lambda_{r}n_{r}%
\]
of $r$ elements of $\mathfrak{s}_{G}\cap N_{G}$. Since $\mathbf{0\notin
Hlb}_{N_{G}}\left(  \sigma_{0}\right)  $, if there were at least one index
$\ j_{\bullet}\in\left\{  1,\ldots,r\right\}  $, for which $\lambda
_{j_{\bullet}}\neq0$. If $\lambda_{j_{\bullet}}=1$ and $\lambda_{j}=0$ for all
$j\in$ $\left\{  1,\ldots,r\right\}  \smallsetminus\left\{  \ j_{\bullet
}\right\}  $, then $n=n_{j_{\bullet}}\in\mathfrak{s}_{G}\cap N_{G}$ which
would contradict our assumption. But even the cases in which either
$\lambda_{j_{\bullet}}=1$ and some other $\lambda_{j}$'s were $\neq0$, or
$\lambda_{j_{\bullet}}\geq2$, would be exluded as impossible because of the
characterization (\ref{Hilbbasis}) of the Hilbert basis $\mathbf{Hlb}_{N_{G}%
}\left(  \sigma_{0}\right)  $ as the set of additively irreducible vectors of
$\sigma_{0}\cap\left(  N_{G}\smallsetminus\left\{  \mathbf{0}\right\}
\right)  $. Hence, $\mathbf{Hlb}_{N_{G}}\left(  \sigma_{0}\right)
\subseteq\mathfrak{s}_{G}\cap N_{G}$.\hfill{}$\square$

\begin{note}
\label{COUNTEREX}(i) For $r=2$ and $r=3,$ condition (\ref{HILBCON}) is
automatically satisfied.\smallskip\ \newline(ii) For $r\geq4$ there is a
plethora of AQS for which (\ref{HILBCON}) is violated. A simple example is the
(non-terminal) CQS $\left(  \mathbb{C}^{4}/G,\left[  \mathbf{0}\right]
\right)  $ of type $\frac{1}{7}\,\left(  1,1,2,3\right)  .$ This singularity
cannot have any crepant, $\mathbb{T}_{N_{G}}$-equivariant resolution, because
setting%
\[
\left\{
\begin{array}
[c]{ccc}%
n_{1}:=\frac{1}{7}\,\left(  1,1,2,3\right)  ^{\intercal}, & n_{2}:=\frac{1}%
{7}\,\left(  2,2,4,6\right)  ^{\intercal}, & n_{3}:=\frac{1}{7}\,\left(
3,3,6,2\right)  ^{\intercal},\\
&  & \\
n_{4}:=\frac{1}{7}\,\left(  4,4,1,5\right)  ^{\intercal}, & n_{5}:=\frac{1}%
{7}\,\left(  5,5,3,1\right)  ^{\intercal}, & n_{6}:=\frac{1}{7}\,\left(
6,6,5,4\right)  ^{\intercal},
\end{array}
\right.
\]
we get
\[
\mathfrak{s}_{G}\cap N_{G}=\{e_{1},e_{2},e_{3},e_{4},n_{1}\}\subsetneqq
\mathbf{Hlb}_{N_{G}}\left(  \sigma_{0}\right)  =\left\{
\begin{array}
[c]{c}%
e_{1},e_{2},e_{3},e_{4},\\
n_{1},n_{2},n_{3},n_{4},n_{5}%
\end{array}
\right\}  .
\]
(iii) In \S \ref{ONETWOPARSER} we shall present certain cyclic quotient
singularity series of arbitrary dimension for which condition (\ref{HILBCON})
turns out to be also \textit{sufficient}. Nevertheless, this is not true
\textit{in general} for $r\geq4.$ As it has been shown in \cite[\S 4.2, pp.
65-66]{Firla} and \cite[Ex. 10, p. 213]{Fir-Zi}, there are exactly
$10$\ four-dimensional Gorenstein cyclic quotient singularities with acting
group order $<100$ which fulfil (\ref{HILBCON}) and possess no crepant,
$\mathbb{T}_{N_{G}}$-equivariant resolutions$.$ Among them, the CQS with the
smallest possible acting group order is that one having the type\emph{\ }%
$\frac{1}{39}\,\left(  1,5,8,25\right)  .$
\end{note}

\section{Non-C.I.'s I: 1- and 2-Parameter CQS-Series\label{ONETWOPARSER}}

\noindent{}Asking whether non-c.i. Gorenstein cyclic quotient singularities of
given type (\ref{typeCQS}) can be resolved as desired, we begin with the
examination of those CQS whose junior simplex contains lattice points living
in a convenient geometric locus (in order to be able to keep track of how the
possible maximal lattice triangulations are built). More precisely, we
consider:\smallskip

\noindent{}$\bullet$ $1$\textit{-parameter CQS-series} $\left(  \mathbb{C}%
^{r}/G,\left[  \mathbf{0}\right]  \right)  ,$ for which the lattice points
belonging to $\left(  \mathfrak{s}_{G}\smallsetminus\left\{  e_{1}%
,\ldots,e_{r}\right\}  \right)  \cap N_{G}$ are \textit{collinear}, so that
the maximal lattice triangulations of the junior simplex $\mathfrak{s}_{G}$
are uniquely determined.\smallskip

\noindent{}$\bullet$ $2$\textit{-parameter CQS-series} $\left(  \mathbb{C}%
^{r}/G,\left[  \mathbf{0}\right]  \right)  ,$ for which the lattice points
belonging to $\left(  \mathfrak{s}_{G}\smallsetminus\left\{  e_{1}%
,\ldots,e_{r}\right\}  \right)  \cap N_{G}$ are \textit{coplanar}, so that
each of the simplices of the required maximal lattice triangulations of
$\mathfrak{s}_{G}$ is to be described as \textit{join} of a lattice polygon
(resp., a lattice segment) with an $\left(  r-4\right)  $-dimensional (resp.,
an $\left(  r-3\right)  $-dimensional) lattice simplex.\smallskip\

\noindent{}(For the somewhat lengthy proofs of Theorems \ref{1PARTHM},
\ref{DESEXC}, \ref{2-PARTHM}, and \ref{2-PARAMARITHM}, see \cite{DHH, DH}.)

\begin{theorem}
[1-parameter CQS]\label{1PARTHM}If $\left(  \mathbb{C}^{r}/G,\left[
\mathbf{0}\right]  \right)  $ is a Gorenstein CQS, such that $r-1$ weights in
its type are equal \emph{(}with $r\geq3),$ then it is isomorphic to the CQS of
type
\begin{equation}
\fbox{$%
\begin{array}
[c]{ccc}
& \dfrac{1}{l}\ (\underset{\left(  r-1\right)  \text{\emph{-times}}%
}{\underbrace{1,1,\ldots,1,1}},\ l-\left(  r-1\right)  ) &
\end{array}
$} \label{montype}%
\end{equation}
with $l=\left\vert G\right\vert \geq r.$ Moreover, we have\emph{:\smallskip
}\newline\emph{(i)} This msc-singularity is isolated if and only if
\emph{gcd}$\left(  l,r-1\right)  =1.\smallskip$\newline\emph{(ii)} There
exists a \emph{unique }maximal, coherent triangulation $\mathcal{T}$ of
$\mathfrak{s}_{G}$ w.r.t. $N_{G}$ inducing a unique crepant $\mathbb{T}%
_{N_{G}}$-equivariant partial projective desingularization \emph{(\ref{DESING}%
). }This is a (full) desingularization (i.e., $\mathcal{T}$ \ is basic w.r.t.
$N_{G}$) if and only if condition \emph{(\ref{HILBCON})} is satisfied. In
particular, \emph{(\ref{HILBCON})} is equivalent to the following:%
\begin{equation}
\fbox{$%
\begin{array}
[c]{lll}
& \text{\emph{\textbf{Either}}}\emph{\ \ }l\equiv0\text{ \emph{mod}}\left(
r-1\right)  \ \ \text{\emph{\textbf{or}\ }}\emph{\ }l\equiv1\text{ \emph{mod}%
}\left(  r-1\right)  . &
\end{array}
$} \label{rescon}%
\end{equation}

\end{theorem}

\begin{theorem}
[Exceptional prime divisors]\label{DESEXC}Suppose that $\left(  \mathbb{C}%
^{r}/G,\left[  \mathbf{0}\right]  \right)  ,$ $r\geq3$\emph{,} is a Gorenstein
CQS of type \smallskip\emph{(\ref{montype}). }If $l$ satisfies
\emph{(\ref{rescon}), }then\emph{\ }the exceptional locus of
\emph{(\ref{DESING})} consists of $%
\genfrac{\lfloor}{\rfloor}{}{}{l}{r-1}%
$ prime divisors $\left\{  D_{j}\left\vert \ 1\leq j\leq%
\genfrac{\lfloor}{\rfloor}{}{1}{l}{r-1}%
\right.  \right\}  $ on $X(N_{G},\widehat{\Delta}_{G}\left(  \mathcal{T}%
\right)  )$, having the following structure\emph{:}\newline%
\[
D_{j}\cong\mathbb{P}(\mathcal{O}_{\mathbb{P}_{\mathbb{C}}^{r-2}}%
\oplus\mathcal{O}_{\mathbb{P}_{\mathbb{C}}^{r-2}}\left(  l-\left(  r-1\right)
j\right)  )\text{ \ \ \ \smallskip\emph{(as}}\emph{\ }\mathbb{P}_{\mathbb{C}%
}^{1}\text{\emph{-bundles}}\emph{\ }\text{\emph{over}}\emph{\ }\mathbb{P}%
_{\mathbb{C}}^{r-2}\emph{)}%
\]
for all$\ j\in\left\{  1,2,\ldots,%
\genfrac{\lfloor}{\rfloor}{}{}{l}{r-1}%
-1\right\}  $, and%
\[
D_{%
\genfrac{\lfloor}{\rfloor}{}{1}{l}{r-1}%
}\cong\left\{
\begin{array}
[c]{lll}%
\mathbb{P}_{\mathbb{C}}^{r-1} & , & \text{\emph{if}\ \ \ }l\equiv
1\ \emph{mod}\left(  r-1\right)  ,\\
\mathbb{P}_{\mathbb{C}}^{r-2}\times\mathbb{C} & , & \text{\emph{if}
\ \ }l\equiv0\ \emph{mod}\left(  r-1\right)  .
\end{array}
\right.
\]

\end{theorem}

\begin{theorem}
[2-parameter CQS]\label{2-PARTHM}Let $\left(  \mathbb{C}^{r}/G,\left[
\mathbf{0}\right]  \right)  $ be a Gorenstein msc-CQS of type
\emph{(\ref{typeCQS})} with $l=\left\vert G\right\vert \geq r\geq3$\emph{, }
for which at least $r-2$ of its defining weights are equal. Then $X\left(
N_{G},\Delta_{G}\right)  =\mathbb{C}^{r}/G$ has crepant, $\mathbb{T}_{N_{G}}%
$-equivariant desingularizations $f_{\mathcal{T}}:X(N_{G},\widehat{\Delta}%
_{G}\left(  \mathcal{T}\right)  )\longrightarrow X\left(  N_{G},\Delta
_{G}\right)  $ if and only if \emph{(\ref{HILBCON}) }is satisfied. Moreover,
at least one of these desingularizations is projective.
\end{theorem}

\noindent{}Conditions equivalent to (\ref{HILBCON}) which can be
\textit{directly} expressed in terms of the defining weights\ occur in the
following case:

\begin{theorem}
[Arithmetic conditions for certain 2-parameter CQS]\label{2-PARAMARITHM}Let
$r$ be an integer $\geq3$ and $l$ an integer $\geq r.$ Write $l-\left(
r-2\right)  =\ a+b,$ where $a,b$ are integers $\geq1.$ Furthermore, set $t:=$
\emph{gcd}$\left(  b,l\right)  =$\emph{\ gcd}$\left(  a+\left(  r-2\right)
,l\right)  ,$ $t^{\prime}:=$\emph{\ gcd}$\left(  a,l\right)  ,$ and consider
$\nu_{1},\nu_{2}\in\mathbb{N}$, such that $\nu_{2}\left(  a+\left(
r-2\right)  \right)  -\nu_{1}l=t.$ Next, define
\[
\overline{p}:=\frac{\nu_{2}\cdot a-\nu_{1}\cdot l}{t^{\prime}},\ \ \ q:=\frac
{l}{t\cdot t^{\prime}},\ \ \ p:=\left[  \overline{p}\right]  _{q}\
\]
and write $q/p$ as regular continued fraction
\[
\frac{q}{p}=\lambda_{1}+\frac{1}{\lambda_{2}+\dfrac{1}{\lambda_{3}+%
\genfrac{.}{.}{0pt}{1}{{}}{\ddots}%
_{%
\begin{array}
[c]{c}%
+\dfrac{1}{\lambda_{\kappa}}%
\end{array}
}}}%
\]
with $\lambda_{i}\geq2$, $\ \forall i,\ 1\leq i\leq\kappa.$ Then, for the
Gorenstein CQS of type
\begin{equation}
\fbox{$%
\begin{array}
[c]{ccc}
& \dfrac{1}{l}\ (\underset{\left(  r-2\right)  \text{\emph{-times}}%
}{\underbrace{1,1,\ldots,1,1}},\ a,b), &
\end{array}
$}\label{2-type}%
\end{equation}
\emph{(\ref{HILBCON})} is equivalent to the following\emph{:}%
\begin{equation}
\fbox{$%
\begin{array}
[c]{ll}%
\text{\textbf{\emph{Either}}} & \text{\emph{gcd}}\left(  a,b,l\right)  =r-2,\\
& \\
\text{\emph{\textbf{or}}} & \left\{
\begin{array}
[c]{l}%
\text{\emph{gcd}}\left(  a,b,l\right)  =1,\ \ \left[  t\right]  _{r-2}=\left[
t^{\prime}\right]  _{r-2}=1,\ \\
\,\ \\
\dfrac{p-\overline{p}}{q}\equiv0\ \text{\emph{mod}}\left(  r-2\right)  ,\\
\\
\lambda_{i}\equiv0\ \text{\emph{mod}}\left(  r-2\right)  ,\forall
i,\ i\in\left(  \left[  2,\kappa-1\right]  \cap2\mathbb{Z}\right)  ,\\
\\
\text{\emph{and} }\lambda_{\kappa}\equiv1\ \text{\emph{mod}}\left(
r-2\right)  ,\text{\emph{\ whenever} }\kappa\text{ \emph{is even}}.
\end{array}
\right.
\end{array}
$}\label{CONDKETT}%
\end{equation}

\end{theorem}

\begin{remark}
\label{WHITEPOINT}(i) The method of building maximal
triangulations $\mathcal{T}$ of $\mathfrak{s}_{G}$ can be roughly
explained by means of Figure \ref{Fig.4} (in which $r=4$).

\begin{figure}[h]
\input{fig2.pstex_t}
\caption{}\label{Fig.4}
\end{figure}

\noindent We consider an arbitrary maximal (necessarily basic) triangulation
of the lattice polygon $\mathfrak{Q}_{G},$ and then we construct $\mathcal{T}$
$\ $by forming the joins of $e_{1}$ and $e_{2}$ with all of its triangles. The
white point belongs to $N_{G},$ and $\mathfrak{Q}_{G}$ itself becomes the
triangle having $e_{3},e_{4},$ and this point as its vertices, if and only if
the first of conditions (\ref{CONDKETT}) is satisfied. In this case, such a
maximal triangulation $\mathcal{T}$ of the entire $\mathfrak{s}_{G}$ is
automatically basic (w.r.t. $N_{G}$). If the white point \textit{does not}
belong to $N_{G},$ then the basicness of such a $\mathcal{T}$ $\ $amounts to
the second of conditions (\ref{CONDKETT}).\smallskip\ \newline(ii) If one of
the conditions (\ref{CONDKETT}) is satisfied, all compactly supported
exceptional prime divisors w.r.t. $f_{\mathcal{T}}$ $\ $\ are the total spaces
of fibrations having basis\emph{\ }$\mathbb{P}_{\mathbb{C}}^{r-3},$%
\emph{\ }and\emph{\ }typical fiber isomorphic either to\emph{\ }%
$\mathbb{P}_{\mathbb{C}}^{1}$\emph{\ }or to a\emph{\ }non-singular compact
toric surface (i.e., to a\emph{\ }$\mathbb{P}_{\mathbb{C}}^{2}$\emph{\ }or to
an\emph{\ }$\mathbb{F}_{\varkappa}=\mathbb{P}(\mathcal{O}_{\mathbb{P}%
_{\mathbb{C}}^{1}}\oplus\mathcal{O}_{\mathbb{P}_{\mathbb{C}}^{1}}\left(
\varkappa\right)  )$\emph{, }probably blown up at finitely many points, cf.
\cite[Thm. 1.28, p. 42]{Oda}).\smallskip\ \newline(iii) For the (rather
tricky) computation of the cohomology group dimensions (\ref{COHDIMFORMULA})
of the underlying space of any crepant desingularization of these cyclic
quotient singularities we refer to \cite[\S 7]{DHH}.
\end{remark}

\begin{examples}
(i) The subseries of non-isolated CQS with defining types%
\[
\tfrac{1}{\left(  \xi+\xi^{\prime}+1\right)  \cdot\left(  r-2\right)
}\,(\underset{\left(  r-2\right)  \text{\emph{-}times}}{\underbrace
{1,1,\ldots,1,1}},\xi\cdot\left(  r-2\right)  ,\xi^{\prime}\cdot\left(
r-2\right)  ),\ \ \xi\emph{,\ }\xi^{\prime}\emph{\ }\in\mathbb{N}%
\emph{,}\smallskip
\]
and\emph{\ }gcd$\left(  \xi,\xi^{\prime}\right)  =1$\emph{, }$r\geq4$\emph{,
}satisfies obviously the first of the conditions (\ref{CONDKETT}%
).\smallskip\ \newline(ii) The subseries of isolated CQS with defining
types\emph{\ }
\[
\tfrac{1}{2\left(  r-1\right)  ^{i}+r-2}\,(\underset{\left(  r-2\right)
\text{\emph{-}times}}{\underbrace{1,1,\ldots,1,1}},\left(  r-1\right)
^{i},\left(  r-1\right)  ^{i})\,
\]
and\emph{\ }$i\in\mathbb{N}$\emph{, }$r\geq4$\emph{, }satisfies
the second of the conditions (\ref{CONDKETT}). \smallskip
\newline
(iii) The example of $4$-dimensional subseries due to Mohri
\cite{Mohri}:%
\[
\tfrac{1}{4\,\xi}\,\left(  1,1,2\,\xi-1,2\,\xi-1\right)  ,\ \ \xi\in
\mathbb{N},
\]
satisfies the second of the conditions (\ref{CONDKETT}) and
contains only isolated singularities. Note that also the single suitably resolvable CQS of type $\tfrac{1}%
{11}\left(  1,1,3,6\right)  \ $\emph{\ }found in \cite{Mohri} belongs to the
subseries of isolated cyclic\emph{\ }quotient singularities of type
\[
\tfrac{1}{4r-5}\underset{\left(  r-2\right)  \text{\emph{-}times}%
}{(\underbrace{1,1,\ldots,1,1}},r-1,2r-2)
\]
satisfying obviously the second of the conditions (\ref{CONDKETT}). Moreover,
there are examples like\emph{\ }$1/28\left(  1,1,1,4,21\right)  $ for
which\emph{\ }$p=0,$ $q=1.$
\end{examples}

\section{Non-C.I.'s II: The GP-Singularity Series\label{GPSERIES}}

\noindent{}Another Gorenstein non-c.i. cyclic quotient singularity series of
particular interest, admitting the required resolutions, is the so-called
\textit{geometric progress singularity series} (GPSS($r;k$), for short, with
type (\ref{kkk})). The purpose of this section is to give a proof of the
following Theorem (appearing as Conjecture 10.2 in \cite{DH}):

\begin{theorem}
\label{GPSSTHM}All Gorenstein CQS $(\mathbb{C}^{r}/G,[\mathbf{0}])$ of type
\begin{equation}
\fbox{$%
\begin{array}
[c]{ccc}
&  & \\
& \dfrac{1}{\left(  \dfrac{k^{r}-1}{k-1}\right)  }\ \,\left(  1,k,k^{2}%
,k^{3},\ldots,k^{r-2},k^{r-1}\right)  & \\
&  &
\end{array}
$} \label{kkk}%
\end{equation}
admit $\mathbb{T}_{N_{G}}$-equivariant projective, crepant resolutions for
\textbf{all }$r\geq3$ and \textbf{all} $k\geq2.$ In particular, for $k=2,$
there is a unique resolution of this sort.
\end{theorem}

\begin{remark}
\label{REMGPSS1}(i) Setting $l:=%
{\textstyle\sum\limits_{i=0}^{r-1}} k^{i}=\tfrac{k^{r}-1}{k-1}$
for the order of $G$ acting on $\mathbb{C}^{r}$, we see that
\[
N_{G}=\mathbb{Z}^{r}+\tfrac{1}{l}\left(  1,k,k^{2},\ldots,k^{r-1}\right)
^{\intercal},\text{ with }\det(N_{G})=\tfrac{1}{l}.
\]
(ii) Since gcd$(k^{i},l)=1,$ $\forall i\in\{0,\ldots,r-1\},$ all members of
the GPSS($r;k$) are \textit{isolated} (and therefore \textit{msc}-)
singularities (cf. Proposition \ref{isol2} and Corollary \ref{ISMSC}).
\end{remark}

\begin{lemma}
\label{LGPSS1}If we denote by $W\left(  r;k\right)  =(w_{ij})_{1\leq i,j\leq
r} $ the $(r\times r)$-matrix with%
\[
w_{ij}:=\left[  k^{i-1}\cdot k^{j-1}\right]  _{l}=\left[  k^{i+j-2}\right]
_{l}%
\]
as its entries, then
\begin{equation}
\left\vert \det(W\left(  r;k\right)  )\right\vert =\left(  k^{r}-1\right)
^{r-1}=l^{r-1}\left(  k-1\right)  ^{r-1}. \label{DETW}%
\end{equation}

\end{lemma}

\noindent{}\textsc{Proof}. Since%
\[
k^{r}=\left(  k-1\right)  l+k^{0}\Longrightarrow k^{r+\beta}=k^{\beta}\left(
k-1\right)  l+k^{\beta},\ \ \forall\beta\in\mathbb{Z}_{\geq0},
\]
we have $w_{ij}=\left[  k^{i+j-2}\right]  _{l}=k^{\left[  i+j-2\right]  _{r}%
}.$ On the other hand, performing the elementary operations%
\[
\left(  i\text{-th row}\right)  \leadsto\left(  i\text{-th row}\right)
-k^{(i-1)}\cdot\left(  \text{first row}\right)  ,\ \forall i\in\{2,3,\ldots
,r\},
\]
we transfer $W\left(  r;k\right)  $ into a matrix of the form
\[
\left(
\begin{array}
[c]{ccccc}%
1 & k & k^{2} & \cdots & k^{r-1}\\
0 & 0 & 0 & \cdots & -k^{r}+1\\
0 & 0 & \cdots & -k^{r}+1 & \star\\
\vdots & \vdots & \rotatebox{70}{$\ddots$} & \star & \star\\
\vdots & \rotatebox{70}{$\ddots$} & \star & \vdots & \vdots\\
-k^{r}+1 & \star & \cdots & \cdots & \star
\end{array}
\right)  .
\]
Hence, $\left\vert \det(W\left(  r;k\right)  )\right\vert $ is given by the
formula (\ref{DETW}).\hfill{}$\square$

\begin{lemma}
\label{LGPSS2}\emph{(i)} Setting $\breve{W}\left(  r;k\right)  :=\left(
\frac{1}{l(k-1)}W\left(  r;k\right)  \right)  ^{-1},$ we have%
\begin{equation}
\breve{W}\left(  r;k\right)  =\left(
\begin{array}
[c]{ccccc}%
-1 & 0 & \cdots & 0 & k\\
0 & 0 & \cdots & k & -1\\
\vdots & 0 & k & -1 & 0\\
0 & \rotatebox{70}{$\ddots$} & \rotatebox{70}{$\ddots$} & \cdots & \vdots\\
k & -1 & 0 & \cdots & 0
\end{array}
\right)  , \label{WHATINVERSE}%
\end{equation}
with%
\begin{equation}
\left\vert \det(\breve{W}\left(  r;k\right)  )\right\vert =l(k-1).
\label{DETWHAT}%
\end{equation}
\emph{(ii) }The regular linear transformation $\Phi:\mathbb{R}^{r}%
\longrightarrow\mathbb{R}^{r},$ with
\begin{equation}
\Phi\left(  \mathbf{x}\right)  :=\breve{W}\left(  r;k\right)  \mathbf{x,\ \ }%
\forall\mathbf{x}\in\mathbb{R}^{r}, \label{PhiW}%
\end{equation}
maps $N_{G}$ onto%
\begin{equation}
\breve{N}_{G}=\left\{  \left(
\lambda_{1},\ldots,\lambda_{r}\right)
^{\intercal}\in\mathbb{Z}^{r}\ \left\vert \
{\textstyle\sum\limits_{i=1}^{r}} \right.
\lambda_{i}\equiv0(\text{\emph{mod} }(k-1))\right\}  ,
\label{EQUALLATTICES}%
\end{equation}
and the junior simplex $\mathfrak{s}_{G}$ onto
\begin{equation}
\fbox{$%
\begin{array}
[c]{lll}
& \mathfrak{\breve{s}}_{G}=\text{\emph{conv}}(\{\left.  \breve{w}%
_{j}\ \right\vert \ 1\leq j\leq r\}), &
\end{array}
$} \label{COLUMNSIMPLEX}%
\end{equation}
where
\[
\breve{w}_{j}:=\left\{
\begin{array}
[c]{ll}%
-e_{1}+ke_{r}, & \text{\emph{if} }j=1,\text{ }\\
-e_{j}+ke_{j-1}, & \text{\emph{if} }j\in\{2,\ldots,r\},
\end{array}
\right.
\]
denotes the $j$-th column vector of the matrix obtained by
$\breve{W}\left( r;k\right)  $ by interchanging its $j$-th with
its $(r+2-j)$-th column for all $j\in\{2,\ldots,r\}$. \emph{(This
rearrangement of the} \emph{index set for enumerating the vertices
of} $\mathfrak{\breve{s}}_{G}$ \emph{will turn out to be
convenient in the subsequent Lemmas.)}
\end{lemma}

\noindent{}\textsc{Proof}. (i) Let $w_{i}$ denote the $i$-th row of $W\left(
r;k\right)  .$ For all $i\in\{1,\ldots,r\}$ we have%
\[
k\cdot w_{ij}-w_{i+1\ j}=k\cdot k^{\left[  i+j-2\right]  _{r}}-k^{\left[
i+j-1\right]  _{r}}=\left\{
\begin{array}
[c]{ll}%
0, &
\begin{array}
[c]{l}%
\text{if }i+j-1\leq r-1\text{ }\\
\text{or }i+j-2\geq r,
\end{array}
\\
k^{r}-1, & \text{if }i+j-1=r.
\end{array}
\right.
\]
Thus, $\breve{W}\left(  r;k\right)  $ is the matrix (\ref{WHATINVERSE})
because
\[
k\cdot w_{i}-w_{i+1}=\left(  k^{r}-1\right)  e_{r+1-i},\ \forall
i\in\{1,\ldots,r-1\},\text{ and }k\cdot w_{r}-w_{1}=\left(  k^{r}-1\right)
e_{1},
\]
and (\ref{DETWHAT}) follows directly from (\ref{DETW}).\smallskip
\ \newline(ii) By definition, the determinant of $\breve{N}_{G}=\Phi(N_{G})$
equals $k-1,$ and
\[
\breve{N}_{G}=\breve{W}\left(  r;k\right)  \mathbb{Z}^{r}+\mathbb{Z}%
(k-1,0,0,\ldots,0,0)^{\intercal}.
\]
$\breve{N}_{G}$ is included into $\left\{  \left(
\lambda_{1},\ldots ,\lambda_{r}\right)
^{\intercal}\in\mathbb{Z}^{r}\ \left\vert \
{\textstyle\sum\limits_{i=1}^{r}} \right.
\lambda_{i}\equiv0(\text{mod }(k-1))\right\}  .$ But also this
lattice has determinant $k-1,$ leading to equality
(\ref{EQUALLATTICES}). (\ref{COLUMNSIMPLEX}) is
obvious.\hfill{}$\square$

\begin{lemma}
\label{LGPSS3}We have
\begin{equation}
\mathfrak{\breve{s}}_{G}\cap\breve{N}_{G}=\left\{  \left(  \lambda_{1}%
,\ldots,\lambda_{r}\right)  ^{\intercal}\in\mathbb{Z}_{\geq0}^{r}\
\left\vert \ {\textstyle\sum\limits_{i=1}^{r}} \right.
\lambda_{i}=k-1\right\}  \cup\{\left.  \breve{w}_{j}\ \right\vert
\ 1\leq j\leq r\}. \label{LATPTGPSS1}%
\end{equation}
In particular,
\begin{equation}
\sharp(\mathfrak{\breve{s}}_{G}\cap\breve{N}_{G})=\tbinom{k+r-2}{r-1}+r.
\label{LATPTGPSS2}%
\end{equation}

\end{lemma}

\noindent{}\textsc{Proof}. Since $\mathfrak{\breve{s}}_{G}\subseteq\left\{
\mathbf{x}=\left(  x_{1},\ldots,x_{r}\right)  ^{\intercal}\in\mathbb{R}%
^{r}\ \left\vert \ {\textstyle\sum\limits_{j=1}^{r}}
\right.  x_{j}=k-1\right\}  ,$ we have%
\[
\mathfrak{\breve{s}}_{G}\cap\breve{N}_{G}=\left\{  \mathbf{x}=\left(
x_{1},\ldots,x_{r}\right)  ^{\intercal}\in\mathbb{Z}^{r}\ \left\vert
\ \mathbf{x}\in\text{conv}(\{\left.  \breve{w}_{j}\ \right\vert \ 1\leq j\leq
r\})\right.  \right\}  .
\]
We first observe that $(k-1)e_{j}\in\mathfrak{\breve{s}}_{G}\cap\breve{N}%
_{G},$ for all $j\in\{1,\ldots,r\}.$ (For instance, $(k-1)e_{1}=%
{\textstyle\sum\limits_{j=1}^{r}} \frac{k^{\left(  j-1\right)
}}{l}\breve{w}_{j}.$ The other inclusions follow by symmetry.)
Next, we consider an
\[
\mathbf{x}=\left(  x_{1},\ldots,x_{r}\right)  ^{\intercal}\in\text{conv}%
(\{\left.  \breve{w}_{j}\ \right\vert \ 1\leq j\leq r\})\cap\left(  \breve
{N}_{G}\mathbb{r}\{\left.  \breve{w}_{j}\ \right\vert \ 1\leq j\leq
r\}\right)  .
\]
This can be written as linear combination
\[
\mathbf{x}=\left(  x_{1},\ldots,x_{r}\right)  ^{\intercal}=%
{\textstyle\sum\limits_{j=1}^{r}} \eta_{j}\breve{w}_{j},\text{ for
suitable }\eta_{j}\text{'s}\in\lbrack0,1).
\]
Since $\mathbf{x}\in\mathbb{Z}^{r}$ we have $x_{j}\geq0,$ for all
$j\in\{1,\ldots,r\}.$ Hence, $\mathbf{x}$ belongs to
\[
\left\{  \left(  x_{1},\ldots,x_{r}\right)
^{\intercal}\in\mathbb{R}_{\geq 0}^{r}\left\vert
{\textstyle\sum\limits_{j=1}^{r}} x_{j}=k-1\right.  \right\}
=\text{ conv}(\{\left.  (k-1)e_{j}\ \right\vert \ 1\leq j\leq
r\}),
\]
and both equalities (\ref{LATPTGPSS1}) and (\ref{LATPTGPSS2}) are
true.\hfill{}$\square$

\begin{lemma}
\label{LGPSS4}Let $\mathbf{s}(\varepsilon_{1},\ldots,\varepsilon_{r}%
)\subset\{\left.  \left(  x_{1},\ldots,x_{r}\right)
^{\intercal}\in \mathbb{R}^{r}\right\vert
{\textstyle\sum\nolimits_{j=1}^{r}}
x_{j}=k-1\}$ denote the simplices%
\[
\mathbf{s}(\varepsilon_{1},\ldots,\varepsilon_{r}):=\text{ \emph{conv}%
}(\left\{  \left.  \varepsilon_{j}u_{j}+\left(  1-\varepsilon_{j}\right)
\breve{w}_{j}\ \right\vert \ 1\leq j\leq r\right\}  ),
\]
where $\varepsilon_{j}\in\{0,1\}$ and $u_{j}:=\left(  k-1\right)
e_{j},\ \forall j\in\{1,\ldots,r\}.$ Then there is a unique triangulation
$\mathfrak{T}(r;k)$ of $\mathfrak{\breve{s}}_{G}$ having $\left\{  \left.
u_{j},\breve{w}_{j}\right\vert 1\leq j\leq r\right\}  $ as its vertex set,
namely%
\begin{equation}
\fbox{$\mathfrak{T}(r;k)=\left\{  \left.  \mathbf{s}(\varepsilon_{1}%
,\ldots,\varepsilon_{r})\ \right\vert \ (\varepsilon_{1},\ldots,\varepsilon
_{r})\in\{0,1\}^{r}\mathbb{r}\{(0,0,\ldots,0,0)\}\right\}  .$}
\label{TRIANGTAU}%
\end{equation}

\end{lemma}

\noindent{}\textsc{Proof}. First note that $\mathbf{s}(0,0,\ldots
,0,0)=\mathfrak{\breve{s}}_{G}.$ Let $T$ be an arbitrary triangulation of the
convex hull of the point set $\left\{  \left.  u_{j},\breve{w}_{j}\right\vert
1\leq j\leq r\right\}  $. If $\mathbf{t}$ is an $(r-1)$-dimensional simplex
belonging to $T,$ then there is no index $j\in\{1,\ldots,r\}$ such that
$\left\{  u_{j},\breve{w}_{j}\right\}  \subset\mathbf{t}.$ Assuming, in the
contrary direction, the existence of such an index $j,$ we would have
$\frac{1}{k}u_{j}+\frac{k-1}{k}\breve{w}_{j}\in\mathbf{t},$ with
\[
\frac{1}{k}u_{j}+\frac{k-1}{k}\breve{w}_{j}=\left\{
\begin{array}
[c]{ll}%
u_{r}, & \text{if }j=1,\text{ }\\
u_{j-1}, & \text{if }j\in\{2,\ldots,r\},
\end{array}
\right.
\]
which would be absurd (because $\mathbf{t}$ could not be an $(r-1)$%
-dimensional simplex of a triangulation). Hence, any triangulation of the
convex hull of $\left\{  \left.  u_{j},\breve{w}_{j}\right\vert 1\leq j\leq
r\right\}  $ must have $\left\{  \left.  \mathbf{s}(\varepsilon_{1}%
,\ldots,\varepsilon_{r})\ \right\vert \ (\varepsilon_{1},\ldots,\varepsilon
_{r})\in\{0,1\}^{r}\mathbb{r}\{(0,0,\ldots,0,0)\}\right\}  $ as maximal
dimensional simplices. In fact, it is easy to verify that the intersection of
any two simplices of this sort is either a face of both or the empty set, and
that
\[
\text{Vol}(\mathbf{s}(\varepsilon_{1},\ldots,\varepsilon_{r}))=\left\{
\begin{array}
[c]{ll}%
\frac{\sqrt{r}}{\left(  r-1\right)  !}\frac{\left\vert \det\left(  \breve
{W}\left(  r;k\right)  \right)  \right\vert }{\det\left(  \breve{N}%
_{G}\right)  }\overset{\text{(\ref{DETWHAT})}}{=}\frac{\sqrt{r}}{\left(
r-1\right)  !}l, & \text{if }(\varepsilon_{1},\ldots,\varepsilon
_{r})=(0,\ldots,0),\text{ }\\
\  & \\
\frac{1}{k-1}\frac{\sqrt{r}}{\left(  r-1\right)  !}\ \left(  k-1\right)  ^{%
{\textstyle\sum\nolimits_{j=1}^{r}}
\varepsilon_{j}}, & \text{if }(\varepsilon_{1},\ldots,\varepsilon_{r}%
)\neq(0,\ldots,0),
\end{array}
\right.
\]
(cf. formula (\ref{vol-reg-s})). Since%
\begin{align*}
\sum_{(\varepsilon_{1},\ldots,\varepsilon_{r})\in\{0,1\}^{r}\mathbb{r}%
\{(0,0,\ldots,0,0)}\text{Vol}(\mathbf{s}(\varepsilon_{1},\ldots,\varepsilon
_{r}))  &  =\tfrac{1}{k-1}\tfrac{\sqrt{r}}{\left(  r-1\right)  !}\left(
\sum_{\rho=1}^{r}\tbinom{r}{\rho}\left(  k-1\right)  ^{\rho}\right)
\smallskip\\
&  =\tfrac{1}{k-1}\tfrac{\sqrt{r}}{\left(  r-1\right)  !}\left(
k^{r}-1\right)  \smallskip\\
&  =\tfrac{\sqrt{r}}{\left(  r-1\right)  !}l=\text{Vol}(\mathfrak{\breve{s}%
}_{G}),
\end{align*}
the support of $\mathfrak{T}(r;k)$ given in (\ref{TRIANGTAU}) equals
$\mathfrak{\breve{s}}_{G},$ and $\mathfrak{T}(r;k)$ is therefore the
\textit{unique} triangulation of $\mathfrak{\breve{s}}_{G}$ having $\left\{
\left.  u_{j},\breve{w}_{j}\right\vert 1\leq j\leq r\right\}  $ as vertex
set.\hfill{}$\square$

\begin{definition}
Let $\widetilde{\Phi}:\mathbb{R}^{r}\longrightarrow\mathbb{R}^{r}$ be the
unimodular transformation%
\[
\widetilde{\Phi}\left(  \mathbf{x}\right)  :=L\,\mathbf{x,\ \ }\forall
\mathbf{x}\in\mathbb{R}^{r},
\]
where%
\[
L:=\left(
\begin{array}
[c]{ccccc}%
1 & 0 & 0 & \cdots & 0\\
0 & 1 & 0 & \cdots & 0\\
\vdots & \vdots & \ddots & \cdots & \vdots\\
0 & \cdots & 0 & 1 & 0\\
1 & 1 & \cdots & 1 & 1
\end{array}
\right)  .
\]
Then $\widetilde{\Phi}$ maps the hyperplane $\{\left.
\mathbf{x}\in \mathbb{R}^{r}\right\vert
{\textstyle\sum\nolimits_{j=1}^{r}} x_{j}=k-1\}$ onto $\left\{
\mathbf{x}\in\mathbb{R}^{r}\left\vert x_{r}=k-1\right.  \right\}
,$ with
\[
\widetilde{\Phi}(u_{j})=\left\{
\begin{array}
[c]{ll}%
(k-1)e_{i}+(k-1)e_{r}, & \text{if }j\in\{1,\ldots,r-1\},\text{ }\\
(k-1)e_{r}\ (=u_{r}), & \text{if }j=r,
\end{array}
\right.
\]
and%
\[
\widetilde{\Phi}(\breve{w}_{j})=\left\{
\begin{array}
[c]{ll}%
-e_{1}+(k-1)e_{r}, & \text{if }j=1,\text{ }\\
-e_{j}+ke_{j-1}+(k-1)e_{r}, & \text{if }j\in\{2,\ldots,r-1\}.\\
ke_{r-1}+(k-1)e_{r}, & \text{if }j=r.
\end{array}
\right.
\]
For all $j\in\{1,\ldots,r\}$ let us use the abbreviations%
\[
\widetilde{u}_{j}:=\widetilde{\Phi}(u_{j}),\ \ \widetilde{w}_{j}%
:=\widetilde{\Phi}(\breve{w}_{j}),\ \ \text{and}\ \ \overline{w}%
_{j}:=\widetilde{w}_{j}-(k-1)e_{r}.
\]
We observe that $\widetilde{\Phi}$ transfers the triangulation $\mathfrak{T}%
(r;k)$ of $\mathfrak{\breve{s}}_{G}$ onto the triangulation%
\[
\widetilde{\Phi}\left(  \mathfrak{T}(r;k)\right)  =\left\{  \left.
\widetilde{\mathbf{s}}(\varepsilon_{1},\ldots,\varepsilon_{r})\ \right\vert
\ (\varepsilon_{1},\ldots,\varepsilon_{r})\in\{0,1\}^{r}\mathbb{r}%
\{(0,0,\ldots,0,0)\}\right\}
\]
of $\widetilde{\mathfrak{s}}_{G}=\widetilde{\Phi}\left(  \mathfrak{\breve{s}%
}_{G}\right)  ,$ where
\[
\widetilde{\mathbf{s}}(\varepsilon_{1},\ldots,\varepsilon_{r}):=\widetilde
{\Phi}(\mathbf{s}(\varepsilon_{1},\ldots,\varepsilon_{r}))=\text{conv}%
(\left\{  \left.  \varepsilon_{j}\widetilde{u}_{j}+\left(  1-\varepsilon
_{j}\right)  \widetilde{w}_{j}\ \right\vert \ 1\leq j\leq r\right\}  ),
\]
and we define the $r$-\textit{dimensional} (!) lattice polytope
\[
\fbox{$%
\begin{array}
[c]{lll}
& P\left(  r;k\right)  :=\text{ conv}(\{\{\widetilde{u}_{j},\overline{w}%
_{j}\ \left\vert \ 1\leq j\leq r\right.  \})\subset\mathbb{R}^{r}. &
\end{array}
$}%
\]

\end{definition}

\begin{lemma}
\label{LGPSS5}The facets of $P\left(  r;k\right)  $ are exactly those
belonging to the set%
\[
\left\{  \left.  \overline{\mathbf{s}}(\varepsilon_{1},\ldots,\varepsilon
_{r})\ \right\vert \ (\varepsilon_{1},\ldots,\varepsilon_{r})\in
\{0,1\}^{r}\right\}  ,
\]
where%
\[
\overline{\mathbf{s}}(\varepsilon_{1},\ldots,\varepsilon_{r}%
):=\text{\emph{conv}}(\left\{  \left.  \varepsilon_{j}\widetilde{u}%
_{j}+\left(  1-\varepsilon_{j}\right)  \overline{w}_{j}\ \right\vert \ 1\leq
j\leq r\right\}  ).
\]

\end{lemma}

\noindent{}\textsc{Proof}. Since $%
{\textstyle\sum\nolimits_{j=1}^{r}} k^{r-j}\overline{w}_{j}=0,$
the origin $\mathbf{0}$ is an interior point of $P\left(
r;k\right)  ,$ and we may assume that the coordinates of each
point $\mathbf{x}=\left(  x_{1},\ldots,x_{r}\right)
^{\intercal}\in P\left( r;k\right)  $ satisfy inequalities of the
form
\begin{equation}
\sum_{j=1}^{r}\eta_{j}x_{j}\leq k-1,\text{ for suitable }r\text{-tuples }%
(\eta_{1},\ldots,\eta_{r})\in\mathbb{R}^{r}. \label{INEQUALP}%
\end{equation}
Since $\widetilde{u}_{j}\in P\left(  r;k\right)  $ (resp., $\overline{w}%
_{j}\in P\left(  r;k\right)  $) $\forall j\in\{1,\ldots,r\},$ valid
inequalities for $P\left(  r;k\right)  $ of the form (\ref{INEQUALP}) must
satisfy%
\begin{equation}
\left\{
\begin{array}
[c]{l}%
\eta_{j}+\eta_{r}\leq1,\forall j\in\{1,\ldots,r-1\},\\
\eta_{r}\leq1,
\end{array}
\right.  \label{INEQ1}%
\end{equation}
and%
\begin{align}
-\eta_{1}  &  \leq k-1,\label{INEQ2}\\
k\eta_{j-1}-\eta_{j}  &  \leq k-1,\ \ \forall j\in\{2,\ldots
,r-1\},\label{INEQ3}\\
k\eta_{r-1}  &  \leq k-1, \label{INEQ4}%
\end{align}
respectively. Let $F$ be a facet of $P\left(  r;k\right)  .$ Assume that the
supporting hyperplane of $F$ is described by an equation $\sum_{j=1}^{r}%
\eta_{j}x_{j}=k-1.$ If there were an index, say $j\in\{2,\ldots,r-1\},$ such
that both $\widetilde{u}_{j}$ and $\overline{w}_{j}$ belong to $F,$ then we
would have%
\begin{equation}
\left.
\begin{array}
[c]{l}%
\eta_{j}+\eta_{r}=1\\
k\eta_{j-1}-\eta_{j}=k-1
\end{array}
\right\}  \Longrightarrow k\eta_{j-1}+\eta_{r}=k. \label{EQ5}%
\end{equation}
Since (\ref{EQ5}) is also valid for $P\left(  r;k\right)  $ itself, all
inequalities (\ref{INEQ1}), (\ref{INEQ2}), (\ref{INEQ3}), (\ref{INEQ4}) would
be satisfied. Hence,%
\[
\left.
\begin{array}
[c]{r}%
\eta_{j-1}+\eta_{r}\leq1\overset{\text{(\ref{EQ5})\smallskip}}{\Longrightarrow
}\eta_{r}\leq0\\
\eta_{j}+\eta_{r}=1
\end{array}
\right\}  \Longrightarrow\eta_{j}\geq1\overset{\text{(\ref{INEQ3})\smallskip}%
}{\Longrightarrow}\eta_{j+1}\geq k\eta_{j}-(k-1)\geq1,
\]
and using the same argument,
\[
\eta_{j+1}\geq1\Longrightarrow\eta_{j+2}\geq1\Longrightarrow\cdots
\Longrightarrow\eta_{r-1}\geq1.
\]
On the other hand, (\ref{INEQ4}) would give $\eta_{r-1}\leq\frac{k-1}{k}<1,$
leading to contradiction. Analogously, by (\ref{INEQ1}), (\ref{INEQ2}),
(\ref{INEQ3}) and (\ref{INEQ4}) one shows that $\widetilde{u}_{j}$ and
$\overline{w}_{j}$ cannot simultaneously belong to $F$ even if $j=1$ or $j=r.$
Thus, $F$ is necessarily of the form%
\[
F=\overline{\mathbf{s}}(\varepsilon_{1},\ldots,\varepsilon_{r}),\text{ for a
suitable }r\text{-tuple }(\varepsilon_{1},\ldots,\varepsilon_{r}%
)\in\{0,1\}^{r}.
\]
It remains to prove that \textit{all} $\overline{\mathbf{s}}(\varepsilon
_{1},\ldots,\varepsilon_{r})$'s are realized as facets of $P\left(
r;k\right)  .$ If $(\varepsilon_{1},\ldots,\varepsilon_{r})$ equals
$(0,0,\ldots,0,0)$ (resp., $(1,1,\ldots,1,1)$), then we get obviously the
\textit{bottom} (resp., the \textit{top}) \textit{facet} of $P\left(
r;k\right)  .$ Next, choose an arbitrary simplex
\[
\overline{\mathbf{s}}(\varepsilon_{1},\ldots,\varepsilon_{r})\text{ with
}(\varepsilon_{1},\ldots,\varepsilon_{r})\in\{0,1\}^{r}\mathbb{r}%
\{(0,0,\ldots,0,0),(1,1,\ldots,1,1)\}
\]
and assume without loss of generality (i.e., up to a permutation of
coordinates) that
\[
\varepsilon_{j}=1,\ \ \forall j\in\{1,\ldots,\rho\},\text{ \ and
\ }\varepsilon_{j}=0,\ \ \forall j\in\{\rho+1,\ldots,r\},\text{ }%
\]
for some $\rho\in\{1,\ldots,r-1\}.$ Defining
\[
\eta_{j}:=\left\{
\begin{array}
[c]{ll}%
\tfrac{k^{r-\rho}-1}{k^{r-\rho}}, & \text{if }j\in\{1,\ldots,\rho-1\},\\
& \\
\frac{1}{k^{r-\rho}}, & \text{if }j=\rho,\\
& \\
\frac{k^{r-j}-1}{k^{r-j}}, & \text{if }j\in\{\rho+1,\ldots,r\},
\end{array}
\right.
\]
one checks easily that $\eta_{j}$'s fulfil (\ref{INEQ1}), (\ref{INEQ2}),
(\ref{INEQ3}) and (\ref{INEQ4}), and that the $r$ affinely independent points
\[
\left\{  \left.  (k-1)e_{j}\!+\!(k-1)e_{r}\right\vert 1\leq j\leq\rho\right\}
,\left\{  \left.  -e_{j}\!+\!(k-1)e_{j-1}\right\vert \rho+1\leq j\leq
r-1\right\}  ,\left\{  ke_{r-1}\right\}  ,
\]
satisfy the equality%
\[
\sum_{j=1}^{r}\eta_{j}x_{j}=k-1.
\]
Hence, their convex hull $\overline{\mathbf{s}}(\varepsilon_{1},..,\varepsilon
_{r})$ constitutes a facet of $P\left(  r;k\right)  ,$ as asserted.\hfill
{}$\square$

\begin{corollary}
\label{CORCOH}The triangulation $\mathfrak{T}(r;k)$ of $\mathfrak{\breve{s}%
}_{G}$ is coherent.
\end{corollary}

\noindent{}\textsc{Proof}. Using Lemma \ref{LGPSS5} and the projection
$\varpi:\mathbb{R}^{r}\longrightarrow\left\{  \mathbf{x}\in\mathbb{R}%
^{r}\ \left\vert \ x_{r}=k-1\right.  \right\}  ,$ with%
\[
\varpi(\mathbf{x}):=\left(  x_{1},x_{2},\ldots,x_{r-1},k-1\right)
^{\intercal},\ \ \forall\mathbf{x}=\left(  x_{1},\ldots,x_{r}\right)
^{\intercal}\in\mathbb{R}^{r},
\]
we see that the set $\left\{  \left.  \overline{\mathbf{s}}(\varepsilon
_{1},\ldots,\varepsilon_{r})\ \right\vert \ (\varepsilon_{1},\ldots
,\varepsilon_{r})\in\{0,1\}^{r}\mathbb{r}\{(0,0,\ldots,0,0)\}\right\}  $
consisting of the\ facets of $P\left(  r;k\right)  $ which belong to its
\textquotedblleft higher envelope\textquotedblright\ is mapped via $\varpi$
onto the triangulation
\[
\varpi\left(  \left\{  \left.  \overline{\mathbf{s}}(\varepsilon_{1}%
,\ldots,\varepsilon_{r})\ \right\vert \ (\varepsilon_{1},\ldots,\varepsilon
_{r})\in\{0,1\}^{r}\mathbb{r}\{(0,0,\ldots,0,0)\}\right\}  \right)
=\widetilde{\Phi}\left(  \mathfrak{T}(r;k)\right)
\]
of $\widetilde{\mathfrak{s}}_{G}$ because%
\[
\varpi\left(  \overline{\mathbf{s}}(\varepsilon_{1},\ldots,\varepsilon
_{r})\right)  =\widetilde{\mathbf{s}}(\varepsilon_{1},\ldots,\varepsilon
_{r}),\ \ \forall(\varepsilon_{1},\ldots,\varepsilon_{r})\in\{0,1\}^{r}%
\mathbb{r}\{(0,0,\ldots,0,0)\}.
\]
Hence, the $\widetilde{\Phi}\left(  \mathfrak{T}(r;k)\right)  $-support
function $\theta:\widetilde{\mathfrak{s}}_{G}\longrightarrow(0,k-1]\subset
\mathbb{R}$ defined by the formula
\[
\theta\left(  \mathbf{x}\right)  :=\max\{\left.  t\in\left[  0,k-1\right]
\mathbb{\ }\right\vert \ \left(  x_{1},x_{2},\ldots,x_{r-1},t\right)
^{\intercal}\in P\left(  r;k\right)  \},
\]
for all $\mathbf{x}=\left(  x_{1},x_{2},\ldots,x_{r-1},k-1\right)
^{\intercal}\in\widetilde{\mathfrak{s}}_{G},$ is strictly upper convex. This
means that%
\[
\theta\circ\left.  \widetilde{\Phi}\right\vert _{\mathfrak{\breve{s}}_{G}%
}:\mathfrak{\breve{s}}_{G}\longrightarrow(0,k-1]\subset\mathbb{R}%
\]
is a strictly upper convex $\mathfrak{T}(r;k)$-support function on
$\mathfrak{\breve{s}}_{G}.$\hfill{}$\square$

\begin{remark}
An alternative proof of Corollary \ref{CORCOH} can be obtained by observing
that $\mathfrak{T}(r;k)$ is actually a \textit{lexicographic triangulation},
and by using the fact that lexicographic triangulations are coherent (see Lee
\cite{Lee}).
\end{remark}

\begin{example}
\label{EXGP}The unique, coherent triangulation $\Phi^{-1}(\mathfrak{T}(3;4))$
(with $\Phi$ as defined in (\ref{PhiW})) of the original junior simplex
$\mathfrak{s}_{G}$ (w.r.t. $N_{G}$) which is induced by $\mathfrak{T}(3;4)$
(for $r=3,$ $k=4$), and has $\left\{  \left.  \Phi^{-1}(u_{j}),e_{j}%
\right\vert 1\leq j\leq3\right\}  $ as its vertex set, is depicted
in Figure \ref{Fig.5}.

\begin{figure}[h]
\epsfig{file=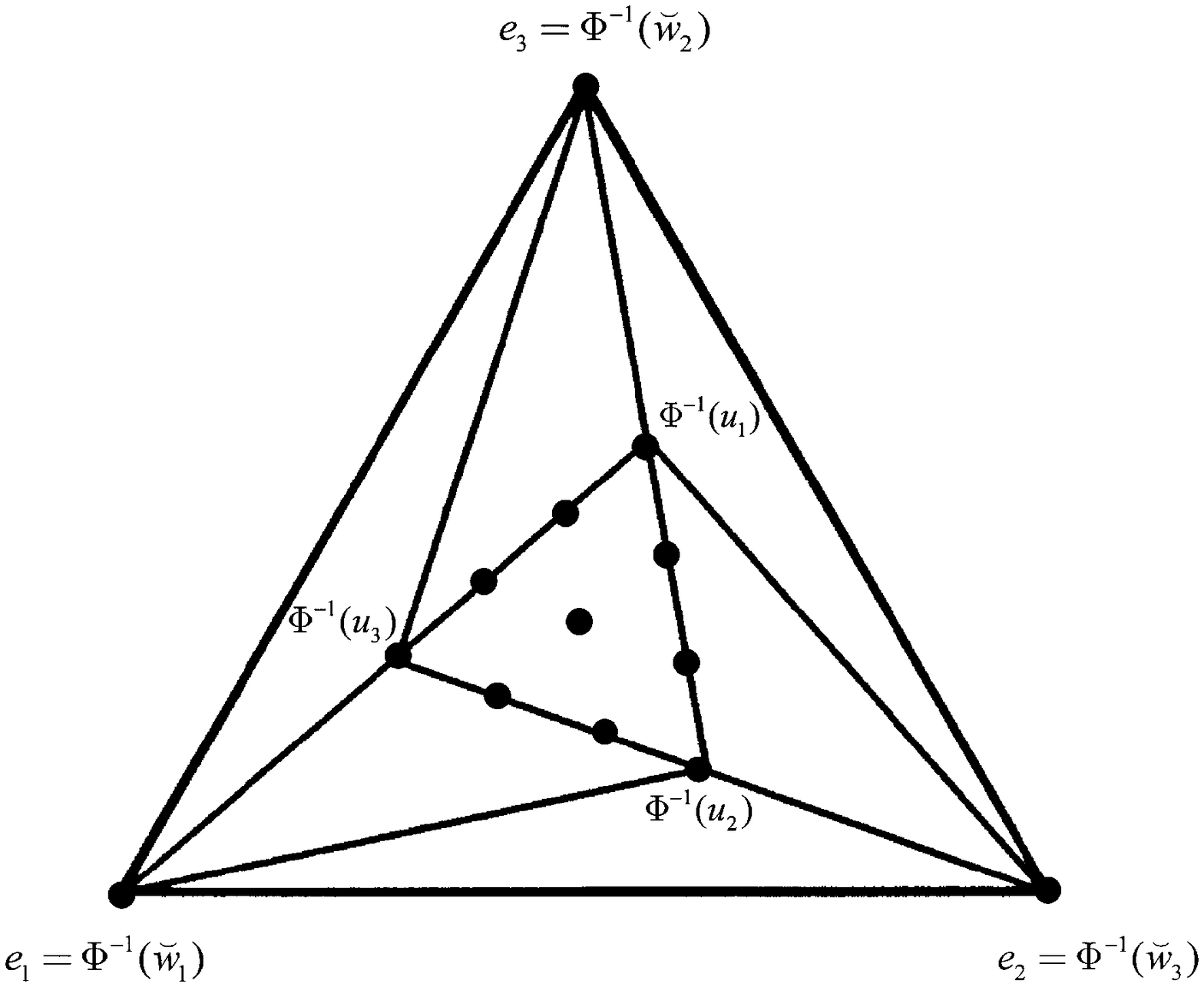, height=8cm, width=10cm}
\caption{}\label{Fig.5}
\end{figure}
\end{example}

\begin{remark}
\label{REMGPSS2}By Remark \ref{REMGPSS1} (ii) and Lemma \ref{LGPSS3}, the
$(r-1)$-dimensional lattice simplex%
\[
\fbox{$%
\begin{array}
[c]{lll}
& \mathfrak{\mathring{s}}_{G}:=\mathbf{s}(1,1,\ldots,1,1)=\text{conv}\left(
\mathfrak{\breve{s}}_{G}\mathbb{r}\{\left.  \breve{w}_{j}\ \right\vert \ 1\leq
j\leq r\}\right)  , &
\end{array}
$}%
\]
(included in the interior of $\mathfrak{\breve{s}}_{G}=$ conv$(\{\left.
\breve{w}_{j}\ \right\vert \ 1\leq j\leq r\})$) contains the $\tbinom
{k+r-2}{r-1}$ non-vertex lattice points of $\mathfrak{\breve{s}}_{G}$. Since
\[
\mathfrak{\mathring{s}}_{G}=\text{conv}(\{\left.  u_{j}\ \right\vert \ 1\leq
j\leq r\})=\left(  k-1\right)  \cdot\text{conv}(\{\left.  e_{j}\ \right\vert
\ 1\leq j\leq r\})
\]
is the dilation of a \textit{basic simplex} (w.r.t. $\breve{N}_{G}$) by the
factor $k-1,$ there is an affine transformation $\mathbb{R}^{r}\longrightarrow
\mathbb{R}^{r-1}\times\{0\}\subset\mathbb{R}^{r}$ whose restriction on the
affine hull aff$(\mathfrak{\breve{s}}_{G})$ of $\mathfrak{\breve{s}}_{G}$ is a
bijection, say $\Upsilon,$ mapping the lattice aff$(\mathfrak{\breve{s}}%
_{G})\cap\breve{N}_{G}$ onto the standard rectangular lattice $\mathbb{Z}%
^{r-1}=\mathbb{Z}^{r-1}\times\{0\}\subset\mathbb{R}^{r-1}\times\{0\},$ and
$\mathfrak{\mathring{s}}_{G}$ onto%
\[
\Upsilon\left(  \mathfrak{\mathring{s}}_{G}\right)  =\left(  k-1\right)
\cdot\text{conv}(\{\mathbf{0},e_{1},e_{1}+e_{2},\ldots,e_{1}+e_{2}%
+\cdots+e_{r-1}\}).
\]
Hence, $\Upsilon^{-1}(\mathbf{T}\left(  r-1;k-1\right)  )$ is a basic coherent
triangulation of $\mathfrak{\mathring{s}}_{G}$ (w.r.t. $\breve{N}_{G}$), where
$\mathbf{T}\left(  r-1;k-1\right)  $ denotes the triangulation of
$\Upsilon\left(  \mathfrak{\mathring{s}}_{G}\right)  $ (w.r.t. $\mathbb{Z}%
^{r-1}$) defined in Example \ref{HVS} (with $d=r-1$).
\end{remark}

\begin{figure}[h]
\epsfig{file=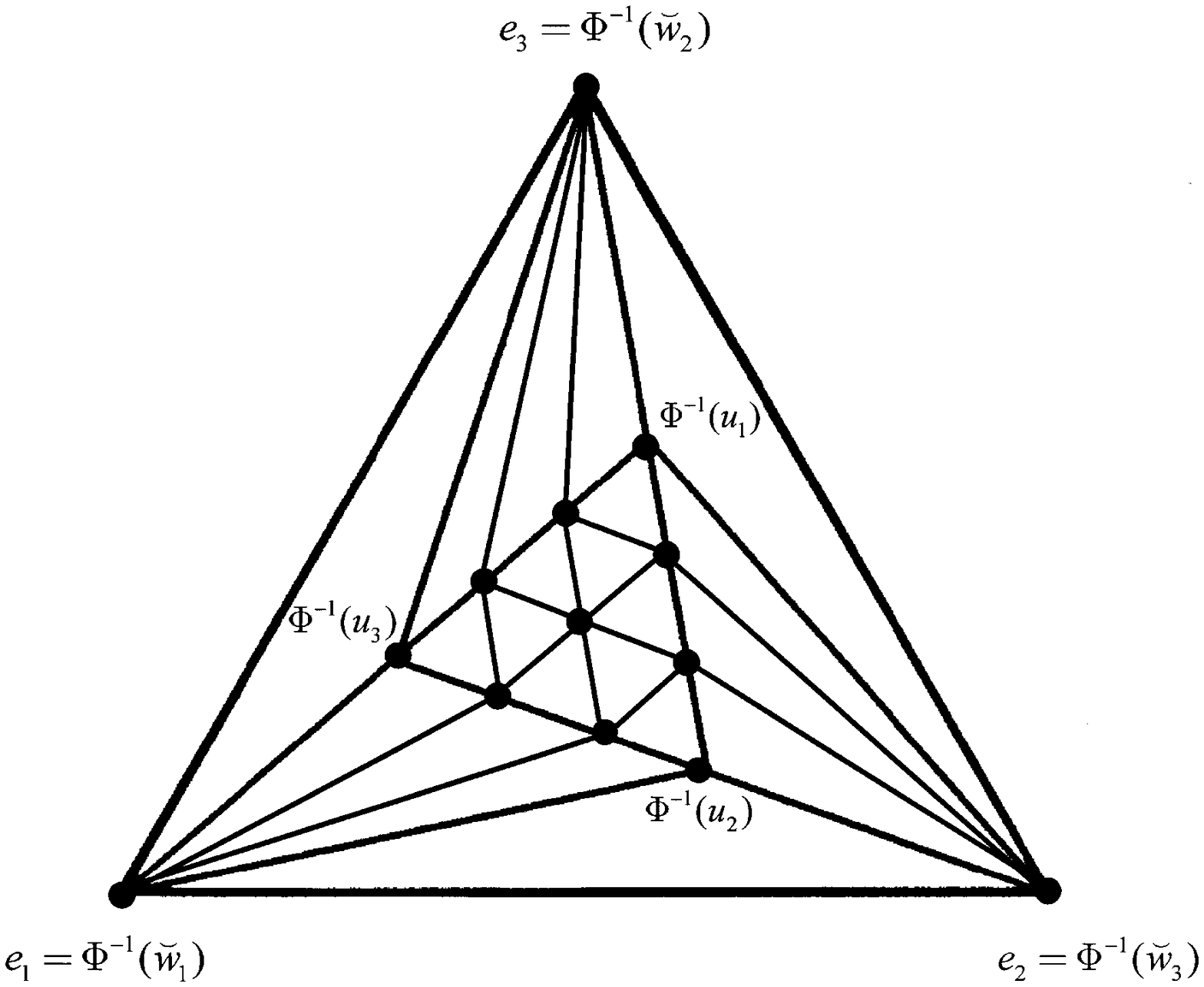,  height=8cm, width=10cm}
\caption{}\label{Fig.6}
\end{figure}

\noindent{}$\rhd$ \textsc{Proof of Theorem }\ref{GPSSTHM}: (i)
\textit{Basicness}. Setting
\[
\mathcal{E}_{\nu_{1},\ldots,\nu_{\rho}}:=\left\{
\begin{array}
[c]{l}%
\text{conv}(\left\{  \breve{w}_{\nu_{1}},\ldots,\breve{w}_{\nu_{\rho}%
}\right\}  )\ast\mathbf{s}\\%
\begin{array}
[c]{l}%
\text{(together with}\\
\text{their faces)}%
\end{array}
\end{array}
\left\vert
\begin{array}
[c]{c}%
\mathbf{s}\text{ an }\left(  r-1-\rho\right)  \text{-dimensional }\\
\text{simplex of }\Upsilon^{-1}(\mathbf{T}\left(  r-1;k-1\right)  )\text{ }\\
\text{belonging to the face}\\
\text{conv}(\{\left.  u_{j}\right\vert j\in\{1,..,r\}\mathbb{r}\{\nu
_{1},..,\nu_{\rho}\}\})\text{ }\\
\text{of the simplex }\mathfrak{\mathring{s}}_{G}\text{ }%
\end{array}
\right.  \right\}  ,
\]
for all $\rho\in\{1,\ldots,r-1\}$ and all index subfamilies $1\leq\nu_{1}%
<\nu_{2}<\cdots<\nu_{\rho}\leq r,$ we define the triangulation
\[
\mathfrak{E}:=\bigcup_{\rho=1}^{r-1}\,\ \bigcup_{1\leq\nu_{1}<\nu_{2}%
<\cdots<\nu_{\rho}\leq r}\mathcal{E}_{\nu_{1},\nu_{2},\ldots,\nu_{\rho}}%
\]
which refines $\left.  \mathfrak{T}(r;k)\right\vert _{\mathfrak{\breve{s}}%
_{G}\mathbb{r}\text{int}(\mathfrak{\mathring{s}}_{G})}$ (with $\mathfrak{T}%
(r;k)$ as given in (\ref{TRIANGTAU})). The set of $(r-1)$-dimensional
simplices of $\mathfrak{E}$ consists of well-defined joins (cf. the proof of
Lemma \ref{LGPSS4}). Glueing $\Upsilon^{-1}(\mathbf{T}\left(  r-1;k-1\right)
)$ and $\mathfrak{E}$ together we obtain a lattice triangulation (w.r.t.
$\breve{N}_{G}$)
\[
\fbox{$%
\begin{array}
[c]{lll}
& \mathcal{T}\left(  r;k\right)  :=\Upsilon^{-1}(\mathbf{T}\left(
r-1;k-1\right)  )\cup\mathfrak{E} &
\end{array}
$}\text{ }%
\]
of the entire simplex $\mathfrak{\breve{s}}_{G}.$ The triangulation
$\Upsilon^{-1}(\mathbf{T}\left(  r-1;k-1\right)  )$ itself is basic. Since for
all $\rho\in\{1,\ldots,r-1\},${\
\[%
\begin{array}
[c]{l}%
\text{aff}_{\mathbb{Z}}(\left\{  \breve{w}_{\nu_{1}},\ldots,\breve{w}%
_{\nu_{\rho}}\right\}  \cup\text{ }\{\left.  u_{j}\ \right\vert \ j\in
\{1,..,r\}\mathbb{r}\{\nu_{1},..,\nu_{\rho}\}\})\medskip\\
=\breve{N}_{G}\cap\text{aff}(\text{conv}(\left\{  \breve{w}_{\nu_{1}}%
,\ldots,\breve{w}_{\nu_{\rho}}\right\}  )\cup\text{conv}(\{\left.
u_{j}\ \right\vert \ j\in\{1,..,r\}\mathbb{r}\{\nu_{1},..,\nu_{\rho}\}\})),
\end{array}
\]
}$\mathfrak{E}$ is basic by \cite[Thm. 3.5, pp. 206-207]{DHZ}. Thus, the
entire $\mathcal{T}\left(  r;k\right)  $ is also basic. (Alternatively, since%
\[
\left\{
\begin{array}
[c]{l}%
\sharp\left\{  (r-1)\text{-dimensional simplices of }\mathfrak{T}%
(r;k)\right\}  =2^{r}-1,\\
\ \\
\sharp\left\{
\begin{array}
[c]{c}%
(r-1)\text{-dimensional simplices }\\
\text{of }\Upsilon^{-1}(\mathbf{T}\left(  r-1;k-1\right)  )
\end{array}
\right\}  =\left(  k-1\right)  ^{r-1},\\
\ \\
\text{and, analogously, for all }\rho\in\{1,\ldots,r-1\},\\
\ \\
\sharp\left\{
\begin{array}
[c]{c}%
\left(  r-1-\rho\right)  \text{-dimensional simplices }\\
\text{of }\Upsilon^{-1}(\mathbf{T}\left(  r-1;k-1\right)  )\text{ belonging to
the face}\\
\text{conv}(\{\left.  u_{j}\ \right\vert \ j\in\{1,..,r\}\mathbb{r}\{\nu
_{1},..,\nu_{\rho}\}\})\text{ of }\mathfrak{\mathring{s}}_{G}%
\end{array}
\right\}  =\left(  k-1\right)  ^{r-\rho-1},
\end{array}
\right.
\]
we get%
\begin{align*}
\sharp\left\{
\begin{array}
[c]{c}%
(r-1)\text{-dimensional }\\
\text{simplices of }\mathcal{T}\left(  r;k\right)
\end{array}
\right\}   &  =\left(  2^{r}-1\right)  +\sum_{\rho=0}^{r-1}\left(  \tbinom
{r}{\rho}\left(  k-1\right)  ^{r-\rho-1}-\tbinom{r}{\rho}\right)  \smallskip\\
&  =\left(  2^{r}-1\right)  -\sum_{\rho=0}^{r-1}\tbinom{r}{\rho}+\tfrac
{1}{k-1}\sum_{\rho=0}^{r-1}\left(  \tbinom{r}{\rho}\left(  k-1\right)
^{r-\rho}\right)  \smallskip\\
&  =\tfrac{1}{k-1}\sum_{\rho=0}^{r-1}\left(  \tbinom{r}{\rho}\left(
k-1\right)  ^{r-\rho}\right)  =\tfrac{(\left(  k-1\right)  +1)^{r}-1}{k-1}=l,
\end{align*}
and $\mathcal{T}\left(  r;k\right)  $ has to be a basic
triangulation of $\mathfrak{\breve{s}}_{G}$ according to Remark
\ref{BASICCR}). We conclude that $\Phi^{-1}\left(
\mathcal{T}\left(  r;k\right)  \right)  $ is a basic triangulation
of the junior simplex $\mathfrak{s}_{G}$ (w.r.t. $N_{G}$), where
$\Phi$ is the regular linear transformation (\ref{PhiW}). The
basic triangulation $\Phi^{-1}\left(  \mathcal{T}\left( 3;4\right)
\right)  $ of $\mathfrak{s}_{G}$ (refining that one of Example
\ref{EXGP}) is shown in Figure \ref{Fig.6}.

\noindent{}(ii) \textit{Coherence}. We define the $\mathcal{T}\left(
r;k\right)  $-support function $\mathbf{\Psi}:\mathfrak{\breve{s}}%
_{G}\longrightarrow\mathbb{R}$ as follows:%
\[
\mathfrak{\breve{s}}_{G}\ni\mathbf{x}\longmapsto\mathbf{\Psi}(\mathbf{x}%
):=\left\{
\begin{array}
[c]{ll}%
\psi_{\text{Heav}}^{(r-1)}(\Upsilon(\mathbf{x})), & \text{if }\mathbf{x}%
\in\mathfrak{\mathring{s}}_{G},\\
& \\
\psi_{\nu_{1},\nu_{2},\ldots,\nu_{\rho}}(\mathbf{x}), & \text{if }\left\{
\begin{array}
[c]{l}%
\mathbf{x}\in\left\vert \mathcal{E}_{\nu_{1},\nu_{2},\ldots,\nu_{\rho}%
}\right\vert ,\\
\text{for some indices}\\
1\leq\nu_{1}<\cdots<\nu_{\rho}\leq r,\\
\text{with }\rho\in\{1,\ldots,r-1\},
\end{array}
\right.
\end{array}
\right.
\]
where $\psi_{\text{Heav}}^{(r-1)}$ is the function defined in \ref{HVS} (with
$d=r-1$), and%
\[
\psi_{\nu_{1},\nu_{2},\ldots,\nu_{\rho}}(\mathbf{x}):=t\cdot(\theta
(\widetilde{\Phi}(\mathbf{x}_{1})))+\left(  1-t\right)  \cdot\psi
_{\text{Heav}}^{(r-1)}(\Upsilon(\mathbf{x}_{2})),
\]
for all $\mathbf{x}=t\mathbf{x}_{1}+(1-t)\mathbf{x}_{2}\in$ conv$(\left\{
\breve{w}_{\nu_{1}},\ldots,\breve{w}_{\nu_{\rho}}\right\}  )\ast\mathbf{s}$
belonging to $\left\vert \mathcal{E}_{\nu_{1},\ldots,\nu_{\rho}}\right\vert ,$
with $t\in\lbrack0,1],$
\[
\mathbf{x}_{1}\in\text{conv}(\left\{  \breve{w}_{\nu_{1}},\ldots,\breve
{w}_{\nu_{\rho}}\right\}  )\text{ and }\mathbf{x}_{2}\in\mathbf{s\in
}\left\vert \left.  \mathcal{T}\left(  r;k\right)  \right\vert _{\text{conv}%
(\{\left.  u_{j}\right\vert j\in\{1,..,r\}\mathbb{r}\{\nu_{1},..,\nu_{\rho
}\}\})}\right\vert \text{ }.
\]
Since $\theta$ (by Corollary \ref{CORCOH}), as well as $\left.  \psi
_{\text{Heav}}^{(r-1)}\right\vert _{\mathfrak{\mathring{s}}_{G}}$ and
$\psi_{\nu_{1},\ldots,\nu_{\rho}}$'s, are strictly upper convex, and the
latter ones coincide on their common domains of linearity, we deduce by the
Patching Lemma\ \ref{PATCH} that also $\mathbf{\Psi}$ is strictly upper
convex. This means that $\mathcal{T}\left(  r;k\right)  $ is a coherent
triangulation of $\mathfrak{\breve{s}}_{G}.$ Consequently, $\Phi^{-1}\left(
\mathcal{T}\left(  r;k\right)  \right)  $ is a coherent triangulation of the
junior simplex $\mathfrak{s}_{G}.\medskip$\newline(iii) \textit{The special
case in which} $k=2.$ In this case, $\mathfrak{\mathring{s}}_{G}$ itself is
already basic (w.r.t. $\breve{N}_{G}$), and the only maximal (and necessarily
basic) triangulation of the junior simplex $\mathfrak{s}_{G}$ (w.r.t. $N_{G}$)
is $\Phi^{-1}\left(  \mathfrak{T}(r;k)\right)  .$ Its uniqueness and coherence
follow from Lemma \ref{LGPSS4} and Corollary \ref{CORCOH}, respectively.\hfill
{}$\square$

\begin{note}
\label{NOTEGPSS}(i) For $k\geq3,$ besides $\Phi^{-1}\left(  \mathcal{T}\left(
r;k\right)  \right)  ,$ there are lots of other basic triangulations of
$\mathfrak{s}_{G},$ due to those of $\Phi^{-1}\left(  \mathfrak{\mathring{s}%
}_{G}\right)  ;$ cf. Remark \ref{REMARKHYPERS} (ii).\smallskip\ \newline(ii)
\textit{Open Problem}: As it was proven recently by Sebestean \cite{Sebestean}
(for $r=4,$ $k=2$), the smooth fourfold obtained by the unique projective
crepant resolution of the (non-symplectic) CQS\ $(\mathbb{C}^{4}%
/G,[\mathbf{0}])$ of type $\frac{1}{15}(1,2,4,8)$ coincides with the Hilbert
scheme $G$-Hilb$(\mathbb{C}^{4})$ of $G$-clusters. It is therefore natural to
ask if this is in general true for \textit{all} the members of the series
GPSS$(r;2)$ (or not) whenever $r\geq5.$
\end{note}

\begin{exercise}
Compute the non-trivial cohomology dimensions (\ref{COHDIMSIG}) of any crepant
resolution space of any member of the geometric progress singularity series
GPSS($r;k$). (\textit{Hint}. Consider $\Phi^{-1}\left(  \mathcal{T}\left(
r;k\right)  \right)  $ as a composite of geometrically more \textquotedblleft
elementary\textquotedblright\ triangulations, and use the inclusion-exclusion
property of lattice point enumerators, combined with the multiplicative
property of the polynomial generating the $\mathbf{h}$-vectors of joins of
triangulations \cite[p. 466]{Chan2}, and with Theorem \ref{BMcMTHM} and
formula (\ref{COHDIMFORMULA}).)
\end{exercise}

\begin{remark}
Concerning the Existence Problem, it is worthwhile stressing the qualitative
difference between the behaviour of the $1$- and $2$-parameter singularity
series discussed in \S \ref{ONETWOPARSER} and that one of the geometric
progress singularity series GPSS($r;k$). The one or two parameters in the
types of the first mentioned singularities have to obay to
\textit{restrictive} arithmetic conditions in order to lead to crepant
resolutions (cf. (\ref{rescon}) and (\ref{CONDKETT})), whereas the parameter
$k\geq2$ in the GP-singularity series is \textit{unconditionally }free in this respect.
\end{remark}

\section{Second Existence Criterion via UBT\label{SECEXCRITERION}}

\noindent{}Let\emph{\ } $(\mathbb{C}^{r}/G,[\mathbf{0}])$ be a Gorenstein AQS
with $l=\left\vert G\right\vert $\emph{ }and $r\geq4$\emph{.} The presence of
\textit{basic} triangulations $\mathcal{T}$ of $\mathfrak{s}_{G}$ (w.r.t.
$N_{G}$) implies the equality $l=(r-1)!$Vol$\left(  \mathfrak{s}_{G}\right)
=\mathfrak{f}_{r-1}(\mathcal{T})$ (by (\ref{VOLEQUALITY})). Bounding the
cardinality $\mathfrak{f}_{r-1}(\mathcal{T})$ of the facets of any such
$\mathcal{T}$ from above by a number depending only on the number of the
available lattice points in $\mathfrak{s}_{G},$ it is possible to obtain a
second \textit{necessary existence condition} which is highly effective and of
purely geometric nature. It comes as no surprise to learn that such a number
involves the cardinality $\mathfrak{f}_{r-1}\left(  \text{CycP}_{r}\left(
\sharp\left(  \mathfrak{s}_{G}\cap N_{G}\right)  \right)  \right)  $ of the
facets of the $r$-dimensional cyclic polytope with $\sharp\left(
\mathfrak{s}_{G}\cap N_{G}\right)  $ vertices, because it reminds you of the
celebrated \textit{UBT }\ref{UBTSS}\textit{ for simplicial spheres}.
Nevertheless, this has first to be suitably modified to be valid for
\textit{simplicial balls} (like $\mathcal{T}$); see Theorem \ref{SECUBB}.
Unfortunately, even if we use the latter upper bound, we do lose some
information whenever our singularity is \textit{non-isolated}, because we are
throwing away a considerable part of the individual contributions of lattice
points which belong to the boundary of $\mathfrak{s}_{G}.$ In fact, our
expectation concerning a general, \textit{tight} upper bound for
$l=\mathfrak{f}_{r-1}(\mathcal{T})$ is expressed in the following:

\begin{conjecture}
\label{MAINCONJ}Let\emph{\ } $(\mathbb{C}^{r}/G,[\mathbf{0}]),$ $r\geq
4$\emph{,} be a Gorenstein AQS with $l=\left\vert G\right\vert \emph{,}$ and
$\mathfrak{s}_{G}$ the corresponding junior simplex. If $\mathfrak{s}_{G}$ has
a basic triangulation, then $l$ has the following upper bound\emph{:}
\begin{equation}
\fbox{$%
\begin{array}
[c]{ccc}
&  & \\
& l\leq\ \mathfrak{f}_{r-1}\left(  \emph{CycP}_{r}\left(  \sharp\left(
\mathfrak{s}_{G}\cap N_{G}\right)  \right)  \right)  -%
{\displaystyle\sum\limits_{k=2}^{r-1}} \ (r-k)\left(  \sharp\left(
\mathfrak{B}_{G}\left(  1,k\right)  \right)
\right)  -1, & \\
& \  &
\end{array}
$} \label{conjformula}%
\end{equation}
with $\sharp\left(  \mathfrak{B}_{G}\left(  1,k\right)  \right)  $\emph{'s} as
given in \emph{(\ref{BG1K}).}
\end{conjecture}

\begin{note}
For the proof of (\ref{conjformula}) it would suffice to show that
UBT-Conjecture \ref{Z-conj} is true. In Theorem \ref{IICR} we prove
(\ref{conjformula}) only for $r=4,$ and give the weaker upper bound for
$r\geq5.$
\end{note}

\begin{theorem}
[Second Necessary Existence Condition]\label{IICR} Let\emph{\ } $(\mathbb{C}%
^{r}/G,[\mathbf{0}]),$ $r\geq4$\emph{,} be a Gorenstein AQS with $l=\left\vert
G\right\vert \emph{,}$ and $\mathfrak{s}_{G}$ the corresponding junior
simplex. If $\mathfrak{s}_{G}$ has a basic triangulation $\mathcal{T},$ then
$l$ has as upper bound
\begin{equation}
\fbox{$%
\begin{array}
[c]{c}%
\\
l\leq\mathfrak{f}_{3}\left(  \emph{CycP}_{4}\left(  \sharp\left(
\mathfrak{s}_{G}\cap N_{G}\right)  \right)  \right)  -2\ \left(  \sharp\left(
\mathfrak{B}_{G}\left(  1,2\right)  \right)  \right)  -\left(  \sharp\left(
\mathfrak{B}_{G}\left(  1,3\right)  \right)  \right)  -1\\
\
\end{array}
$} \label{fircr2}%
\end{equation}
for $r=4$\emph{,} and
\begin{equation}
\fbox{$%
\begin{array}
[c]{ccc}
&  & \\
& l\leq\ \mathfrak{f}_{r-1}\left(  \emph{CycP}_{r}\left(  \sharp\left(
\mathfrak{s}_{G}\cap N_{G}\right)  \right)  \right)  -%
{\displaystyle\sum\limits_{k=2}^{r-1}} \ \left(  \sharp\left(
\mathfrak{B}_{G}\left(  1,k\right)  \right)  \right)
-1 & \\
& \  &
\end{array}
$} \label{criter2}%
\end{equation}
for $r\geq5$\emph{. (The number} $\sharp\left(  \mathfrak{s}_{G}\cap
N_{G}\right)  $ \emph{can be calculated by the formulae (\ref{DR-FORM}) and
(\ref{DR-TRANSF}) given in Appendix \ref{LATPJS}. The numbers} $\sharp\left(
\mathfrak{B}_{G}\left(  1,k\right)  \right)  $, $2\leq k\leq r-1$,
\emph{occuring in (\ref{fircr2}), (\ref{criter2}),} \emph{are computable
}either\emph{ by (\ref{BG1K}) and (\ref{AB-FORM}) }or\emph{ by using
(\ref{BG1K}) and then counting the lattice points lying in the relative
interior of each of the }$\left(  k-1\right)  $\emph{-dimensional faces of
}$\mathfrak{s}_{G}$\emph{, by (\ref{RECIPRLOW}), (\ref{DR-FORM}) and
(\ref{DR-TRANSF}),} \emph{applied for these faces instead for }$\mathfrak{s}%
_{G}$ \emph{itself).}
\end{theorem}

\noindent{}\textsc{Proof}\textit{. }Using the notation of Appendix \ref{APPA},
apply (\ref{UP-LO2}) (to get (\ref{criter2})\emph{ }for $r\geq5$) just by
setting $d=r-1,\ \mathcal{S}=\mathcal{T},$ $\mathfrak{b}=\sharp\left(
\mathfrak{s}_{G}\cap N_{G}\right)  ,$ $\mathfrak{b}^{\prime}=%
{\textstyle\sum\nolimits_{k=2}^{r-1}}
\ \left(  \sharp\left(  \mathfrak{B}_{G}\left(  1,k\right)  \right)  \right)
+r.$ Correspondingly, to get (\ref{fircr2}) for $r=4,$ apply Theorem
\ref{BIGTH} by setting $d=3$, $\mathbf{s}=\mathfrak{s}_{G},$ $\mathcal{S}%
=\mathcal{T},$ $b_{1}=4$, and $b_{k}=\sharp\left(  \mathfrak{B}_{G}\left(
1,k\right)  \right)  $, for $k\in\{2,3\}$. Of course, for the desingularizing
space $X(N_{G},\widehat{\Delta_{G}}\left(  \mathcal{T}\right)  )$ of $X\left(
N_{G},\Delta_{G}\right)  =\mathbb{C}^{r}/G$ being induced by $\mathcal{T}$, we
have $l=\mathfrak{f}_{r-1}\left(  \mathcal{T}\right)  $ by (\ref{VOLEQUALITY}%
).\hfill{}$\square$

\begin{corollary}
\label{TET-COR}Let \emph{\ }$(\mathbb{C}^{4}/G,[\mathbf{0}])$ be a Gorenstein
cyclic quotient msc-singularity of type $\frac{1}{l}\left(  \alpha_{1}%
,\alpha_{2},\alpha_{3},\alpha_{4}\right)  .$ Then the inequality
\emph{(\ref{fircr2})} can be written as\emph{\ }follows\emph{:}
\begin{equation}
\fbox{$%
\begin{array}
[c]{c}%
\\
\!\!l\leq\!\tfrac{\sharp\left(  \mathfrak{s}_{G}\cap N_{G}\right)
(\sharp\left(  \mathfrak{s}_{G}\cap N_{G}\right)  -3)}{2}-%
{\displaystyle\sum\limits_{i=1}^{4}}
\frac{\text{\emph{gcd}}\left(  \alpha_{i},l\right)  }{2}\ \!-%
{\displaystyle\sum\limits_{1\leq i<j\leq4}}
\ \!\text{\emph{gcd}}\left(  \alpha_{i},\alpha_{j},l\right)  +7,\!\!\\
\
\end{array}
$} \label{tetrah}%
\end{equation}
\emph{\smallskip}where $\sharp\left(  \mathfrak{s}_{G}\cap N_{G}\right)  $ is
known by the formulae \emph{(\ref{MODPO1}), (\ref{MODPO2}), (\ref{MODPO3}).}
\end{corollary}

\noindent\ \noindent{}\textsc{Proof}. Obviously,
\[
\sharp\left(  \mathfrak{B}_{G}\left(  1,2\right)  \right)  =\sum_{1\leq\nu
_{1}<\nu_{2}\leq4}\ \sharp\left(  \mathfrak{B}_{G}\left(  1,2;\nu_{1},\nu
_{2}\right)  \right)  .
\]
If for any pair of indices $\nu_{1},\nu_{2}$, with $1\leq\nu_{1}<\nu_{2}\leq
4$, we define $\left\{  \nu_{3},\nu_{4}\right\}  $ to be the complement set
$\left\{  1,2,3,4\right\}  \smallsetminus\left\{  \nu_{1},\nu_{2}\right\}  $,
then
\begin{equation}
\sharp\left(  \mathfrak{B}_{G}\left(  1,2;\nu_{1},\nu_{2}\right)  \right)
=\text{gcd}\left(  \alpha_{\nu_{3}},\alpha_{\nu_{4}},l\right)  -1 \label{N-12}%
\end{equation}
by (\ref{GCD-FORM}). Analogously,
\[
\sharp\left(  \mathfrak{B}_{G}\left(  1,3\right)  \right)  =\sum_{1\leq\nu
_{1}<\nu_{2}<\nu_{3}\leq4}\ \sharp\left(  \mathfrak{B}_{G}\left(  1,3;\nu
_{1},\nu_{2},\nu_{3}\right)  \right)  ,
\]
and if for any triple of indices $\nu_{1},\nu_{2},\nu_{3}$, with $1\leq\nu
_{1}<\nu_{2}<\nu_{3}\leq4$, $\left\{  \nu_{4}\right\}  $ denotes the
complement set $\left\{  1,2,3,4\right\}  \smallsetminus\left\{  \nu_{1}%
,\nu_{2},\nu_{3}\right\}  $, then the number of the interior points of each
$2$-face of the junior tetrahedron $\mathfrak{s}_{G}$ equals
\begin{align}
\sharp\left(  \mathfrak{B}_{G}\left(  1,3;\nu_{1},\nu_{2},\nu_{3}\right)
\right)   &  =\frac{1}{2}\ \left[  \text{gcd}\left(  \alpha_{\nu_{4}%
},l\right)  -1-\left(  \sum_{j=1}^{3}\text{\ gcd}\left(  \alpha_{\nu_{4}%
},\alpha_{\nu_{j}},l\right)  -1\right)  \right] \nonumber\\
&  \ \ \nonumber\\
&  =\frac{1}{2}\ \left[  \text{gcd}\left(  \alpha_{\nu_{4}},l\right)
-\sum_{j=1}^{3}\ \text{gcd}\left(  \alpha_{\nu_{4}},\alpha_{\nu_{j}},l\right)
\right]  +1 \label{N123}%
\end{align}
by the refined Ping-Pong Lemma \ref{REF p-p} (see Remark \ref{pprem}).
Substituting (\ref{N-12}), (\ref{N123}) into formula (\ref{fircr2}) we obtain
(\ref{tetrah}).\hfill{}$\square$

\begin{example}
Let \emph{\ }$(\mathbb{C}^{4}/G,[\mathbf{0}])$ be the CQS of type $\frac
{1}{12}\left(  \alpha_{1},\alpha_{2},\alpha_{3},\alpha_{4}\right)  ,$ with
$\alpha_{1}=\alpha_{2}=2,$ $\alpha_{3}=3,$ and $\alpha_{4}=5.$ Since
$\sharp\left(  \mathfrak{s}_{G}\cap N_{G}\right)  =7,$ and the right-hand side
of (\ref{tetrah}) equals%
\[
\tfrac{7(7-3)}{2}-\tfrac{1}{2}\ \left[
{\textstyle\sum\limits_{i=1}^{4}}
\ \text{gcd}\left(  \alpha_{i},l\right)  \right]  -%
{\textstyle\sum\limits_{1\leq i<j\leq4}} \ \text{gcd}\left(
\alpha_{i},\alpha_{j},l\right)  +7=14-4=10,
\]
it does not admit any crepant resolution because $10<12=l$.
\end{example}

\begin{remark}
[Comparison of the two Existence Criteria]Which of the necessary
conditions (\ref{HILBCON}) and (\ref{fircr2})-(\ref{criter2})
given in Theorems \ref{KILLER} and \ref{IICR}, respectively,\ is
better? The answer to this question depends on how one would like
to interpret the adjective \textquotedblleft
better\textquotedblright. Undoubtedly, (\ref{HILBCON})
\textquotedblleft kills\textquotedblright\ more candidates for
having crepant resolutions. For instance, for the CQS
$(\mathbb{C}^{4}/G,[\mathbf{0}])$ of type $\frac{1}{9}(1,2,3,3)$
(\ref{tetrah}) holds as equality but
$\frac{1}{9}(5,1,6,6)^{\intercal}\in\mathbf{Hlb}_{N_{G}}\left(
\sigma _{0}\right)  \mathbb{r}(\mathfrak{s}_{G}\cap N_{G}).$
(Hence, this CQS does not have any crepant resolution.) On the
other hand, in view of Theorem \ref{SEBOTHM}, the determination of
the Hilbert basis is a time-consuming procedure compared with the
lattice point enumeration of the junior simplex (in particular, in
high dimensions and for acting groups with big orders).
\end{remark}

\begin{exercise}
For the Gorenstein CQS of type $\frac{1}{12}(1,3,3,5)$ show that
(\ref{tetrah}) holds as strict inequality, though $\frac{1}{12}%
(3,9,9,3)^{\intercal}\in\mathbf{Hlb}_{N_{G}}\left(  \sigma_{0}\right)
\mathbb{r}(\mathfrak{s}_{G}\cap N_{G}).$
\end{exercise}

\section{Sketching an Auxiliary Algorithm\label{ALGORITHM}}

\begin{figure}[h]
\epsfig{file=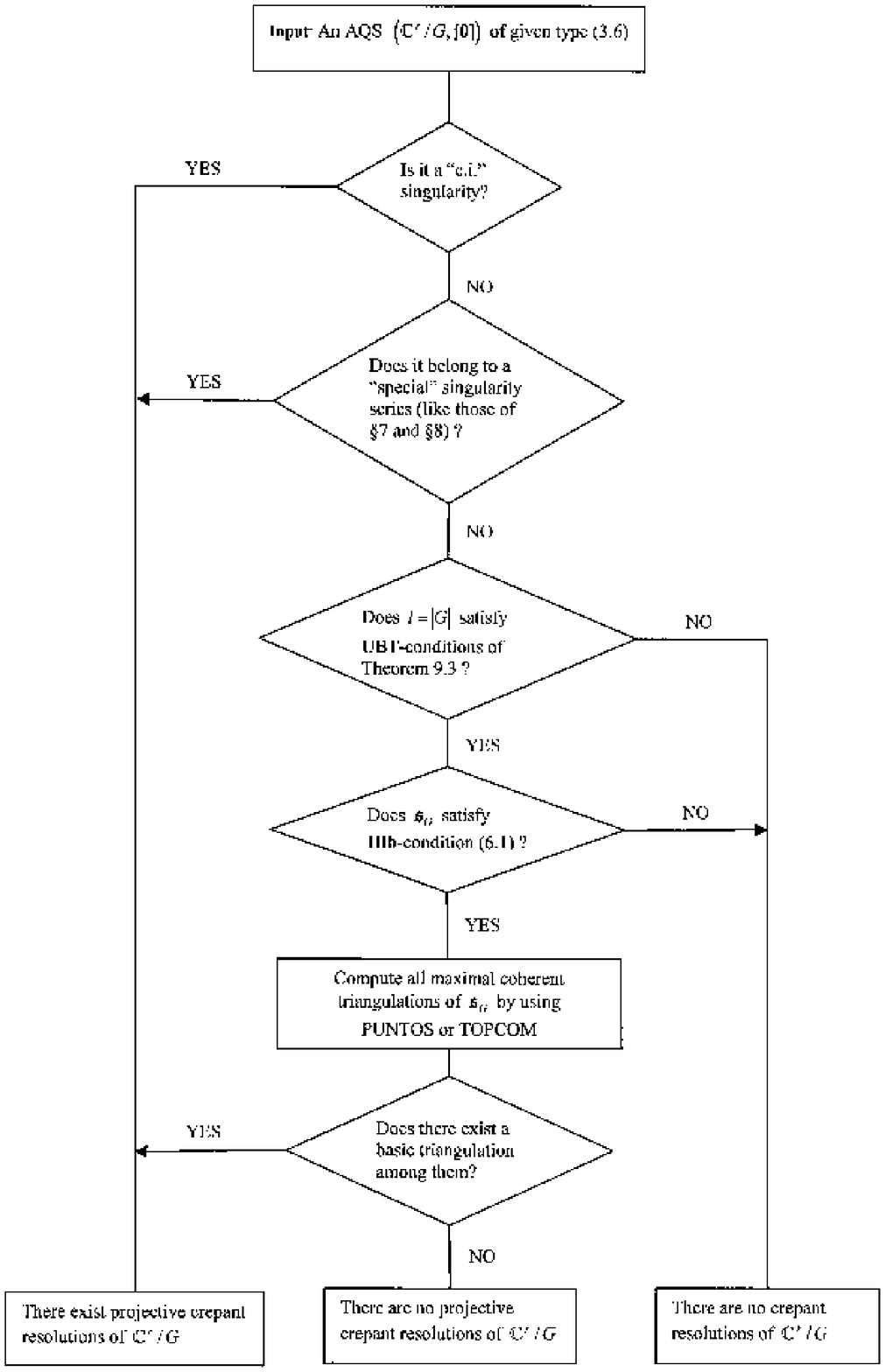, height=18cm, width=12.5cm}
\caption{}\label{Fig.7}
\end{figure}

\noindent{}Taking into account what we have discussed so far, it
is possible, for given AQS $(\mathbb{C}^{r}/G,[\mathbf{0}])$ of
type (\ref{typeAQS}), to outline an algorithm in order to examine
whether it admits the desired resolutions, but at the cost of
\textit{increasing computational complexity} (in the consecutive
steps). More precisely, the auxiliary algorithm we have in mind
(summarized in Figure \ref{Fig.7}) is built up as follows:\medskip

\noindent{}$\vartriangleright$ \textsc{Step }$1$. If $\mathfrak{s}_{G}$ is
lattice equivalent to a Watanabe simplex, then $\mathbb{C}^{r}/G$ admits of
projective crepant desingularizations according to Theorem \ref{DHZTHM}. If
not, we go to Step $2$.\medskip

\noindent{}$\vartriangleright$ \textsc{Step }$2$. If $(\mathbb{C}%
^{r}/G,[\mathbf{0}])$ is non-c.i but belongs to \textquotedblleft
special\textquotedblright\ singularity series (like those $1$- and
$2$-parameter series of \S \ref{ONETWOPARSER} having weights satisfying
conditions (\ref{rescon}), (\ref{CONDKETT}), or even the entire
GP-singularities series of \S \ref{GPSERIES}), which have projective crepant
resolutions \textit{by construction}, we stop; otherwise we proceed. (To
continue increasing our stock of \textquotedblleft special\textquotedblright%
\ singularity series of this kind would be a real challenge for future
work.)\medskip

\noindent{}$\vartriangleright$ \textsc{Step }$3$. We \textit{count the lattice
points} of the junior simplex $\mathfrak{s}_{G}$ involved in (\ref{fircr2}),
resp., (\ref{criter2}), by the formulae given in Appendix \ref{LATPJS}, and
then we check if these inequalities for $l=\left\vert G\right\vert $ are valid
or not. We proceed only if (\ref{fircr2}) (for $r=4$), resp., (\ref{criter2})
(for $r\geq5$), are indeed valid; otherwise $\mathbb{C}^{r}/G$ does not admit
any crepant desingularization by Theorem \ref{IICR}.\medskip

\noindent{}$\vartriangleright$ \textsc{Step }$4$. We determine the
Hilbert basis $\mathbf{Hlb}_{N_{G}}\left(  \sigma_{0}\right)  $
(see Remark \ref{HBREM}), and control if it satisfies condition
(\ref{HILBCON}). We proceed to the next (final) step only if
(\ref{HILBCON}) is satisfied; otherwise the quotient space
$\mathbb{C}^{r}/G$ does not have any crepant desingularization by
Theorem \ref{KILLER}.\medskip

\noindent{}$\vartriangleright$ \textsc{Step }$5$. If $(\mathbb{C}%
^{r}/G,[\mathbf{0}])$ happens to pass all the above tests without stop in the
one or the other stage, we have to find out \textit{all the junior lattice
points}
\[
\left\{  \left.  n_{g}\in N_{G}\text{ (as in (\ref{LPTREPRG}))}\right\vert
\text{ age}(g)=1\right\}  =\mathfrak{s}_{G}\cap N_{G}%
\]
(not just their cardinality!), to run \texttt{Puntos}
\cite{DELOERA2} or \texttt{TOPCOM} \cite{RAMBAU} for the point
configuration $\mathcal{V}=\mathfrak{s}_{G}\cap N_{G}$ in order to
specify the coherent lattice triangulations of $\mathfrak{s}_{G},$
to separate the \textit{maximal} ones$,$ and then to count the
number of $(r-1)$-dimensional simplices in each of them; cf. Note
\ref{DELOERANOTE}. Projective crepant desingularizations of
$\mathbb{C}^{r}/G$ are present as long as this number equals $l$
for at least one of them. (\textquotedblleft
Sporadic\textquotedblright\ counterexamples, like the isolated CQS
of type $\frac{1}{39}\,\left(  1,5,8,25\right)  $ mentioned in
Remark \ref{COUNTEREX} (iii), indicate why Step $5$ or any similar
computer-assisted procedure seems -as yet- to be unavoidable.)

\newpage
\newpage

\appendix

\section{Triangulations and Upper Bound Theorems\label{APPA}}

\noindent{}Triangulations (as geometric simplicial subdivisions of polytopes
or polytopal complexes) are treated in the classical framework of the
categorial inclusions:
\[
\left\{
\begin{array}
[c]{c}%
\text{geometric }\\
\text{simplicial complexes}%
\end{array}
\right\}  \subset\left\{
\begin{array}
[c]{c}%
\text{polytopal }\\
\text{complexes}%
\end{array}
\right\}  \subset\left\{
\begin{array}
[c]{c}%
``\text{regular cell }\\
\text{complexes\textquotedblright}\\
\text{(i.e. regular finite }\\
\text{CW complexes)}%
\end{array}
\right\}  .\smallskip
\]

\noindent$\bullet$ \textbf{Notation}. The symbol \textquotedblleft$\approx
$\textquotedblright\ between two topological spaces indicates the existence of
an homeomorphism from the one onto the other. A topological space $X$ is
called a \textit{sphere} (resp., a \textit{ball}) if $X\approx\mathbb{S}^{k}$
(resp., $X\approx\mathbb{B}^{k}$), for some $k$, where $\mathbb{S}^{k}$ and
$\mathbb{B}^{k}$ denote the standard $k$-dimensional sphere $\mathbb{S}%
^{k}=\partial\mathbb{B}^{k+1}=\left\{  \mathbf{x}\in\mathbb{R}^{k+1}\left\vert
\ \left\Vert \mathbf{x}\right\Vert =1\right.  \right\}  $ and the standard
$k$-dimensional ball $\mathbb{B}^{k}=\left\{  \mathbf{x}\in\mathbb{R}%
^{k}\left\vert \ \left\Vert \mathbf{x}\right\Vert \leq1\right.  \right\}  $,
respectively. If $k$ is assumed to be fixed, then we simply say that such an
$X$ is a $k$-\textit{sphere }(resp. a $k$-\textit{ball}). \textit{Polytopes}
will be always convex, defined as in \cite[Lecture $1$]{Ziegler}.\medskip

\noindent$\bullet$ \textbf{Regular cell complexes}. A \textit{regular cell
complex }$\mathcal{K}$ is a finite collection of balls $\mathbf{c}$ in a
Hausdorff topological space $\left\vert \mathcal{K}\right\vert =\bigcup
\left\{  \mathbf{c}\ \left\vert \ \mathbf{c}\in\mathcal{K}\right.  \right\}  $
such that\newline(i) the relative interiors int$\left(  \mathbf{c}\right)  $
of all $\mathbf{c}$'s partition $\left\vert \mathcal{K}\right\vert $, i.e.,
each element of $\left\vert \mathcal{K}\right\vert $ lies in exactly one
int$\left(  \mathbf{c}\right)  $, and\newline(ii) the relative boundary
$\partial\mathbf{c}$ of every $\mathbf{c}\in\mathcal{K}$ is a union of some
members of $\mathcal{K}$.\smallskip\newline The balls $\mathbf{c}%
\in\mathcal{K}$ are called the \textit{closed cells }of $\mathcal{K}$ and
their interiors int$\left(  \mathbf{c}\right)  $ the \textit{open cells }of
$\mathcal{K}$. $\left\vert \mathcal{K}\right\vert $ is called the
\textit{underlying space} (or \textit{the support}) of $\mathcal{K}$. The
dimension of a (closed) cell $\mathbf{c}$, for which $\mathbf{c}%
\approx\mathbb{B}^{k}$, is defined to be $k$. (Such a cell is particularly
called a $k$-\textit{cell}). If $\mathbf{c}_{1},\mathbf{c}_{2}\in\mathcal{K}$
and $\mathbf{c}_{1}\subset\mathbf{c}_{2}$, then $\mathbf{c}_{1}$ is said to be
a \textit{face }of $\mathbf{c}_{2}$. (We use the notation: $\mathbf{c}%
_{1}\prec\mathbf{c}_{2}$). $0$-and $1$-cells are called \textit{vertices }and
\textit{edges}, respectively. $\mathcal{K}$ is defined to be \textit{pure }if
all maximal cells have the same dimension. $\mathcal{K}^{\prime}$ is a
\textit{subcomplex} of $\mathcal{K}$ if $\mathbf{c}\in\mathcal{K}^{\prime}$
implies that every face of $\mathbf{c}$ belongs to $\mathcal{K}^{\prime}$.
Note that a regular cell complex is homeomorphic to the order complex of its
face poset.\medskip

\noindent$\bullet$ \textbf{Polytopal complexes}. A \textit{polytopal complex}
$\mathcal{S}$ consists of a finite family of polytopes in $\mathbb{R}^{d}$
such that\smallskip\newline(i) if $P\in\mathcal{S}$ and $F\prec P$, then
$F\in\mathcal{S}$, and\newline(ii) if $P_{1},P_{2}\in\mathcal{S}$ have
non-empty intersection, then $P_{1}\cap P_{2}\prec P_{1}$, and $P_{1}\cap
P_{2}\prec P_{2}$.\smallskip\newline Since every polytope is topologically a
ball, a polytopal complex $\mathcal{S}$ is a regular cell complex whose
(closed) cells (called also \textit{faces}) are the participating polytopes,
and whose underlying space $\left\vert \mathcal{S}\right\vert $ is the union
of these polytopes. (If $\mathcal{S}$ is a polytopal complex, we denote by
vert$\left(  \mathcal{S}\right)  $ the set of its vertices. If $\mathcal{S}$
is, in addition, pure, we call the dim$(\mathcal{S})$-faces \textit{facets} of
$\mathcal{S}.$)\medskip\

\noindent$\bullet$ \textbf{Geometric simplicial complexes}. A
\textit{geometric simplicial complex }is by definition a polytopal complex all
of whose (closed) cells are simplices. We frequently denote the simplices of
such an $\mathcal{S}$ by $F$ or $\mathbf{s}$ instead of $\mathbf{c}$. If
$\left\vert \mathcal{S}\right\vert \approx\mathbb{S}^{k}$ (resp., if
$\left\vert \mathcal{S}\right\vert \approx\mathbb{B}^{k}$), then $\mathcal{S}$
is called a \textit{simplicial }$k$-\textit{sphere }(resp., a
\textit{simplicial }$k$-\textit{ball}).

\begin{example}
\label{Boundarycomplex}\noindent Every $d$-polytope $P$ together with all of
its faces forms a polytopal $d$-complex $\mathcal{S}_{P}$. For a $d$-polytope
$P$ the \textit{boundary complex }$\mathcal{S}_{\partial P}$ of $P$ is defined
to be the $\left(  d-1\right)  $-dimensional polytopal complex consisting of
the proper faces of $P$ together with $\varnothing$ and having support
$\left\vert \mathcal{S}_{\partial P}\right\vert =\partial P$. The
\textit{facets} of $P$ are defined to be the facets of $\mathcal{S}_{\partial
P}.$ Obviously, $\mathcal{S}_{\partial P}$ is a geometric pure simplicial
complex (in fact, a simplicial $\left(  d-1\right)  $-sphere) if and only if
$P$ is a simplicial polytope.
\end{example}

\noindent$\bullet$ \textbf{Abstract simplicial complexes}. Geometric
simplicial complexes can be obtained as \textit{geometric realizations }of
\textquotedblleft abstract\textquotedblright\ simplicial (finite) complexes.
An \textit{abstract simplicial }(\textit{finite}) \textit{complex
}$\mathcal{S}\left(  \mathcal{V}\right)  $ with vertex set $\mathcal{V}$ is a
finite collection of subsets $F$ of $\mathcal{V}$ having the properties
:\smallskip\newline(i) If $v\in\mathcal{V}$, then $\left\{  v\right\}
\in\mathcal{S}\left(  \mathcal{V}\right)  $, and (ii) if $F\in\mathcal{S}%
\left(  \mathcal{V}\right)  $ and $F^{\prime}\subset F$, then $F^{\prime}%
\in\mathcal{S}\left(  \mathcal{V}\right)  $.\smallskip\newline The elements
$F\in\mathcal{S}\left(  \mathcal{V}\right)  $ are called \textit{abstract
simplices }or \textit{faces}. For an $F\in\mathcal{S}\left(  \mathcal{V}%
\right)  $ one defines dim$\left(  F\right)  :=\sharp\left(  F\right)  -1$ and
dim$\left(  \mathcal{S}\left(  \mathcal{V}\right)  \right)  :=$ max$\left\{
\text{dim}\left(  F\right)  \left\vert F\in\mathcal{S}\left(  \mathcal{V}%
\right)  \right.  \right\}  $ as the \textit{dimension }of $\mathcal{S}\left(
\mathcal{V}\right)  .$ (If dim$\left(  F\right)  =k+1$, then $F$ is said to be
an \textit{abstract }$k$-\textit{simplex }or a $k$-\textit{face}). A
\textit{subcomplex} of $\mathcal{S}\left(  \mathcal{V}\right)  $ is an
abstract simplicial complex contained in $\mathcal{S}\left(  \mathcal{V}%
\right)  $ whose vertex-set is a subset of $\mathcal{V}$. (Sometimes, for $F$
a subset of $\mathcal{V}$, one denotes the abstract simplex with vertex set
$F$ by $2^{F}$.)

\begin{definition}
Let $\mathcal{S}\left(  \mathcal{V}\right)  $ be an abstract simplicial
complex with vertex set $\mathcal{V}$ and $\iota:\mathcal{V}\rightarrow
\mathbb{R}^{d}$ an injective map, such that\smallskip\newline(i) the elements
of the images $\iota\left(  F\right)  $ are affinely independent for all
$F\in\mathcal{S}\left(  \mathcal{V}\right)  $, and\smallskip\newline(ii)
int$\left(  \text{conv}\left(  \iota\left(  F\right)  \right)  \right)  \cap$
int$\left(  \text{conv}\left(  \iota\left(  F^{\prime}\right)  \right)
\right)  =\varnothing$, for all $F,F^{\prime}\in\mathcal{S}\left(
\mathcal{V}\right)  $, $F\neq F^{\prime}$.\smallskip\newline Then
$\bigcup_{F\in\mathcal{S}\left(  \mathcal{V}\right)  }\ $int$\left(
\text{conv}\left(  \iota\left(  F\right)  \right)  \right)  $ is called a
\textit{geometric realization} of $\mathcal{S}\left(  \mathcal{V}\right)  $
w.r.t. $\iota$.
\end{definition}

\noindent{}Geometric realizations always exist, the underlying spaces of any
two geometric realizations of $\mathcal{S}\left(  \mathcal{V}\right)  $ are
homeomorphic to each other, and therefore any \textit{support }$\left\vert
\mathcal{S}\left(  \mathcal{V}\right)  \right\vert $ \textquotedblleft
realizing\textquotedblright\ $\mathcal{S}\left(  \mathcal{V}\right)  $ is
well-defined\ in the topological category. In many cases, we shall denote both
geometric and abstract simplicial complexes by the letter $\mathcal{S}$. When,
for some reason, our intention is to stress what kind of complexes is meant
(if this is not clear from the context), and our starting point is a geometric
simplicial complex $\mathcal{S}$, we denote by $\mathcal{S}^{\text{abs}}$ the
corresponding abstract simplicial complex; and conversely, when our starting
point is an abstract simplicial complex $\mathcal{S}$, we consider a fixed
realization $\left\vert \mathcal{S}\right\vert $. Moreover, we mostly use
$\mathcal{V}$ and vert$\left(  \mathcal{S}\right)  $ (for the vertex set)
interchangeably.\medskip\

\noindent{}$\bullet$ $\mathbf{f}$- \textbf{and} $\mathbf{h}$-\textbf{vectors}.
The\textbf{ }$\mathbf{f}$-\textit{vector} $\mathbf{f}\left(  \mathcal{S}%
\right)  =\left(  \mathfrak{f}_{-1}\left(  \mathcal{S}\right)  ,\mathfrak{f}%
_{0}\left(  \mathcal{S}\right)  ,\ldots,\mathfrak{f}_{d}\left(  \mathcal{S}%
\right)  \right)  $ of a $d$-dimensional geometric or abstract simplicial
complex $\mathcal{S}$ is defined by setting
\[
\mathfrak{f}_{-1}\left(  \mathcal{S}\right)  :=-1,\text{ and\ }\mathfrak{f}%
_{i}\left(  \mathcal{S}\right)  :=\sharp\{i\text{-dimensional faces of
}\mathcal{S}\},\ \forall i\in\{0,\ldots,d\}.
\]
The\emph{ }$\mathbf{h}$-\textit{vector } $\mathbf{h}\left(  \mathcal{S}%
\right)  =\left(  \mathfrak{h}_{0}\left(  \mathcal{S}\right)  ,\mathfrak{h}%
_{1}\left(  \mathcal{S}\right)  ,\ldots,\mathfrak{h}_{d+1}\left(
\mathcal{S}\right)  \right)  $ of $\mathcal{S}$ is defined by the equation%
\begin{equation}
\mathbf{h}\left(  \mathcal{S};t\right)  =\left(  1-t\right)  ^{d}%
\ \mathbf{f}(\mathcal{S};\tfrac{t}{1-t}), \label{HFEQ}%
\end{equation}
where\emph{ }$\mathbf{f}\left(  \mathcal{S};t\right)  :=%
{\textstyle\sum\limits_{i=0}^{d+1}} \mathfrak{f}_{i-1}\left(
\mathcal{S}\right)  t^{i}\in\mathbb{Z}_{\geq
0\,}\left[  t\right]  \,\,$and $\mathbf{h}\left(  \mathcal{S};t\right)  :=%
{\textstyle\sum\limits_{i=0}^{d+1}} \mathfrak{h}_{i}\left(
\mathcal{S}\right)  t^{i}\in\mathbb{Z}\left[
t\right]  .$ Note, in particular, that (\ref{HFEQ}) gives%
\begin{equation}
\mathfrak{f}_{j-1}\left(  \mathcal{S}\right)  =%
{\textstyle\sum\limits_{i=0}^{j}} \ \ \tbinom{d-i}{d-j}\
\mathfrak{h}_{i}\left(  \mathcal{S}\right)  ,\ \forall
j,\ 0\leq j\leq d+1. \label{F-H-FORMULA}%
\end{equation}
The $\mathbf{f}$-vector $\mathbf{f}\left(  P\right)  $ of a simplicial
polytope $P$ is by definition the $\mathbf{f}$-vector $\mathbf{f}\left(
\mathcal{S}_{\partial P}\right)  $ of $\mathcal{S}_{\partial P}$ (as in
Example \ref{Boundarycomplex}).\medskip\

\noindent{}$\bullet$ \textbf{UBT for the facets of simplicial balls}. We
denote by CycP$_{d}\left(  k\right)  $ the \textit{cyclic} $d$%
-\textit{polytope} with $k$ vertices. As it is known, the number of its facets
equals
\begin{equation}
\mathfrak{f}_{d-1}\left(  \text{CycP}_{d}\left(  k\right)  \right)
=\tbinom{k-\left\lceil \frac{d}{2}\right\rceil }{\left\lfloor \frac{d}%
{2}\right\rfloor }+\tbinom{k-1-\left\lceil \frac{d-1}{2}\right\rceil
}{\left\lfloor \frac{d-1}{2}\right\rfloor }. \label{FDCYCLIC}%
\end{equation}

\noindent This is due to Gale's evenness condition and to the fact that
CycP$_{d}\left(  k\right)  $ is $\left\lfloor d/2\right\rfloor $-neighbourly
(cf. \cite[p. 24]{Ziegler}). Let us first recall the classical UBT and LBT for
simplicial \textit{spheres}, and then explain how one obtains an UBT for the
facets of simplicial \textit{balls}.

\begin{theorem}
[{Upper Bound Theorem for Simplicial Spheres, \cite[II.4.5]{Stanley2}}%
]\label{UBTSS}\negthinspace\negthinspace The $\mathbf{f}$-vector coordinates
of a simplicial $\left(  d-1\right)  $-sphere $\mathcal{S}$ with
$\mathfrak{f}_{0}\left(  \mathcal{S}\right)  =k$ vertices satisfy the
following inequalities\emph{:}
\[
\mathfrak{f}_{i}\left(  \mathcal{S}\right)  \leq\mathfrak{f}_{i}\left(
\emph{CycP}_{d}\left(  k\right)  \right)  \ ,\ \forall\ i\ ,\ 0\leq i\leq
d-1.
\]

\end{theorem}

\begin{theorem}
[Lower Bound Theorem for Simplicial Spheres, \cite{KALAI}]\label{LBTH}The
$\mathbf{h}$-vector coordinates of a simplicial $\left(  d-1\right)  $-sphere
$\mathcal{S}$ with $\mathfrak{f}_{0}\left(  \mathcal{S}\right)  =k$ vertices
satisfy the following inequalities\emph{:}
\begin{equation}
\mathfrak{h}_{1}\left(  P\right)  =k-d\leq\mathfrak{h}_{i}\left(  P\right)
\ ,\ \forall\ i\ ,\ 2\leq i\leq d. \label{LBTFORMULA}%
\end{equation}

\end{theorem}

\begin{theorem}
[UBT for the Facets of Simplicial Balls, \cite{DAIS}]\label{SECUBB}Let
$\mathcal{S}$ be a simplicial $d$-ball with $\mathfrak{f}_{0}\left(
\mathcal{S}\right)  =\mathfrak{b}$ vertices. Suppose that $\mathfrak{f}%
_{0}\left(  \partial\mathcal{S}\right)  =\mathfrak{b}^{\prime}$.
Then\emph{:\medskip}
\begin{equation}
\mathfrak{f}_{d}\left(  \mathcal{S}\right)  \leq\mathfrak{f}_{d}\left(
\emph{CycP}_{d+1}\left(  \mathfrak{b}\right)  \right)  -\left(  \mathfrak{b}%
^{\prime}-d\right)  . \label{UP-LO2}%
\end{equation}

\end{theorem}

\noindent{}\textsc{Sketch of Proof}. Introduce the auxiliary vector
$\widetilde{\mathbf{h}}\left(  \mathcal{S}\right)  =(\widetilde{\mathfrak{h}%
}_{0}\left(  \mathcal{S}\right)  ,\ldots,\widetilde{\mathfrak{h}}_{d+1}\left(
\mathcal{S}\right)  )$ with{\small
\[
\widetilde{\mathfrak{h}}_{i}\left(  \mathcal{S}\right)  :=\left\{
\begin{array}
[c]{lll}%
\mathfrak{h}_{i}\left(  \mathcal{S}\right)  & , & \text{for \ \ \ }0\leq
i\leq\left\lfloor \frac{d+1}{2}\right\rfloor ,\\
\  &  & \\
\mathfrak{h}_{i}\left(  \mathcal{S}\right)  -\left(  \mathfrak{h}_{d-i}\left(
\partial\mathcal{S}\right)  -\mathfrak{h}_{d+1-i}\left(  \partial
\mathcal{S}\right)  \right)  & , & \text{for \ \ \ }\left\lfloor \frac{d+1}%
{2}\right\rfloor +1\leq i\leq d+1.
\end{array}
\right.
\]
}By \cite[Thm. 4.3, p. 136]{Schenzel} we see that%
\[
\mathfrak{h}_{i}\left(  \mathcal{S}\right)  \leq\mathfrak{h}_{i}\left(
\text{CycP}_{d+1}\left(  \mathfrak{b}\right)  \right)  =\tbinom{\mathfrak{b}%
-d+i-2}{i},\,\,\,\forall i,\,\,\,0\leq i\leq\left\lfloor \tfrac{d+1}%
{2}\right\rfloor .
\]
On the other hand, using Stanley's \textquotedblleft$\mathfrak{h}$ of
$\partial$\textquotedblright-Lemma \cite[Lemma 2.3, p. 253]{Stanley5}, i.e.,
\[
\mathfrak{h}_{i-1}\left(  \partial\mathcal{S}\right)  -\mathfrak{h}_{i}\left(
\partial\mathcal{S}\right)  =\mathfrak{h}_{\left(  d+1\right)  -i}\left(
\mathcal{S}\right)  -\mathfrak{h}_{i}\left(  \mathcal{S}\right)  ,\ \ \forall
i,\ 0\leq i\leq d+1,\smallskip
\]
we get $\widetilde{\mathfrak{h}}_{i}\left(  \mathcal{S}\right)  =\widetilde
{\mathfrak{h}}_{\left(  d+1\right)  -i}\left(  \mathcal{S}\right)
,\,\,\forall i\in\{0,1,\ldots,d+1\},$ and therefore
\begin{equation}
\widetilde{\mathfrak{h}}_{i}\left(  \mathcal{S}\right)  \leq\mathfrak{h}%
_{i}\left(  \text{CycP}_{d+1}\left(  \mathfrak{b}\right)  \right)
,\,\,\,\forall i,\,\,\,0\leq i\leq d+1. \label{INEQHAHAT}%
\end{equation}
Passing to the $\mathbf{f}$-vector, and using Dehn-Sommerville relations for
$\mathbf{h}(\partial\mathcal{S}),$ we verify easily via (\ref{INEQHAHAT})
that
\begin{equation}
\mathfrak{f}_{i}\left(  \mathcal{S}\right)  \leq\mathfrak{f}_{i}\left(
\text{CycP}_{d+1}\left(  \mathfrak{b}\right)  \right)  -\sum\limits_{j=d-i}%
^{\left\lfloor \frac{d}{2}\right\rfloor \smallskip}\ \tbinom{j}{d-i}\ \left(
\mathfrak{h}_{j}\left(  \partial\mathcal{S}\right)  -\mathfrak{h}_{j-1}\left(
\partial\mathcal{S}\right)  \right)  . \label{WEAKUBTSB}%
\end{equation}
For $i=d$, (\ref{WEAKUBTSB}) gives{\small
\begin{align*}
\mathfrak{f}_{d}\left(  \mathcal{S}\right)   &  \leq\mathfrak{f}_{d}\left(
\text{CycP}_{d+1}\left(  \mathfrak{b}\right)  \right)  -\sum\limits_{j=0}%
^{\left\lfloor \frac{d}{2}\right\rfloor \smallskip}\ \left(  \mathfrak{h}%
_{j}\left(  \partial\mathcal{S}\right)  -\mathfrak{h}_{j-1}\left(
\partial\mathcal{S}\right)  \right) \\
&  =\mathfrak{f}_{d}\left(  \text{CycP}_{d+1}\left(  \mathfrak{b}\right)
\right)  -\mathfrak{h}_{\left\lfloor \frac{d}{2}\right\rfloor }\left(
\partial\mathcal{S}\right)  \leq\mathfrak{f}_{d}\left(  \text{CycP}%
_{d+1}\left(  \mathfrak{b}\right)  \right)  -\mathfrak{h}_{1}\left(
\partial\mathcal{S}\right)  ,
\end{align*}
}where the latter inequality comes from (\ref{LBTFORMULA}) for the simplicial
sphere $\partial\mathcal{S}$. Now obviously, $\mathfrak{h}_{1}\left(
\partial\mathcal{S}\right)  =\mathfrak{b}^{\prime}-d$.\hfill$\square\medskip$

\noindent$\bullet$ \textbf{Subdivisions}.\textsf{ }There are several kinds of
\textit{subdivisions }of simplicial and polytopal complexes which correspond
to the different distinctive features of the partitioning objects and depend
on the way they have to fit together. (Here we follow Stanley's terminology
from \cite{Stanley4, Stanley2}.)

\begin{definition}
Let $\mathcal{S}$ be an abstract simplicial complex.\emph{\ }A
\textit{topological simplicial subdivision }of $\mathcal{S}$\emph{\ }is a pair
$\left(  \mathcal{S}^{\prime},\varphi\right)  $ consisting of an
abstract\emph{\ }simplicial complex $\mathcal{S}^{\prime}$ together with a map
$\varphi:\mathcal{S}^{\prime}\rightarrow\mathcal{S}$ satisfying the following
conditions:\smallskip\newline(i) For every $F\in\mathcal{S}$\emph{,
}the\emph{\ }\textit{restriction} $\mathcal{S}_{F}^{\prime}:=\varphi
^{-1}\left(  2^{F}\right)  $ of $\mathcal{S}^{\prime}$ to $F$ is a subcomplex
of\emph{\ }$\mathcal{S}^{\prime}$ having a geometric realization $\left\vert
\mathcal{S}_{F}\right\vert $ which is a $\left(  \text{dim}\left(  F\right)
\right)  $-ball\emph{.\smallskip}\newline(ii) For each\emph{\ }$F^{\prime}%
\in\mathcal{S}^{\prime}$, $F=\varphi\left(  F^{\prime}\right)  \in\mathcal{S}$
if and only if $F^{\prime}$ is an \textit{interior face }of\emph{\ }%
$\mathcal{S}_{F}^{\prime}$\emph{, }i.e., if and only if $\left\vert F^{\prime
}\right\vert $ is not contained in the boundary $\partial\left\vert
\mathcal{S}_{F}^{\prime}\right\vert =\left\vert \mathcal{S}_{F}^{\prime
}\right\vert \smallsetminus$ int$\left(  \left\vert \mathcal{S}_{F}^{\prime
}\right\vert \right)  $.\smallskip\newline(Often one calls $\mathcal{S}%
^{\prime}$ a topological simplicial subdivision of $\mathcal{S}$ and omits
$\varphi$ if it is self-evident from the context.)\emph{\ }
\end{definition}

\begin{definition}
\label{QUASI-SUB}Let $\mathcal{S}$ be an abstract simplicial complex. A
topological simplicial subdivision\textit{\ }$\left(  \mathcal{S}^{\prime
},\varphi\right)  $ of $S$\ is called\emph{\ }\textit{quasigeometric} if for
every face $F^{\prime}$ of $\mathcal{S}^{\prime}$ there does not exist a face
$F\in\mathcal{S}$ for which\smallskip\ \emph{\ }\newline(i) dim$\left(
F\right)  <$ dim$\left(  F^{\prime}\right)  $ and\smallskip\ \newline(ii) each
vertex $v$ of $F^{\prime}$ lies on some subset of $F$ (depending on $v$).
\end{definition}

\begin{definition}
\label{GEOMTRIANG}Let $\mathcal{S}$ be a geometric simplicial complex
and\emph{\ }$\mathcal{S}^{\text{abs}}$\emph{\ }the corresponding abstract
simplicial complex. A geometric simplicial complex\emph{\ }$\mathcal{S}%
^{\prime}$ is called a \textit{geometric simplicial subdivision }or a
\textit{geometric triangulation }or simply a \textit{triangulation }of
$\mathcal{S}$ (and, respectively, $\mathcal{S}^{\prime\,\text{\emph{abs}}}$ a
\textit{geometric simplicial subdivision }or\emph{\ }\textit{a geometric
triangulation }of $\mathcal{S}^{\text{\emph{abs}}}$) if\emph{\ }$\left\vert
\mathcal{S}\right\vert =\left\vert \mathcal{S}^{\prime}\right\vert $ and every
simplex in $\mathcal{S}^{\prime}$ is contained in some simplex in
$\mathcal{S}.$
\end{definition}

\begin{definition}
\label{POLSUBDIV}Let $\mathcal{S}$ be a polytopal complex. A polytopal
complex\emph{\ }$\mathcal{S}^{\prime}$ is called a \textit{polytopal
subdivision } of $\mathcal{S}$\emph{\ }if $\left\vert \mathcal{S}\right\vert
=\left\vert \mathcal{S}^{\prime}\right\vert $ and every polytope in
$\mathcal{S}^{\prime}$ is contained in some polytope in $\mathcal{S}.$ If
$\mathcal{S}^{\prime}$ is a simplicial complex, then we again say that
$\mathcal{S}^{\prime}$ is a \textit{geometric simplicial subdivision }or
\textit{a geometric triangulation }or simply a \textit{triangulation
}\emph{of} $\mathcal{S}$.
\end{definition}

\begin{note}
Quasigeometric simplicial subdivisions are geometric (because of the affine
independence of the vertices of a geometric simplex) but the converse is not
always true. Moreover, not every topological simplicial subdivision is
quasigeometric. For counterexamples we refer to \cite[p. 468]{Chan2} and
\cite[p. 814]{Stanley4}\emph{.}
\end{note}

\noindent$\bullet$ \textbf{Working in the PL-category}.\textsf{ }For the proof
of the Upper Bound Theorem \ref{BIGTH} we shall make use of notions and
propositions from \textquotedblleft PL-topology\textquotedblright. Working in
the PL-category, i.e., in the category of simplicial complexes with piecewise
linear maps (as morphisms), we can take full advantage of the fact that most
of the theoretic arguments do not depend on specific geometric realizations
and that many topological operations with PL balls or PL spheres (like
starring, linking, gluing etc.) produce again other PL balls or PL spheres
(i.e., something which is by no means true in general within the \textit{usual
}topological category). Standard references for the \textquotedblleft
PL-topology\textquotedblright\ are the books of Glaser \cite{Glaser}, Hudson
\cite{Hud}, and Rourke \& Sanderson \cite{Rou-Sa}.

\begin{definition}
Let $\mathcal{S}_{1}$\emph{\ }and $\mathcal{S}_{2}$ be two abstract simplicial
complexes. A \textit{simplicial map }$\varphi:\mathcal{S}_{1}\rightarrow
\mathcal{S}_{2}$ is a function\emph{\ }$\varphi:$ vert$\left(  \mathcal{S}%
_{1}\right)  \rightarrow$ vert$\left(  \mathcal{S}_{2}\right)  $, such that
whenever $\left\{  v_{0},\ldots,v_{k}\right\}  $\ is an (abstract) simplex,
then $\left\{  \varphi\left(  v_{0}\right)  ,\ldots,\varphi\left(
v_{k}\right)  \right\}  $ is an (abstract) simplex too. If such a $\varphi$
is, in addition, a homeomorphism\emph{\ }(i.e. $\left\vert \varphi\right\vert
\left(  \left\vert \mathcal{S}_{1}\right\vert \right)  \approx\left\vert
\mathcal{S}_{2}\right\vert $), then\emph{\ }$\varphi$\emph{\ }is called a
\textit{simplicial homeomorphism}\emph{. }Passing to geometric realizations, a
simplicial map $\left\vert \varphi\right\vert :\left\vert \mathcal{S}%
_{1}\right\vert \rightarrow\left\vert \mathcal{S}_{2}\right\vert $%
\emph{\ }carries the vertices\emph{\textit{\ }}of\emph{\ }$\mathcal{S}_{1}$ to
the vertices of $\mathcal{S}_{2}$ and the geometrically realized simplices of
$\mathcal{S}_{1}$ \textit{linearly }onto those of\emph{\ }$\mathcal{S}_{2}$.
(For $\left\vert \varphi\right\vert $ to be \textit{linear } means that for
each point $v\in\left\vert \mathcal{S}_{1}\right\vert $\emph{,} which is
uniquely expressible as a convex linear combination $v=\sum_{i=0}^{k}%
\lambda_{i}v_{i}$ by the \textit{barycentric coordinates }$\lambda_{0}%
,\ldots,\lambda_{k}$ with\emph{\ }respect to an ambient simplex conv$\left(
\left\{  v_{0},\ldots,v_{k}\right\}  \right)  $\emph{, }one has $\left\vert
\varphi\right\vert \left(  v\right)  =\sum_{i=0}^{k}\lambda_{i}\,\left\vert
\varphi\right\vert \left(  v_{i}\right)  $).
\end{definition}

\begin{definition}
Let\emph{\ }$\mathcal{S}_{1}$ and $\mathcal{S}_{2}$ be two abstract simplicial
complexes. A map $\varphi:\mathcal{S}_{1}\rightarrow\mathcal{S}_{2}$%
\emph{\ }is called \textit{piecewise linear }(or PL \textit{map}) if for
some\emph{\ }geometric realizations of $\mathcal{S}_{1}$ and\emph{\ }%
$\mathcal{S}_{2}$\emph{, }the corresponding map $\left\vert \varphi\right\vert
:\left\vert \mathcal{S}_{1}\right\vert \rightarrow\left\vert \mathcal{S}%
_{2}\right\vert $ satisfies anyone of the following equivalent
conditions:\emph{\smallskip}\newline(i) There exist geometric triangulations
$\mathcal{S}_{1}^{\prime}$, $\mathcal{S}_{2}^{\prime}$ of the complexes
$\mathcal{S}_{1}$ and $\mathcal{S}_{2}$\emph{, }respectively, relative to
which $\left\vert \varphi\right\vert :\left\vert \mathcal{S}_{1}^{\prime
}\right\vert \rightarrow\left\vert \mathcal{S}_{2}^{\prime}\right\vert $ is
simplicial.\smallskip\newline(ii) There is a geometric triangulation
$\mathcal{S}_{1}^{\prime}$ \emph{of }$\mathcal{S}_{1}$\emph{, }relative to
which $\left\vert \varphi\right\vert :\left\vert \mathcal{S}_{1}^{\prime
}\right\vert \rightarrow\left\vert \mathcal{S}_{2}\right\vert $ is
linear.\smallskip\newline(It can be shown that this definition depends neither
on the particular choice of the geometric realizations of $S_{1}$ and $S_{2}%
$\ nor on the particular choice of the geometric triangulations in (i), (ii).
Note that a simplicial map is a PL map but the converse is not always true).
\end{definition}

\begin{definition}
(i) Two abstract simplicial complexes $\mathcal{S}_{1}$\emph{,} $\mathcal{S}%
_{2}$ are called \textit{PL} \textit{homeomorphic }(denoted by $\mathcal{S}%
_{1}\underset{\text{PL}}{\approx}\mathcal{S}_{2}$) if there is a PL
map\emph{\ }$\varphi:\mathcal{S}_{1}\rightarrow\mathcal{S}_{2}$ which is also
a homeomorphism\emph{.\smallskip}\newline(ii) In particular, an abstract
$k$-dimensional simplicial complex $\mathcal{S}$ is called a\emph{\ }%
(\textit{simplicial}) PL $k$-\textit{ball }(resp., PL $k$-\textit{sphere}) if
$\mathcal{S}$ is PL homeomorphic to\emph{\ }the\emph{\ }$k$-simplex (resp., to
the boundary of the\emph{\ }$\left(  k+1\right)  $-simplex).\smallskip
\newline[Geometric triangulations of topological spheres (resp., balls) are
not necessarily\emph{\ }PL spheres (resp., PL balls). It is well-known, for
instance, that all\emph{\ }geometric triangulations of a\emph{\ }$k$-sphere
are PL\emph{\ }$k$-spheres for\emph{\ }$k\leq3$, whereas there exist non-PL
geometric triangulations of a $k$-sphere for\emph{\ }$k\geq5$.]\smallskip
\newline(iii) A regular cell complex\emph{\ }$\mathcal{S}$ is called a PL
$k$-\textit{ball }(resp., a PL $k$-\textit{sphere}) if the (simplicial) order
complex of its face poset is a simplicial PL\emph{\ }$k$-ball\emph{\ }(resp.,
$k$-sphere).\emph{\ }
\end{definition}

\noindent$\bullet$ \textbf{Joins, stars and links}. Let $\mathcal{S}$ be an
abstract simplicial complex, $v\in$ vert$\left(  \mathcal{S}\right)  $, and
$F$ a face of $\mathcal{S}$. For $v\notin$ vert$\left(  F\right)  $, $v\ast F$
is defined to be the simplex with vertex set vert$\left(  F\right)
\cup\left\{  v\right\}  $, i.e. the so-called \textit{join} of $v$ with $F$.
In general, if $\mathcal{S}_{1}$, $\mathcal{S}_{2}$ are two abstract
simplicial complexes on disjoint vertex sets $\mathcal{V}_{1}$, $\mathcal{V}%
_{2}$, the \textit{join} of $\mathcal{S}_{1}$ and $\mathcal{S}_{2}$ is defined
to be the simplicial complex
\[
\mathcal{S}_{1}\ast\mathcal{S}_{2}:=\left\{  F\in\mathcal{V}_{1}%
\cup\mathcal{V}_{2}\ \left\vert \ F\cap\mathcal{V}_{1}\in\mathcal{S}_{1}\text{
and }F\cap\mathcal{V}_{2}\in\mathcal{S}_{2}\right.  \right\}
\]
For $w\notin$ vert$\left(  \mathcal{S}\right)  $, $w\ast\mathcal{S}$ is
nothing but the simplicial complex whose faces are $\left\{  \varnothing
\right\}  \cup\left\{  w\ast F\ \left\vert \ F\in\mathcal{S}\right.  \right\}
\cup\mathcal{S}$, i.e. the \textit{cone }(with apex $w$) over $\mathcal{S}$.
If $w^{\prime}\notin\mathcal{S}$ and $w^{\prime}$ is different from $w$, then
the double joining
\[
\left\{  w,w^{\prime}\right\}  \ast\mathcal{S}:=w\ast\left(  w^{\prime}%
\ast\mathcal{S}\right)
\]
is the \textit{suspension} of $\mathcal{S}$ w.r.t. the additional vertices
$w$, $w^{\prime}$.\medskip\newline$\bullet$ For $v\in$ vert$\left(
\mathcal{S}\right)  $, let
\[
\left\{
\begin{array}
[c]{l}%
\text{star}_{v}\left(  \mathcal{S}\right)  :=\left\{  F\in\mathcal{S\ }%
\left\vert \ v\in\text{vert}\left(  F\right)  \right.  \right\}  ,\smallskip\\
\\
\text{ast}_{v}\left(  \mathcal{S}\right)  :=\left\{  F\in\mathcal{S\ }%
\left\vert \ v\notin\text{vert}\left(  F\right)  \right.  \right\}
,\smallskip\\
\\
\overline{\text{star}_{v}}\left(  \mathcal{S}\right)  :=\left\{  \text{
}F^{\prime}\in\mathcal{S\ }\left\vert
\begin{array}
[c]{c}%
\ F^{\prime}\text{ faces of all }F\in\mathcal{S}\text{ }\\
\text{for which }v\in\text{vert}\left(  F\right)
\end{array}
\right.  \right\}  ,\smallskip\smallskip\\
\\
\overline{\text{ast}_{v}}\left(  \mathcal{S}\right)  :=\left\{  \text{
}F^{\prime}\in\mathcal{S\ }\left\vert
\begin{array}
[c]{c}%
\ F^{\prime}\text{ faces of all }F\in\mathcal{S}\text{ }\\
\text{for which }v\notin\text{vert}\left(  F\right)
\end{array}
\right.  \right\}  ,\smallskip\smallskip\\
\\
\text{link}_{v}\left(  \mathcal{S}\right)  :=\left\{  F\in\mathcal{S\ }%
\left\vert \ v\notin\text{vert}\left(  F\right)  ,\ v\ast F\in\mathcal{S}%
\right.  \right\}  ,\smallskip
\end{array}
\right.
\]
denote the \textit{star}, the \textit{antistar}, the \textit{closed}
\textit{star}, the \textit{closed antistar,} and the \textit{link }of $v$ in
$\mathcal{S}$, respectively. The last three form subcomplexes of $\mathcal{S}$
and are related as follows:%
\[
\overline{\text{star}_{v}}\left(  \mathcal{S}\right)  \cap\overline
{\text{ast}_{v}}\left(  \mathcal{S}\right)  =\text{ link}_{v}\left(
\mathcal{S}\right)  ,\text{ and\smallskip\ }\overline{\text{star}_{v}}\left(
\mathcal{S}\right)  =v\ast\text{link}_{v}\left(  \mathcal{S}\right)  \text{.}%
\]

\begin{proposition}
\label{PL-prop}For \emph{simplicial} \emph{PL }spheres and \emph{PL }balls
$\mathcal{S}$ we have the following implications\emph{:}\smallskip
\smallskip\newline$%
\begin{array}
[c]{llll}%
\text{\emph{(i)}} & \mathcal{S}\text{ is a \emph{PL }}k\text{-ball} &
\Longrightarrow & \partial\mathcal{S}\text{ is a \emph{PL }}\left(
k-1\right)  \text{-sphere\smallskip\smallskip}\\
&  &  & \\
\text{\emph{(ii)}} & \left.
\begin{array}
[c]{l}%
\mathcal{S}\text{ is a \emph{PL }}k\text{-ball}\\
\text{\emph{(}or a \emph{PL }}k\text{-sphere\emph{)}}\\
\text{and }w\notin\mathcal{S}%
\end{array}
\right\}  & \Longrightarrow & \left\{
\begin{array}
[c]{l}%
\text{the cone }w\ast\mathcal{S}\text{ is a \emph{PL }}\left(  k+1\right)
\text{-ball }\\
\text{with boundary }=\mathcal{S}\cup\left(  w\ast\partial\mathcal{S}\right)
\end{array}
\right.  \smallskip\smallskip\\
&  &  & \\
\text{\emph{(iii)}} & \mathcal{S}\text{ is a \emph{PL }}k\text{-sphere} &
\Longrightarrow & \text{\emph{link}}_{v}\left(  \mathcal{S}\right)  \text{ is
a \emph{PL }}\left(  k-1\right)  \text{-sphere\smallskip\smallskip}\\
&  &  & \\
\text{\emph{(iv)}} & \mathcal{S}\text{ is a \emph{PL }}k\text{-sphere} &
\Longrightarrow & \text{\smallskip\smallskip}\overline{\text{\emph{star}}_{v}%
}\left(  \mathcal{S}\right)  \text{ is a \emph{PL }}k\text{-ball}\\
&  &  & \\
\text{\emph{(v)}} & \mathcal{S}\text{ is a \emph{PL }}k\text{-sphere} &
\Longrightarrow & \overline{\text{\emph{ast}}_{v}}\left(  \mathcal{S}\right)
\text{ is a \emph{PL }}k\text{-ball\smallskip\smallskip}\\
&  &  & \\
\emph{(vi)} & \left.
\begin{array}
[c]{l}%
\mathcal{S}\text{ is a \emph{PL }}k\text{-sphere,}\\
v\in\text{\emph{vert}}\left(  \mathcal{S}\right)
\end{array}
\right\}  & \Longrightarrow & \left\{
\begin{array}
[c]{l}%
v\ast\overline{\text{\emph{ast}}_{v}}\left(  \mathcal{S}\right)  \text{ is a
\emph{PL }}\left(  k+1\right)  \text{-ball }\\
\text{with boundary }=\mathcal{S}%
\end{array}
\right.  \smallskip\smallskip
\end{array}
$
\end{proposition}

\noindent\textsc{Proof}\textit{. }(i) This is obvious because the boundary of
any simplicial subdivision of $\mathcal{S}$ is the restriction of this
subdivision to the boundaries of the participating simplices.\smallskip
\newline(ii) Since, in general, the join of two abstract simplicial complexes
is PL homeomorphic to the join of the images of these complexes under any PL
homeomorphisms (see e.g. \cite[ II.5, p. 22, and II.17, p. 41]{Glaser}), we
have $w\ast\mathcal{S}=$ $w\ast\left(  \text{a }k\text{-simplex}\right)  $ (or
$=w\ast\left(  \text{the boundary of a }\left(  k+1\right)  \text{-simplex}%
\right)  $).\smallskip\newline(iii) See \cite[proof of the Corollary 1.16, p.
24]{Hud}.\smallskip\newline(iv) By (iii) and (ii) we get%
\[
\overline{\text{star}_{v}}\left(  \mathcal{S}\right)  =v\ast\text{link}%
_{v}\left(  \mathcal{S}\right)  \underset{\text{PL}}{\approx}v\ast\left(
\text{the boundary of a }k\text{-simplex}\right)  \underset{\text{PL}}%
{\approx}\left(  \text{a }k\text{-simplex}\right)  \text{.}%
\]
\newline(v) Since $\overline{\text{ast}_{v}}\left(  \mathcal{S}\right)
=\mathcal{S}\smallsetminus$ star$_{v}\left(  \mathcal{S}\right)  $ and
star$_{v}\left(  \mathcal{S}\right)  $ has a PL $k$-ball (and consequently a
stellar $k$-ball) as closure (by (iv)), $\overline{\text{ast}_{v}}\left(
\mathcal{S}\right)  $ is a stellar $k$-ball according to \cite[II.15, p.
37]{Glaser}. From the equivalence of stellar- and PL-homeomorphism-property
(cf. \cite[II.17, p. 41]{Glaser}) we conclude that $\overline{\text{ast}_{v}%
}\left(  \mathcal{S}\right)  $ is indeed a PL $k$-ball.\smallskip\newline(vi)
That $v\ast\overline{\text{ast}_{v}}\left(  \mathcal{S}\right)  $ is a
PL\emph{\ }$\left(  k+1\right)  $-ball follows from (ii). Its boundary equals
\[
\overline{\text{ast}_{v}}\left(  \mathcal{S}\right)  \cup\left(  v\ast
\partial\left(  \overline{\text{ast}_{v}}\left(  \mathcal{S}\right)  \right)
\right)  =\overline{\text{ast}_{v}}\left(  \mathcal{S}\right)  \cup
\overline{\text{star}_{v}}\left(  \mathcal{S}\right)  =\mathcal{S}%
\]
and the proof is completed.\hfill{}$\square$

\begin{corollary}
\label{AST-LINK}Let $\mathcal{S}$ be a \emph{(}simplicial\emph{)} \emph{PL
}$k$-ball, $k\geq2$, $v\in$ \emph{vert}$\left(  \partial\mathcal{S}\right)  $
a vertex of its boundary, and $w\in\emph{vert}\left(  \text{\emph{link}}%
_{v}\left(  \mathcal{S}\right)  \cap\partial\mathcal{S}\right)  $ a boundary
vertex of its link. Then the suspension
\[
\mathcal{S}_{v,w}:=\left\{  v,w\right\}  \ast\overline{\text{\emph{ast}}_{w}%
}\left(  \text{\emph{link}}_{v}\left(  \partial\mathcal{S}\right)  \right)
\]
is a \emph{(}simplicial\emph{) PL }$k$-ball and \smallskip\smallskip
$\overline{\text{\emph{star}}_{v}}\left(  \partial\mathcal{S}\right)  $ is a
subcomplex of the boundary $\partial\mathcal{S}_{v,w}.$
\end{corollary}

\noindent\textsc{Proof}\textit{. }By Proposition \ref{PL-prop},
\[%
\begin{array}
[c]{ll}%
\partial\mathcal{S}\text{ is a PL\emph{\ }}\left(  k-1\right)  \text{-sphere}
& \ \ \text{(using (i)),}\\
\text{link}_{v}\left(  \partial\mathcal{S}\right)  \text{is a PL\emph{\ }%
}\left(  k-2\right)  \text{-sphere} & \ \ \text{(using (iii)),\smallskip}\\
\overline{\text{ast}_{w}}\left(  \text{link}_{v}\left(  \partial
\mathcal{S}\right)  \right)  \text{ is a PL\emph{\ }}\left(  k-2\right)
\text{-ball} & \ \ \text{(using (v)),}\\
w\ast\overline{\text{ast}_{w}}\left(  \text{link}_{v}\left(  \partial
\mathcal{S}\right)  \right)  \text{ is a PL\emph{\ }}\left(  k-1\right)
\text{-ball} & \ \ \text{(using (ii)).}%
\end{array}
\]
Hence, the first claim for $\mathcal{S}_{v,w}=v\ast(w\ast\overline
{\text{ast}_{w}}\left(  \text{link}_{v}\left(  \partial\mathcal{S}\right)
\right)  )$ is true by (ii). The verification of the second claim is a
consequence of the fact that $w\ast\overline{\text{ast}_{w}}\left(
\text{link}_{v}\left(  \partial\mathcal{S}\right)  \right)  $ is a
PL\emph{\ }$\left(  k-1\right)  $-ball whose boundary is the link$_{v}\left(
\partial\mathcal{S}\right)  $ (by using (vi)).\hfill{}$\square$

\begin{definition}
\label{IND-SUB}Let $\mathcal{S}$ be an abstract simplicial complex. A
subcomplex $\mathcal{S}^{\prime}$ of $\mathcal{S}$ is called \textit{induced
}if for any face $F\in\mathcal{S}$\emph{, }vert$\left(  F\right)
\subset\mathcal{\ }$vert$\left(  \mathcal{S}^{\prime}\right)  $
implies\emph{\ }$F\in\mathcal{S}^{\prime}$\emph{\ }(i.e., if the vertex set of
a face lies in the subcomplex, then the whole face lies in the subcomplex).
\end{definition}

\begin{definition}
\label{IND-TR}Let\emph{\ }$\mathbf{s}$ be an abstract $d$-simplex (considered
as simplicial complex consisting of itself together with all of its faces). An
abstract simplicial complex $\mathcal{S}$ is called an \textit{induced
simplicial subdivision}\emph{\ }or an\emph{\ }\textit{induced triangulation
}of\textit{\emph{\ }}$\mathbf{s}$ \textit{\ }if there is a PL
homeomorphism\emph{\ }$\varphi:\mathbf{s}\rightarrow\mathcal{S}$ such that for
every face $F$ of $\mathbf{s}$\emph{, }the image\emph{\ }$\varphi\left(
F\right)  $ is an induced subcomplex of $\mathcal{S}$.
\end{definition}

\begin{remark}
If $\mathbf{s}$ is an abstract $d$-simplex\ (viewed as simplicial complex),
then every induced simplicial subdivision\emph{\ }$\mathcal{S}$\emph{\ }%
of\emph{\ }$\mathbf{s}$\emph{\ }is quasigeometric\emph{\ }(see \ref{QUASI-SUB}).
\end{remark}

\begin{proposition}
[Gluing PL Balls]\label{VERKLEBEN}Let $\mathcal{S}$, $\mathcal{S}^{\prime}$ be
two simplicial \emph{PL }$k$-balls and $\mathcal{S}^{\prime\prime}:=$
$\mathcal{S}\cap\mathcal{S}^{\prime}$. Suppose that $\mathcal{S}^{\prime
\prime}\subset$ $\partial\mathcal{S}\cap\partial\mathcal{S}^{\prime
}.\smallskip$\newline\emph{(i) }If $\mathcal{S}^{\prime\prime}$ is a
\emph{(}simplicial\emph{)} \emph{PL }$\left(  k-1\right)  $-ball, then the
regular cell complex $\mathcal{S}%
{\textstyle\bigcup\nolimits_{\mathcal{S}^{\prime\prime}}}
\mathcal{S}^{\prime}$ which is obtained by gluing $\mathcal{S}$
with
$\mathcal{S}^{\prime}$ along $\mathcal{S}^{\prime\prime}$ is a \emph{PL }%
$k$-ball.\smallskip\newline\emph{(ii) }If $\mathcal{S}%
{\textstyle\bigcup\nolimits_{\mathcal{S}^{\prime\prime}}}
\mathcal{S}^{\prime}$ is a subcomplex of both $\mathcal{S}$, $\mathcal{S}%
^{\prime}$, and for at least one of $\mathcal{S}$, $\mathcal{S}^{\prime}$, the
subcomplex $\mathcal{S}^{\prime\prime}$ is \emph{induced}, then the glued
regular cell complex $\mathcal{S}%
{\textstyle\bigcup\nolimits_{\mathcal{S}^{\prime\prime}}}
\mathcal{S}^{\prime}$ is a \emph{simplicial }complex.
\end{proposition}

\noindent\textsc{Proof}\textit{. }(i) It follows directly from \cite[Corollary
1.28, p. 39]{Hud}.\smallskip\newline(ii) Suppose that $\mathcal{S}%
^{\prime\prime}$ is an induced subcomplex of $\mathcal{S}$. If $\mathcal{S}%
{\textstyle\bigcup\nolimits_{\mathcal{S}^{\prime\prime}}}
\mathcal{S}^{\prime}$ were a non-simplicial complex, then,
w.l.o.g. we may assume that there were a face $F$ of $\mathcal{S}$
and a face $F^{\prime}$ of $\mathcal{S}^{\prime}$, such that
$F\cap F^{\prime}$ is not a \textit{single simplex} (considered
itself as simplicial complex together with all its
faces). But the vertex set of $F\cap F^{\prime}$ $\subset\mathcal{S}%
^{\prime\prime}$ is contained in $F$. Thus, there is one single face
$\widetilde{F}$ of $F$ with vert$(\widetilde{F})=$ vert$\left(  F\cap
F^{\prime}\right)  .$ On the other hand, $\widetilde{F}\in\mathcal{S}%
^{\prime\prime}$, because $\mathcal{S}^{\prime\prime}$ contains all the
vertices of $\widetilde{F}$. But this would mean that $F\cap F^{\prime
}=\widetilde{F}$, which would lead to contradiction.\hfill{}$\square$

\begin{conjecture}
[UBC for the facets of geometric simplex triangulations]\label{Z-conj}Let
$\mathbf{s}$ be an $d$-dimensional simplex \emph{(}considered as an abstract
simplicial complex\emph{)}, and let $\mathcal{V}\subset\mathbf{s}$ be a finite
set of $\sharp\left(  \mathcal{V}\right)  =\mathfrak{b}$ points in
$\mathbf{s}$, so that $b_{k}$ of them are contained in the relative interiors
of the $\left(  k-1\right)  $-dimensional faces of $\mathbf{s}$, with%
\[
b_{1}=d+1,\ \ \mathfrak{b}^{\prime}:=\mathfrak{b}-b_{d+1}\ \ \text{and}%
\ \ \mathfrak{b}=\sum_{k=1}^{d+1}b_{k}>d+1.
\]
Then a geometric triangulation $\mathcal{S}$ of $\mathbf{s}$ with vertex set
$\mathcal{V}$ \ has not more than
\[%
\begin{array}
[c]{l}%
\mathfrak{f}_{d}\left(  \emph{CycP}_{d+1}\left(  \mathfrak{b}\right)  \right)
-%
{\displaystyle\sum\limits_{k=2}^{d}}
\ \left(  d-\left(  k-1\right)  \right)  \ b_{k}-1\\
\ \ \\
\overset{\text{\emph{(\ref{FDCYCLIC})}}}{=}\ \left(
\begin{array}
[c]{c}%
\mathfrak{b}-\left\lceil \tfrac{d+1}{2}\right\rceil \\
\ \ \\
\left\lfloor \tfrac{d+1}{2}\right\rfloor
\end{array}
\right)  +\left(
\begin{array}
[c]{c}%
\mathfrak{b}-1-\left\lceil \tfrac{d}{2}\right\rceil \\
\ \ \\
\left\lfloor \tfrac{d}{2}\right\rfloor
\end{array}
\right)  -%
{\textstyle\sum\limits_{k=2}^{d-1}} \ \left(  d-k\right)  \
b_{k}-\left(  \mathfrak{b}^{\prime}-\left( d+1\right)  \right)  -1
\end{array}
\]
facets \emph{(}$=d$-faces\emph{).}
\end{conjecture}

\begin{remark}
(i) Conjecture \ref{Z-conj} is true in dimension $d=1$\emph{ }because we have
$\mathfrak{f}_{1}\left(  \text{CycP}_{2}\left(  \mathfrak{b}\right)  \right)
=\mathfrak{b}$, and it follows from Euler's polyhedron formula for\emph{\ }%
$d=2$\emph{\ }(where\emph{\ }$\mathfrak{f}_{2}\left(  \text{CycP}_{3}\left(
\mathfrak{b}\right)  \right)  =2\mathfrak{b}-4$). Thus, the first
\textquotedblleft interesting\textquotedblright\ case coming into question is
that for $d=3$ (where $\mathfrak{f}_{3}\left(  \text{CycP}_{4}\left(
\mathfrak{b}\right)  \right)  =\frac{1}{2}\mathfrak{b}\left(  \mathfrak{b}%
-3\right)  $), corresponding to triangulations of the tetrahedron
using\emph{\ }$b_{2}$\emph{\ }additional vertices on its edges\emph{,} $b_{3}$
extra\emph{\ }vertices in its $2$-faces, and $b_{4}$ more vertices in its
relative interior.\smallskip\newline(ii) \ref{Z-conj} is also true (and tight)
in the case in which all \textquotedblleft additional\textquotedblright%
\ vertices lie in the relative interior of\emph{\ }$\mathbf{s}$, that is, if
$b_{2}=b_{3}=\cdots=b_{d}=0$\emph{, }because the upper bound (\ref{UP-LO2})
can be written as
\[
\mathfrak{f}_{d}\left(  \text{CycP}_{d+1}\left(  \mathfrak{b}\right)  \right)
-\left(  b_{2}+\cdots+b_{d}\right)  -1=\mathfrak{f}_{d}\left(  \text{CycP}%
_{d+1}\left(  \mathfrak{b}\right)  \right)  -\left(  \mathfrak{b}^{\prime
}-\left(  d+1\right)  \right)  -1\ .
\]
But whenever there\emph{\ }\textit{are }additional vertices on the boundary,
the above \textquotedblleft new\textquotedblright\ upper bound would obviously
improve (\ref{UP-LO2}) by\emph{\ }subtracting\emph{\ }the \textquotedblleft
extra\textquotedblright\ summand\emph{\ }$%
{\textstyle\sum\limits_{k=2}^{d-1}} \ \left(  d-k\right)  \
b_{k}.$ This would correspond to a better estimation
of a part of a \textquotedblleft missing correction term\textquotedblright%
\ involving\emph{\ }$\mathbf{h}$-vector\emph{\ }components of\emph{\ }%
$\partial\mathcal{S}.\smallskip$\newline(iii) We believe that the
\textquotedblleft right\textquotedblright\ setting for a proof of Conjecture
\ref{Z-conj} is provided by Stanley's theory \cite{Stanley4} of\emph{\ }%
\textquotedblleft local $\mathbf{h}$-vectors\textquotedblright: If
$\mathcal{S}$ is a topological simplicial subdivision of an abstract
$d$-simplex $\mathbf{s},$ and\emph{ $\mathcal{V}$ \ }the vertex set of
$\mathbf{s},$ then the \textit{local} $\mathbf{h}$-\textit{vector}
$\ell_{\mathcal{V}}\left(  \mathcal{S}\right)  :=\left(  \ell_{0}\left(
\mathcal{S}\right)  ,\ldots,\ell_{d+1}\left(  \mathcal{S}\right)  \right)  $
of $\mathcal{S}$\emph{ }is defined by expanding%
\[
\ell_{\mathcal{V}}\left(  \mathcal{S};t\right)  :=%
{\displaystyle\sum\limits_{\mathcal{W}\subseteq\mathcal{V}}} \
\left(  -1\right)  ^{\#\left(
\mathcal{V}\smallsetminus\mathcal{W}\right) }\ \mathbf{h}\left(
\mathcal{S}_{\mathcal{W}};t\right)  \in\mathbb{Z}[t]
\]
w.r.t. $t,$%
\[
\ell_{\mathcal{V}}\left(  \mathcal{S};t\right)  =\ell_{0}\left(
\mathcal{S}\right)  +\ell_{1}\left(  \mathcal{S}\right)  t+\cdots+\ell
_{d}\left(  \mathcal{S}\right)  t^{d}+\ell_{d+1}\left(  \mathcal{S}\right)
t^{d+1},
\]
and has the following properties: \textbf{a)} \textit{Reciprocity}:
$\!\ell_{i}\left(  \mathcal{S}\right)  =\ell_{d-i}\left(  \mathcal{S}\right)
,\forall i\in\{0,..,d+1\},$ \newline\textbf{b)} \textit{Positivity}: $\ell
_{i}\left(  \mathcal{S}\right)  \geq0,\ \forall i\in\{0,\ldots,d+1\},$
whenever is a quasigeometric subdivision of $\mathbf{s}$ (in the sense of
\ref{QUASI-SUB}), and \textbf{c)} \textit{Locality}:
\[
\mathbf{h}\left(  \mathcal{S}^{\prime};t\right)  =%
{\displaystyle\sum\limits_{F\in\mathcal{S}}}
\ \mathbf{\ell}_{F}\left(  \mathcal{S}_{F}^{\prime};t\right)  \ \mathbf{h}%
\left(  \text{link}_{F}\left(  \mathcal{S}\right)  ;t\right)  ,
\]
with $\mathcal{S}^{\prime}$ a topological subdivision of a pure abstract
simplicial $d$-complex $\mathcal{S},$ and%
\[
\text{link}_{F}\left(  \mathcal{S}\right)  :=\left\{  F^{\prime}%
\in\mathcal{S\ }\left\vert \ F\cup F^{\prime}\in\mathcal{S}\text{\emph{,}
}F\cap F^{\prime}=\varnothing\right.  \right\}  .
\]
However, even for $d=3$\emph{, }this is not \textquotedblleft for
free\textquotedblright\ in our case: we would need to establish
the following upper bound for the difference\emph{\ }between the
second and the first coordinate of the local
$\mathbf{h}$-vector\emph{\ }$\ell_{\mathcal{V}}\left(
\mathcal{S}\right)  $ of\emph{\ }$\mathcal{S}$\text{:}
$\ell_{2}\left( \mathcal{S}\right)  -\ell_{1}\left(
\mathcal{S}\right)  \leq\binom {\mathfrak{h}_{1}\left(
\mathcal{S}\right)  }{2},$\emph{ }whereas from the results in
\cite{Stanley4,Chan1,Chan2} we only get
\[
\ell_{2}\left(  \mathcal{S}\right)  \leq\mathfrak{h}_{2}\left(  \mathcal{S}%
\right)  \leq\binom{\mathfrak{h}_{1}\left(  \mathcal{S}\right)  +1}%
{2}\,\,\,\,\,\,\left(  \text{i.e., \thinspace\thinspace\thinspace}\ell
_{2}\left(  \mathcal{S}\right)  -\mathfrak{h}_{1}\left(  \mathcal{S}\right)
\leq\binom{\mathfrak{h}_{1}\left(  \mathcal{S}\right)  }{2}\right)
\]
which is weaker than that we would like to have whenever there are boundary
vertices. In any case, to proceed along these lines depends certainly on a
deep understanding of how \textbf{a)}-\textbf{c)} could be applicable to our
specific situation. Here we restrict ourselves to present another proof for
dimension $d=3$ by passing to triangulations in the category of PL
subdivisions of $\mathbf{s}$ within which the gluing of balls is able to work
without essential obstructions and leads to the desired result.
\end{remark}

\begin{theorem}
\label{BIGTH}For $d=3$, Conjecture \emph{\ref{Z-conj}}
holds\emph{\ }in greater generality\emph{: }any \emph{induced}
triangulation $\mathcal{S}$ \ \emph{(cf. \ref{IND-TR}) }of a
tetrahedron $\mathbf{s}$ \emph{(}considered as abstract simplicial
complex\emph{) }using $b_{2}+b_{3}+b_{4}$ additional vertices
within the edges$/2$-faces$/$interior of $\mathbf{s}$ possesses at
most $\mathfrak{f}_{3}\left(  \emph{CycP}_{4}\left(
\mathfrak{b}\right)  \right) -2b_{2}-b_{3}-1$ facets.
\end{theorem}

\noindent\textsc{Proof}\textit{. }If $b_{2}=b_{3}=0$, then we have nothing to
show. The proof will use induction on $b_{2}+b_{3}$. For fixed $b_{2},b_{3}$,
with $b_{2}+b_{3}>0$, assume that all induced triangulations of tetrahedra
whose number of vertices lying in the relative interior of their edges and of
their $2$-faces is $<b_{2}+b_{3}$ enjoy the desired property. We shall
distinguish two cases.\medskip\newline$\rhd$ \textsc{First case}\textit{.
}Suppose that $b_{3}>0$, i.e., that there are vertices of $\mathcal{S}$ in the
relative interior of at least one $2$-face, say $F$, of the tetrahedron
$\mathbf{s}$. For every vertex $v$ of $F$ (in the simplicial PL $2$-sphere
$\partial\mathcal{S}$) the closed star $\overline{\text{star}_{v}}\left(
\partial\mathcal{S}\right)  $ of $v$ in $\partial\mathcal{S}$ is a simplicial
PL $2$-ball (that is, a simplicial PL disc\footnote{We may freely identify all
simplicial discs which will occur in the arguments of our proof with the
planar graphs consisting only of their vertices and edges (i.e. with their
$1$-skeletons).}), cf. Proposition \ref{PL-prop} (i) and (iv). $\overline
{\text{star}_{v}}\left(  \partial\mathcal{S}\right)  $ is not necessarily an
induced subcomplex of $\partial\mathcal{S}$ (cf. \ref{IND-SUB}). This first
difficulty will be removed as follows.\ \medskip\newline$\blacktriangleright$
\textsc{Claim.}\textit{\ Considering as }\textquotedblleft%
\textit{starting-point}\textquotedblright\textit{\ an arbitrary vertex }%
$v_{0}$ \textit{of the relative interior of }$F$\textit{, we can determine
another vertex }$v_{\mu}$ (\textit{also lying in the relative}
\textit{interior of} $F$)\textit{, such that} $\overline{\text{star}_{v_{\mu}%
}}\left(  \partial\mathcal{S}\right)  $ \textit{forms an} induced
\textit{subcomplex of} $\partial\mathcal{S}$.\ \medskip

\noindent$\blacktriangleright$ \textsc{Proof of the claim.} At first define
\[
\mathfrak{U}\left(  v_{0}\right)  :=\left\{
\begin{array}
[c]{c}%
\text{all simplicial PL discs }\mathbf{D}\text{ with vertex sets belonging }\\
\text{exclusively to }F\text{, such that }v_{0}\text{ lies in their}\\
\text{relative interior, and all }w\in\text{vert}\left(  \partial
\mathbf{D}\right)  \text{ are adjacent to }v_{0}%
\end{array}
\right\}  ,
\]
and
\[
\mathfrak{U}\left(  v_{0};\lambda\right)  :=\left\{  \mathbf{D}\in
\mathfrak{U}\left(  v_{0}\right)  \ \left\vert \ \sharp\left(  \text{vert}%
\left(  \mathbf{D}\right)  \right)  =\lambda\right.  \right\}  ,\smallskip
\]%
\[
\lambda_{0}:=\text{ max}\left\{  \lambda\ \left\vert \ 4\leq\lambda\leq
\sharp\left(  \text{vert}\left(  \partial\mathcal{S}\left\vert _{F}\right.
\right)  \right)  ,\text{ such that }\mathfrak{U}\left(  v_{0};\lambda\right)
\neq\varnothing\right.  \right\}  .\smallskip
\]
Fix a simplicial PL disc $\mathbf{D}_{0}\in\mathfrak{U}\left(  v_{0}%
;\lambda_{0}\right)  $. If $\mathbf{D}_{0}$ is not an induced subcomplex of
$\partial\mathcal{S}$, choose a vertex $v_{1}\neq v_{0}$ of the relative
interior of $\mathbf{D}_{0}$. After that, define analogously\smallskip%
\[
\mathfrak{U}\left(  v_{1}\right)  :=\left\{
\begin{array}
[c]{c}%
\text{all simplicial PL discs }\mathbf{D}\text{ with vertex sets belonging }\\
\text{exclusively to }\mathbf{D}_{0}\text{, such that }v_{1}\text{ lies in
their }\\
\text{relative interior, and all }w\in\text{vert}\left(  \partial
\mathbf{D}\right)  \text{ are adjacent to }v_{1}%
\end{array}
\right\}  ,
\]
and
\[
\mathfrak{U}\left(  v_{1};\lambda\right)  :=\left\{  \mathbf{D}\in
\mathfrak{U}\left(  v_{1}\right)  \ \left\vert \ \sharp\left(  \text{vert}%
\left(  \mathbf{D}\right)  \right)  =\lambda\right.  \right\}  ,\smallskip
\]%
\[
\lambda_{1}:=\text{ max}\left\{  \lambda\ \left\vert \ 4\leq\lambda\leq
\sharp\left(  \text{vert}\left(  \mathbf{D}_{0}\right)  \right)  ,\text{ such
that }\mathfrak{U}\left(  v_{1};\lambda\right)  \neq\varnothing\right.
\right\}  .\smallskip
\]
Fix again a simplicial PL disc $\mathbf{D}_{1}\in\mathfrak{U}\left(
v_{1};\lambda_{1}\right)  $. If $\mathbf{D}_{1}$ is not an induced subcomplex
of $\partial\mathcal{S}$, choose a vertex $v_{2}\neq v_{1}$ of the relative
interior of $\mathbf{D}_{1}$, and repeat this construction for $v_{2}$ ...
etc. We shall call the occuring numbers $\lambda_{0},\lambda_{1},\lambda
_{2},\ldots$ \textit{ring sizes }of $v_{0},v_{1},v_{2},\ldots$ \textit{with
respect to} $v_{0}$.\smallskip\newline We have:\medskip\newline(i)
$\overline{\text{star}_{v_{i}}}\left(  \partial\mathcal{S}\right)
\subseteq\mathbf{D}_{i-1}$, for all $i$, $i=1,2,\ldots\smallskip$ (This is
immediate by definition).\newline(ii) $\lambda_{0}>\lambda_{1}>\lambda
_{2}>\cdots>\lambda_{i-1}>\lambda_{i}>\cdots\smallskip$ \ (Since $v_{i}$ is
contained in the triangle formed by $v_{i-1}$ together with two neighbours
whose connecting edge does not belong to $\overline{\text{star}_{v_{i-1}}%
}\left(  \partial\mathcal{S}\right)  ,$ we have $\lambda_{i-1}>\lambda_{i}%
$.)\newline(iii) This procedure leads (after $\mu$ steps) to a vertex $v_{\mu
}$, such that all members of $\mathfrak{U}\left(  v_{\mu}\right)  $ are PL
homeomorphic to the closed star $\overline{\text{star}_{v_{\mu}}}\left(
\partial\mathcal{S}\right)  $ (i.e., they have no vertices in their relative
interiors besides $v_{\mu}$ itself), and are induced subcomplexes of
$\partial\mathcal{S}$ (and hence of $\mathcal{S}$); moreover, $\lambda_{\mu
}-1$ equals the (graph-theoretic) \textit{degree}\footnote{The \textit{degree}
of a vertex in a graph (without loops) is defined to be the number of the
edges containing this vertex.}\textit{\ }deg$\left(  v_{\mu}\right)  $ of
$v_{\mu}$ within $\overline{\text{star}_{v_{\mu}}}\left(  \partial
\mathcal{S}\right)  $. (This is clear because all $\lambda_{i}$'s are integers
$\geq4$).

\begin{figure}[h]
\epsfig{file=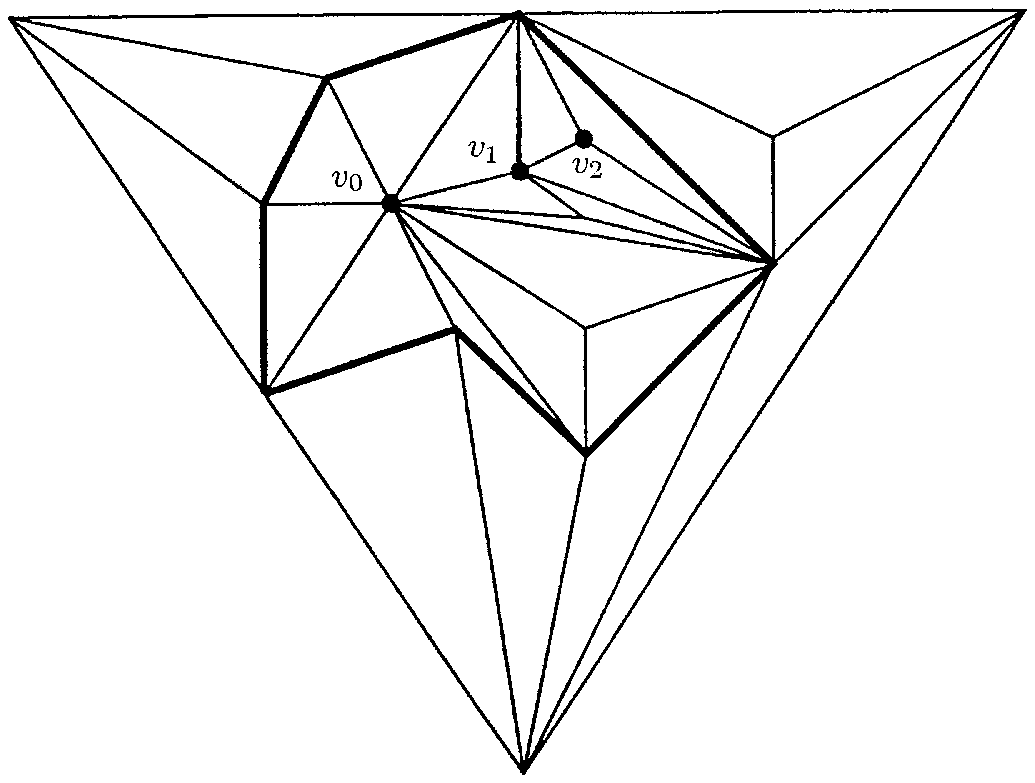, height=5cm, width=7cm}
\caption{}\label{Fig.8}
\end{figure}

Figure \ref{Fig.8} illustrates the above construction for an
example in which $\mu=2$ and $v_{0},v_{1},v_{2}$ have ring sizes
(w.r.t. $v_{0})$ $\lambda_{0}=12$, $\lambda_{1}=6$ and
$\lambda_{2}=4$, respectively.

\begin{figure}[h]
\epsfig{file=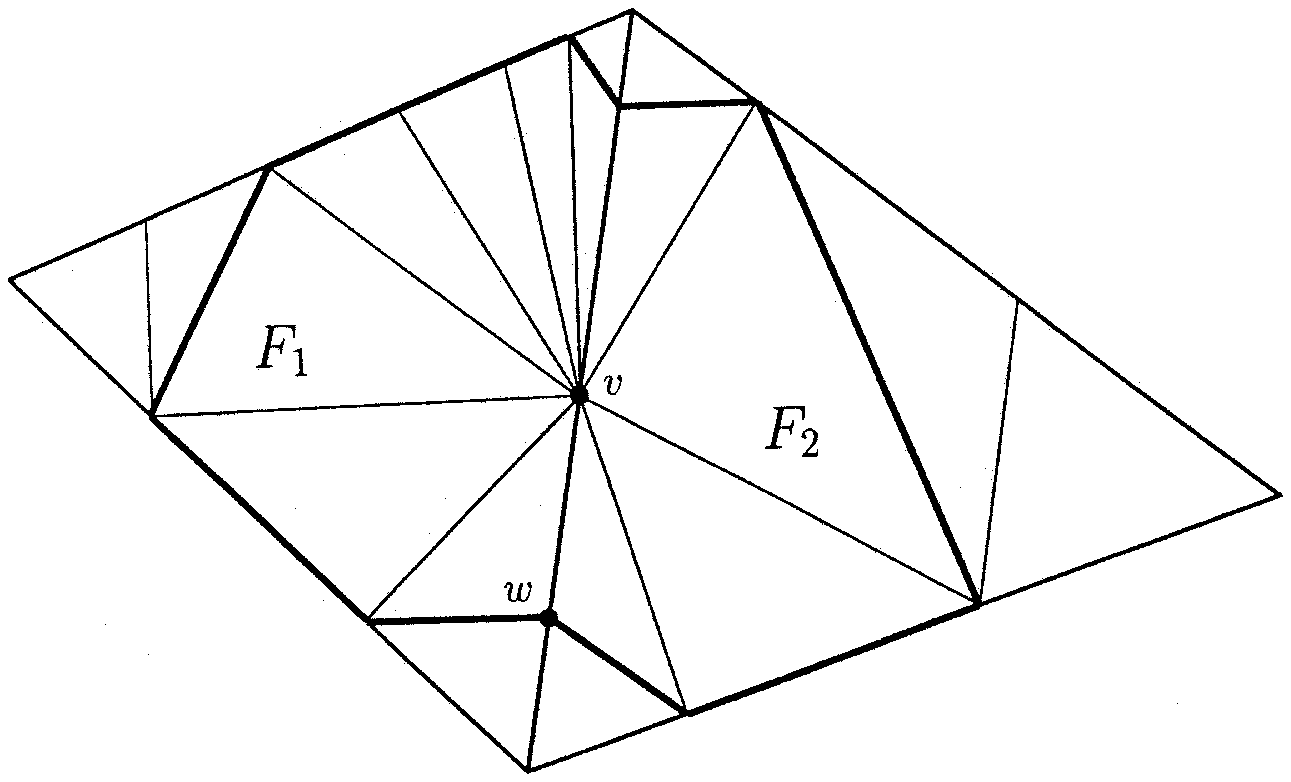, height=5cm, width=9cm}
\caption{}\label{Fig.9}
\end{figure}

\noindent$\rhd$ \textsc{Proof of theorem for the 1st case (continued)}. From
now on, let $v$ denote an (always existing, as verified above) vertex of the
relative interior of $F$, such that $\overline{\text{star}_{v}}\left(
\partial\mathcal{S}\right)  $ is an induced subcomplex of the simplicial PL
$3$-ball $\mathcal{S}$. Take a $w\in$ vert$\left(  F\right)  $ which is
adjacent to $v$ (and hence lying on the boundary of $\overline{\text{star}%
_{v}}\left(  \partial\mathcal{S}\right)  $). Then link$_{v}\left(
\partial\mathcal{S}\right)  $ is a (graph-theoretic) \textit{circuit}%
\footnote{In our special case this is synonymous to a simple closed path with
the number ($=$ \textit{size}) of its edges being equal to the number of its
vertices.}\textit{\ }with\textit{\ size} $=$ deg$\left(  v\right)  $, its
closed antistar $\overline{\text{ast}_{w}}\left(  \text{link}_{v}\left(
\partial\mathcal{S}\right)  \right)  $ (obtained by deleting $w$) is a path of
deg$\left(  v\right)  -2$ edges, and its join $\mathcal{S}_{v,w}:=\left\{
v,w\right\}  \ast\overline{\text{ast}_{w}}\left(  \text{link}_{v}\left(
\partial\mathcal{S}\right)  \right)  $ with the edge $\left\{  v,w\right\}  $
is a simplicial PL $3-$ball consisting of deg$\left(  v\right)  -2>0$
tetrahedra (by Corollary \ref{AST-LINK}). Gluing $\mathcal{S}$ and
$\mathcal{S}_{v,w}$ along their intersection $\mathcal{S}\cap\mathcal{S}%
_{v,w}=$ $\overline{\text{star}_{v}}\left(  \partial\mathcal{S}\right)  $ we
obtain a \textit{simplicial }PL $3$-ball
\[
\mathcal{S}^{\prime}:=\mathcal{S}%
{\displaystyle\bigcup\limits_{\overline{\text{star}_{v}}\left(
\partial \mathcal{S}\right)  }}
\mathcal{S}_{v,w}%
\]
according to Proposition \ref{VERKLEBEN} (i), (ii). On the other hand,
\[
\mathfrak{f}_{3}\left(  \mathcal{S}^{\prime}\right)  =\mathfrak{f}_{3}\left(
\mathcal{S}\right)  +2>\mathfrak{f}_{3}\left(  \mathcal{S}\right)  ,
\]
and the number of the vertices of $\mathcal{S}^{\prime}$ lying in
the relative interior of the $2$-faces of $\mathbf{s}$ equals
$b_{3}-1$. So we are done by induction.\medskip\newline$\rhd$
\textsc{Second case}\textit{. }Assume that $b_{3}=0$ but
$b_{2}>0$, and let $v$ denote a subdivision vertex on an edge
$F_{1}\cap F_{2}$, for $F_{1},F_{2}$ two $2-$faces of the original
tetrahedron $\mathbf{s}$ (which are adjacent to $v$). Consider a
vertex $w$ adjacent to $v$ on the same edge $F_{1}\cap F_{2}$ ($w$
may be a vertex of the original tetrahedron). In this case
$\overline{\text{star}_{v}}\left(  \partial \mathcal{S}\left\vert
_{F_{i}}\right.  \right)  $ is obviously an induced subcomplex of
$\mathcal{S}$ for $i=1,2$ (see Figure \ref{Fig.9}).

As above one proves that%
\[
\mathcal{S}_{v,w}^{\left[  i\right]  }:=\left\{  v,w\right\}  \ast
\overline{\text{ast}_{w}}\left(  \text{link}_{v}\left(  \partial
\mathcal{S}\left\vert _{F_{i}}\right.  \right)  \right)
\]
is a simplicial PL $3$-ball for $i=1,2$; gluing $\mathcal{S}$ and
$\mathcal{S}_{v,w}^{\left[  i\right]  }$ along $\mathcal{S}\cap\mathcal{S}%
_{v,w}^{\left[  i\right]  }=$ $\overline{\text{star}_{v}}\left(
\partial\mathcal{S}\left\vert _{F_{i}}\right.  \right)  $ we obtain a new
\textit{simplicial} PL $3$-ball
\[
\mathcal{S}^{\prime}:=\left(  \mathcal{S}%
{\displaystyle\bigcup\limits_{\overline{\text{star}_{v}}\left(
\partial \mathcal{S}\left\vert _{F_{1}}\right.  \right)  }}
\mathcal{S}_{v,w}^{\left[  1\right]  }\right)  \
{\displaystyle\bigcup}
\ \left(  \mathcal{S}%
{\displaystyle\bigcup\limits_{\overline{\text{star}_{v}}\left(
\partial \mathcal{S}\left\vert _{F_{2}}\right.  \right)  }}
\mathcal{S}_{v,w}^{\left[  2\right]  }\right)
\]
with $\mathfrak{f}_{3}\left(  \mathcal{S}^{\prime}\right)  \geq\mathfrak{f}%
_{3}\left(  \mathcal{S}\right)  +2>\mathfrak{f}_{3}\left(  \mathcal{S}\right)
$, and the number of the vertices of $\mathcal{S}^{\prime}$ lying in the
relative interior of the edges of $\mathbf{s}$ equals $b_{2}-1$. Thus, the
proof is finished by induction.\hfill{}$\square$

\section{Coherent Triangulations and Secondary Polytopes\label{APPCOH}}

\noindent{}\textit{Coherent triangulations} of simplices are exactly what one
needs to express the \textquotedblleft projectivity
condition\textquotedblright\ for crepant birational morphisms desingularizing
(partially or fully) Gorenstein AQS in the language of geometric combinatorics.

\begin{definition}
\label{CohTriang}A triangulation \emph{$\mathcal{T}$} of a polytope $P$ is
called \textit{coherent} (or \textit{regular}) if there exists \textit{a
strictly upper convex} $\mathcal{T}$\emph{-}\textit{support function}
$\psi:\left\vert \mathcal{T}\right\vert \rightarrow\mathbb{R}$, i.e., a
piecewise-linear real function defined on the underlying\emph{\ }%
space\emph{\ }$\left\vert \mathcal{T}\right\vert $, for which
\[
\psi(t\ \mathbf{x}+\left(  1-t\right)  \mathbf{\ y)}\geq t\ \psi\left(
\mathbf{x}\right)  +\left(  1-t\right)  \ \psi\left(  \mathbf{y}\right)
,\text{\emph{\ }for all \ }\mathbf{x},\mathbf{y}\in\left\vert \mathcal{T}%
\right\vert ,\emph{\ }t\in\left[  0,1\right]  ,
\]
so that its domains of linearity are the simplices of $\mathcal{T}$
\emph{\ }having maximal dimension.
\end{definition}

\begin{note}
\label{NOTCOHEX}(i) For a polytope $P$ one can always construct coherent
triangulations $\mathcal{T}$ with given vertex set (e.g., lexicographic,
reverse lexicographic etc.; cf. \cite[\S $2$-\S $4$]{Lee} and \cite[Ch. $8$%
]{Sturmfels}).\smallskip\ \newline(ii) Already in dimension $2$
there are lots of examples of non-coherent triangulations. Figure
\ref{Fig.10} shows two triangulations of a triangle with the same
vertex set. Triangulation (a) is coherent (in fact, affinely
isomorphic to the triangulation given in Figure \ref{Fig.3}),
whereas (b) is non-coherent. (Its \textquotedblleft mirror
image\textquotedblright\ is a different triangulation with the
same vertex set; though, the number of the simplices belonging to
the star of each vertex remains invariant. This is enough to prove
non-coherence; compare \cite[Example 2.4, p. 161]{BFS}.)

\begin{figure}[h]
\epsfig{file=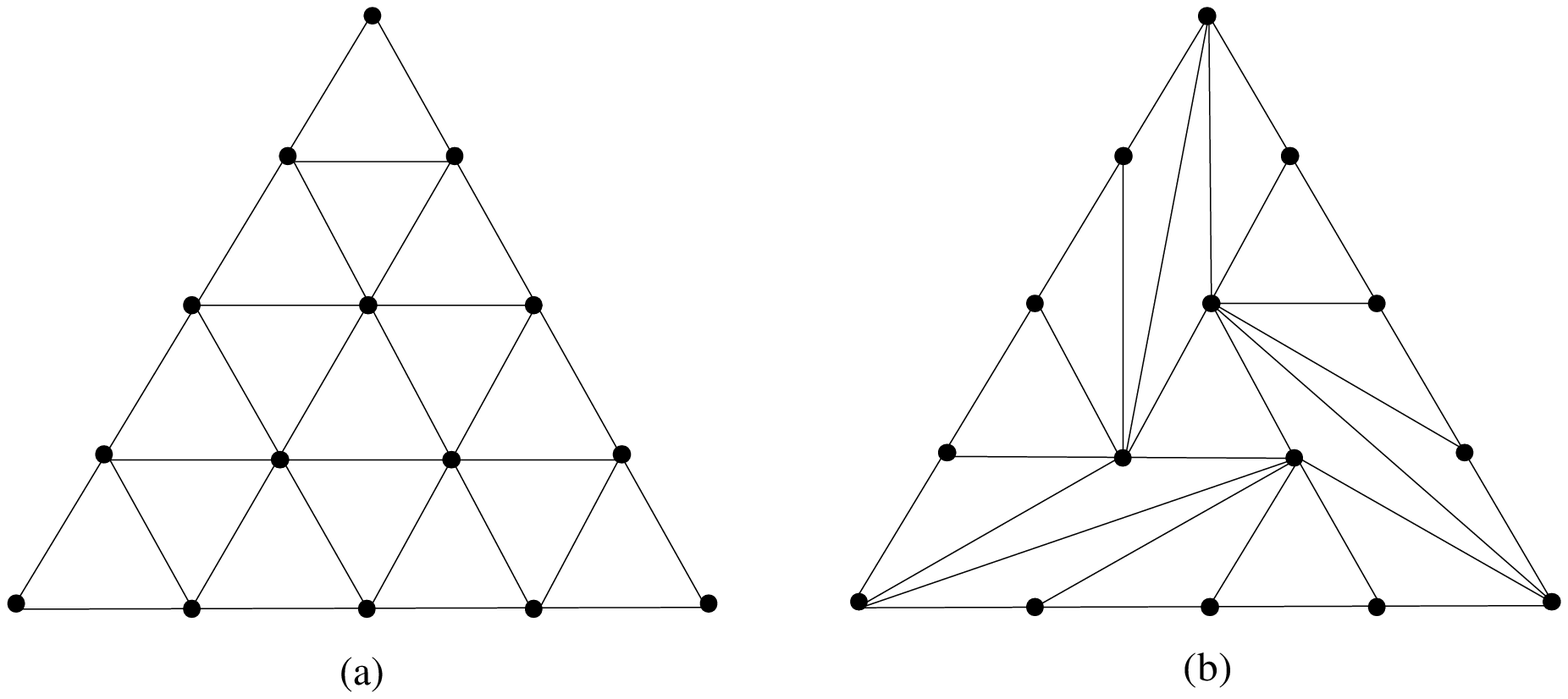, height=5cm, width=10.5cm}
\caption{}\label{Fig.10}
\end{figure}
\end{note}

\noindent{}Next Lemma is used essentially in the proof of Theorem
\ref{GPSSTHM}.

\begin{lemma}
[{Patching Lemma, \cite[2.2.2, pp. 143-145]{BGT}, \cite[5.12, p. 115]{KKMS}}%
]\label{PATCH}Let $P$ be a polytope, $\mathcal{T}=\left\{  \mathbf{s}%
_{i}\left\vert i\in I\right.  \right\}  \ $\emph{(}with $I$ a finite
set\emph{) }a coherent triangulation of $P$, and $\mathcal{T}_{i}=\left\{
\mathbf{s}_{i,j}\left\vert j\in J_{i}\right.  \right\}  \ $\emph{(}$J_{i}$
finite, for all $i\in I$\emph{) }a coherent triangulation of $\mathbf{s}_{i}$,
for all $i\in I$. If $\psi_{i}:\left\vert \mathcal{T}_{i}\right\vert
\rightarrow\mathbb{R}$ denote strictly upper convex $\mathcal{T}_{i}$-support
functions, such that
\[
\psi_{i}\left\vert _{\mathbf{s}_{i}\cap\mathbf{s}_{i^{\prime}}}\right.
=\psi_{i^{\prime}}\left\vert _{\mathbf{s}_{i}\cap\mathbf{s}_{i^{\prime}}%
}\right.
\]
for all $\left(  i,i^{\prime}\right)  \in I\times I$, then $\widetilde
{\mathcal{T}}:=\left\{  \mathbf{s}_{i,j}\left\vert \ j\in J_{i}%
\ \text{\emph{and \ }}i\in I\right.  \right\}  $ forms a \emph{coherent
triangulation }of the initial polytope $P$\emph{\ (}because the above
$\psi_{i}$'s can be canonically \emph{\textquotedblleft}patched
together\emph{\textquotedblright\ }to construct a strictly\textit{ }upper
convex $\widetilde{\mathcal{T}}$\emph{-}support function $\psi$\emph{).}
\end{lemma}

\noindent{}$\bullet$ \textbf{Secondary polytope}. For any finite set of
points\emph{ $\mathcal{V}$ \ }in $\mathbb{R}^{d}$\emph{,} \textit{all
triangulations} $\mathcal{T}$ of a polytope $P=$ conv$\left(  \mathcal{V}%
\right)  $ with vert$\left(  \mathcal{T}\right)  \subseteq\mathcal{V}$ $\ $are
parametrized by the vertices of a \textquotedblleft gigantic\textquotedblright%
\ polytope\emph{ }Un$\left(  \mathcal{V}\right)  $\emph{, }the
so-called \textit{universal polytope} of\emph{ }$P$\emph{\ }(see
\cite[\S 3]{BFS} and \cite[\S 1-\S 4]{DELHSS})\emph{. }Un$\left(
\mathcal{V}\right)  $ projects onto a polytope Sec$\left(
\mathcal{V}\right)  $ whose vertices\emph{ }parametrize only the
coherent\emph{ }$\mathcal{T}$\emph{'}s.

\begin{definition}
\label{DEFSECP}If \emph{$\mathcal{V}$ }$=\{a_{1},\ldots,a_{k}\}$ and $P$ is
$d$-dimensional$,$ the \textit{secondary polytope} Sec$\left(  \mathcal{V}%
\right)  $ of $P$ is the $(k-d-1)$-dimensional polytope
\[
\text{Sec}\left(  \mathcal{V}\right)  :=\text{conv}\left(  \left\{
\mathbf{v}_{\mathcal{T}}\ \left\vert \ \mathcal{T}\text{\ a triangulation of
}P\text{ with vert}\left(  \mathcal{T}\right)  \subseteq\mathcal{V}\right.
\right\}  \right)  \subset\mathbb{R}^{k},
\]
defined as the convex hull of the points
\[
\mathbf{v}_{\mathcal{T}}:=\sum_{i=1}^{k}\left[  \sum_{j=1}^{\nu}\left\{
\text{Vol}\left(  \mathbf{s}_{j}\right)  :a_{i}\in\mathbf{s}_{j}\right\}
\right]  \cdot e_{i},
\]
where $\left\{  \mathbf{s}_{1},\ldots,\mathbf{s}_{\nu}\right\}  $ denotes an
enumeration of the $d$-simplices of $\mathcal{T}$ and $\{e_{1},\ldots,e_{k}\}$
the standard unit vector basis of $\mathbb{R}^{k}.$ (For the main concepts of
the theory of secondary polytopes we refer to \cite{BFS}, \cite{DERAS},
\cite[Ch. 7]{GKZ}, \cite{OP}, and \cite[Lecture $9$]{Ziegler}.)
\end{definition}

\begin{note}
\label{DELOERANOTE}Sec$\left(  \mathcal{V}\right)  $ is in most of
the cases also considerably \textquotedblleft
big\textquotedblright. \ In practice, working with examples for
which $k\ $is relatively small, the vertices of\emph{ }Sec$\left(
\mathcal{V}\right)  $ can be easily determined by using De Loera's
\texttt{Puntos}\footnote{As it is pointed out in \cite[Thm.
1.5.12]{DELOERA1}, the computation of all coherent triangulations
of a $d$-polytope $P=$ conv$(\mathcal{V})$ with \texttt{Puntos}
requires $O(dRk^{d+\left\lfloor d/2\right\rfloor
}L(k-d-1,(k-d)k^{\left\lfloor d/2\right\rfloor }),\eta)$
arithmetic operations. The symbol $R$ denotes the number of
coherent triangulations of the point configuration $\mathcal{V}$
(i.e., $R=\mathfrak{f}_{0}($Sec$\left( \mathcal{V}\right)  )$),
$\eta$ is the
size of the matrix encoding $\mathcal{V},$ and $L(\beta_{1},\beta_{2}%
,\beta_{3})$ denotes the number of arithmetic operations used to
solve a linear system of inequalities of $\beta_{1}$ variables,
$\beta_{2}$ constraints and input size $\beta_{3}.$}\
\cite{DELOERA2} or Rambau's
\texttt{TOPCOM}\footnote{\texttt{TOPCOM} is much faster (as its
written in $C^{++}$) and has more functionalities than
\texttt{Puntos}.} \ \cite{RAMBAU}.
\end{note}

\noindent{}$\bullet$ \textbf{Circuits and bistellar flips}. \textquotedblleft
Flipping\textquotedblright\ along circuits in a given finite set of points
$\mathcal{V}$ $\subset\mathbb{R}^{d}$ is an elementary geometric procedure
which enable us to move from one vertex of Sec$\left(  \mathcal{V}\right)  $
to another.

\begin{definition}
[Circuits]\label{DEFCIRCUITS}A non-empty subset $\mathcal{C}\subset
\mathcal{V}$ is called a \textit{circuit} if any $\mathcal{C}^{\prime
}\subsetneqq\mathcal{C}$ is affinely independent but $\mathcal{C}$ itself is
affinely dependent\emph{. }Up to a real scalar multiple, there is a unique
real affine dependence relation among the elements of a circuit\emph{
}$\mathcal{C}$\emph{. }In fact, one can decompose\emph{ }$\mathcal{C}$ into
those\emph{ }elements\emph{ }$\mathcal{C}_{+}$ occuring with positive
coefficients in this affine linear relation, and\emph{ }$\mathcal{C}%
_{-}:=\mathcal{C}\smallsetminus\mathcal{C}_{+},$ so that
\[
\text{int}\left(  \text{conv}\left(  \mathcal{C}_{+}\right)  \right)
\cap\text{int}\left(  \text{conv}\left(  \mathcal{C}_{-}\right)  \right)
\neq\varnothing\ .
\]

\end{definition}

\begin{lemma}
[{\cite[1.2.1]{DELOERA1}}]\label{CIRCUITTR}Every circuit $\mathcal{C}$\ has
only two\emph{ }triangulations, namely
\[
\mathcal{Y}_{+}\left(  \mathcal{C}\right)  :=\left\{  \text{\emph{conv}%
}(\mathcal{C}\smallsetminus\left\{  v\right\}  )\left\vert v\in\mathcal{C}%
_{+}\right.  \right\}  \ \ \text{\emph{and\ }\ \ }\mathcal{Y}_{-}\left(
\mathcal{C}\right)  :=\left\{  \text{\emph{conv}}(\mathcal{C}\smallsetminus
\left\{  v\right\}  )\left\vert v\in\mathcal{C}_{-}\right.  \right\}  .
\]

\end{lemma}

\begin{definition}
Let $\mathcal{C}$ be a circuit and $\mathcal{T}$ a triangulation of $P=$
conv$\left(  \mathcal{V}\right)  $ with vert$\left(  \mathcal{T}\right)
\subseteq\mathcal{V}$\emph{. }Suppose that there is a sign $\odot\in\{+,-\}$
such that the following conditions are satisfied:\newline(i) The
triangulation\emph{ }$\mathcal{Y}_{\odot}\left(  \mathcal{C}\right)  $\emph{
}is a simplicial subcomplex of\emph{ }$\mathcal{T}$.\newline(ii)\emph{ }The
links of all maximal-dimensional simplices of $\mathcal{Y}_{\odot}\left(
\mathcal{C}\right)  $ within $\mathcal{T}$ coincide,\emph{ }i.e.,\emph{ }they
form the same simplicial subcomplex, say $\mathcal{T}^{\left[  \mathcal{C}%
\right]  }$\emph{, }of $\mathcal{T}$\emph{.\smallskip\ }\newline Then we say
that $\mathcal{T}$ \ is said to be \textit{supported on} $\mathcal{C}.$
\end{definition}

\begin{definition}
[Bistellar Flips]\label{DEFFLIP}Let $\mathcal{C}\subset\mathcal{V}$\emph{\ }be
a circuit and $\mathcal{T}$ a triangulation of\emph{ }$P$\emph{\ }which is
supported on $\mathcal{C}$\emph{. }The triangulation%
\[
\text{FL}_{\mathcal{C}}\left(  \mathcal{T}\right)  :=\left\{  \text{all faces
of\emph{ }}\mathcal{T}\smallsetminus(\mathcal{Y}_{\odot}\left(  \mathcal{C}%
\right)  \ast\mathcal{T}^{\left[  \mathcal{C}\right]  })\text{\emph{\ }%
}\right\}  \cup\left\{  \text{all faces of\emph{ }}\mathcal{Y}_{\boxtimes
}\left(  \mathcal{C}\right)  \ast\mathcal{T}^{\left[  \mathcal{C}\right]
}\text{\emph{\ }}\right\}
\]
of\emph{ }$P,$ where $\{\boxtimes\}=\{+-\}\mathbb{r}\{\odot\},$ is called
the\textit{ bistellar flip}\footnote{Other alternative names used in the
literature are: \textit{elementary transformation, modification, geometric
bistellar operation, surgery with respect to }$Z$ etc.\smallskip}\textit{ of}
$\mathcal{T}$ \textit{along }$\mathcal{C}$.
\end{definition}

\begin{theorem}
[{Flipping Property, \cite[Thm. 2.10, p. 233]{GKZ}}]\label{FLIPPING}Let
$\mathcal{V}$ $\subset\mathbb{R}^{d}$ be a finite set of points and $P=$
\emph{conv}$\left(  \mathcal{V}\right)  .$ For two \emph{(}different\emph{)
}coherent triangulations $\mathcal{T},\mathcal{T}^{\prime}$ with vertex set
$\mathcal{V}$, the vertices $\mathbf{v}_{\mathcal{T}},\mathbf{v}%
_{\mathcal{T}^{\prime}}$ of the secondary polytope $\emph{Sec}\left(
\mathcal{V}\right)  $ corresponding to $\mathcal{T}$ and $\mathcal{T}^{\prime
}$ are joined by an edge of $\emph{Sec}\left(  \mathcal{V}\right)  $ if and
only if there is a circuit $\mathcal{C}$ of $\mathcal{V}$ on which both
$\mathcal{T}$ and $\mathcal{T}^{\prime}$ are supported, and $\mathcal{T}%
^{\prime}=$ \emph{FL}$_{\mathcal{C}}\left(  \mathcal{T}\right)  $, i.e.,
$\mathcal{T}^{\prime}$ is the bistellar flip of $\mathcal{T}$ along
$\mathcal{C}$.
\end{theorem}

\section{Ehrhart Polynomials and Lattice Triangulations\label{APPB}}

\noindent A polytope $P\subset\mathbb{R}^{d}$ is a \textit{lattice polytope}
w.r.t. a lattice $N\subset\mathbb{R}^{d}$ if aff$(P)\cap N\neq\varnothing$ and
vert$(P)\subset N.$ We denote by $N_{P}$ the sublattice of $N$ generated (as
subgroup) by aff$(P)\cap N,$ and by\footnote{In the main text, working with
the junior simplex $\mathfrak{s}_{G}$ (or similar $(r-1)$-dimensional lattice
simplices) we often write Vol$(\mathfrak{s}_{G})$ instead of Vol$_{(N_{G}%
)_{\mathfrak{s}_{G}}}(\mathfrak{s}_{G})$ (by abuse of notation), but it is
always clear from the context what we mean in each case.} Vol$_{N_{P}}(P):=$
Vol$(P)/\det(N_{P})$ its \textit{relative volume}.\smallskip

Two lattice polytopes $P_{i}\subset\mathbb{R}^{d}$ w.r.t. $N,$ $i=1,2,$ are
called \textit{lattice equivalent} to each other if there exists an affine map
$\Phi:\mathbb{R}^{d}\rightarrow\mathbb{R}^{d}$ such that its restriction
$\left.  \Phi\right\vert _{\text{aff}(P_{1})}:$ aff$(P_{1})\rightarrow$
aff$(P_{2})$ is a bijection mapping $P_{1}$ onto the (necessarily
equidimensional) $P_{2},$ every $j$-dimensional face of $P_{1}$ onto a
$j$-dimensional face of $P_{2}$, for all $j=0,1,\ldots,\dim(P_{1})=\dim
(P_{2}),$ and $N_{P_{1}}$ onto $N_{P_{2}}.$ (If rank$(N)=\dim(P_{1}%
)=\dim(P_{2}),$ then these $\Phi$'s are exactly the \textit{affine integral
transformations}, composed of unimodular transformations and lattice
translations.)\medskip

\noindent$\bullet$ \textbf{Ehrhart polynomials}. If\emph{\ }$P\subset
\mathbb{R}^{d^{\prime}}$ is a lattice $d$-dimensional polytope\emph{\ }w.r.t.
$N\subset\mathbb{R}^{d^{\prime}},$ $d\leq d^{\prime},$ and $\nu P=\left\{  \nu
x\in\mathbb{R}^{d^{\prime}}~\left\vert ~x\in P\right.  \right\}  $\emph{\ }the
$\nu$ times dilated polytope $P$ (for $\nu\in\mathbb{N}$), then we denote the
enumerating function of its lattice points by%
\[
\mathbf{Ehr}_{N}\left(  P,\nu\right)  :=\sharp\,\left(  \nu P\cap N\right)  .
\]

\begin{theorem}
[{\cite[\S 3.3]{Beck-Robins}, \cite[28.3]{Hibi}, \cite[4.6.28]{Stanley3}}%
]\label{Ehr-Pol}$\mathbf{Ehr}_{N}\left(  P,\nu\right)  $ can be expressed as a
polynomial
\begin{equation}
\mathbf{Ehr}_{N}\left(  P,\nu\right)  =\mathbf{a}_{0}\left(  P\right)
+\mathbf{a}_{1}\left(  P\right)  \ \nu+\cdots+\mathbf{a}_{d-1}\left(
P\right)  \ \nu^{d-1}+\mathbf{a}_{d}\left(  P\right)  \ \nu^{d}\ \in
\ \mathbb{Q}\left[  \nu\right]  \label{EHR-POL}%
\end{equation}
of degree $d$ \emph{(}the so-called \emph{Ehrhart polynomial }of $P$\emph{).}
\end{theorem}

\begin{note}
(i) To find the number of lattice points of the interior of $\nu P$ one uses
the reciprocity law:%
\begin{equation}
\sharp\left(  \text{int}(\nu P)\cap N\right)  =\left(  -1\right)
^{d}\mathbf{Ehr}_{N}\left(  P,-\nu\right)  . \label{RECIPRLOW}%
\end{equation}
(ii) As it is well-known (cf. \cite[\S 3.4-3.5]{Beck-Robins}), the first, the
last but one and the last coefficient of\emph{\ } $\mathbf{Ehr}_{N}\left(
P,\nu\right)  $ are equal to $\mathbf{a}_{0}\left(  P\right)  =1$,
\begin{equation}
\mathbf{a}_{d-1}\left(  P\right)  =\frac{1}{2}\ \sum_{\text{facets }F\prec
P}\ \text{Vol}_{N_{F}}\left(  F\right)  ,\ \ \text{and~\ ~ }\mathbf{a}%
_{d-1}\left(  P\right)  =\text{Vol}_{N_{P}}\left(  P\right)  ,~
\label{vorle-le}%
\end{equation}
respectively. (For the remaining coefficients for simplices, see below Theorem
\ref{DR-thm})\emph{.\smallskip\ }\newline(iii) $\mathbf{a}_{i}(P_{1}%
)=\mathbf{a}_{i}(P_{2}),$ for all $i\in\{0,1,\ldots,\dim(P_{1})=\dim
(P_{2})\},$ whenever $P_{1}$ and $P_{2}$ are lattice equivalent$.$
\end{note}

{}\noindent{}$\bullet$ \textbf{Ehrhart series}. W.l.o.g.\footnote{If
$d<d^{\prime},$ then we work with $N_{P}$ instead of $N.$} we shall henceforth
assume that $d=d^{\prime}.$ Let
\[
\sigma_{P}:=\left\{  \left(  \lambda,\lambda p\right)  \in\mathbb{R}%
\oplus\mathbb{R}^{d}\ \left\vert \ p\in P,\ \lambda\in\mathbb{R}_{\geq
0}\right.  \right\}
\]
be the $\left(  d+1\right)  $-dimensional rational s.c.p. cone supporting $P$
within $\mathbb{R}^{d+1},$ and let $\mathfrak{R}_{P}$ denote the subring of
$\mathbb{C}[\mathfrak{x}_{0},\mathfrak{x}_{1}^{\pm1},\ldots,\mathfrak{x}%
_{d}^{\pm1}]$ spanned over $\mathbb{C}$ by all Laurent monomials of the form
$\mathfrak{x}_{0}^{\nu}\,\mathfrak{x}_{1}^{\mu_{1}}\cdots\mathfrak{x}_{d}%
^{\mu_{d}},$ where $\nu\in\mathbb{Z}_{\geq0}$ and $(\mu_{1},\ldots,\mu_{d}%
)\in\nu P\cap N.$
$\mathfrak{R}_{P}=\mathbb{C}[\sigma_{P}\cap\mathbb{Z}^{d+1}]$ is
graded by setting deg$(\mathfrak{x}_{0}^{\nu
}\,\mathfrak{x}_{1}^{\mu_{1}}\cdots\mathfrak{x}_{d}^{\mu_{d}}):=\nu.$
Hence, the so-called \textit{Ehrhart power series}
\[
\mathbf{Ehr}_{N}\left(  P;t\right)  :=1+\sum_{\nu=1}^{\infty}\ \mathbf{Ehr}%
_{N}\left(  P,\nu\right)  \ t^{\nu}\in\mathbb{Q}_{\,}\left[  \!\left[
t\right]  \!\right]
\]
\textit{\ }of $P$ is the Hilbert series of the graded normal semigroup
ring\smallskip\ $\mathfrak{R}_{P}=%
{\textstyle\bigoplus\limits_{\nu\geq0}} (\mathfrak{R}_{P})_{\nu},$
and can be therefore written in the form
\[
\mathbf{Ehr}_{N}\left(  P;t\right)  =\frac{\ \mathfrak{h}_{0}^{\ast}\left(
P\right)  +\ \mathfrak{h}_{1}^{\ast}\left(  P\right)  \ t+\cdots
+\ \mathfrak{h}_{d-1}^{\ast}\left(  P\right)  \ t^{d-1}+\ \mathfrak{h}%
_{d}^{\ast}\left(  P\right)  \ t^{d}}{\left(  1-t\right)  ^{d+1}},
\]
having the \textit{Eulerian} $\mathbf{Ehr}_{N}\left(  P,\cdot\right)
$-\textit{polynomial} as its numerator (see \cite[\S 4.3]{Stanley3}).

\begin{definition}
\label{DEFHSTARVECTOR}$\mathbf{h}^{\ast}(P):=(\mathfrak{h}_{0}^{\ast}\left(
P\right)  ,\mathfrak{h}_{1}^{\ast}\left(  P\right)  ,\ldots,\mathfrak{h}%
_{d}^{\ast}\left(  P\right)  )\in\mathbb{Z}^{d+1}$ is called the
$\mathbf{h}^{\ast}$-\textit{vector}\footnote{We adopt here Stanley's notation
from \cite{Stanley5}, but it appears in the literature also under different
names (e.g., as \textquotedblleft$\psi$-vector\textquotedblright\ in
\cite[\S 4]{BD}, as \textquotedblleft$\delta$-vector\textquotedblright\ in
\cite[Ch. IX]{Hibi} etc).} of $P.$
\end{definition}

\begin{theorem}
[{cf. \cite[\S 28]{Hibi}, \cite{Stanley-Decomp, Stanley5}}]\emph{(i)}
$\mathfrak{h}_{j}^{\ast}\left(  P\right)  \geq0,\forall j\in\{0,1,\ldots
,d\}.\smallskip$ \newline\emph{(ii)} $\mathfrak{h}_{0}^{\ast}\left(  P\right)
=1,\ \mathfrak{h}_{1}^{\ast}\left(  P\right)  =\mathbf{Ehr}_{N}\left(
P,1\right)  -\left(  d+1\right)  $\emph{,} \emph{and\ }$\mathfrak{h}_{d}%
^{\ast}\left(  P\right)  =\#\,\left(  \text{\emph{int}}\left(  P\right)  \cap
N\right)  .\smallskip$ \newline\emph{(iii)} The sum of all coordinates of the
$\mathbf{h}^{\ast}$-vector of $P$ equals%
\begin{equation}%
{\textstyle\sum\limits_{j=0}^{d}} \mathfrak{h}_{j}^{\ast}\left(
P\right)  =d!\,\text{\emph{Vol}}(P).
\label{SUMOFHASTARS}%
\end{equation}

\end{theorem}

\begin{proposition}
Each coordinate of the $\mathbf{h}^{\ast}$-vector of $P$ can be expressed as
integer linear combination of the coefficients of its Ehrhart polynomial
\emph{(\ref{EHR-POL})} as follows\emph{:}%
\begin{equation}
\mathfrak{h}_{i}^{\ast}(P)=%
{\displaystyle\sum\limits_{j=0}^{d}} \left(
{\displaystyle\sum\limits_{\kappa=0}^{i}}
(-1)^{\kappa}\tbinom{d+1}{\kappa}\left(  i-\kappa\right)
^{j}\right)
\mathbf{a}_{j}(P),\ \ \forall i\in\{0,1,\ldots,d\}. \label{H-A-FORMULA}%
\end{equation}

\end{proposition}

\noindent{}\textsc{Proof}. We expand
\[
\mathbf{Ehr}_{N}\left(  P,\nu\right)  =\text{dim}_{\mathbb{C}}(\mathfrak{R}%
_{P})_{\nu}=\dfrac{1}{\nu!}\left[  \tfrac{d^{\nu}}{dx^{\nu}}\!\left.  \left(
\frac{%
{\textstyle\sum\limits_{\iota=0}^{d}} \
\mathfrak{h}_{\iota}^{\ast}\left(  P\right)  \ x^{\iota}}{\left(
1-x\right)  ^{d+1}}\right)  \right\vert _{x=0}\right]  \smallskip
\]%
\begin{align*}
&  =\sum_{\iota=0}^{d}\ \mathfrak{h}_{\iota}^{\ast}\left(  P\right)
\ \tbinom{d+1+\left(  \nu-\iota\right)  -1}{\nu-\iota}=\sum_{\iota=0}%
^{d}\ \mathfrak{h}_{\iota}^{\ast}\left(  P\right)  \ \tbinom{\nu-\iota+d}%
{d}\smallskip\\
&  =\frac{1}{d!}\ \left[  \sum_{j=0}^{d}\mathfrak{h}_{j}^{\ast}\left(
P\right)  \ \right]  \ \nu^{d}+\cdots\,\cdots+1,
\end{align*}
as polynomial in the variable $\nu$ and compare coefficients.\hfill{}%
$\square\medskip$

\noindent{}$\bullet$ \textbf{Lattice triangulations}. These are one of our
main tools in \S \ref{CREPANTAQS} and thereafter.

\begin{definition}
[Lattice subdivisions and triangulations]\label{LATTR}A \textit{lattice
subdivision} $\mathcal{S}$ $\ $of a lattice polytope\footnote{Lattice
subdivisions, lattice triangulations and basic triangulations\textit{ }of the
undelying point set of \textit{simplicial complexes} whose vertices are
lattice points are defined similarly. (We use this generalization only for
$\mathbf{T}_{\mathbb{H}_{d}}$ in Example \ref{HVS}.)} $P$ is a polytopal
subdivision of\emph{ }$P$, such that the set vert$\left(  \mathcal{S}\right)
$\emph{\ }of\emph{ }the vertices of\emph{ }$\mathcal{S}$\emph{\ \ }belongs to
the reference lattice, and vert$\left(  P\right)  \subseteq$\emph{\ }%
vert$\left(  \mathcal{S}\right)  $. A \textit{lattice triangulation} of a
lattice polytope $P$ is a lattice subdivision of\emph{ }$P$\emph{\ }which, in
addition, is a triangulation\ (in the sense of \ref{GEOMTRIANG}).
\end{definition}

\begin{definition}
[Basic triangulations]\label{BASICTR}(i) A lattice polytope\emph{ }is called
\textit{elementary} if the lattice points belonging to it are exactly\emph{
}its\emph{ }vertices. A lattice simplex $\mathbf{s}$ is said to be
\textit{basic} (or\emph{ }\textit{unimodular})\textit{ }if its vertices
constitute a part of a $\mathbb{Z}$-basis of the reference lattice (that is,
if\emph{ }its relative volume equals $1/\dim(\mathbf{s})!$)$.\smallskip
$\newline(ii) A lattice triangulation\emph{ }$\mathcal{T}$ $\ $of a lattice
polytope $P$ is defined to be \textit{basic} if it consists only of elementary
basic simplices.
\end{definition}

\begin{note}
\label{EXLATTTR}(i) Obviously, all basic simplices are elementary. On the
other hand, all elementary triangles are basic, but in dimensions $d\geq
3$\ there exist lots of elementary simplices which are non-basic. For
instance, the so-called \textit{Reeve's simplices} \cite{Reeve}:%
\[
\text{RS}(k):=\text{ conv}(\{\mathbf{0},e_{1},e_{2},\ldots,e_{d-1}%
,(1,1,\ldots,1,k)^{\intercal}\})\subset\mathbb{R}^{d}%
\]
are elementary but non-basic (w.r.t. $\mathbb{Z}^{d}$) for $d\geq3$ and
$k\geq2$ because%
\[
d!\text{Vol}(\text{RS}(k))=\left\vert \det(e_{1},e_{2},\ldots,e_{d-1}%
,(1,1,\ldots,1,k)^{\intercal})\right\vert =k\neq1.
\]
(ii) \textquotedblleft Basicness\textquotedblright\ is a property preserved by
lattice equivalence.\smallskip\ \newline(iii) \textit{Open problem}: Hibi and
Ohsugi \cite{Hi-O} discovered a $9$-dimensional\emph{ }$0/1$-polytope\emph{
}(with\emph{ }$15$ vertices) having basic triangulations, but none of\emph{
}whose coherent triangulations is basic. It is not known if such
high-dimensional \textquotedblleft pathological
counterexamples\textquotedblright\ can be also found in the class of lattice
\textit{simplices}.
\end{note}

\begin{theorem}
[{Betke \& McMullen \cite[Thm. 2]{BetM}, Stanley~\cite[Corollary
2.5]{Stanley-Decomp}}]\label{BMcMTHM} If $\mathcal{T}$ is a lattice
triangulation of a $d$-dimensional lattice polytope $P$,
then\emph{\footnote{Note that $\mathfrak{h}_{d+1}\left(  \mathcal{T}\right)
=0$ because $P$ is homeomorphic to a $d$-ball.}}%
\begin{equation}
\mathfrak{h}_{j}^{\ast}\left(  P\right)  \geq\mathfrak{h}_{j}\left(
\mathcal{T}\right)  ,\ \forall j\in\{0,1,\ldots,d\}, \label{HAHASTAR}%
\end{equation}
and $\mathcal{T}$ is basic if and only if \emph{(\ref{HAHASTAR})} hold
\emph{(}simultaneously\emph{)} as equations.
\end{theorem}

\section{Counting the Lattice Points of the Junior Simplex\label{LATPJS}}

{}\noindent{}In this Appendix we explain how one can count, for a given
Gorenstein AQS $(\mathbb{C}^{r}/G,[\mathbf{0}])$ of type (\ref{typeAQS}), the
number of lattice points of the junior simplex $\mathfrak{s}_{G}$ (w.r.t.
$N_{G}$) by making use of Mordell-Pommersheim and Diaz-Robins
formulae.\medskip

\noindent{}$\bullet$ \textbf{Mordell-Pommersheim formula}. By (\ref{vorle-le})
the first coefficient of the Ehrhart polynomial (\ref{EHR-POL}), whose
expression (as rational linear combination of relative volumes of faces of
$P$) turns out to be relatively \textquotedblleft difficult\textquotedblright%
,\ is $\mathbf{a}_{1}\left(  \mathbf{s}\right)  ,$ arising already in the case
in which $d=3$ and $P=\mathbf{s}$ is a lattice tetrahedron. In fact, by
Theorem \ref{MORPO}, the corresponding formula (\ref{MO-PO2}) for
$\mathbf{a}_{1}\left(  \mathbf{s}\right)  $ involves \textit{Dedekind measures
}of the \textquotedblleft dihedral angles\textquotedblright\ of $\mathbf{s}.$

\begin{definition}
(i) If $x\in\mathbb{Q}$, we define
\[
\left(  \!\left(  x\right)  \!\right)  :=\left\{
\begin{array}
[c]{lll}%
x-\left\lfloor x\right\rfloor -\frac{1}{2}, & \text{if} & x\notin\mathbb{Z},\\
0, & \text{if} & x\in\mathbb{Z}.
\end{array}
\right.
\]
(ii) Let $p,q$ be two integers with $q>0$ and gcd$\left(  p,q\right)  =1$. The
\textit{Dedekind sum} DS$\left(  p,q\right)  $ of $p$ and $q$ is defined to
be
\[
\text{DS}\left(  p,q\right)  :=\sum_{j=1}^{q-1}\ \left(  \!\!\!\left(
\frac{j}{q}\right)  \!\!\!\right)  \ \left(  \!\!\!\left(  \frac{pj}%
{q}\right)  \!\!\!\right)  .
\]

\end{definition}

\begin{theorem}
[Mordell-Pommersheim formula, \cite{Mordell, Pommersheim}]\label{MORPO}%
Let\emph{\ }$N\subset\mathbb{R}^{3}$ be a lattice of rank $3,$ $n_{1}%
,n_{2},n_{3},n_{4}\in N$ four affinely independent points, and $\mathbf{s}$
the tetrahedron $\mathbf{s}=\emph{conv}\left(  \left\{  n_{1},n_{2}%
,n_{3},n_{4}\right\}  \right)  $. Denote by
\[
E_{i,j}:=\emph{conv}\left(  \left\{  n_{i},n_{j}\right\}  \right)
\ ,\ \ \ \forall\ i,j\ ,\ \text{ }1\leq i<j\leq4,
\]
its six edges, and by
\[
F_{i}:=\emph{conv}\left(  \left\{  n_{1},n_{2},n_{3},n_{4}\right\}
\smallsetminus\left\{  n_{i}\right\}  \right)  \ ,\ \ \ \forall i,\ 1\leq
i\leq4\ ,
\]
its four facets. For any fixed pair $i$, $j$, $1\leq i<j\leq4$, let
$F_{i^{\prime}}$, $F_{j^{\prime}}$ be the two facets of $\mathbf{s}$
containing its edge $E_{i,j}$ \emph{(}with $\left\{  i^{\prime},j^{\prime
}\right\}  =\left\{  1,2,3,4\right\}  \smallsetminus\left\{  i,j\right\}
,\ i^{\prime}<j^{\prime}$\emph{)},\ $\widetilde{n}_{i^{\prime}},\widetilde
{n}_{j^{\prime}}$ the images of $n_{i^{\prime}},n_{j^{\prime}}$ in
$N(E_{i,j}):=N/N_{E_{i,j}}$ under the canonical projection map $N\rightarrow
N(E_{i,j})$, and $\widetilde{\widetilde{n}}_{i^{\prime}}$, $\widetilde
{\widetilde{n}}_{j^{\prime}\text{ }}$the primitive vectors of $\emph{conv}%
\left(  \left\{  \mathbf{0},n_{i}\right\}  \right)  $ and $\emph{conv}\left(
\left\{  \mathbf{0},n_{j}\right\}  \right)  $, respectively, lying within the
lattice $N(E_{i,j})$. Set
\[
\mathbf{q}_{i,j}:=\frac{\text{\emph{det}}(\widetilde{\widetilde{n}}%
_{i^{\prime}},\widetilde{\widetilde{n}}_{j^{\prime}\text{ }})}%
{\text{\emph{det}}\left(  N(E_{i,j})\right)  }.
\]
Moreover, choosing a\thinspace$\ \mathbb{Z}$-basis, say $\{\widetilde
{\widetilde{n}}_{i^{\prime}},\widehat{n}\}$, of $N(E_{i,j})$, and expressing
$\widetilde{\widetilde{n}}_{j^{\prime}\text{ }}$ as an integer linear
combination of its elements in the form $\widetilde{\widetilde{n}}_{j^{\prime
}\text{ }}=\lambda\cdot\widetilde{\widetilde{n}}_{i^{\prime}}+\mathbf{q}%
_{i,j}\cdot\widehat{n},$ define $\mathbf{p}_{i,j}:=\left[  \lambda\right]
_{\mathbf{q}_{i,j}}.$ Then the number of lattice points of $\mathbf{s}$ is
given by the formula
\begin{equation}
\sharp\left(  \mathbf{s}\cap N\right)  =\mathbf{Ehr}_{N}\left(  \mathbf{s}%
,1\right)  =1+\mathbf{a}_{1}\left(  \mathbf{s}\right)  +\frac{1}{2}~\sum
_{i=1}^{4}~\text{\emph{Vol}}_{N_{F_{i}}}\left(  F_{i}\right)
+\text{\emph{Vol}}_{N_{\mathbf{s}}}\left(  \mathbf{s}\right)  , \label{MO-PO1}%
\end{equation}
and%
\begin{equation}
\mathbf{a}_{1}\left(  \mathbf{s}\right)  =%
{\displaystyle\sum\limits_{1\leq i<j\leq4}} \left(
\dfrac{1}{36\mathbf{q}_{i,j}}~\mathfrak{Y}_{i,i^{\prime}}+\text{
}\widetilde{\text{\emph{DS}}}(\mathbf{p}_{i,j},\mathbf{q}_{i,j})\right)
\,\text{\emph{Vol}}_{N_{E_{i,j}}}\left(  E_{i,j}\right)  , \label{MO-PO2}%
\end{equation}
where%
\[
\mathfrak{Y}_{i,i^{\prime}}:=\frac{\text{\emph{Vol}}_{N_{F_{i}}}\left(
F_{j}\,\right)  }{\text{\emph{Vol}}_{N_{F_{i^{\prime}}}}(F_{j^{\prime}}%
)}+\frac{\text{\emph{Vol}}_{N_{F_{i^{\prime}}}}(F_{j^{\prime}})}%
{\text{\emph{Vol}}_{N_{F_{i}}}\left(  F_{j}\,\right)  },
\]
and
\[
\text{ }\widetilde{\text{\emph{DS}}}(\mathbf{p}_{i,j},\mathbf{q}_{i,j}%
):=\frac{1}{4}-\text{\emph{DS}}(\mathbf{p}_{i,j},\mathbf{q}_{i,j})
\]
denotes the so-called \emph{\textquotedblleft Dedekind
measure\textquotedblright\ }of $\mathbf{p}_{i,j}$ and $\mathbf{q}_{i,j}$.
\end{theorem}

\begin{corollary}
\label{mod-MP}Let \emph{\ }$(\mathbb{C}^{4}/G,[\mathbf{0}])$ be a
$4$-dimensional Gorenstein CQS of type $\frac{1}{l}\left(  \alpha_{1}%
,\alpha_{2},\alpha_{3},\alpha_{4}\right)  $\emph{,} and $\mathfrak{s}_{G}=$
\emph{conv}$\left(  \left\{  e_{1},e_{2},e_{3},e_{4}\right\}  \right)  $ the
corresponding junior tetrahedron. For any $i$, $1\leq i\leq4$, consider an
\emph{(}arbitrary\emph{)} representation of \emph{gcd}$\left(  \alpha
_{i},l\right)  $ as integer linear combination
\[
\text{\emph{gcd}}\left(  \alpha_{i},l\right)  =\gamma_{i}\cdot\alpha
_{i}+\widehat{\gamma}_{i}\cdot l.
\]
Moreover, for any fixed pair $i$,$j$, $1\leq i<j\leq4$, let $\left\{
i^{\prime},j^{\prime}\right\}  =\left\{  1,2,3,4\right\}  \smallsetminus
\left\{  i,j\right\}  $ denote the complement \emph{(}with $i^{\prime
}<j^{\prime}$\emph{)}. Then the number of lattice points of $\mathfrak{s}_{G}$
is given by the formula
\[
\sharp\left(  \mathfrak{s}_{G}\cap N_{G}\right)  =\sharp\left\{  \lambda
\in\mathbb{Z\ }\left\vert \ 1\leq\lambda\leq l-1\ \ \text{\emph{with}\ \ }%
\sum_{i=1}^{4}\ \left[  \lambda\alpha_{i}\right]  _{l}=l\right.  \right\}  +4
\]%
\begin{equation}
\ \ \ =\mathbf{Ehr}_{N_{G}}\left(  \mathfrak{s}_{G},1\right)  =1+\mathbf{a}%
_{1}\left(  \mathfrak{s}_{G}\right)  +\mathbf{a}_{2}\left(  \mathfrak{s}%
_{G}\right)  +\mathbf{a}_{3}\left(  \mathfrak{s}_{G}\right)  \label{MODPO1}%
\end{equation}
with
\[
\mathbf{a}_{1}\left(  \mathfrak{s}_{G}\right)  =\sum_{1\leq i<j\leq4}~\left[
\tfrac{(\text{\emph{gcd}}(\alpha_{i^{\prime}},l))^{2}+(\text{\emph{gcd(}%
}\alpha_{j^{\prime}},l))^{2}}{36~l}~+\text{ }\widetilde{\text{\emph{DS}}%
}(\mathbf{p}_{i,j},\mathbf{q}_{i,j})\right]  \cdot\text{\emph{gcd}}%
(\alpha_{i^{\prime}},\alpha_{j^{\prime}},l)
\]%
\begin{equation}
\mathbf{=}~\frac{1}{12\,l}\left(  \sum_{i=1}^{4}\ \text{\emph{gcd}}\left(
\alpha_{i},l\right)  \right)  +\sum_{1\leq i<j\leq4}\ \left[  \widetilde
{\text{\emph{DS}}}\left(  \mathbf{p}_{i,j},\mathbf{q}_{i,j}\right)
\cdot\text{\emph{gcd}}\left(  \alpha_{i^{\prime}},\alpha_{j^{\prime}%
},l\right)  \right]  \label{MODPO2}%
\end{equation}
and
\begin{equation}
\mathbf{a}_{2}\left(  \mathfrak{s}_{G}\right)  =\frac{1}{4}\ \sum_{i=1}%
^{4}\ \text{\emph{gcd}}\left(  \alpha_{i},l\right)  \ ,\ \ \mathbf{a}%
_{3}\left(  \mathfrak{s}_{G}\right)  =\frac{l}{6}, \label{MODPO3}%
\end{equation}
where now
\[
\mathbf{q}_{i,j}=\frac{l\cdot\text{\emph{gcd}}(\alpha_{i^{\prime}}%
,\alpha_{j^{\prime}},l)}{\text{\emph{gcd}}\left(  \alpha_{i^{\prime}%
},l\right)  \cdot\text{\emph{gcd}}(\alpha_{j^{\prime}},l)},\ \ \ \mathbf{p}%
_{i,j}=\left[  \frac{\left(  -\gamma_{i^{\prime}}\right)  \cdot\alpha
_{i^{\prime}}}{\text{\emph{gcd}}\left(  \alpha_{i^{\prime}},l\right)
}\right]  _{\mathbf{q}_{i,j}}.
\]

\end{corollary}

\noindent Since our lattice $N_{G}=\sum_{i=1}^{4}\,\mathbb{Z\,}e_{i}%
+\mathbb{Z\,}\frac{1}{l}\left(  \alpha_{1},\alpha_{2},\alpha_{3},\alpha
_{4}\right)  ^{\intercal}$ is \textquotedblleft skew\textquotedblright\ and of
rank $4$, to apply Mordell-Pommersheim formula we have to modify appropriately
our \textquotedblleft lattice data\textquotedblright. As we shall see below,
it is more convenient to consider $l\,\mathfrak{s}_{G}$ and $l\,N_{G}$ instead
of $\mathfrak{s}_{G}$ and $N_{G}$ (to avoid bothersome denominators), to work
with the fundamental half-open parallelotope of the cone supporting
$l\,\mathfrak{s}_{G},$ and to evaluate (after that) the relative volumes of
the simplex faces by passing to the intersection with the $l$-times dilated
affine hyperplane $l\,\mathcal{H}_{1}$ of level $1$.

\begin{lemma}
\label{DIST-HY}Let $N\subset\mathbb{R}^{d}$ be a lattice and $\mathcal{H}$ a
$\left(  d-1\right)  $-dimensional linear hyperplane in $\mathbb{R}^{d}$, so
that $\overline{N}:=N\cap\mathcal{H}$ is a lattice of rank $d-1$. Then there
exists an element $n\in N$ such that $N=\overline{N}+\mathbb{Z\,}n$, and
\begin{equation}
\text{\emph{det}}\left(  N\right)  =\left(
\begin{array}
[c]{c}%
\text{\emph{euclidean distance}}\\
\text{\emph{between }}\mathbf{0\ }\text{\emph{and }}n+\mathcal{H}%
\end{array}
\right)  \cdot\text{\emph{det}}\left(  \overline{N}\right)  . \label{dist-H}%
\end{equation}

\end{lemma}

\noindent\textsc{Proof}\textit{.} Let $\left\{  n_{1},\ldots,n_{d-1}%
,n_{d}\right\}  $ be a $\mathbb{Z}$-basis of $N$ such that $\left\{
n_{1},\ldots,n_{d-1}\right\}  $ is a $\mathbb{Z}$-basis of the lattice
$\overline{N}$, $\ M:=$ Hom$_{\mathbb{Z}}\left(  N,\mathbb{Z}\right)  $ the
dual lattice of $N$, and let $v$ denote the normal vector of $\mathcal{H}$.
Then $\frac{v}{\left\langle v,n_{d}\right\rangle }\in M$ because
\[
\left\langle \frac{v}{\left\langle v,n_{d}\right\rangle },n_{i}\right\rangle
=0,\ \forall i,\ 1\leq i\leq d-1,\text{ \ and \ \ }\left\langle \frac
{v}{\left\langle v,n_{d}\right\rangle },n_{d}\right\rangle =1.
\]
Furthermore, it is primitive and by definition we have
\[
\left(
\begin{array}
[c]{c}%
\text{euclidean\ distance}\\
\text{between\emph{\ \ }}\mathbf{0\ \ }\text{and\ \emph{\ }}n_{d}+\mathcal{H}%
\end{array}
\right)  =\frac{\left\langle v,n_{d}\right\rangle }{\left\Vert v\right\Vert
}=\frac{1}{\left\Vert \frac{v}{\left\langle v,n_{d}\right\rangle }\right\Vert
}.
\]
Thus, setting $n=n_{d}$ and applying \cite[Corollary of p. 25]{Cassels}, we
get (\ref{dist-H}).\hfill$\square$

\begin{lemma}
[{\cite[Ch. 21, 2.E, formula (12) on p. 453]{Co-Sl}}]The \emph{(}%
standard\emph{)} volume of a \emph{regular} $d$-dimensional simplex
$\mathbf{s}$ of edge length $\sqrt{2}$ equals
\begin{equation}
\text{\emph{Vol}}\left(  \mathbf{s}\right)  =\dfrac{\sqrt{d+1}}{d!}\,\,
\label{vol-reg-s}%
\end{equation}

\end{lemma}

\noindent\textsc{Proof of Corollary }\ref{mod-MP}. Let $l\,\mathfrak{s}_{G}$
in $\left(  l\,N_{G}\right)  _{\mathbb{R}}$ be the $l$-times dilated junior
lattice simplex with $l\,N_{G}=\sum_{i=1}^{4}\,\mathbb{Z\,}\left(  l\,e_{i}\right)  +\mathbb{Z\,}%
\left(  \alpha_{1},\alpha_{2},\alpha_{3},\alpha_{4}\right)
^{\intercal},$ and
\[
Q:=\mathbf{Par}\left(  \mathbb{R}_{\geq0}\,\left(  l\,e_{1}\right)
+\mathbb{R}_{\geq0}\,\left(  l\,e_{2}\right)  +\mathbb{R}_{\geq0}\,\left(
l\,e_{3}\right)  +\mathbb{R}_{\geq0}\,\left(  l\,e_{4}\right)  \right)
\cap\left(  l\,N_{G}\right)
\]%
\begin{align*}
&  =\left\{  n\mathbf{\in}\left(  l\,N_{G}\right)  \ \left\vert \ n\mathbf{=}%
\sum_{i=1}^{4}\ \varepsilon_{i}\ \left(  l\,e_{i}\right)  ,\ \text{with\ }%
0\leq\varepsilon_{i}<1,\ \forall i,\ 1\leq i\leq4\right.  \right\}  \medskip\\
&  =\left\{  \left(  \left[  j\,\alpha_{1}\right]  _{l},\left[  j\,\alpha
_{2}\right]  _{l},\left[  j\,\alpha_{3}\right]  _{l},\left[  j\,\alpha
_{4}\right]  _{l}\right)  ^{\intercal}\ \left\vert \ j\in\left\{
0,1,\ldots,l-1\right\}  \right.  \right\}  .
\end{align*}
Defining
\[
Q_{i}:=\left\{  n=\left(  n^{\left(  1\right)  },n^{\left(  2\right)
},n^{\left(  3\right)  },n^{\left(  4\right)  }\right)  ^{\intercal}\in
Q\ \left\vert \ n^{\left(  i\right)  }=0\right.  \right\}  ,\ \forall
i,\ 1\leq i\leq4,
\]
it is easy to verify that
\begin{equation}
\left\{
\begin{array}
[c]{l}%
\sharp\left(  Q_{i}\right)  =\text{ gcd}\left(  \alpha_{i},l\right)
,\ \forall i,\ 1\leq i\leq4\ ,\text{ and}\smallskip\\
\sharp\left(  Q_{i}\cap Q_{j}\right)  =\text{ gcd}\left(  \alpha_{i}%
,\alpha_{j},l\right)  ,\ \forall i,j,\ 1\leq i<j\leq4
\end{array}
\right\}  \ \label{QUI}%
\end{equation}
To apply \ref{MORPO} it suffices to consider $l\,\mathfrak{s}_{G}=$
conv$\left(  \left\{  \left(  l\,e_{1}\right)  ,\,\left(  l\,e_{2}\right)
,\,\left(  l\,e_{3}\right)  ,\,\left(  l\,e_{4}\right)  \right\}  \right)  $
w.r.t. the lattice $\overline{l\,N_{G}}$ $:=l\,N_{G}\cap l\,\mathcal{H}_{1}$
of rank $3.$ Note that the euclidean distance between $\mathbf{0}$ and
\smallskip\smallskip$le_{1}+l\,\mathcal{H}_{1}$ is equal to $\frac{l}{2}$ .
Thus, by (\ref{dist-H}) and (\ref{VOLFORMULA}), we deduce that
\[
\text{det}\left(  \overline{l\,N_{G}}\right)  \cdot\frac{l}{2}=\text{
det}\left(  l\,N_{G}\right)  ,\ \ \tfrac{\text{det}\left(  l\,\mathbb{Z}%
^{4}\right)  }{\text{det}\left(  l\,N_{G}\right)  }=\tfrac{l^{4}}%
{\text{det}\left(  l\,N_{G}\right)  }=\sharp\left(  Q\right)  =l
\]%
\begin{equation}
\Rightarrow\text{det}\left(  l\,N_{G}\right)  =\text{det}(\left(
l\,N_{G}\right)  _{l\,\mathfrak{s}_{G}})=l^{3},\ \text{det}\left(
\overline{l\,N_{G}}\right)  =\text{det}(\left(  \overline{l\,N_{G}}\right)
_{l\,\mathfrak{s}_{G}})=2\,l^{2}.\ \label{DET-BER}%
\end{equation}
\newline$\rhd$ \textsc{First step}\textit{. }Let $F_{i}=$ conv$\left(
\left\{  \left(  l\,e_{1}\right)  ,\,\left(  l\,e_{2}\right)  ,\,\left(
l\,e_{3}\right)  ,\,\left(  l\,e_{4}\right)  \right\}  \smallsetminus\left\{
\left(  l\,e_{i}\right)  \right\}  \right)  $, $1\leq i\leq4$, be the facets
and $E_{i,j}=$conv$\left(  \left\{  \left(  l\,e_{i}\right)  ,\,\left(
l\,e_{j}\right)  \right\}  \right)  ,\ 1\leq i<j\leq4,$ the edges of the
tetrahedron $l\,\mathfrak{s}_{G}$. How does one compute their relative
volumes\thinspace with respect to $\overline{l\,N_{G}}$ ? By (\ref{vol-reg-s})
the standard volumes are the following:
\begin{equation}
\text{Vol}\left(  l\,\mathfrak{s}_{G}\right)  =\frac{2\,l^{3}}{3!}=\frac
{l^{3}}{3}\ ,\ \ \text{Vol}\left(  F_{i}\right)  =\frac{\sqrt{3}}{2}%
\,l^{2}\ ,\ \ \text{Vol}\left(  E_{i,j}\right)  =\sqrt{2}\,l.
\label{STAND-VOL}%
\end{equation}
On the other hand, applying again Lemma \ref{DIST-HY} and (\ref{VOLFORMULA})
for the facets and the edges of $l\,\mathfrak{s}_{G}$ we get similarly:

\[
\text{det}\left(  \text{aff}\left(  F_{i}\right)  \cap\overline{l\,N_{G}%
}\right)  \cdot\tfrac{l}{\sqrt{3}}=\text{ det}\left(
\text{lin}\left( F_{i}\right)  \cap l\,N_{G}\right)  ,\ \
\]%
\[
\dfrac{\text{det}(%
{\textstyle\sum\nolimits_{j\in\{1,2,3,4\}\mathbb{r}\{i\}}}
\mathbb{Z\,}\left(  l\,e_{j}\right)  )}{\text{det}\left(
\text{lin}\left( F_{i}\right)  \cap l\,N_{G}\right)
}=\dfrac{l^{3}}{\text{det}\left( \text{lin}\left(  F_{i}\right)
\cap l\,N_{G}\right)  }=\sharp\left( Q_{i}\right)
\overset{\text{(\ref{QUI})}}{=}\text{ gcd}\left(  \alpha
_{i},l\right)  \
\]%
\begin{equation}
\Longrightarrow\text{det}\left(  \text{lin}\left(  F_{i}\right)
\cap l\,N_{G}\right)  =\tfrac{l^{3}}{\text{gcd}\left(
\alpha_{i},l\right)
}\ ,\text{ det}\left(  \text{aff}\left(  F_{i}\right)  \cap\overline{l\,N_{G}%
}\right)  =\tfrac{\sqrt{3}\,l^{2}}{\text{gcd}\left(
\alpha_{i},l\right)
},\ \label{DET-BE2}%
\end{equation}
and
\[
\text{det}\left(  \text{aff}\left(  E_{i,j}\right)  \cap\overline{l\,N_{G}%
}\right)  \cdot\tfrac{l}{\sqrt{2}}=\text{det}\left(
\text{lin}\left( E_{i,j}\right)  \cap l\,N_{G}\right)  \ ,
\]
$\ $%
\begin{align*}
\dfrac{\text{det}\left(  \mathbb{Z\,}\left(
l\,e_{i^{\prime}}\right) +\mathbb{Z\,}\left(
l\,e_{j^{\prime}}\right)  \right)  }{\text{det}\left(
\text{lin}\left(  E_{i,j}\right)  \cap l\,N_{G}\right)  } &  =\dfrac{l^{2}%
}{\text{det}\left(  \text{lin}\left(  E_{i,j}\right)  \cap
l\,N_{G}\right)
}\\
&  =\sharp\left(  Q_{i^{\prime}}\cap Q_{j^{\prime}}\right)
\overset {\text{(\ref{QUI})}}{=}\text{ gcd}\left(
\alpha_{i^{\prime}},\alpha _{j^{\prime}},l\right)
\end{align*}%
\begin{equation}
\Longrightarrow\left\{
\begin{array}
[c]{l}%
\text{det}\left(  \text{lin}\left(  E_{i,j}\right)  \cap
l\,N_{G}\right)
=\tfrac{l^{2}}{\text{gcd}(\alpha_{i^{\prime}},\alpha_{j^{\prime}},l)},\text{
and}\\
\text{det}\left(  \text{aff}\left(  E_{i,j}\right)  \cap\overline{l\,N_{G}%
}\right)
=\tfrac{\sqrt{2}\,l}{\text{gcd}(\alpha_{i^{\prime}},\alpha
_{j^{\prime}},l)}%
\end{array}
\right\}  \ .\ \label{DET-BE3}%
\end{equation}
Combining (\ref{STAND-VOL}) with (\ref{DET-BER}), (\ref{DET-BE2}),
(\ref{DET-BE3}), we finally obtain
\begin{equation}
\text{Vol}_{(\overline{l\,N_{G}})_{l\,\mathfrak{s}_{G}}}\left(
l\,\mathfrak{s}_{G}\right)  =\tfrac{\text{Vol}\left(  l\,\mathfrak{s}%
_{G}\right)  }{\text{det}(\left(  \overline{l\,N_{G}}\right)
_{l\,\mathfrak{s}_{G}})}=\tfrac{l}{6}\ ,\label{REVOL1}%
\end{equation}%
\begin{equation}
\text{Vol}_{\left(  \overline{l\,N_{G}}\right)  _{F_{i}}}\left(
F_{i}\right) =\tfrac{\text{Vol}\left(  F_{i}\right)
}{\text{det}\left(  \text{aff}\left( F_{i}\right)
\cap\overline{l\,N_{G}}\right)  }=\tfrac{\text{gcd}\left(
\alpha_{i},l\right)  }{2},\ \label{REVOL2}%
\end{equation}%
\begin{equation}
\text{Vol}_{(\overline{l\,N_{G}})_{E_{i,j}}}\left(  E_{i,j}\right)
=\tfrac{\text{Vol}\left(  E_{i,j}\right)
}{\text{det}(\text{aff}\left(
E_{i,j}\right)  \cap\overline{l\,N_{G}})}=\text{gcd}(\alpha_{i^{\prime}%
},\alpha_{j^{\prime}},l).\label{REVOL3}%
\end{equation}
$\rhd$ \textsc{Second step}\textit{. }The procedure of the determination of
$\mathbf{q}_{i,j}$'s and $\mathbf{p}_{i,j}$'s (whose Dedekind measures lead to
the evaluation of the contributions of the \textquotedblleft dihedral
angles\textquotedblright\ to the counting of lattice points) is a little bit
more complicated. To simplify it, we shall this time transform
$l\,\mathfrak{s}_{G}$. For all indices $i,j,\ 1\leq i<j\leq4$, we define an
integer translation
\[
\mathbf{s}_{i,j}:=\text{ }l\,\mathfrak{s}_{G}-l\,e_{i}=\text{ conv}\left(
\left\{  \mathbf{0},\,l\left(  e_{j^{\prime}}-e_{i}\right)  ,\,l\left(
e_{i^{\prime}}-e_{i}\right)  ,\,l\left(  e_{j}-e_{i}\right)  \right\}
\right)
\]
and
\[
N_{i,j}:=\overline{l\,N_{G}}-l\,e_{i}%
\]%
\[
=\mathbb{Z\,}l\,(e_{j^{\prime}}-e_{i})+\mathbb{Z\,}l\,\left(  e_{i^{\prime}%
}-e_{i}\right)  +\mathbb{Z\,}l\mathbb{\,}\left(  e_{j}-e_{i}\right)
+\mathbb{Z\,}l\,(\left(  \alpha_{1},\alpha_{2},\alpha_{3},\alpha_{4}\right)
^{\intercal}-(%
{\textstyle\sum\limits_{\iota=1}^{4}} \,\alpha_{\iota})\,e_{i})
\]
Furthermore, we define a unimodular transformation $\Phi_{i,j}:\left(
N_{i,j}\right)  _{\mathbb{R}}\longrightarrow\left(  N_{i,j}\right)
_{\mathbb{R}}$ by
\[
\Phi_{i,j}\left(  e_{j}\right)  =e_{i}\,,\,\,\,\Phi_{i,j}\left(  e_{i}\right)
=e_{i}-e_{j}\,,\,\,\,\Phi_{i,j}\left(  e_{i^{\prime}}\right)  =e_{i}%
+e_{i^{\prime}}\,\,\,\,\text{and\thinspace\thinspace\thinspace\thinspace}%
\Phi_{i,j}\left(  e_{j}\right)  =e_{i}+e_{j^{\prime}}\,,
\]
and linear extension. $\Phi_{i,j}$ transfers $\mathbf{s}_{i,j}$ onto
\[
\Phi_{i,j}\left(  \mathbf{s}_{i,j}\right)  =\text{conv}\left(  \left\{
\begin{array}
[c]{c}%
\mathbf{0\ },\ \ \Phi_{i,j}\,\left(  l\left(  e_{j^{\prime}}-e_{i}\right)
\right)  =l(e_{j}+e_{j^{\prime}}),\smallskip\\
\Phi_{i,j}\,\left(  l\left(  e_{i^{\prime}}-e_{i}\right)  \right)  =l\left(
e_{i^{\prime}}+e_{j}\right)  ,\,\Phi_{i,j}\left(  l\left(  e_{j}-e_{i}\right)
\right)  =le_{j}%
\end{array}
\right\}  \right)  ,
\]
the lattice $N_{i,j}$ onto
\[
\widehat{N_{i,j}}:=\Phi_{i,j}\left(  N_{i,j}\right)
\]%
\[
=\mathbb{Z\,}l\,(e_{j^{\prime}}+e_{j})+\mathbb{Z\,}l\,e_{j}+\mathbb{Z\,}%
l\mathbb{\,}\left(  e_{i^{\prime}}+e_{j}\right)
+\mathbb{Z\,}\left(  \left(
{\textstyle\sum\limits_{\iota\in\{1,2,3,4\}\mathbb{r}\{i\}}}
\,\alpha_{\iota}\right)  e_{j}+\alpha_{i^{\prime}}\,e_{i^{\prime}}%
+\alpha_{j^{\prime}}\,e_{j^{\prime}}\right)
\]%
\[
=\mathbb{Z\,}l\,e_{j}+\mathbb{Z\,}l\,e_{j^{\prime}}+\mathbb{Z\,}%
l\mathbb{\,}e_{i^{\prime}}+\mathbb{Z\,}\left(  -\,\alpha_{i}\,\,e_{j}%
+\alpha_{i^{\prime}}\,e_{i^{\prime}}+\alpha_{j^{\prime}}\,e_{j^{\prime}%
}\right)  \ ,
\]
the edges $E_{i,j}$ onto $\Phi_{i,j}\left(  E_{i,j}\right)  =$ conv$\left(
\left\{  \mathbf{0},l\mathbb{\,}e_{j}\right\}  \right)  $ and the facets
$F_{i^{\prime}}$ and $F_{j^{\prime}}$ containing $E_{i,j}$ onto
\[
\Phi_{i,j}\left(  F_{i^{\prime}}\right)  =\text{conv}\left(  \left\{
\mathbf{0},\,l\mathbb{\,}e_{j},\,l\left(  e_{i^{\prime}}+e_{j}\right)
\right\}  \right)  \,\,\text{and\thinspace\thinspace}\Phi_{i,j}(F_{j^{\prime}%
})=\text{conv}(\{\mathbf{0},\,l\mathbb{\,}e_{j},\,l(e_{j}+e_{j^{\prime}})\}),
\]
respectively. Since gcd$(-\,\alpha_{i},\alpha_{i^{\prime}},\alpha_{j^{\prime}%
})=$ gcd$(\alpha_{i},\alpha_{i^{\prime}},\alpha_{j^{\prime}})=1$, it is
det$(\widehat{N_{i,j}})=l^{2}$. \ On the other hand,
\[
\Phi_{i,j}\left(  E_{i,j}\right)  \cap\widehat{N_{i,j}}=\text{conv}\left(
\left\{  \mathbf{0},l\mathbb{\,}e_{j}\right\}  \right)  \cap\widehat{N_{i,j}%
}=\sharp\left(  Q_{i^{\prime}}\cap Q_{j^{\prime}}\right)  =\text{ gcd}\left(
\alpha_{i^{\prime}},\alpha_{j^{\prime}},l\right)  ,
\]
i.e. $\frac{l\ e_{j}}{\text{gcd}(\alpha_{i^{\prime}},\alpha_{j^{\prime}}%
,l)}\ $ is a primitive vector. Setting $\widetilde{N_{i,j}}:=\widehat{N_{i,j}%
}/(\mathbb{Z\ }\tfrac{l\ e_{j}}{\text{gcd}(\alpha_{i^{\prime}},\alpha
_{j^{\prime}},l)})$ we get%
\[
\text{det}(\widetilde{N_{i,j}})\cdot\frac{l}{\text{gcd}\left(  \alpha
_{i^{\prime}},\alpha_{j^{\prime}},l\right)  }=\text{det}(\widehat{N_{i,j}%
})\Rightarrow\text{det}(\widetilde{N_{i,j}})=l\cdot\text{gcd}(\alpha
_{i^{\prime}},\alpha_{j^{\prime}},l).
\]
If we denote by $\widetilde{n}_{i^{\prime}},\widetilde{n}_{j^{\prime}}$ the
images of $\,l\left(  e_{i^{\prime}}+e_{j}\right)  ,\ l(e_{j}+e_{j^{\prime}})$
under the canonical epimorphism $\widehat{N_{i,j}}\rightarrow\widetilde
{N_{i,j}}$ (i.e. $\widetilde{n}_{i^{\prime}}=le_{i^{\prime}}$,\ $\widetilde
{n}_{j^{\prime}}=le_{j^{\prime}}$ within $\widetilde{N_{i,j}}$), and by
$\widetilde{\widetilde{n}}_{i^{\prime}}$, $\widetilde{\widetilde{n}%
}_{j^{\prime}\text{ }}$the primitive vectors of conv$\left(  \left\{
\mathbf{0},\widetilde{n}_{i^{\prime}}\right\}  \right)  $ and
conv$(\{\mathbf{0},\widetilde{n}_{j^{\prime}}\})$ w.r.t. $\widetilde{N_{i,j}}%
$, then
\[
\left\Vert \widetilde{\widetilde{n}}_{i^{\prime}}\right\Vert \cdot\frac
{l}{\text{gcd}\left(  \alpha_{i^{\prime}},\alpha_{j^{\prime}},l\right)
}=\text{ det}\left(  \text{lin}\left(  F_{i^{\prime}}\right)  \cap
\widehat{N_{i,j}}\right)
\]%
\[
=\frac{l^{2}}{\sharp\left(  \text{conv}\left(  \left\{  l\mathbb{\,}%
e_{j},\,l\left(  e_{i^{\prime}}+e_{j}\right)  \right\}  \right)  \cap
\widehat{N_{i,j}}\right)  }=\frac{l^{2}}{\text{gcd}\left(  \alpha_{j^{\prime}%
},l\right)  }.
\]
Consequently,
\[
\widetilde{\widetilde{n}}_{i^{\prime}}=\left(  \tfrac{l\cdot\text{gcd}%
(\alpha_{i^{\prime}},\alpha_{j^{\prime}},l)}{\text{gcd}(\alpha_{j^{\prime}%
},l)}\right)  e_{i^{\prime}}\ \ \ \left(  \text{and analogously,
\ \ }\widetilde{\widetilde{n}}_{j^{\prime}}=\left(  \tfrac{l\cdot
\text{gcd}\left(  \alpha_{i^{\prime}},\alpha_{j^{\prime}},l\right)
}{\text{gcd}\left(  \alpha_{i^{\prime}},l\right)  }\right)  e_{j^{\prime}%
}\right)  .
\]
This means that
\[
\mathbf{q}_{i,j}=\frac{\left\Vert \widetilde{\widetilde{n}}_{i^{\prime}%
}\right\Vert \cdot\left\Vert \widetilde{\widetilde{n}}_{j^{\prime}}\right\Vert
}{\text{det}\left(  \widetilde{N_{i,j}}\right)  }=\frac{l\cdot\text{gcd}%
\left(  \alpha_{i^{\prime}},\alpha_{j^{\prime}},l\right)  }{\text{gcd}\left(
\alpha_{i^{\prime}},l\right)  \cdot\text{gcd}\left(  \alpha_{j^{\prime}%
},l\right)  }\ \ .
\]
Since
\[
\widetilde{N_{i,j}}=\mathbb{Z\,}l\,e_{j^{\prime}}+\mathbb{Z\,}l\mathbb{\,}%
e_{i^{\prime}}+\mathbb{Z\,}\left(  \alpha_{i^{\prime}}\,e_{i^{\prime}}%
+\alpha_{j^{\prime}}\,e_{j^{\prime}}\right)  ,
\]
the vector $\gamma_{i^{\prime}}\cdot\left(  \alpha_{i^{\prime}}\,e_{i^{\prime
}}\right)  +$gcd$\left(  \alpha_{j^{\prime}},l\right)  \,e_{j^{\prime}}$
belongs to $\widetilde{N_{i,j}}$, and from
\[
\text{det}(\mathbb{Z\ }\widetilde{\widetilde{n}}_{i^{\prime}}+\mathbb{Z\ (}%
\gamma_{i^{\prime}}\cdot\left(  \alpha_{i^{\prime}}\,e_{i^{\prime}}\right)
+\text{gcd}(\alpha_{j^{\prime}},l)\,e_{j^{\prime}}))=l\cdot\text{gcd}%
(\alpha_{i^{\prime}},\alpha_{j^{\prime}},l)=\text{det}(\widetilde{N_{i,j}})
\]
we conclude that $\mathbb{\ \{}\widetilde{\widetilde{n}}_{i^{\prime}}%
,\gamma_{i^{\prime}}\cdot\left(  \alpha_{i^{\prime}}\,e_{i^{\prime}}\right)
+$gcd$(\alpha_{j^{\prime}},l)\,e_{j^{\prime}}\}$ is a $\mathbb{Z}$-basis of
$\widetilde{N_{i,j}}$ with
\[
\widetilde{\widetilde{n}}_{j^{\prime}}=\left(  \tfrac{l\cdot\text{gcd}\left(
\alpha_{i^{\prime}},\alpha_{j^{\prime}},l\right)  }{\text{gcd}\left(
\alpha_{i^{\prime}},l\right)  }\right)  \ e_{j^{\prime}}\smallskip=\left(
\tfrac{\left(  -\gamma_{i^{\prime}}\right)  \cdot\alpha_{i^{\prime}}%
}{\text{gcd}\left(  \alpha_{i^{\prime}},l\right)  }\right)  \widetilde
{\widetilde{n}}_{i^{\prime}}+\mathbf{q}_{i,j}\cdot\mathbb{(}\gamma_{i^{\prime
}}\cdot\left(  \alpha_{i^{\prime}}\,e_{i^{\prime}}\right)  +\text{gcd}\left(
\alpha_{j^{\prime}},l\right)  \,e_{j^{\prime}}).
\]
Hence, $\mathbf{p}_{i,j}=\left[  \frac{\left(
-\gamma_{i^{\prime}}\right)
\cdot\alpha_{i^{\prime}}}{\text{gcd}\left(
\alpha_{i^{\prime}},l\right) }\right]  _{\mathbf{q}_{i,j}},$ and
the proof is completed after the substitution of (\ref{REVOL1}),
(\ref{REVOL2}), (\ref{REVOL3}) into (\ref{MO-PO1}),
(\ref{MO-PO2}).\hfill {}$\square$

\begin{remark}
One can analogously compute the $\mathbf{a}_{i}\left(  \mathfrak{s}%
_{G}\right)  $'s in the case in which the acting group $G$ is abelian (not
necessarily cyclic), again by lattice transforming and by Mordell-Pommersheim
formula. In these more complicated expressions the greatest common divisors
are replaced by denumerants of restricted weighted vectorial partitions.
\end{remark}

\noindent$\bullet$ \textbf{Diaz-Robins formula}. To present
\textit{Diaz-Robins formula}, by means of which one computes the coefficients
of the Ehrhart polynomial of a lattice simplex of \textit{arbitrary}
dimension, let us first recall the notion of \textit{Hermite normal form}
which will enable us to choose a convenient coordinate system for the simplex vertices.

\begin{theorem}
[{\cite[II.2 and II.3, pp. 15-18]{Newman}}]\label{HNF}For a given integer (or
rational) $\left(  d\times d^{\prime}\right)  $-matrix $A$ of full row rank,
there is a\emph{\ }unimodular matrix $U\in$ \emph{GL}$\left(  d^{\prime
},\mathbb{Z}\right)  $, such that $AU$ is lower-triangular with positive
diagonal elements. Each off-diagonal element of $AU$ is non-negative and
strictly less than the diagonal element in its column. We say that $AU$ is in
\emph{Hermite normal form}$.$ If $\det(A)\neq0,$ then $U$ is uniquely
determined, and $\mathbf{HNF}\left(  A\right)  :=AU$ is \emph{the Hermite
normal form of} $A$.
\end{theorem}

\noindent Now let $\mathbf{s\subset\ }\mathbb{R}^{d}$ be a $d$-dimensional
simplex whose vertices belong to (the standard lattice) $\mathbb{Z}^{d}$.
W.l.o.g., we may assume that $\mathbf{s}=$ conv$\left(  \left\{
\mathbf{0},n_{1},\ldots,n_{d-1},n_{d}\right\}  \right)  .$ The matrix $\left(
n_{1},\ldots,n_{d-1},n_{d}\right)  $ is by Theorem \ref{HNF} left-equivalent
to
\[
\mathbf{\Lambda}_{\mathbf{s}}:=\mathbf{HNF}\left(  \left(  n_{1},n_{2}%
,\ldots,n_{d-1},n_{d}\right)  \right)  =\left(
\begin{array}
[c]{cccc}%
\lambda_{11} & 0 & \cdots & 0\\
\lambda_{21} & \lambda_{22} & \cdots & 0\\
\vdots & \vdots & \ddots & \vdots\\
\lambda_{d1} & \lambda_{d2} & \cdots & \lambda_{dd}%
\end{array}
\right)  .
\]
Consider
\[
\widetilde{\mathbf{\Lambda}_{\mathbf{s}}}:=\left(
\begin{array}
[c]{c}%
\mathbf{\Lambda}_{\mathbf{s}}\ \ 0\\
\mathbf{1}%
\end{array}
\right)  \ ,\ \ \ \ \text{ with\ \ \ \ }\mathbf{1}:=(\underset{\left(
d+1\right)  -\text{times}}{\underbrace{1,1,\ldots,1,1}}).
\]
Define
\[
\varpi_{j}:=\prod_{1\leq i<j}\ \lambda_{ii\ },\ \ \forall j,\ \ 1\leq j\leq
d\ ,\text{\ \ and\ \ \ }\varpi_{d+1}:=\varpi_{d},
\]
\[
\mathfrak{G}_{\mathbf{s}}:=\left(  \mathbb{Z\ }/\
\varpi_{1}\mathbb{Z}\right) \times\left(  \mathbb{Z\ }/\
\varpi_{2}\mathbb{Z}\right)  \times\cdots \times\left(  \mathbb{Z\
}/\ \varpi_{d}\mathbb{Z}\right).
\]
Moreover, for each group element $\mathfrak{g}=\left(  \mathfrak{g}%
_{1},\mathfrak{g}_{2},\ldots,\mathfrak{g}_{d}\right)  \in\mathfrak{G}%
_{\mathbf{s}}$, set:
\[
\boldsymbol{\varepsilon}_{j}\left(  \mathfrak{g}\right)  :=\left\langle
\left(  0,\mathfrak{g}_{1},\mathfrak{g}_{2},\ldots,\mathfrak{g}_{d}\right)
\text{, }(j\text{-th column of }\widetilde{\mathbf{\Lambda}_{\mathbf{s}}%
})\right\rangle .
\]

\begin{theorem}
[{Diaz-Robins formula; cf. \cite{Diaz-Rob1, Diaz-Rob2} and \cite[\S 5]{Chen}}%
]\label{DR-thm}The Ehrhart polynomial
\[
\mathbf{Ehr}_{N}\left(  \mathbf{s},\nu\right)  =\sum_{i=0}^{d}\mathbf{a}%
_{i}\left(  \mathbf{s}\right)  \ \nu^{i}%
\]
of $\mathbf{s}$\ has the following exponential generating function\emph{:}$\ $%
\begin{equation}%
{\textstyle\sum\limits_{\nu=0}^{\infty}} \mathbf{Ehr}_{N}\left(
\mathbf{s},\nu\right)  e^{-2\pi x\nu}=\tfrac
{1}{2^{d+1}\ \left\vert \mathfrak{G}_{\mathbf{s}}\right\vert }%
{\textstyle\sum\limits_{\mathfrak{g}\in\mathfrak{G}_{\mathbf{s}}}}
{\textstyle\prod\limits_{j=1}^{d+1}}
(1\!+\!\text{\emph{coth}}(\tfrac{\pi}{\varpi_{j}}(x\!+\!\sqrt{-1}%
\boldsymbol{\varepsilon}_{j}\left(\mathfrak{g}\right))))
\label{DR-FORM}%
\end{equation}
and $\mathbf{a}_{i}\left(  \mathbf{s}\right)$ equals the
coefficient of $\frac{1}{x^{i+1}}$ in the Laurent expansion at
$x=0$ of
\[
\frac{\pi^{i+1}}{i!\ 2^{d-i}\ \left\vert \mathfrak{G}_{\mathbf{s}}\right\vert
}\ \sum_{\mathfrak{g}\in\mathfrak{G}_{\mathbf{s}}}\ \left[  \prod_{j=1}%
^{d+1}\ (1+\text{\emph{coth}}(\frac{\pi}{\varpi_{j}}(x+\sqrt{-1}%
\ \boldsymbol{\varepsilon}_{j}\left(  \mathfrak{g}\right)  )))\right]  .
\]

\end{theorem}

\noindent{}Let us now describe how can one apply Diaz-Robins formula to the
case of the simplex we are interested in. Let \emph{\ }$(\mathbb{C}%
^{r}/G,[\mathbf{0}]),$ $r\geq4,$ be a Gorenstein AQS of type (\ref{typeAQS}),
and $\mathfrak{s}_{G}=$ conv$\left(  \left\{  e_{1},e_{2},\ldots
,e_{r}\right\}  \right)  $ the corresponding junior simplex. As our lattice
$N_{G}$ is \textquotedblleft skew\textquotedblright, to count
\[
\sharp\left(  \mathfrak{s}_{G}\cap N_{G}\right)  =\sharp\left\{  \left.
\left(  j_{1},..,j_{\kappa}\right)  \in%
{\textstyle\prod\limits_{\mu=1}^{\kappa}} \left(  \left\{
0,1,..,q_{\mu}\right\}  \right)  \right\vert
{\textstyle\sum\limits_{i=1}^{r}} \delta_{i}\left(
j_{1},..,j_{\kappa}\right)  =\text{exp}\left(  G\right) \right\}
+r
\]
we have to transform\emph{\ }$\mathfrak{s}_{G}$ onto another appropriate simplex.

\begin{corollary}
\label{COR-DR}There exists a lattice simplex $\overline{\mathfrak{s}}%
_{G}^{\,\prime\prime}$ w.r.t. $\mathbb{Z}^{r-1}$ such that%
\begin{equation}
\sharp\left(  \overline{\mathfrak{s}}_{G}^{\,\prime\prime}\cap\mathbb{Z}%
^{r-1}\right)  =\emph{\ }\sharp\left(  \mathfrak{s}_{G}\cap N_{G}\right)  ,
\label{DR-TRANSF}%
\end{equation}
and therefore $\sharp\left(  \mathfrak{s}_{G}\cap N_{G}\right)  $ can be
computed by applying formula \emph{(\ref{DR-FORM})} for $\overline
{\mathfrak{s}}_{G}^{\,\prime\prime}.$
\end{corollary}

\noindent{}\textsc{Proof}. This will be done in three steps.\smallskip
\smallskip\newline$\rhd$ \textsc{First step.}\textit{\ }We first perform a
translation in order to insert the zero point as a vertex, and
define $\overline{\mathfrak{s}}_{G}:=\mathfrak{s}_{G}-e_{1}=$
conv$\left(  \left\{
\mathbf{0},e_{2}-e_{1},\ldots,e_{r}-e_{1}\right\}  \right)  $ with
vertex set belonging to the lattice
$\overline{N}_{G}:=%
{\textstyle\sum\limits_{i=2}^{r}}
\,\mathbb{Z\,}\left(  e_{i}-e_{1}\right)  +%
{\textstyle\sum\limits_{\mu=1}^{\kappa}}
\,\mathbb{Z\,}\tfrac{1}{q_{\mu}}\left(  \alpha_{\mu}-\vartheta_{\mu}%
\,e_{1}\right),$ where $\alpha_{\mu}:=\left(
\alpha_{\mu,1},\ldots,\alpha_{\mu,r}\right) ^{\intercal}$ and
$\vartheta_{\mu}:=\sum_{i=1}^{r}\alpha_{\mu,i}$\emph{.
}Obviously,%
\[
\sharp\left(  \mathfrak{s}_{G}\cap N_{G}\right)  =\sharp\left(  \overline
{\mathfrak{s}}_{G}\cap\overline{N}_{G}\right)  .
\]
$\rhd$\emph{\ }\textsc{Second step.}\textit{\ }We define a unimodular
transformation\emph{\ }$\Phi:\mathbb{R}^{r}\longrightarrow$ $\mathbb{R}^{r}$
by
\[
\Phi\left(  e_{1}\right)  :=e_{1}-e_{2},\,\ \,\Phi\left(  e_{2}\right)
:=e_{1},\ \,\Phi\left(  e_{j}\right)  :=e_{1}+e_{j},\,\ \,\forall j,\ \,3\leq
j\leq r,
\]
and linear extension. $\Phi$ maps $\overline{\mathfrak{s}}_{G}$ onto
\[
\overline{\mathfrak{s}}_{G}^{\,\prime}:=\Phi\left(  \overline{\mathfrak{s}%
}_{G}\right)  =\text{conv}\left(  \left\{  \mathbf{0},e_{2},e_{2}+e_{3}%
,e_{2}+e_{4},\ldots,e_{2}+e_{r}\right\}  \right)
\]
and $\overline{N}_{G}$ onto
\[
\overline{N}_{G}^{\,\prime}:=\Phi\left(  \overline{N}_{G}\right)
=\mathbb{Z\ }e_{2}+\mathbb{Z\ }\left(  e_{2}+e_{3}\right)  +\cdots
+\mathbb{Z\ }\left(  e_{2}+e_{r}\right)  +\sum_{\mu=1}^{\kappa}\,\mathbb{Z\,}%
\frac{\alpha_{\mu}^{\prime}}{q_{\mu}},
\]
with
\[
\alpha_{\mu}^{\prime}:=\Phi\left(  \alpha_{\mu}-\vartheta_{\mu}\,e_{1}\right)
=\Phi\left(  \alpha_{\mu}\right)  -\vartheta_{\mu}\ \Phi\left(  e_{1}\right)
\]%
\begin{align*}
=\left(  \vartheta_{\mu}-\alpha_{\mu,1}\right)  e_{2}+\sum_{j=3}%
^{r}\ \alpha_{\mu,j}\ e_{j}\ .
\end{align*}
Thus,
\[
\overline{N}_{G}^{\,\prime}=\mathbb{Z\ }e_{2}+\mathbb{Z\ }e_{3}+\cdots
+\mathbb{Z\ }e_{r}+\sum_{\mu=1}^{\kappa}\,\mathbb{Z\,}\frac{\alpha_{\mu
}^{\prime}}{q_{\mu}}%
\]
and we can consider both $\overline{\mathfrak{s}}_{G}^{\,\prime}$ and
$\overline{N}_{G}^{\,\prime}$ within $\mathbb{R}^{r-1}$ by identifying it with
the set $\left\{  \left(  x_{1},\ldots,x_{r}\right)  \in\mathbb{R}%
^{r}\ \left\vert \ x_{1}=0\right.  \right\}  $\emph{. }Moreover,\emph{\ }%
$\sharp\left(  \overline{\mathfrak{s}}_{G}^{\,\prime}\cap\overline{N}%
_{G}^{\,\prime}\right)  =\sharp\left(  \overline{\mathfrak{s}}_{G}%
\cap\overline{N}_{G}\right)  $\emph{.\smallskip}\newline$\rhd$ \textsc{Third
step.}\textit{\ }Define $A$ to be the rational $\left(  r-1\right)
\times\left(  r-1+\kappa\right)  $-matrix formed by the column vectors which
generate $\overline{N}_{G}^{\,\prime}$:%
\[
A:=(e_{2},e_{3},\ldots,e_{r},\tfrac{1}{q_{1}}\ \alpha_{1}^{\prime}%
,\ldots,\tfrac{1}{q_{\kappa}}\ \alpha_{\kappa}^{\prime}).
\]
According to Theorem \ref{HNF} (with $d=$ $r-1$, $d^{\prime}=r-1+\kappa$), we
can determine a unimodular matrix $U$ such that \ $AU=(R\ \ 0),$ where $R$ is
a rational, non-singular $\left(  r-1\right)  \times\left(  r-1\right)
$-matrix being in Hermite normal form. Hence, if we set\emph{\ }
\begin{align*}
\overline{\mathfrak{s}}_{G}^{\,\prime\prime}  &  :=R^{-1}\overline
{\mathfrak{s}}_{G}^{\,\prime}\\
&  =\text{conv}\left(  \left\{  \mathbf{0,\ }R^{-1}\ e_{2},R^{-1}\ \left(
e_{2}+e_{3}\right)  ,R^{-1}\ \left(  e_{2}+e_{4}\right)  ,\ldots
,R^{-1}\ \left(  e_{2}+e_{r}\right)  \right\}  \right)  ,
\end{align*}
the lattice $\overline{N}_{G}^{\,\prime}$ is transformed onto $\overline
{N}_{G}^{\,^{\prime\prime}}:=R^{-1}\ \overline{N}_{G}^{\,\prime}%
=\mathbb{Z}^{r-1},$\ i.e., onto the \textit{standard} lattice of rank
$r-1$\emph{. }But then we are done because\emph{\ }$\sharp\left(
\overline{\mathfrak{s}}_{G}^{\,\prime\prime}\cap\mathbb{Z}^{r-1}\right)
=\emph{\ }\sharp(\overline{\mathfrak{s}}_{G}^{\,\prime}\cap\overline{N}%
_{G}^{\,\prime}),$ and we can apply Theorem \ref{DR-thm} for the simplex
$\overline{\mathfrak{s}}_{G}^{\,\prime\prime}\subset\mathbb{Z}^{r-1}$\ (with
$d=$ $r-1$).\hfill{}$\square$

\end{document}